\documentclass[12pt]{amsbook}
\allowdisplaybreaks
\usepackage{amssymb}
\usepackage{graphicx}
\usepackage{fullpage}
\usepackage{enumitem}
\usepackage[pdfencoding=auto,colorlinks,linkcolor=,citecolor=]{hyperref}
\usepackage{cite}
\usepackage{dsfont}
\usepackage{eucal}
\usepackage[all,cmtip]{xy}

\newcommand{\Aut}{\mathsf{Aut}}
\newcommand{\bbar}{\mathsf{bar}}

\newcommand{\Ch}{\mathsf{Ch}}
\newcommand{\Cliffq}{\mathsf{Cliff}(\qq)}
\newcommand{\cobar}{\mathsf{cobar}}

\newcommand{\Cotor}{\mathsf{Co tor}}
\newcommand{\D}{{\rm D}}
\newcommand{\kD}{\kk\D}
\newcommand{\Db}{\sfD^{\mathsf{b}}}

\newcommand{\Dsg}{\sfD_{\mathsf{sg}}}
\newcommand{\Dcsg}{\sfD_{\mathsf{csg}}}
\newcommand{\End}{\mathsf{End}}
\newcommand{\iEnd}{\mathcal{E}\mathit{nd}}
\newcommand{\ep}{\varepsilon}
\newcommand{\Ext}{\mathsf{Ext}}

\newcommand{\half}{{\textstyle\frac{1}{2}}}
\newcommand{\halfs}{{\scriptstyle\frac{1}{2}}}
\newcommand{\HH}{H\!H}
\newcommand{\HHinf}{\overset{\raisebox{-.4em}{$\scriptscriptstyle\infty$}}{\HH}{}}
\newcommand{\Hd}{\partial}
\newcommand{\HHd}{\delta}
\newcommand{\Hom}{\mathsf{Hom}}
\newcommand{\iHom}{\mathcal{H}\mathit{om}}

\newcommand{\id}{\mathsf{id}}
\newcommand{\ik}{\mathsf{ik}}
\newcommand{\Inj}{\mathsf{Inj}}
\newcommand{\kG}{{\kk G}}
\newcommand{\KInj}{\mathsf{KInj}}
\newcommand{\Mat}{\mathsf{Mat}}
\newcommand{\MCM}{\mathsf{MCM}}
\newcommand{\om}{{\boldsymbol{\omega}}}
\newcommand{\omb}{{\bar\om}}

\newcommand{\op}{{\mathsf{op}}}
\newcommand{\phat}{{}^{^\wedge}_p}
\newcommand{\twohat}{{}^{^\wedge}_2}
\newcommand{\threehat}{{}^{^\wedge}_3}
\newcommand{\ellhat}{{}^{^\wedge}_\ell}
\newcommand{\Qhat}{{}^{^\wedge}_\bQ}
\newcommand{\Proj}{\mathsf{Proj}}
\newcommand{\Q}{{\rm Q}}
\newcommand{\kQ}{{\kk\Q}}
\newcommand{\res}{\mathop{\rm res}}
\newcommand{\SD}{{\rm SD}}
\newcommand{\kSD}{{\kk\SD}}
\newcommand{\Simp}{\mathsf{Simp}}

\newcommand{\Sol}{\mathsf{Sol}}
\newcommand{\Soc}{\mathsf{Soc}}

\newcommand{\Sq}{\mathsf{Sq}}
\newcommand{\StMod}{\mathsf{StMod}}
\newcommand{\stmod}{\mathsf{stmod}}

\newcommand{\Thick}{\mathsf{Thick}}

\newcommand{\Tor}{\mathsf{Tor}}

\newcommand{\bfa}{\mathbf a}

\newcommand{\bC}{\mathbb C}
\newcommand{\bF}{\mathbb F}

\newcommand{\bQ}{\mathbb Q}

\newcommand{\bZ}{\mathbb Z}
\newcommand{\bfx}{\mathbf x}

\newcommand{\cA}{\mathcal A}
\newcommand{\cC}{\mathcal C}

\newcommand{\cF}{\mathcal F}

\newcommand{\cL}{\mathcal L}

\newcommand{\fa}{\mathfrak{a}}
\newcommand{\fb}{\mathfrak{b}}
\newcommand{\fm}{\mathfrak{m}}
\newcommand{\fp}{\mathfrak{p}}
\newcommand{\kk}{\mathsf{k}}
\newcommand{\sfD}{\mathsf{D}}
\newcommand{\sfK}{\mathsf{K}}
\newcommand{\MM}{\mathsf{M}}
\newcommand{\NN}{\mathsf{N}}
\newcommand{\qq}{\mathsf{q}}

\newcommand{\xx}{\hat{x}}
\newcommand{\yy}{\hat{y}}
\newcommand{\zz}{\hat{z}}

\newcommand{\blue}{\color{blue}}

\renewcommand{\le}{\leqslant}
\renewcommand{\ge}{\geqslant}

\makeindex
\numberwithin{equation}{section}

\theoremstyle{plain}
\newtheorem{lemma}[equation]{Lemma}
\newtheorem{theorem}[equation]{Theorem}
\newtheorem{proposition}[equation]{Proposition}
\newtheorem{corollary}[equation]{Corollary}

\theoremstyle{definition}

\newtheorem{definition}[equation]{Definition}
\newtheorem{example}[equation]{Example}
\newtheorem{notation}[equation]{Notation}

\theoremstyle{remark} 
\newtheorem{remark}[equation]{Remark} 
\newtheorem{remarks}[equation]{Remarks}

\newtheorem{case}{Case}[section]

\iftrue
\author{David J. Benson} 
\address{Institute of Mathematics \\ 
University of Aberdeen \\ 
Aberdeen AB24 3UE \\ 
United Kingdom}
\fi

\title{Classifying spaces of finite groups \\ 
of tame representation type}
\subjclass{Primary: 20J06, Secondary: 16E45, 55P35, 55P60, 55S30}
\keywords{
$A_\infty$ algebras,
cohomology of groups, 
cosingularity categories,
dihedral groups,
generalised quaternion groups,
Hochschild cohomology,
loop spaces,  
Massey products, 
$p$-completed classifying space,  
projective special linear group, 
projective special unitary group,
semidihedral groups,
singularity categories
}

\begin{document}

\begin{abstract}
Thanks to the work of Karin Erdmann, we know a great deal about the
representation theory of blocks of finite groups with tame
representation type. Our purpose here is to examine the $p$-completed 
classifying spaces of these blocks and their loop spaces. We pay
special attention to the $A_\infty$ algebra structures, and
singularity and cosingularity categories.
\end{abstract}

\maketitle

\tableofcontents

\chapter*{Preface}

The purpose of this document is to describe the $A_\infty$ algebra
structures on the cohomology of the classifying space $H^*BG$, and on the homology
on the loop space of the $p$-completed classifying space $H_*\Omega
BG\twohat$, when $G$ is a finite 
group with dihedral, semidihedral, or generalised quaternion Sylow
$2$-subgroups. These are the groups whose group algebras have tame
representation type.  In some cases, we are able to completely
describe the singularity and cosingularity categories of
these, while in others this seems to be more difficult.

Part of the point of this work is to describe a variety of
computational techniques for approaching questions such as these.
The main tool in the dihedral and semidihedral case, as it was in
the cyclic case~\cite{Benson/Greenlees:2021a,Benson/Greenlees:2023a}, 
is to put a grading on the basic algebra of the
principal block. This gives rise to a double grading on group
cohomology, and a triple grading on the Hochschild cohomology of the
group cohomology.
This technique gives us no information in the
generalised quaternion case, but an explicit computation involving
minimal resolutions comes to our rescue in this case.

One curious outcome of this work is that if $G$ has semidihedral or
generalised quaternion Sylow $2$-subgroups, and no normal subgroup of
index two, then $H^*BG$ is formal, meaning that it is quasi-isomorphic
to the ring $H^*BG$ with higher multiplications equal to zero,
see Theorems~\ref{th:SD1-formal} and \ref{th:Q1-formal}. 
The same happens in one of the other cases with semidihedral Sylow
$2$-subgroups, see Theorem~\ref{th:SD2-formal}. 
This also happens
for finite groups with elementary abelian Sylow $2$-subgroups in
characteristic two, but necessary and sufficient
conditions for this to occur are not known.\bigskip

\noindent
{\bf Acknowledgements.}
The investigations described in this work grew out of work with John 
Greenlees~\cite{Benson/Greenlees:2021a,Benson/Greenlees:2023a} in which we
described the situation for finite groups with cyclic Sylow
$p$-subgroups. I would like to take the opportunity to thank 
him for his influence on this work, which would not have been
carried out without his input and encouragement. My thanks go to
the University of Warwick (visits supported by EPSRC grant EP/P031080/1) and the 
Isaac Newton Institute in Cambridge 
(programme \emph{`Groups, representations and applications: new
perspectives’}, supported by EPSRC grant EP/R014604/1), both of whose hospitality allowed
me extensive discussions with Greenlees in 2019, early 2020, and 2022,
leading to the work from which this grew. I thank Srikanth Iyengar for
patiently explaining various pieces of commutative algebra to me,
John Greenlees,
Manuel Rivera and Ran Levi for helping me with the homotopy theory, and 
David Craven for correspondence about the structure of tame blocks and
various errors in the literature. I have flagged these under the index
entry ``errors''. I also thank the referees for helpful feedback.

Finally, this work would probably never have happened without the
confinement imposed by the Covid-19 pandemic, but there's frankly no way I'm
going to thank this curs\`ed virus.

\chapter{Introduction and background}

\section{Motivation}

The reader may wonder why the $p$-completed classifying space of a
finite group and its
loop space are interesting objects of study, and why we are interested
in their singularity and cosingularity categories. The answer comes
from Koszul duality,\index{Koszul duality} as we now explain.

First a bit of background on Koszul duality for finite dimensional
algebras. If $A$ is a finite dimensional algebra over a field $\kk$, with Jacobson
radical\index{Jacobson radical} $J$, then the graded $\kk$-algebra
$E=\Ext^*_A(A/J,A/J)$ naturally
comes with the structure of $A_\infty$-algebra, well defined up to
quasi-isomorphism, and $A/J$ is an $A_\infty$ $A$-$E$-bimodule.
It induces an equivalence of bounded derived categories $\Db(A)\simeq\Db(E)$.
In particular, the algebra $A$ can be reconstructed from $E$ as
the endomorphisms in $\Db(E)$ of the
object $A/J$; see Example~\ref{eg:A-from-Ext}.
If $G$ is a finite $p$-group, $\kk$ has characteristic $p$, and $A$ is the group algebra $\kG$ then
the algebra $E$ is the group cohomology $H^*(G,\kk)$. It follows that
$\kG$ can be
recovered as the derived endomorphisms over
the $A_\infty$ algebra $H^*(G,\kk)$ of the
residue field $\kk$. However, for a more general
finite group, $H^*(G,\kk)$ only ``sees'' the trivial module, and not the
other simple $kG$-modules. Different finite groups may have
quasi-isomorphic cohomology rings as $A_\infty$ algebras, and so $\kG$
cannot be recovered from $H^*(G,\kk)$. Instead, what is recovered this way is
the homology of the loopspace of the $p$-completed classifying space
$C_*\Omega BG\phat$, again as an $A_\infty$ algebra.\medskip

In more detail, 
let $G$ be a finite group and $\kk$ a field of
characteristic $p$. The stable module category\index{stable!module category}
$\StMod(\kG)$ is a
much studied triangulated category\index{triangulated category} with a
symmetric monoidal tensor product making it a tensor triangulated
category with tensor identity the trivial module $\kk$. See for
example Chapter~1 of Benson, Iyengar and Krause~\cite{Benson/Iyengar/Krause:2012a}.
The graded endomorphism ring of $\kk$ is the 
Tate cohomology\index{Tate cohomology} ring 
$\hat H^*(G,\kk)$. This is usually not Noetherian, and very often most
of the products in negative degree are zero, see 
Benson and Carlson~\cite{Benson/Carlson:1992a}. The part in
non-negative degrees is the ordinary cohomology ring $H^*(G,\kk)$,
which is a Noetherian\index{Noetherian cohomology}
graded commutative $\kk$-algebra whose
Krull dimension\index{Krull dimension} is
equal to the $p$-rank of $G$ (Corollary~7.8 of
Quillen~\cite{Quillen:1971b}). The category $\StMod(\kG)$ sits inside
a slightly larger category $\KInj(\kG)$,\index{KInj@$\KInj(\kG)$} the
homotopy category of complexes of injective $\kG$-modules, see Benson
and Krause~\cite{Benson/Krause:2008a}.  The tensor identity is the
injective resolution $\ik$ of the trivial module, and the endomorphism ring of
$\ik$ in $\KInj(\kG)$ is the cohomology ring $H^*(G,\kk)$, 
which is isomorphic to the cohomology of the classifying space, $H^*(BG;\kk)$.
The compact objects in $\KInj(\kG)$
form a triangulated category equivalent to the bounded derived
category $\Db(\kG)$.\index{bounded derived category}
The cohomology $H^*(G,\kk)$ 
should be regarded not just as a $\kk$-algebra, because it carries a
great deal of extra structure. This makes it an $A_\infty$-algebra,
which captures the higher associativity structure of the
endomorphisms, see Section~\ref{se:Ainfinity}. It is even an
$E_\infty$-algebra, capturing all the higher commutativity.

If $G$ is a finite $p$-group then $\kG$ may be recovered from the
$A_\infty$ algebra $H^*(BG;\kk)$ as the endomorphisms of  the trivial
module. This sets up a Koszul duality between $\kG$ and
$H^*(BG;\kk)$. Each is the endomorphism ring of the tensor identity for
the other. If we regard $H^*(BG;\kk)$ as an $E_\infty$ algebra, we can do
better, and recover the Hopf algebra structure on $\kG$. This then
enables us to locate the group $G$ itself inside $\kG$ as the
grouplike elements of this Hopf algebra, namely the elements $x$ that
are sent to $x\otimes x$ by the comultiplication. Note that the
group algebra $\kG$ of a finite $p$-group $G$ does not
determine the group $G$, even up to isomorphism, 
because of the recent negative solution of the modular isomorphism
problem, by Garc{\'\i}a-Lucas, Margolis and
del R{\'\i}o~\cite{Garcia-Lucas/Margolis/delRio:2022a}.%
\index{modular isomorphism problem}

As mentioned above,
if $G$ is not a $p$-group then the $A_\infty$ algebra 
$H^*(BG;\kk)$ only ``sees''
the part of $\KInj(\kG)$ generated by $\ik$. This is all that can be
reconstructed from $H^*(BG;\kk)$. If we take the endomorphism ring of
$\kk$ over the $A_\infty$ algebra $H^*(BG;\kk)$, what we recover is not
$\kG$ but $H_*(\Omega BG\phat;\kk)$, the homology of the loop space of the 
Bousfield--Kan $p$-completion\index{Bousfield--Kan $p$-completion} 
of $BG$, again as an $A_\infty$ algebra. So this may be regarded as a
sort of double Koszul dual of $\kG$. If we again take Koszul dual, and
look at the endomorphism ring of $\kk$ over the $A_\infty$ algebra
$H_*(\Omega BG\phat;\kk)$, we get back to $H^*(BG;\kk)$. This duality 
is investigated extensively in Dwyer, Greenlees and
Iyengar~\cite{Dwyer/Greenlees/Iyengar:2006a}. 

The upshot of this is that $H_*(\Omega BG\phat;\kk)$ acts as a sort of
homotopical proxy for $\kG$. Its singularity category is closely
related to the singularity category of $\kG$, which is the small stable
module category $\stmod(\kG)$ consisting of the finitely generated
$\kG$-modules. 

Here is a diagram of the Koszul dualities discussed above.
\[ \xymatrix@=2cm{\kG\ar@{~>}[r] & H^*(BG;\kk)\ar@{~>}@<2ex>[d]\\
&H_*(\Omega BG\phat;\kk)\ar@{~>}@<2ex>[u]} \]

The next point to make is that the Bousfield--Kan
$p$-completion of $BG$ only sees $G$ through the eyes of fusion. If
we're given $BG$, we can recover $G$ as $\pi_1(BG)$. But
$\pi_1(BG\phat)$ is the largest $p$-group quotient of $G$, so for
example if $G$ has no normal subgroup of index $p$ then 
$\pi_1(BG\phat)$ is trivial. Knowing $BG\phat$ is equivalent to
knowing a Sylow $p$-subgroup, and the fusion system and linking system
of the group $G$. So once we $p$-complete, we may as well be talking
about fusion systems and linking systems, rather than finite
groups, see Section~\ref{se:fusion}. This has the advantage that there
are \emph{exotic} fusion systems that don't come from finite groups. A
remarkable theorem of Chermak shows that given a fusion system, there
exists a linking system lying over it, and it is unique, see
Theorem~\ref{th:Chermak}. This has the advantage that a fusion system
on a finite $p$-group is enough information to determine a
$p$-completed classifying space. For example, a block of a finite
group algebra has an associated fusion system on its defect group,
using fusion of subpairs, so a block has a $p$-completed classifying
space. However, it is not at all clear that this the cochains on this
classifying space are the endomorphism $A_\infty$ algebra of any
object in $\KInj$ of the block.\medskip

This paper is a sequel to Benson and
Greenlees~\cite{Benson/Greenlees:2021a,Benson/Greenlees:2023a}. 
In those papers, we examined the
structure of $H^*BG$ and $H_*\Omega BG\phat$ as $A_\infty$ 
algebras, for finite groups
$G$ with cyclic Sylow $p$-subgroups,\index{Sylow subgroup!cyclic} 
over a field $\kk$ of
characteristic $p$, and elucidated the structure of
their singularity and cosingularity categories.
Here, we examine the case of
blocks with dihedral, semidihedral, or (generalised) quaternion 
defect groups in characteristic
two. The reason for the interest in this class of blocks is that these
are the ones where the group algebra has tame representation type, see
Section~\ref{se:tame}, so we are studying the next interesting case
after the cyclic defect case. 
We have a great deal of information about the module categories, from
the work of Erdmann~\cite{Erdmann:1977a,
Erdmann:1979a,
Erdmann:1987a,
Erdmann:1988a,
Erdmann:1988b,
Erdmann:1988c,
Erdmann:1988d,
Erdmann:1990a,
Erdmann:1990b,
Erdmann:1990c,
Erdmann:1992a}, 
and this allows us to do some very explicit computations.

In the final chapter we give a small glimpse
into what happens for blocks of wild representation type. There we run
into a dichotomy first discovered by Ran Levi in his thesis~\cite{Levi:1995a}, between
polynomial growth and (semi-)exponential growth of the homology of $\Omega
BG\phat$, see Section~\ref{se:poly-exp}. We discuss examples of both
types, and end with some open problems.

\section{Summary of results}

Tame blocks of finite groups have 
dihedral, semidihedral, or (generalised) quaternion
defect groups. We discuss these three cases in turn 
in the next three chapters. The results 
are given in detail in the introductory sections of the three
chapters, dealing with the three types of defect groups.
The ring structure on the homology of the loop space
$H_*\Omega BG\twohat$ was first computed by
Levi~\cite{Levi:1995a} in these cases.
The following table shows where the discussions of the various 
cases can be found. Note that for tame type, the homotopy type of $BG$
is determined by the number of conjugacy classes of elements of order
two and four; see Sections~\ref{se:dihedral-Sylow},
\ref{se:semidihedral-Sylow}, and
\ref{se:generalised-quaternion-Sylow}. For definition of the
Hochschild cohomology $\HHinf^*\fa$ of an $A_\infty$ algebra $\fa$,
see Section~\ref{se:HH}, and for the isomorphism $\HHinf H^*BG\cong
\HHinf^*H_*\Omega BG\phat$ see Theorem~\ref{th:HHOmega}.\bigskip

\begin{center}
\renewcommand{\arraystretch}{1.2}
\begin{tabular}{|c|c|c|c|c|c|c|}
\hline
Sylow& ccls of elts
&\multicolumn{2}{c|}{$\HH^*$ of}
&$\HHinf^*H^*BG$&\multicolumn{2}{c|}{$A_\infty$ structure of} 
\\ \cline{3-4}\cline{6-7}
$p$-subgroup&of order&$H^*BG$&$H_*\Omega BG\twohat$&$\cong$
&$H^*BG$&$H_*\Omega BG\twohat$\\
& 2\qquad 4 &&&$\HHinf^*H_*\Omega BG\twohat$&& \\ \hline
Cyclic & ---\quad\ --- &\cite{Benson/Greenlees:2023a} &
\cite{Benson/Greenlees:2023a}&\cite{Benson/Greenlees:2023a}&
\cite{Benson/Greenlees:2021a}&\cite{Benson/Greenlees:2021a} \\ \hline
& 1\qquad 1&\ref{pr:HHHBG-D1}&\ref{pr:HHHOmegaBG-D1}
&\ref{rk:D1-HHCBG}&\ref{th:AinftyHBG-D1}&\ref{th:AinftyHOmegaBG-D1} \\
Dihedral& 2\qquad 1&\ref{pr:HHHBG-D2}&
\ref{pr:HHHOmegaBG-D2}&\ref{rk:D2-HHCBG}&
\ref{th:AinftyHBG-D2}&\ref{th:AinftyHOmegaBG-D2} \\
& 3\qquad 1&\ref{th:HHHBD}&
\cite{Siegel/Witherspoon:1999a}, \ref{th:HHkD}&
\cite{Siegel/Witherspoon:1999a}, \ref{th:HHkD}&
\ref{th:AinftyHBG}&\ref{se:OmegaBDtwohat} \\ \hline
&1\qquad 1&\ref{th:HHHBG-SD1}&\ref{th:HOmegaBGtwohat-SD1}
&\ref{th:SD1-HHCBG}&\ref{th:SD1-formal}&\ref{th:SD1-m3} \\ 
Semidihedral&2\qquad 1&\ref{th:HHHBG-SD2}&
\ref{th:SD2-HHHOmega}&\ref{th:SD2-HHCBG}&
\ref{th:SD2-formal}&\ref{th:SD2-m3} \\
&1\qquad 2&&&&(\ref{se:SD3})&(\ref{se:SD3}) \\
&2\qquad 2&&\cite{Generalov/Nikulin:2020a}&
\cite{Generalov/Nikulin:2020a}&&\ref{se:OmegaBSDtwohat} \\
\hline
Generalised&1\qquad 1&\ref{co:Q1-formal}&\ref{co:Q1-formal}&
\ref{co:Q1-formal}&\ref{th:Q1-formal}&\ref{co:Q1-formal} \\
quaternion&1\qquad  2&\ref{th:HHHBQ2classes}&&&\ref{th:Qu-non-formal}& \\
&1\qquad 3&\ref{th:Q1-HHHBG}&\cite{Hayami:2006a}&
\cite{Hayami:2006a}&&\ref{se:OmegaBQtwohat} \\ 
\hline
\end{tabular}\bigskip
\end{center}

For ease of reference, we list here the cohomology algebras 
in the various cases of Sylow subgroups, and the
loop space homology in those cases where the group is not
$p$-nilpotent,\index{p-nilpotent@$p$-nilpotent finite group}
because if $G$ is $p$-nilpotent then $H_*\Omega BG\phat \cong
\kG/O^p(G)$ (see Notation~\ref{no:gt}).
The degrees are written homologically, followed by the
degrees coming from the internal grading in the 
cases where they exist.\bigskip\pagebreak[3]

\noindent
{\bf Cohomology algebras} $H^*BG\cong \Ext^*_{\kG}(\kk,\kk)$\bigskip

\noindent
Cyclic, order $p^n$, inertial index $q|(p-1)$:
\[ \kk[x] \otimes \Lambda(t),\qquad |x|=-(2q,p^n),\qquad |t|=-(2q-1,p^n-(p^n-1)/q).\smallskip \]

\noindent
Klein four group, one class of involutions:
\[ \kk[\xi,\eta,t]/(\xi\eta+t^3),\qquad |\xi|=|\eta|=-(3,3),\qquad |t|=-(2,2).\smallskip \]

\noindent
Klein four group, three classes of involutions:
\[ \kk[x,y], \qquad |x|=|y|=-(1,1).\smallskip \]

\noindent
Dihedral, order $4q$ ($q\ge 2$), one class of involutions:
\[ \kk[\xi,\eta,t]/(\xi\eta),\quad |\xi|=-(3,q+1,q),
\quad |\eta|=-(3,q,q+1),
\quad |t|=-(2,q,q).\smallskip \]

\noindent
Dihedral, order $4q$ ($q\ge 2$), two classes of involutions:
\[ \kk[\xi,y,t]/(\xi y),\quad |\xi|=-(3,q+1,q),
\quad |y|=-(1,0,1), \quad |t|=-(2,q,q).\smallskip \]

\noindent
Dihedral, order $4q$ ($q\ge 2$), three classes of involutions:
\[ \kk[x,y,t]/(xy),\quad |x|=-(1,1,0),\quad
|y|=-(1,0,1), \quad |t|=-(2,q,q).\smallskip \]

\noindent
Semidihedral, order $8q$ ($q\ge 2$), one class of involutions, one of
elements of order four:
\[ \kk[x,y,z]/(x^2y+z^2), \qquad |x|=-(3,q+1),
\qquad |y|=-(4,4q),\qquad |z|=-(5,3q+1).\smallskip \]

\noindent
Semidihedral, order $8q$ ($q\ge 2$), two classes of involutions, one
of elements of order four:
\[ \kk[x,y,z]/(x^2y+z^2), \qquad |x|=-(1,1-q),
\qquad |y|=-(4,4q),\qquad |z|=-(3,q+1).\smallskip \]

\noindent
Semidihedral, order $8q$ ($q\ge 2$), one class of involutions, two of
elements of order four:
\[ \kk[y,z,w,v]/(y^3,vy,yz,v^2+z^2w),\qquad
|y|=-1,\quad
|z|=-3,\quad
|w|=-4,\quad
|v|=-5.\smallskip \]

\noindent
Semidihedral, order $8q$ ($q\ge 2$), two classes of involutions, two
of elements of order four:
\[ \kk[x,y,z,w]/(xy,y^3,yz,z^2+x^2w),\qquad |x|=|y|=-1,\qquad |z|=-3, 
\qquad  |w|=-4.\smallskip \]

\noindent
Quaternion or generalised quaternion of order $8q$, one class of elements of order four:
\[ \kk[z]\otimes\Lambda(y),\qquad
|z|=-4,\qquad |y|=-3.\smallskip \]

\noindent
Generalised quaternion of order $8q$, two classes of elements of order four:
\[ \kk[y,z]/(y^4), \qquad |y|=-1, \qquad
|z|=-4.\smallskip \]

\noindent
Quaternion of order $8$, three classes of elements of order four:
\[ \kk[u,v,z]/(u^2+uv+v^2,u^2v+uv^2), \qquad
|u|=|v|=-1, \qquad
|z|=-4.\smallskip \]

\noindent
Generalised quaternion of order $8q$ ($q\ge 2$), three classes of
elements of order four:
\[ \kk[x,y,z]/(xy,x^3+y^3),\qquad
|x|=|y|=-1,\qquad
|z|=-4.\smallskip \]

\bigskip\pagebreak[3]

\noindent
{\bf Loop space homology} $H_*\Omega BG\phat$\bigskip

\noindent
Cyclic, order $p^n$, inertial index $q|(p-1)$ ($q\ge 2$):
\[ \kk[\tau]\otimes \Lambda(\xi),\qquad
|\tau|=(2q-2,p^n-(p^n-1)/q), \qquad
|\xi|=(2q-1,p^n). \]
Exception, $p^n=3$, $q=2$:\qquad $\kk[\tau,\xi]/(\xi^2+\tau^3)$,\qquad
$|\tau|=(2,2)$,\qquad $|\xi|=(3,3)$.\medskip

\noindent
Dihedral, order $4q$ ($q\ge 1$), one class of involutions:
\[ \Lambda(\tau) \otimes \kk\langle \alpha,\beta\mid
\alpha^2=0,\beta^2=0\rangle,\quad
|\tau|=(1,q,q),\quad
|\alpha|=(2,q+1,q),\quad
|\beta|=(2,q,q+1).\smallskip \]

\noindent
Dihedral, order $4q$ ($q\ge 2$), two classes of involutions:
\[ \Lambda(\tau)\otimes \kk\langle \alpha,Y\mid
\alpha^2=0,Y^2=0\rangle,\quad
|\tau|=(1,q,q),\quad 
|\alpha|=(2,q+1,q),\quad
|Y|=(0,0,1).\smallskip \]

\noindent
Semidihedral, order $8q$ ($q\ge 2$), one class of involutions, one of
elements of order four:
\[ \Lambda(\hat x,\hat y) \otimes \kk[\hat z],\qquad
|\hat x|=(2,q+1),\qquad
|\hat y|=(3,4q),\qquad
|\hat z|=(4,3q+1).\smallskip \]

\noindent
Semidihedral, order $8q$ ($q\ge 2$), two classes of involutions, one of
elements of order four:
\[ \Lambda(\hat x,\hat y)\otimes \kk[\hat z],\qquad
|\hat x|=(0,1-q),\qquad
|\hat y|=(3,4q),\qquad
|\hat z|=(2,q+1).\smallskip \]

\noindent
Semidihedral, order $8q$ ($q\ge 2$), one class of involutions, two of
elements of order four:
\[ \Lambda(\eta)\otimes 
\kk\langle \hat y,\hat z\mid \hat y^2=\hat z^2=0\rangle,\qquad
|\eta|=1,\qquad
|\hat y|=0,\qquad
|\hat z|=2.\smallskip \]

\noindent
Quaternion or generalised quaternion, one class of elements of order
four:
\[ \Lambda(\hat z) \otimes \kk[\hat y], \qquad
|\hat z|=3,\qquad |\hat y|= 2.\smallskip \]

\noindent
Quaternion or generalised quaternion, two classes of elements of order
four:
\[ \Lambda(\hat y,\hat z)\otimes \kk[\eta], \qquad
|\hat y|=0, \qquad
|\hat z|=3,\qquad
|\eta|=2.\smallskip \]

\section{Notation and conventions}

In this chapter, we give some of the background, and set the stage for  
this project. We begin with the notations and conventions used throughout.\bigskip

\begin{notation}\label{no:gt}
 We use the following standard group theoretic
notations. For a finite group $G$, we write 
$O_p(G)$\index{OpG@$O_p(G)$, $O_{p'}(G)$, $O(G)$, $O^p(G)$, $O^{p'}(G)$} 
for the largest normal $p$-subgroup of $G$ and 
$O_{p'}(G)$
for the largest normal $p'$-subgroup, i.e., the largest normal
subgroup of order not divisible by $p$. When $p=2$, we write
$O(G)$ for $O_{2'}(G)$, the largest normal odd
order subgroup of $G$.

We write $O^p(G)$
for the smallest normal subgroup of $G$ for which
the quotient is a $p$-group, and $O^{p'}(G)$
for the smallest normal
subgroup for which the quotient is a $p'$-group.
\end{notation}

\begin{notation}\label{no:GLetc}
We write $\Gamma L(n,p^m)$,\index{GammaL@$\Gamma L(n,p^m)$} 
$GL(n,p^m)$,\index{GL@$GL(n,p^m)$} and
$SL(n,p^m)$\index{SL@$SL(n,p^m)$} 
for the groups of semi-linear automorphisms, linear
automorphisms, and special (i.e., determinant one) linear
automorphisms of a vector space of dimension $n$ over $\bF_{p^m}$.
We write $P\Gamma L(n,p^m)$,\index{PGammaL@$P\Gamma L(n,p^m)$}
$PGL(n,p^m)$,\index{PGL@$PGL(n,p^m)$} and 
$PSL(n,p^m)$\index{PSL@$PSL(n,p^m)$} for the
corresponding groups of 
projective transformations,\index{projective!transformations} 
namely the quotients of these groups by
the subgroups of scalar transformations. Similarly, we write
$\Gamma U(n,p^m)$\index{GammaU@$\Gamma U(n,p^m)$},
$GU(n,p^m)$\index{GU@$GU(n,p^m)$}, and
$SU(n,p^m)$\index{SU@$SU(n,p^m)$} for the corresponding groups of 
unitary\index{unitary group}
transformations of a vector space of dimension $n$ over $\bF_{p^{2m}}$
with respect to the conjugation given by the 
Frobenius automorphism\index{Frobenius!automorphism} of 
order two of the field. We write
$P\Gamma U(n,p^m)$\index{PGammaU@$P\Gamma U(n,p^m)$},
$PGU(n,p^m)$\index{PGU@$PGU(n,p^m)$}, and
$PSU(n,p^m)$\index{PSU@$PSU(n,p^m)$} 
for the quotient by the subgroup of scalar transformations.
Closely related groups $SL^\pm(2,p^m)$, $SU^\pm(2,p^m)$ will be 
defined in Section~\ref{se:semidihedral-Sylow}, and their
isoclinic groups $SL^\circ(2,p^m)$ and $SU^\circ(2,p^m)$ will be
defined in Section~\ref{se:isoclinism}. The group $PGL^*(2,p^{2m})$ is
also defined in Section~\ref{se:semidihedral-Sylow}.
\end{notation}

\begin{notation}\label{no:C*}
All chains, cochains, homology and cohomology will have coefficients
in a field $\kk$. So when we write $C^*X$, $H^*X$, $C_*X$, $H_*X$, we
mean $C^*(X;\kk)$, $H^*(X;\kk)$, $C_*(X;\kk)$, $H_*(X;\kk)$
respectively. Since we are interested in both homology and cohomology,
we shall write all degrees homologically, so that for example
cohomology elements are given negative degrees. We write $\tilde H_*$ for
reduced homology, which is the same as ordinary homology in positive
degrees, but in degree zero it is the kernel of the augmentation map
$H_0\to\bZ$.

We shall frequently use the fact that $C^*BG$, regarded as a
differential graded algebra, is quasi-isomorphic to the
differential graded algebra of endomorphisms
of a projective resolution of the trivial module, $\End_{\kG}(P_*)$.
Such a quasi-isomorphism induces an isomorphism in cohomology 
$H^*BG\cong\Ext^*_{\kG}(\kk,\kk)$. The link between $C^*BG$ and $\End_\kG(P_*)$ is given
by the Rothenberg--Steenrod construction,\index{Rothenberg--Steenrod construction} 
as explained for example in Rothenberg and Steenrod~\cite{Rothenberg/Steenrod:1965a}, or
Section~4 of Benson and Krause~\cite{Benson/Krause:2008a}.
\end{notation}

\section{\texorpdfstring{$A_\infty$ algebras and quasi-isomorphisms}
{A∞ algebras and quasi-isomorphisms}}\label{se:Ainfinity}

The concept of $A_\infty$ algebra was introduced by  
Stasheff~\cite{Stasheff:1963a,Stasheff:1963b}, and further information
can be found in 
Boardman and Vogt~\cite{Boardman/Vogt:1973a},
Kadeishvili~\cite{Kadeishvili:1982a},
Keller~\cite{Keller:2001a,Keller:2002a,Keller:2006b},
Section~9.2 of Loday and Vallette~\cite{Loday/Vallette:2012a},
Chapter~7 of Witherspoon~\cite{Witherspoon:2019a},
and 
the theses of Lef\`evre-Hasegawa~\cite{Lefevre-Hasegawa:2003a} and
Prout\'e~\cite{Proute:2011a}. The basic idea is to
formalise the notion of being ``associative up to all higher
homotopies''. The remarkable thing is that this notion turns out to be
expressible in a natural way, entirely algebraically.

Recall that an $A_{\infty}$ algebra\index{Ainfinity@$A_\infty$ algebra} 
over a field $\kk$ is a $\bZ$-graded vector space
$\fa$ with graded maps $m_n\colon \fa^{\otimes n}\rightarrow \fa$ of degree $n-2$ for
$n\ge 1$ satisfying 
\begin{equation}\label{eq:Ainfty} 
\sum_{r+s+t=n}(-1)^{r+st}m_{r+1+t}(\id^{\otimes r}\otimes m_s \otimes
\id^{\otimes t})=0 
\end{equation}
for $n\ge 1$. 

\begin{remark}
The cases of~\eqref{eq:Ainfty} with $n=1$, $2$, $3$ are as follows.
\begin{align*}
m_1 m_1 &= 0,\\
m_1m_2 &=m_2(m_1\otimes \id + \id \otimes m_1),\\
m_2(\id \otimes m_2 - m_2\otimes \id)&=m_1m_3 
+m_3(m_1\otimes \id \otimes \id + \id \otimes m_1 \otimes \id 
+ \id \otimes \id \otimes m_1).
\end{align*}
The map $m_1$ is therefore a differential. The map
$m_2$ is not necessarily associative, but it satisfies the Leibniz
rule, making $m_1$ a derivation with
respect to $m_2$. Furthermore, $m_2$
induces an associative product on $H_*\fa$. 
If $m_1=0$, so that $\fa=H_*\fa$, then
$m_3$ can be regarded as a Massey triple
product\index{Massey product} on $\fa$. 
\end{remark}

\begin{example}
A DG algebra\index{DG!algebra} 
(differential graded algebra)\index{differential!graded algebra} 
$\fa$ can be regarded as an $A_\infty$
algebra with $m_1$ the differential, $m_2$ the product,
and $m_i=0$ for $i>2$. The identities~\eqref{eq:Ainfty} are easy to
verify in this case. The identities for $n=1$ and $2$ say that $m_1$
is a differential obeying the Leibniz rule. The right hand side for
$n=3$ is zero, so the equation says that $m_2$ is associative. For
$n\ge 4$, the equation says that $0=0$.
\end{example}

A \emph{morphism} of $A_\infty$ algebras\index{morphism of $A_\infty$ algebras} 
$f\colon \fa \to \fa'$ consists of graded maps
$f_n \colon \fa^{\otimes n} \to \fa'$ of degree $n-1$ satisfying
\begin{equation}\label{eq:finfty} 
\sum_{r+s+t=n} (-1)^{r+st}
f_{r+1+t}(\id^{\otimes r}\otimes m_s\otimes \id^{\otimes t})=
\sum_{\substack{1\le r\le n\\i_1+\dots+i_r=n}} (-1)^\sigma m'_r(f_{i_1}\otimes \dots \otimes f_{i_r})  
\end{equation}
where in the sum on right hand side, 
$\sigma=\sum_{j=1}^{r-1}(r-j)(i_j-1)$ and $m'_r$ are the operations
in $\fa'$.  The \emph{composition}\index{composition of $A_\infty$ morphisms} 
of two morphisms $f\colon\fa\to\fa'$ and
$g\colon\fa'\to\fa''$ is given by
\[ (g\circ f)_n=\sum_{\substack{1\le r\le n\\i_1+\dots+i_r=n}}
(-1)^\sigma g_r(f_{i_1}\otimes \cdots\otimes f_{i_r}) \]
where $\sigma$ is as above.

\begin{remark}
The first two cases of~\eqref{eq:finfty} are as follows.
\begin{align*}
f_1m_1&=m'_1f_1,\\
f_1m_2&=m'_2(f_1\otimes f_1)+m'_1f_2
+f_2(m_1\otimes \id +\id \otimes m_1).
\end{align*}
Thus $f_1$ is a map of complexes with respect to the differential
$m_1$, and commutes with the product $m_2$ up to a homotopy given by
$f_2$. 
\end{remark}

\begin{definition}
The morphism $f$ is said to  
be a \emph{quasi-isomorphism}\index{quasi-isomorphism} 
of $A_\infty$ algebras  
if $f_1$ induces an isomorphism in homology with respect to the  
differentials $m_1,m'_1$.  
 
We say that an $A_\infty$ algebra $\fa$ is 
\emph{formal}\index{formal $A_\infty$ algebra}
if it is quasi-isomorphic to an $A_\infty$ algebra with $m_i=0$ for
$i\ne 2$. 
\end{definition}

\begin{theorem}\label{th:Ainfinityqi}
A quasi-isomorphism of $A_\infty$ algebras is a homotopy equivalence.
\end{theorem}
\begin{proof}
See Section~4 of~\cite{Keller:2001a}, Propostion~2.4.1.1 and
Corollaire~2.4.2.2 of~\cite{Lefevre-Hasegawa:2003a}.  
\end{proof}

\begin{example}
If $\fa$ is a DG algebra then it is formal 
as an $A_\infty$ algebra if and only if it is formal as a DG algebra;
namely if and only if there are quasi-isomorphisms of DG algebras 
$\fa\leftarrow \fa' \rightarrow H_*\fa$, where $H_*\fa$ is regarded as 
a DG algebra with zero differential. 
\end{example}

A theorem of Kadeishvili~\cite{Kadeishvili:1982a} (see also
Keller~\cite{Keller:2001a,Keller:2002a},
Merkulov~\cite{Merkulov:1999a}, Petersen~\cite{Petersen:2020a}) may be stated as
follows. 

\begin{theorem}\label{th:Kadeishvili}\index{Kadeishvili's Theorem}
Suppose that we are given an $A_\infty$
algebra $\fa$ over a field $\kk$. Let $Z_*(\fa)$ be the cocycles,
$B_*(\fa)$ be the coboundaries, and $H_*(\fa)=Z_*(\fa)/B_*(\fa)$,
with respect to the differential $m_1$.
Choose a vector space splitting $f_1\colon H_*(\fa) \to Z_*(\fa)
\subseteq \fa$
of the quotient. Then the homology $H_*\fa$ has an $A_\infty$ structure
with $m_1=0$ and $m_2$ the multiplication on $H_*\fa$ induced by the
multiplication on $\fa$,  and $f_1$ extends to a quasi-isomorphism of $A_\infty$
algebras $f\colon H_*\fa \to \fa$.

If $\fa$ happens to carry auxiliary gradings respected by the maps $m_i$
then $H_*\fa$ inherits the grading, and the $A_\infty$
algebra structure maps $m_i$ on $H_*\fa$ and the quasi-isomorphism $f$ can be
chosen to respect the gradings.
\end{theorem}
\begin{proof}
The idea of the proof is an inductive procedure which goes as follows.
The morphisms 
\[ m_2(f_1\otimes f_1),\ f_1m_2 \colon H_*\fa\otimes H_*\fa \to \fa \]
are homotopic. In the presence of a grading on $\fa$, these maps preserve the
grading, so a homogeneous homotopy can be chosen. 
This is a morphism $f_2 \colon H_*\fa \otimes H_*\fa \to \fa$ of degree $(1,0)$
such that 
\[ f_1m_2 = m_1f_2 + m_2(f_1\otimes f_1). \]
Next, consider the map
\[ m_2(f_1\otimes f_2 - f_2 \otimes f_1)  + f_2(1\otimes
  m_2-m_2\otimes 1) \colon (H_*\fa)^{\otimes 3} \to \fa \]
of degree $(1,0)$. Composing with $m_1\colon \fa \to \fa$, a short
calculation shows that we get zero. So we can add something in the
image of $f_1$ to get a coboundary. Thus there exist maps $m_3\colon
(H_*\fa)^{\otimes 3} \to H_*\fa$ of degree $(1,0)$ and $f_3 \colon (H_*\fa)^{\otimes 3} \to \fa$
of degree $(2,0)$ such that 
\[ f_1m_3 - m_1f_3 = m_2(f_1\otimes f_2-f_2\otimes f_1)
+f_2(1\otimes m_2-m_2\otimes 1). \]
Continuing this way, we obtain maps $m_i$ of degree $(i-2,0)$ giving
an $A_\infty$ structure on $H_*\fa$, and
$f_i$ of degree $(i-1,0)$ giving a quasi-isomorphism $H_*\fa \to \fa$.
Note that equation~\eqref{eq:finfty} simplifies slightly because
$m_1=0$ on $H_*\fa$. If $\fa$ is a DG algebra then
it simplifies further, because $m_i=0$ for $i>2$ on $\fa$. But we can't
just write down formulas for the $m_i$ and $f_i$, because at each
stage we are choosing a splitting for a surjection of (graded) vector
spaces. Continuing the proof in the above fashion is messy, and a
cleaner way to do it may be found in Kadeishvili~\cite{Kadeishvili:1982a} and
Section~3 of Petersen~\cite{Petersen:2020a}.
\end{proof}

\begin{definition}\label{def:minimal}
We say that an $A_\infty$ algebra $\fa$ is
\emph{minimal}\index{minimal!$A_\infty$ algebra} if $m_1=0$.
Given any $A_\infty$ algebra $\fa$,
a minimal $A_\infty$ algebra quasi-isomorphic to $\fa$ is called a 
\emph{minimal model}\index{minimal!model} for $\fa$. Kadeishvili's
Theorem~\ref{th:Kadeishvili} says that every $A_\infty$ algebra $\fa$
has a minimal model, and it is $H_*\fa$ endowed with a suitable
$A_\infty$ structure. 
If $\fa$ carries an internal grading preserved
by the structure maps then the structure maps on a minimal model may
be taken to preserve the internal grading.
\end{definition}

\begin{example}\label{eg:A-from-Ext}
Let $A$ be a finite dimensional algebra over an algebraically closed
field $\kk$, and let $J$ be the Jacobson radical\index{Jacobson radical}
of $A$. Let $P_*$ be a projective resolution of $A/J$ over $A$. Then
the endomorphism ring $\fa=\End_A(P_*)$ is a differential graded algebra.
Its homology is $H_*\fa=\Ext^{-*}_A(A/J,A/J)$, and this becomes an $A_\infty$
algebra by Kadeishvili's Theorem~\ref{th:Kadeishvili}. The basic algebra of $A$ can be
recovered as a quiver with relations from $H_*\fa$ by a recipe described in
Section~2 of Keller~\cite{Keller:2002a}. This recipe only uses the
part of the $A_\infty$ structure that involves
$\Ext^n_A(A/J,A/J)$ for $n=0$, $1$ and $2$, so the rest is determined by this. In this
sense, the $\Ext$ algebra is determined by generators in degree one
and $A_\infty$ relations in
degree two.

In more detail, the recipe described in~\cite{Keller:2002a} is as
follows. Let us suppose that $A$ is a basic algebra, and therefore
given by a quiver with relations, $\kk Q/I$. Let $J_Q$ be the ideal in
$kQ$ generated by the arrows, so that $I\subseteq J_Q^2$.
Then $A/J\cong kQ/J_Q$, and $\Ext^1_A(A/J,A/J)$ is dual to
$J_Q/J_Q^2\cong J/J^2$. The quiver algebra $\kk Q$ is
the tensor algebra over $\kk Q/J_Q$ on $J_Q/J_Q^2$. Then $\Ext^2_A(A/J,A/J)$ is
the dual of $I/(IJ_Q+J_QI)$, and
the maps
\[ m_i\colon \Ext^1_A(A/J,A/J)^{\otimes i} \to \Ext^2_A(A/J,A/J) \]
dualise to give a map
\[ I/(IJ_Q+J_QI) \to \prod_{i=1}^\infty (J_Q/J_Q^2)^{\otimes i}. \]
Since $A$ is finite dimensional, all high enough powers of $J_Q/J_Q^2$ are
generated by the image of this map, and the corresponding
ideal in
\[ \kk Q=\kk Q/J_Q\oplus \bigoplus_{i=1}^\infty (J_Q/J_Q^2)^{\otimes i} \]
determined by the duals of the $m_i$ is exactly $I$.
\end{example}

\begin{example}\label{eg:Einfinity}
Let $G$ be a finite group and $\kk$ a field. Let $P_*$ be a projective
resolution of $\kk$ over $\kG$. Then the endomorphism ring
$\fa=\End_\kG(P_*)$ is a differential graded algebra. Its homology is $H_*\fa=
H^{-*}(G,\kk)=\Ext^{-*}_{\kG}(\kk,\kk)$ is an $A_\infty$ algebra. However,
the algebra $\kG$ can only be recovered from this if $\kk$ is the only
simple $\kG$-module, namely if $G$ is a finite $p$-group. If $G$ is not
a finite $p$-group then all that can be recovered from $H^*(G,\kk)$ is a
Sylow $p$-subgroup and the fusion system on it. See
Section~\ref{se:fusion} for further details. Even recovering the Sylow
$p$-subgroup requires us to regard $H^*(G,\kk)$ as an $E_\infty$
algebra, a structure which encodes all the higher commutativities as
well as the higher associativities. To see this, consider the case
of a finite $p$-group. Two finite $p$-groups have quasi-isomorphic
cohomology rings as $A_\infty$ algebras if and only if their group
algebras are isomorphic, whereas they have isomorphic cohomology rings
as $E_\infty$ algebras if and only if the groups are isomorphic.
The modular isomorphism problem asks whether
a finite $p$-group can be recovered from its modular group
algebra. This was recently settled in the negative by
Garc{\'\i}a-Lucas, Margolis and del
R{\'\i}o~\cite{Garcia-Lucas/Margolis/delRio:2022a}, who found an
example of a pair of non-isomorphic groups of order $2^9$ with isomorphic $2$-modular
group algebras.\index{modular isomorphism problem}
\end{example}

\begin{remark}\label{rk:formal}\index{formal $A_\infty$ algebra}
In the inductive procedure described in the proof of
Theorem~\ref{th:Kadeishvili}, if it happens that 
\[ m_2(f_1\otimes f_2 - f_2 \otimes f_1)+f_2(1\otimes m_2-m_2\otimes 1) = 0 \]
then it follows that we may take $m_i\colon (H_*\fa)^{\otimes i} \to
H_*\fa$ and $f_i \colon (H_*\fa)^{\otimes i} \to \fa$ to be zero for all
$i\ge 3$. In this case, we deduce that $\fa$ is formal. We shall apply
this in Section~\ref{se:Q1} the case of finite groups with 
a generalised quaternion Sylow
$2$-subgroup and no normal subgroup of index two.
\end{remark}

\section{Derived categories}%
\index{derived!category}\label{se:D}

We shall be working a lot with derived categories and quotients of
them, so we briefly recall the definition and basic
properties. Singularity and cosingularity categories will be discussed
in Section~\ref{se:sing}. Recall from Notation~\ref{no:C*} that we use
homological gradings for all our complexes, so that the differential
decreases the degree by one. Background can be found in 
Chapter~6 of Christensen, Foxby and Holm~\cite{Christensen/Foxby/Holm:2024a}, 
Holm, J\o rgensen and Rouquier~\cite{Holm/Jorgensen/Rouquier:2010a},
Positselski~\cite{Positselski:2011a}, 
Verdier~\cite{Verdier:1996a},
Weibel~\cite{Weibel:1994a},
Yekutieli~\cite{Yekutieli:2020a} and many other places.

If $R$ is a ring, we denote by 
$\Ch(R)$\index{Ch@$\Ch(R)$, $\Ch^{\mathsf{b}}(R)$, $\Ch^+(R)$, $\Ch^-(R)$} 
the category whose objects are the chain complexes of
$R$-modules and whose morphisms are the 
\emph{chain maps},\index{chain map} namely the
degree preserving maps that
commute with the differential. This is an abelian category. The full
subcategory $\Ch^b(R)$ consists of \emph{bounded} complexes, namely
those where all but a finite number of the terms are zero. We also
have $\Ch^+(R)$, the complexes whose terms in large enough negative
degrees are zero, and $\Ch^-(R)$, the complexes whose terms in large
enough positive degrees are zero.

Two maps $f,g\colon X\to Y$ in $\Ch(R)$ are \emph{homotopic}\index{homotopic maps}
if there is a chain homotopy $s$, namely a map that increases degree by one,
and satisfying $f-g=ds+sd$. This is an equivalence relation, and
composition of chain maps is well defined on homotopy classes of maps.
Then $\sfK(R)$\index{K@$\sfK(R)$,  $\sfK^{\mathsf{b}}(R)$, $\sfK^+(R)$, $\sfK^-(R)$} 
is defined to be the category whose objects
are the chain complexes of $R$-modules, and whose arrows are the
homotopy classes of maps.  This is no longer an abelian category, but
it is a \emph{triangulated category}\index{triangulated category}. 
The axioms for a triangulated
category may be found in Definition~10.2.1 of~\cite{Weibel:1994a}.
Since homotopic maps of chain
complexes induce the same map on homology, we can regard homology as a
functor $H_*$ from $\sfK(R)$ to graded $R$-modules.We say that a map $f$ in $\sfK(R)$
is a \emph{quasi-isomorphism}\index{quasi-isomorphism} if $H_*(f)$ is
an isomorphism. The \emph{derived category} $\sfD(R)$ is
obtained from $\sfK(R)$ by adjoining inverses for
quasi-isomorphisms. So a morphism from $X$ to $Y$ in $\sfD(R)$ is a
diagram $X\xleftarrow{\sim} X' \to Y$ where $X\xleftarrow{\sim}X'$ is
a quasi-isomorphism and $X'\to Y$ is a morphism in $\sfK(R)$. For
details of this localisation process, see Section~10.3 of~\cite{Weibel:1994a}.
Inside $\sfK(R)$, we define full triangulated subcategories
$\sfK^{\mathsf{b}}(R)$, $\sfK^+(R)$, $\sfK^-(R)$
corresponding to
$\Ch^{\mathsf{b}}(R)$, $\Ch^+(R)$, $\Ch^-(R)$. We further define
$\sfK^{+,\mathsf{b}}(R)$, $\sfK^{-,\mathsf{b}}(R)$, to be the objects $X$ in
$\sfK^+(R)$, respectively $\sfK^-(R)$ for which $H_*(X)$ is bounded.

Replacing a bounded complex by a quasi-isomorphic bounded below complex
of projectives induces an equivalence $\sfK^{+,\mathsf{b}}(\Proj(R))\to
\Db(R)$. Similarly, replacing a bounded complex by a quasi-isomorphic
bounded above complex of injectives induces an equivalence
$\sfK^{-,\mathsf{b}}(\Inj(R))\to\Db(R)$. 
For unbounded complexes, it's a bit more complicated, as explained in
Chapter~5 of~\cite{Christensen/Foxby/Holm:2024a}. We say that a
complex $P$ is \emph{semi-projective}\index{semi-projective} if it has the
lifting property with respect to surjective quasi-isomorphisms. In
other words, given a map $P\to Y$ and a surjective quasi-isomorphism
$X\xrightarrow{\sim} Y$ as in the diagram
\[ \xymatrix{&P\ar[d]\ar@{.>}[dl]\\X\ar[r]^\sim& Y\ar[r]&0} \]
the dotted arrow exists making the diagram commute. Dually, a complex
$I$ is \emph{semi-injective}\index{semi-injective} if it has the
extension property with respect to injective quasi-isomorphisms.
\[ \xymatrix{0\ar[r]&Y\ar[r]^\sim\ar[d]&X\ar@{.>}[dl]\\&I} \]
Given any complex $X$, there exists a 
\emph{semi-projective resolution}\index{semi-projective!resolution} of
$X$, namely a semi-projective complex $P$ and a quasi-isomorphism
$P\to X$, and a \emph{semi-injective resolution}\index{semi-injective!resolution}
of $X$, namely a semi-injective complex $I$ and a quasi-isomorphism
$X\to I$. See Theorems~5.2.15 and~5.3.19
of~\cite{Christensen/Foxby/Holm:2024a}. 
Taking resolutions in this fashion induces equivalences 
$\sfK\,\mathsf{Semi\text{-}Proj}(R)\to\sfD(R)$ and
$\sfK\,\mathsf{Semi\text{-}Inj}(R)\to\sfD(R)$, see Proposition~6.4.3
of~\cite{Christensen/Foxby/Holm:2024a}. 

Next, consider a DG (differential graded)
algebra $A$. DG $A$-modules are already complexes, so
rather than considering 
complexes of DG modules, we let $\sfK(A)$ be the
homotopy category of DG $A$-modules. The case where
$A$ consists of an algebra in degree zero with zero differential, a DG
$A$-module is just a complex of $A$-modules, so this agrees with the
previous definition. As before, the derived category is obtained by
adjoining inverses for the quasi-isomorphisms.

Finally,
let $\fa$ be an $A_\infty$ algebra over $\kk$. The 
\emph{derived category}\index{derived!category} 
$\sfD(\fa)$ has as its objects the
$A_\infty$ modules over $\fa$ and as arrows the homotopy classes of
$A_\infty$ morphisms. In this category, by Theorem~\ref{th:Ainfinityqi}, quasi-isomorphisms
automatically have inverses, so the localisation step is
unnecessary. 
This generalises the case where $m_i=0$ for $i>2$, which is the case
of DG algebras and DG modules discussed above. In particular, a DG
module quasi-isomorphism may not have a DG module homotopy inverse,
but it always has an $A_\infty$ module homotopy inverse.

\section{Hochschild cohomology}%
\index{Hochschild cohomology!of an $A_\infty$ algebra}%
\index{HH@$\HHinf^*$}\label{se:HH}

For background on Hochschild cohomology of algebras, see Chapter~IX of
Cartan and Eilenberg~\cite{Cartan/Eilenberg:1956a}, Chapter~9 of
Weibel~\cite{Weibel:1994a}, and the book of
Witherspoon~\cite{Witherspoon:2019a}. In this section, we describe the
Hochschild cohomology of an augmented $A_\infty$ algebra $\fa$ and its relation to the
Hochschild cohomology of $H_*\fa$.

The Hochschild cohomology of $\fa$ as an augmented
$A_\infty$ algebra is
described as follows (see \S3 of Getzler and
Jones~\cite{Getzler/Jones:1990a}, and Definition~12.6 of 
Stasheff~\cite{Stasheff:1970a}). Let $\Sigma\bar\fa$ be the graded vector
space given by $(\Sigma\bar\fa)_i=\bar\fa_{i-1}$, where $\bar\fa$ is
the kernel of the augmentation $\fa\to\kk$, and write $[a_1|\dots|a_n]$
for the element corresponding to $a_1\otimes\dots\otimes a_n$ in
$(\Sigma\bar\fa)^{\otimes n}$. Then given a right $\fa$-module $M$ and a
left $\fa$-module $N$, 
the bar complex\index{bar!complex}\index{bar@$\bbar_{*,*}(M,\fa,N)$}
\[ \bbar_{*,*}(M,\fa,N)
=\bigoplus_{n\ge
  0}M\otimes(\Sigma\bar\fa)^{\otimes n}\otimes N \]
has a differential defined by 
\begin{equation}\label{eq:AinfinityHC*}
\begin{split}
\Hd(x\otimes& [a_1|\dots|a_n]\otimes y)=
\sum_{j=0}^n\pm m_{j+1}(x,a_1,\dots,a_j)\otimes 
[a_{j+1}|\dots|a_n]\otimes y \\
&{}+\sum_{0\le i+j\le n}\pm x\otimes [a_1|\dots|a_i|
m_j(a_{i+1},\dots,a_{i+j})|a_{i+j+1}|\dots|a_n] \otimes y\\
&\qquad+\sum_{j=0}^n\pm x\otimes[a_1|\dots|a_{n-j}]
\otimes m_{j+1}(a_{n-j+1},\dots,a_n,y). 
\end{split}
\end{equation}
In some sense the signs are given by the usual sign rules, but this
takes some interpretation.
It is explained in
Keller~\cite{Keller:2002a} how to compute these signs by looking at the
reduced tensor coalgebra of the suspension. The signs are also
explicitly described in Section~2 of Berglund and
B\"orjeson~\cite{Berglund/Borjeson:2020a}. Fortunately, we are mostly 
working in characteristic two, which allows us to ignore the signs.
There are two gradings on the bar complex. The first is the
homological grading coming from the number of bars, and the second is
the internal grading given by adding the degrees of the elements
involved. We write $\bbar_*(M,\fa,N)$ for the total complex where the
homological and internal gradings are added.

The \emph{bar resolution}\index{bar!resolution} of $\fa$ as an
$\fa$-$\fa$-bimodule is $\bbar_*(\fa,\fa,\fa)$. So if
$M$ is an $\fa$-$\fa$-bimodule, we have Hochschild cochains
$\Hom_{\fa,\fa}(\bbar_{*}(\fa,\fa,\fa),M)$. If $\fa$ is augmented, the augmentation
$\fa\to\kk$ gives an isomorphism
\[ \Hom_{\fa,\fa}(\bbar_{*}(\fa,\fa,\fa),M)
\cong  \Hom_\kk(\bbar_{*}(\kk,\fa,\kk),M) \] 

\begin{definition}\label{def:bar-construction}
We call $\bbar_{*,*}(\kk,\fa,\kk)$ or its total complex $\bbar_*(\kk,\fa,\kk)$ the \emph{bar
  construction}\index{bar!construction} on $\fa$, and abbreviate it to
$\bbar_{*,*}(\fa)$, respectively $\bbar_*(\fa)$.\index{bar@$\bbar_*(\fa)$}
\end{definition}

The
differential on $\Hom_\kk(\bbar_*(\fa),M)$ is induced from~\eqref{eq:AinfinityHC*}:
\begin{equation}\label{eq:AinfinityHH*}
\begin{split}
(\HHd f)[a_1|\dots|a_n]=\\
m_1f[a_1|&\dots|a_n]+\sum_{0\le i+j\le n}
\pm
f[a_1|\dots|a_i|m_j(a_{i+1},\dots,a_{i+j})|a_{i+j+1}|\dots|a_n], 
\end{split}
\end{equation}
see also Section~1 of Roitzheim and Whitehouse~\cite{Roitzheim/Whitehouse:2011a}.
The homology of this complex is the \emph{Hochschild cohomology}
$\HHinf^*(\fa,M)$. We write $\HHinf^*\fa$ for 
$\HHinf^*(\fa,\fa)$. We also have Hochschild homology $\HHinf_*\fa$,
which is the homology of the complex
\[ \bbar_*(\fa,\fa,\fa)\otimes_{\fa,\fa}\fa \cong
  \bbar_*(\fa)\otimes_\kk\fa. \]

If the only non-vanishing $m_i$ is $m_2$ then this differential on the
bar resolution agrees with the
classical notion of Hochschild cohomology of an algebra. If $\fa$ is a
DG algebra,\index{DG!algebra}\index{differential!graded algebra}
so the only non-vanishing $m_i$ are $m_1$ and $m_2$,
we can think of this as the usual Hochschild cohomology of a DG algebra. 

\begin{theorem}\label{th:HHHA}
If $\fa$ is an $A_\infty$ algebra, there 
is a conditionally convergent spectral sequence 
\begin{equation}\label{eq:HHHA} 
\HH^*H_*\fa\Rightarrow \HHinf^*\fa, 
\end{equation}
where the left hand side is the Hochschild cohomology of the
homology algebra $H_*\fa$, not taking into account any higher
structure, and the right hand side is the Hochschild cohomology of
$\fa$ as defined above. If some $E^r$ has the property that each
graded piece is finite dimensional then the spectral sequence is
strongly convergent.
\end{theorem}
\begin{proof}
The filtration of the complex $\Hom_{\fa,\fa}(\bbar_{*}(\fa,\fa,\fa),\fa)$ 
coming from the double complex $\bbar_{*,*}(\fa,\fa,\fa)$ is the
filtration by number of bars. This filtration 
gives rise to the spectral sequence. The
differentials are determined by the maps $m_i$.  
See Section~5 of~\cite{Benson/Greenlees:2023a} for details, and see
Boardman~\cite{Boardman:1999a} for a general discussion of conditional and
strong convergence.
\end{proof}

\begin{corollary}\label{co:AinfinityqiHH}
If $f\colon\fa\to\fa'$ is a quasi-isomorphism of $A_\infty$ algebras,
then there is an induced isomorphism
$\HHinf^*\fa\cong\HHinf^*\fa'$. This isomorphism preserves cup product
and Gerstenhaber bracket.
\end{corollary}
\begin{proof}
By Theorem~\ref{th:Ainfinityqi}, $f$ has a homotopy inverse
$g\colon\fa'\to\fa$. We then have maps of conditionally convergent spectral sequences
\[ \xymatrix{\HH^*H^*(\fa,\fa)\ar@<.5ex>[r]^{f_*}\ar@{=>}[d]&
\HH^*H^*(\fa,\fa')\ar@<.5ex>[l]^{g_*}\ar@<.5ex>[r]^{g^*}\ar@{=>}[d]&
\HH^*H^*(\fa',\fa')\ar@<.5ex>[l]^{f^*}\ar@{=>}[d]\\
\HHinf^*(\fa,\fa)\ar@<.5ex>[r]^{f_*}&
\HHinf^*(\fa,\fa') \ar@<.5ex>[l]^{g_*}\ar@<.5ex>[r]^{g^*}&
\HHinf^*(\fa',\fa') \ar@<.5ex>[l]^{f^*}} \]
The maps of $E^2$ pages are isomorphisms, and hence so are the maps of
targets (see for example Theorem~5.2
of~\cite{Benson/Greenlees:2023a}). The maps of targets are induced by
maps of Hochschild 
complexes, so they preserve all higher structure on Hochschild cohomology
such as cup product and Gerstenhaber bracket.
\end{proof}

\begin{remark}
If $\fa$ is an $A_\infty$ algebra, then $H_*\fa$ may either be
regarded as an ordinary associative algebra whose Hochschild
cohomology is written $\HH^*H_*\fa$, or as an $A_\infty$ algebra via
Theorem~\ref{th:Kadeishvili}, in which case we write
$\HHinf^*H_*\fa$. The spectral sequence~\eqref{eq:HHHA} relates these
two notions, and takes the form $\HH^*H_*\fa\Rightarrow\HHinf^*H_*\fa$.
\end{remark}

Kadeishvili~\cite{Kadeishvili:1980a} discusses the relationship 
between $A_\infty$ structure and Hochschild cohomology, and obtains 
the following.

\begin{proposition}\label{pr:HH}
Suppose that the action of $m_i$ on $H_*\fa$ is zero
for $i=1$ and $2<i<n$. Then $m_n\colon H_*\fa^{\otimes n}\to H_*\fa$ 
is a Hochschild $n$-cocycle\index{Hochschild cohomology} 
on $H_*\fa$, of internal (homological) degree $n-2$.

If $f\colon H_*\fa \to H_*\fa$ is a quasi-isomorphism to another
$A_\infty$ structure $m'_i$ on $H_*\fa$ satisfying the same assumptions,
with $f_1$ equal to the identity, then $m_n-m'_n$ is a Hochschild
coboundary, and all such occur this way.
Thus valid choices for $m_n$ form a well defined 
class in degree $(-n,n-2)$ in $\HH^*H_*\fa$.

If $\fa$ is graded and the $A_\infty$ structure preserves internal degree then
$m_n$ represents a Hochschild class of degree $(-n,n-2,0)$ on $H_*\fa$.
Thus if, regarding $\fa$ as a doubly graded algebra and ignoring the 
$m_n$ with $n>2$, the Hochschild cohomology ring 
$\HH^*H_*\fa$ has no non-zero elements of degree $(-n,n-2,0)$ with $n>2$, 
then the $A_\infty$ structure on $\fa$ is formal.
\end{proposition}
\begin{proof}
Equation \eqref{eq:Ainfty} implies that
\[ -m_2(\id \otimes m_n) + \sum_{r=0}^{n-1}(-1)^r m_n(\id^{\otimes r} \otimes m_2 \otimes
  \id^{\otimes(n-r-1)})+ (-1)^{n} m_2(m_n\otimes \id)=0 \] 
and so $m_n$ is a Hochschild cocycle on the cohomology ring. If
$f\colon H_*\fa \to H_*\fa$ is a quasi-isomomorphism with $f_1$ equal
to the identity then equation~\ref{eq:finfty} implies that
\begin{multline*} 
f_1m_n + \sum_{r=0}^{n-2}(-1)^rf_{n-1}(\id^{\otimes r} \otimes m_2 \otimes
  \id^{\otimes(n-r-2)}) \\= m'_2(f_1 \otimes f_{n-1}) +
  (-1)^{n}m'_2(f_{n-1}\otimes f_1) + m'_n(f_1\otimes\cdots\otimes  f_1). 
\end{multline*}
Now $f_1$ is the identity and $m'_2 = m_2$, so this becomes
\[ m'_n-m_n = -m_2(\id \otimes f_{n-1}) +
  \sum_{r=0}^{n-2}(-1)^rf_{n-1}(\id^{\otimes r} \otimes m_2 \otimes
  \id^{\otimes(n-r-2)})+(-1)^{n-1}m_2(f_{n-1}\otimes\id). \]
The right hand side is the formula for the Hochschild coboundary of
$f_{n-1}$, which may be taken to be any $(n-1)$-cochain.

The last part follows by an easy
inductive argument on $n$, beginning with $n=3$.
\end{proof}

\begin{definition}
We say that an $A_\infty$ algebra $\fa$ is 
\emph{intrinsically formal}\index{intrinsically formal} if given
another $A_\infty$ algebra $\fa'$ and an isomorphism of associative algebras
$H_*\fa \cong H_*\fa'$ there is a quasi-isomorphism $\fa\to\fa'$
inducing it. 
\end{definition}

\begin{remark}
Clearly an intrinsically formal $A_\infty$ algebra is formal.
If $\fa$ carries an internal grading then the isomorphism
$H_*\fa\cong H_*\fa'$ is required to preserve the induced internal
grading. So a graded $A_\infty$ algebra may be intrinsically formal
while the corresponding ungraded algebra is not. 

Proposition~\ref{pr:HH} implies that if there are no non-zero
classes of degree $(-n,n-2)$ in $\HH^*H_*\fa$ with $n>2$ then $\fa$ is
intrinsically formal. If $\fa$ carries an internal grading then we
only require that there are no non-zero classes of degree
$(-n,n-2,0)$, which is a weaker condition.
\end{remark}

\begin{theorem}\label{th:Massey}
Let $x_1,\dots,x_n$ ($n\ge 3$) be elements of the homology of a 
DG algebra $\fa$, and
suppose that the Massey product\index{Massey product} 
$\langle x_1,\dots,x_n\rangle$ is non-empty. 
Consider an $A_\infty$ structure on $H_*\fa$ given by Kadeishvili's
theorem, and suppose that $m_i=0$ for $2<i\le n-2$. Then 
\[ \ep m_n(x_1,\dots,x_n)\in \langle x_1,\dots,x_n\rangle, \]
where $\ep=(-1)^{\sum_{j=1}^{n-1}(n-j)|x_j|}$.
\end{theorem}
\begin{proof}
This is described in Theorem~3.1 of Lu, Palmieri, Wu and
Zhang~\cite{Lu/Palmieri/Wu/Zhang:2009a},
and corrected in Theorem~3.2 of Buijs, Moreno-Fern\'andez and
Murillo~\cite{Buijs/MorenoFernandez/Murillo:2020a}. 
\end{proof}

\section{The Gerstenhaber circle product}%
\index{Gerstenhaber!circle product}\label{se:circle}

Let $A$ be an associative
$\kk$-algebra and $M$ an $A$-$A$-bimodule. 
Gerstenhaber~\cite{Gerstenhaber:1963a}, introduced a
circle product on Hochschild cocycles of $A$, that is related to
$A_\infty$ structure, as we now explain. Recall from
Section~\ref{se:HH} that $\bbar_{*,*}(A)$ is the bar resolution of $A$. If $f\colon \bbar_{m,*}(A)\to M$ is an
$m$-cochain and $g\colon \bbar_{n,*}(A)\to A$ is an $n$-cochain, we define
$f\circ_i g\colon \bbar_{m+n,*}(A)\to M$ to be the $(m+n-1)$-cochain given on the basis by
\[ f\circ_i g\,[a_1|\dots|a_i|b_1|\dots|b_n|a_{i+1}|\dots|a_m]=
f[a_1|\dots|a_i|\,g[b_1|\dots|b_n]\,|a_{i+1}|\dots|a_m]. \]
We then define the 
\emph{circle product}\index{circle product}\footnote{Beware 
that the degrees in Section~5 of~\cite{Gerstenhaber:1963a}
differ by one from the Hochschild degrees, so that the circle product
is additive on degrees. This is explained in
Section~7, on the lower part of page 279 of~\cite{Gerstenhaber:1963a}.}
\[ f\circ g = \sum_{i=0}^m(-1)^{(|g|-1)i}f\circ_i g. \]

\begin{example}
The statement $m_2\circ m_2=0$ is equivalent to the associativity of
multiplication, because
\[ (m_2 \circ m_2) [a_1|a_2|a_3] = m_2(m_2(a_1,a_2),a_3) -
  m_2(a_1,m_2(a_2,a_3)) = (a_1a_2)a_3-a_1(a_2a_3). \]
\end{example}

\begin{theorem}
The circle product on $\HH^*(A,A)$ satisfies
\begin{equation*}
  (h\circ f)\circ g - (-1)^{(|f|-1)(|g|-1)}(h\circ g)\circ f
  = h \circ (f\circ g - (-1)^{(|f|-1)(|g|-1)}g\circ f),
\end{equation*}
and is related to the differential, the 
cup product,\index{cup product} and 
Gerstenhaber bracket\index{Gerstenhaber!bracket}
on cochains by the formulas
\begin{align} 
\HHd f&=(-1)^{|f|-1}m_2\circ f - f\circ m_2\label{eq:circ1} \\
f \cup g &= (m_2\circ_0 f) \circ_{|f|-1} g\label{eq:circ2} \\
f\circ \HHd g - \HHd(f\circ g) &{}+ (-1)^{|g|-1}\HHd f\circ g 
= (-1)^{|g|-1}(g\cup f - (-1)^{|f||g|}f\cup g)\label{eq:circ3} \\
[f,g] &=f\circ g - (-1)^{(|f|-1)(|g|-1)}g\circ f.\label{eq:circ4}
\end{align}
where $m_2$ is the multiplication in $A$.
\end{theorem}
\begin{proof}
This is proved in~\cite{Gerstenhaber:1963a}.
\end{proof}

In terms of the circle product, equation~\eqref{eq:Ainfty} can be
rewritten as 
\begin{equation}\label{eq:Ainfty-circ} 
\sum_{i+j=n+1}(-1)^im_i\circ m_j = 0. 
\end{equation}
So suppose, for example, that $\fa$ is an $A_\infty$ algebra with
$m_1=0$, and $m_i=0$ for $2<i<n$. 
We saw in Proposition~\ref{pr:HH}
that $m_n$ is a Hochschild $n$-cocycle on $A=H_*\fa$. 
We can see this easily using this formulation, since the condition
reduces to $m_2 \circ m_n +(-1)^n m_n \circ m_2=0$, or equivalently
using~\eqref{eq:circ1}, $\HHd m_n = 0$.
It also follows from this formulation that under these circumstances, for
$n<i<2n-2$ the condition is again that $m_i$ should be a Hochschild
$n$-cocycle, and if these are all coboundaries then they can be
rechosen to be zero. Then
the condition for $m_{2n-2}$ is
\[ m_2\circ m_{2n-2} + (-1)^nm_n\circ m_n + (-1)^{2n-2}m_{2n-2}\circ
  m_2 =0, \]
which can be rewritten using~\eqref{eq:circ1} as
\[ \HHd m_{2n-2} = (-1)^nm_n \circ m_n. \]
Continuing this way, we obtain the following, which will help
understand what is going on in Section~\ref{se:AinftyHBD}.

\begin{proposition}\label{pr:circle}
Let $n,t\ge 2$, and 
let $\fa$ be an $A_\infty$ algebra, such that
that for $i<t$ we have $m_i=0$ unless $i$ is congruent to $2$ modulo
$n-2$. Then
\begin{enumerate}[label={\rm(\arabic*)}]
\item If $t$ is not congruent to $2$ modulo $n-2$ then $m_t$ is a
  Hochschild cocycle.
\item If $t=s(n-2)+2$ then $m_t$ satisfies the coboundary condition
\begin{equation*} 
\HHd m_{s(n-2)+2} =  \sum_{i=1}^{s-1} (-1)^{in}m_{i(n-2)+2} \circ m_{(s-i)(n-2)+2}. 
\end{equation*}
\item
Suppose that $\HH^iH^*\fa=0$ for $i$ not congruent to $2$ modulo
$n-2$. Then an $A_\infty$ structure on $H^*\fa$ quasi-isomorphic to
that on $\fa$ may be chosen with $m_i=0$ unless $i$ is congruent to
$2$ modulo $n-2$.\qed
\end{enumerate}
\end{proposition}

\begin{remark}
There is one more piece of structure on Hochschild cohomology that is
related to the circle product. Over a field $k$ of characteristic $p$, 
there is a linear $p$-power operation defined on odd degree elements, 
which together with 
the Gerstenhaber bracket makes the shift of Hochschild cohomology
$\Sigma\HH^*(A)$ into a graded
restricted Lie algebra.\index{restricted Lie algebra}%
\index{Lie algebra, restricted}
So if $x\in\HH^{2n+1}(A)$ then $x^{[p]}\in\HH^{2np+1}(A)$. If $f\colon\bbar_{*,*}(A)\to
A$ is a cocycle representing $x$ then the circle product of $p$
copies of $f$, $((f\circ f)\circ\cdots)\circ f$, represents $x^{[p]}$.
For details see Tourtchine~\cite{Tourtchine:2006a}, 
Briggs and Rubio y Degrassi~\cite{Briggs/Rubio-y-Degrassi:2022a}. The latter
paper makes use of this to show that in characteristic $p$, the $p$-restricted Lie algebra
structure on Hochschild cohomology is invariant under 
stable equivalences of Morita type between self-injective
algebras.\index{stable!equivalence of Morita type}
\end{remark}

\section{\texorpdfstring{Bousfield--Kan $p$-completion}
{Bousfield--Kan p-completion}}\index{Bousfield--Kan $p$-completion}%
\index{completion}\index{p-completion@$p$-completion}

We shall use the $p$-completion of Bousfield and
Kan~\cite{Bousfield/Kan:1972a}, namely the completion with respect to
the field $\bF_p$ of $p$ elements.\footnote{Warning: when \cite{Bousfield/Kan:1972a}
writes $Z_p$, this refers to the field of $p$ elements, not the
$p$-adic integers.} 
We write $X\phat$ for 
the $p$-completion of a space $X$. This comes with a natural map $X\to
X\phat$, and has the properties listed in
Theorem~\ref{th:Bousfield-Kan} below. We begin with a definition.

\begin{definition}\label{def:nilpotent}
A space $X$ is \emph{nilpotent} if $\pi_1(X)$ is a nilpotent group and 
the action of $\pi_1(X)$ on $\pi_i(X)$, $i\ge 2$, is nilpotent. That 
is, there is a series 
\[ \pi_i(X)=H_1\ge H_2 \ge \dots \ge H_{n_i}=0 \]
of $\pi_1(X)$-invariant subgroups with the property that for $1\le 
j\le n_i-1$, the induced action of $\pi_1(X)$ on $H_j/H_{j+1}$ is trivial. 
\end{definition}

\begin{theorem}\label{th:Bousfield-Kan}
The Bousfield--Kan $p$-completion has the following properties.
\begin{enumerate}[label={\rm(\roman*)}]
\item A map of spaces $X\to Y$ induces a mod $p$ reduced homology
equivalence $\tilde H_*X \to \tilde H_*Y$ (Notation~\ref{no:C*}) if and only if it induces a weak
homotopy equivalence between the completions $X\phat \to Y\phat$. \label{it:Bousfield-Kan/whe}
\item A space $X$ is said to be $\bF_p$-\emph{good}, or
$p$-\emph{good},\index{p-good@$p$-good}  if the map $\tilde H_*X \to
\tilde H_*X\phat$ is an isomorphism, otherwise $X$ is
$\bF_p$-\emph{bad}, or $p$-\emph{bad}.\index{p-bad@$p$-bad} 
$X$ is said to be $\bF_p$-\emph{complete}, 
or $p$-\emph{complete},\index{p-complete@$p$-complete} if $X
\to X\phat$ is a weak homotopy equivalence. The following are equivalent:
{\rm (a)} $X$ is $p$-good, {\rm (b)} $X\phat$ is $p$-complete, {\rm (c)}
$X\phat$ is $p$-good. Thus if $X$ is $p$-bad, then however
many times we complete it, it remains $p$-bad. \label{it:Bousfield-Kan/good}
\item If $X$ is connected and 
$\pi_1X$ is finite then $X$ is $p$-good for all primes $p$.
In this case, we have $\pi_1X\phat\cong \pi_1X/O^p(\pi_1X)$
(Notation~\ref{no:gt}). \label{it:Bousfield-Kan/pi1} 
\item If $X$ is connected and $\pi_iX$ is finite for all $i\ge 1$ then $X\phat$ is a
nilpotent space whose homotopy groups $\pi_i(X\phat)$ are all finite $p$-groups.\label{it:Bousfield-Kan/nilpotent}
\end{enumerate}
\end{theorem}
\begin{proof}
Parts \ref{it:Bousfield-Kan/whe} and \ref{it:Bousfield-Kan/good} are proved in Section~I.I.5
of~\cite{Bousfield/Kan:1972a}, while part \ref{it:Bousfield-Kan/pi1} is proved in 
Proposition~I.VII.5.1
of~\cite{Bousfield/Kan:1972a}. Part~\ref{it:Bousfield-Kan/nilpotent}
is proved in Proposition~I.VII.4.3 of~\cite{Bousfield/Kan:1972a}.
\end{proof}

\begin{corollary}\label{co:BGphat}
If $G$ is a finite group then
the classifying space\index{classifying space} 
$BG$\index{BG@$BG\phat$} is
$p$-good, its completion $BG\phat$ is a $p$-complete, 
nilpotent space\index{nilpotent space}
and $\pi_1(BG\phat)\cong G/O^p(G)$.\index{pi1@$\pi_1(BG\phat)$}
The map $BG\to BG\phat$ induces a mod $p$ homology equivalence
$H_*BG\cong H_*BG\phat$ and a mod $p$ cohomology equivelance
$H^*BG\phat\cong H^*BG$.
The space $BG$ is already
$p$-complete if and only if $G$ is a finite $p$-group.

The following are equivalent.
\begin{enumerate}[label={\rm(\arabic*)}]
\item $BG\phat$ is simply connected.\index{simply connected}
\item $\Omega BG\phat$ is connected.
\item $G$ has no normal subgroup of index $p$.
\end{enumerate}
\end{corollary}
\begin{proof}
The statements that $BG$ is $p$-good and $\pi_1BG\phat\cong G/O^p(G)$ follow from
Theorem~\ref{th:Bousfield-Kan}\,\ref{it:Bousfield-Kan/pi1}, because $\pi_1(BG)\cong
G$. Then by
Theorem~\ref{th:Bousfield-Kan}\,\ref{it:Bousfield-Kan/good}, $H_*BG\to
H_*BG\phat$ is an isomorphism and
$BG\phat$ is $p$-complete. By the Universal Coefficient Theorem,
$H^*BG\phat\to H^*BG$ is also an isomorphism.
By
Theorem~\ref{th:Bousfield-Kan}\,\ref{it:Bousfield-Kan/nilpotent} $BG\phat$ is
nilpotent. 
The equivalence of
(1) and (2) follows from the general isomorphism
$\pi_i(X)\cong\pi_{i-1}(\Omega X)$. For the equivalence of (1) and
(3), since $\pi_1BG\phat\cong G/O^p(G)$, $BG\phat$ is simply connected
if and only if $G=O^p(G)$, which happens if and only if $G$ has no normal
subgroup of index $p$.
\end{proof}

\begin{remark}\label{rk:EMSS}
In general, for a nilpotent space\index{nilpotent space} 
$X$ with finite dimensional homology groups, the
Eilenberg--Moore spectral sequence%
\index{Eilenberg--Moore spectral sequence} is a spectral sequence of
Hopf algebras\index{Hopf algebra} converging to the homology
of the loop space:\index{loop space}
\begin{equation}\label{eq:Cotor}\index{Cotor}
\Ext^{*,*}_{H^*X}(\kk,\kk) \cong \Cotor_{*,*}^{H_*X}(\kk,\kk)
\Rightarrow H_*\Omega X 
\end{equation}
For the simply connected case, see Eilenberg--Moore~\cite{Eilenberg/Moore:1966a}, 
Smith~\cite{Smith:1967a,Smith:1969a,Smith:1970b}. For $p$-complete
spaces, the Eilenberg--Moore spectral sequences is discussed in 
Dwyer~\cite{Dwyer:1975b}.

The isomorphism between Ext and Cotor in~\eqref{eq:Cotor} follows from
the following lemma (see also Lemma~A1.1.6 of Ravenel~\cite{Ravenel:1986a}).

\begin{lemma}\label{le:Cotor}
Let $C$ be a graded coalgebra of finite
type over a field $\kk$. Then the graded dual $C^*$ is a graded
algebra of finite type. The doubly graded cobar complex
$\cobar_{*,*}(C,C,C)$ is dual to the doubly graded bar
complex $\bbar_{*,*}(C^*,C^*,C^*)$, and  
\[ \Cotor_{*,*}^C(\kk,\kk)\cong \Ext^{*,*}_{C^*}(\kk,\kk). \]
\end{lemma}
\begin{proof}
For graded vector spaces of finite type (i.e., each graded piece has
finite dimension), duality is exact and the double dual is naturally
isomorphic to the original.
Since $C$ has finite type, the doubly graded cobar complex $\cobar_{*,*}(C,C,C)$
is a free bicomodule resolution of $C$ of finite type, so
$Q=\kk\square_C\cobar_{*,*}(C,C,C)$ is a free right comodule resolution of $\kk$
of finite type, where $\square_C$ is the cotensor product over $C$. The dual is
$P=\Hom(Q,\kk)=\bbar_{*,*}(C^*,C^*,C^*)\otimes_{C^*}\kk$, a left module
resolution of $\kk$ over $C^*$ of finite type, and the double dual is again $Q$.

So we have
\begin{align*}
\Cotor_{*,*}^C(M,\kk) &\cong H_{*,*}(Q\,\square_C\, \kk) \\
&\cong H_{*,*}(\Hom_\kk(P,\kk)\,\square_C\, \kk) \\
&\cong H_{*,*}(\Hom_C(P, \kk)) \\
&\cong \Ext^{*,*}_{C^*}(\kk,\kk).
\qedhere
\end{align*}
\end{proof}

If the finite group $G$ is not 
$p$-nilpotent,\index{p-nilpotent@$p$-nilpotent finite group} then $BG$ is not a nilpotent
space, because $\pi_1BG\cong G$. However, by
Corollary~\ref{co:BGphat}, $BG\phat$ is a nilpotent space, and $BG\to
BG\phat$ induces an isomorphism
$H^*BG\phat\cong H^*BG$. So from~\eqref{eq:Cotor} we get a spectral sequence 
\[ \Ext^{*,*}_{H^*BG}(\kk,\kk) \Rightarrow H_*\Omega BG\phat. \]
This expresses the cochain level statement for the DG algebra of (derived)
endomorphisms\index{End@$\iEnd_{C^*BG}(\kk)$}
\[ \iEnd_{C^*BG}(\kk) \simeq C_*\Omega BG\phat. \]
The proof of this is similar to the proof of Lemma~\ref{le:Cotor},
dualising the cobar complex.\index{cobar complex}
It follows that we have a functor 
$\iHom_{C^*BG}(\kk,-)$\index{Hom@$\iHom_{C^*BG}(\kk,-)$}
from the derived category of $C^*BG$-modules to the derived category of $C_*\Omega BG\phat$-modules.

In the other direction, for any path connected space $X$ 
the Rothenberg--Steenrod construction~\cite{Rothenberg/Steenrod:1965a}%
\index{Rothenberg--Steenrod construction} gives 
\[ \iEnd_{C_*\Omega X}(\kk) \simeq C^*X. \]
We can apply this either to $X=BG$ to obtain $\iEnd_{\kk G}(\kk)
\simeq C^*BG$ or to $X=BG\phat$ to obtain 
\[ \iEnd_{C_*\Omega BG\phat}(\kk)\simeq C^*BG\phat \simeq C^*BG. \]
In the latter case, we get a functor 
$\iHom_{C_*\Omega BG\phat}(\kk,-)$\index{Hom@$\iHom_{C_*\Omega BG\phat}(\kk,-)$}
from $C_*\Omega BG\phat$-modules to $C^*BG$-modules.
\end{remark}

\begin{lemma}[The fibre lemma]\index{fibre lemma}\label{le:fibre}
Suppose that $F \to E \to B$ is a fibration sequence, and that
the action of $\pi_1(B)$ on $H_iF$ is nilpotent for all $i\ge 0$. That
is, there exists a series
\[ H_iF = H_1\ge H_2\ge \cdots \ge H_{n_i}=0 \]
of $\pi_1(B)$-invariant subgroups with the property that for $1\le
j\le n_i-1$, the induced action of $\pi_1(B)$ on $H_j/H_{j+1}$ is
trivial (cf.\ Definition~\ref{def:nilpotent}). 
Then $E\phat \to B\phat$ is a fibration, with fibre homotopy
equivalent to $F\phat$.
\end{lemma}
\begin{proof}
Recalling that our convention is that $H_*F$ denotes mod $p$
homology, this is the case $R=\bF_p$ of the Mod-$R$ Fibre Lemma~II.5.1 of 
Bousfield and Kan~\cite{Bousfield/Kan:1972a}.
\end{proof}

\begin{proposition}\label{pr:fibration-phat}
Let $G$ be a finite group, and embed $G$ in a finite unitary group
$G \to U(n)$. Then there are fibration sequences
\begin{enumerate}[label={\rm(\roman*)}]
\item $U(n)/G \to BG \to BU(n)$, \label{it:fibration-phat/uncompleted}
and 
\item $(U(n)/G)\phat \to BG\phat \to BU(n)\phat$. \label{it:fibration-phat/completed}
\end{enumerate}
\end{proposition}
\begin{proof}
\ref{it:fibration-phat/uncompleted}
Let $EU(n)$ be the complex Stiefel variety (or any contractible space
on which $U(n)$ acts freely). Then we can use $EU(n)$ for $EG$,
and the required fibration is
\[ U(n)/G \to EU(n)/G \to EU(n)/U(n). \]

\ref{it:fibration-phat/completed} Since $\pi_1(BU(n))$ is trivial, this follows from
Lemma~\ref{le:fibre} and part \ref{it:fibration-phat/uncompleted}.
\end{proof}

\section{Classifying spaces and fusion systems}\label{se:fusion}

Let $G$ be a finite group and $\kk$ be a field of characteristic
$p$. Let $EG$ be a contractible space with a free $G$-action, and let
$BG$ be the quotient $EG/G$. Since $C_*EG$ is a free resolution of
$\kk$ as a $\kG$-module,
the cohomology ring $H^*BG$ is isomorphic to the group cohomology
$H^*(G,\kk)=\Ext^*_{\kG}(\kk,\kk)$.\index{Ext@$\Ext^*_{\kG}(\kk,\kk)$}
Furthermore, if $P_*$ is any projective resolution\index{projective!resolution} of $\kk$
as a $\kG$-module, then the DG algebra\index{DG!algebra}%
\index{differential!graded algebra}
$\Hom_{\kG}(P_*,P_*)$ is quasi-isomorphic to $C^*BG$.

By a theorem of Cartan and Eilenberg (Theorem~XII.10.1
of~\cite{Cartan/Eilenberg:1956a}, the cohomology $H^*BG$ only
depends on the Sylow $p$-subgroup $D$ of $G$ and the fusion system 
on it defined by conjugation in $G$. This is defined as
follows.\index{Sylow subgroup}

\begin{definition}
Let $D$ be a Sylow $p$-subgroup of a finite group $G$.
For subgroups $H$ and $K$ of $D$, we define $\Hom_G(H,K)$ to
be the set of group homomorphisms from $H$ to $K$ that are induced by
conjugation by some element of $G$, $\{\phi\colon H\to K\mid \exists
g\in G\ \forall h\in H\ \phi(h)=ghg^{-1}\}$.
The \emph{fusion category}\index{fusion!category} 
of $G$ over $D$ is the category $\cF_D(G)$ whose
objects are the subgroups of $D$, and whose morphisms are given by
$\Hom_G(H,K)$. The \emph{fusion system}\index{fusion!system} 
of $G$ over $D$ consists of
$D$ together with the fusion category.
\end{definition}

Abstract fusions systems are studied in the books of Puig~\cite{Puig:2009a},
Aschbacher, Kessar and
Oliver~\cite{Aschbacher/Kessar/Oliver:2011a} and
Craven~\cite{Craven:2011b}. There is a set of axioms, devised by
Puig (his terminology is 
\emph{Frobenius $P$-categories}\index{Frobenius!$P$-category}). 
These can be found in~\cite[Chapters~2 and~21]{Puig:2009a},
\cite[Definition~2.2]{Aschbacher/Kessar/Oliver:2011a},
\cite[Definition~4.1]{Craven:2011b}.
Not every fusion system comes from a finite group in the above
way. We shall assume that the saturation axiom\index{saturation axiom} 
is part of the definition.

\begin{definition}
Let $D$ be a Sylow $p$-subgroup of a finite group $G$.
A subgroup $H$ of $D$ is \emph{$p$-centric}\index{centric subgroup} 
in $G$ if $Z(H)$ is a
Sylow $p$-subgroup of $C_G(H)$. This is equivalent to saying that
$C_G(H)=Z(H)\times O_{p'}C_G(H)$ (Notation~\ref{no:gt}).
The \emph{centric linking system}%
\index{centric linking system}\index{linking system}
of $G$ over $D$ is the category
$\cL_D(G)$
whose objects are the subgroups of $D$ that are $p$-centric in $G$,
 and whose morphisms are the quotient of $\{g\in G\mid gHg^{-1}\le
 K\}$ by the action of $O_{p'}C_G(H)$. There is an obvious functor
 $\cL_D(G) \to \cF_D(G)$.

Again, there is a set of axioms for an abstract centric linking system $\cL$ over a
fusion system on a $p$-group $D$. 
A \emph{$p$-local finite group}\index{p-local@$p$-local finite group}
consists of a finite $p$-group $D$ together with a fusion system $\cF$
over $D$ and a centric linking system $\cL\to \cF$.
\end{definition}

Given a $p$-local finite group $(D,\cF,\cL)$, its 
\emph{classifying space}\index{classifying space} $|\cL|$ 
is defined to be the nerve
of the category $\cL$. It is a $p$-good space (Proposition~1.12
of Broto, Levi and Oliver~\cite{Broto/Levi/Oliver:2003a}). 
The $p$-local finite group $(D,\cF,\cL)$
can be recovered from the homotopy type of $|\cL|\phat$
(Theorem~7.4 of~\cite{Broto/Levi/Oliver:2003a}). Thus the information
contained in the space $|\cL|\phat$ is equivalent to that contained in the
$p$-local finite group $(D,\cF,\cL)$.
Each can be recovered from the other. Note, however, that
$(D,\cF,\cL)$ cannot be recovered from the mod $p$ cohomology
$H^*|\cL|\simeq H^*|\cL|\phat$ as an $A_\infty$ algebra, but rather
the full $E_\infty$ structure is needed, for the reasons explained in
Example~\ref{eg:Einfinity}.

\begin{theorem}\label{th:BLO1}
Let $(D,\cF,\cL)$ be a $p$-local finite group coming from a finite
group $G$. 
Then the natural map $|\cL_D(G)| \to BG$ is a mod $p$ cohomology
equivalence, and so induces a homotopy equivalence $|\cL_D(G)|\phat
\to BG\phat$.
\end{theorem}
\begin{proof}
This is Proposition~1.1 of Broto, Levi and Oliver~\cite{Broto/Levi/Oliver:2003b}.
\end{proof}

\begin{theorem}\label{th:BLO2}
For a $p$-local finite group $(D,\cF,\cL)$, the cohomology
$H^*|\cL|$ is isomorphic to the ring of 
stable elements\index{stable!elements} in $H^*BD$ with respect to the
fusion system $\cF$, in
the sense of Cartan and Eilenberg, Theorem~XII.10.1
of~\cite{Cartan/Eilenberg:1956a}. 
\end{theorem}
\begin{proof}
This is Theorem~B of~\cite{Broto/Levi/Oliver:2003a}.
\end{proof}

\begin{theorem}\label{th:Oliver}
Suppose that $G$ and $H$ are finite groups, and there is a fusion
preserving isomorphism from the Sylow $p$-subgroup $D$ of $G$ to that of
$H$. Then there is a homotopy equivalence $BG\phat \to BH\phat$.
\end{theorem}
\begin{proof}
For $p=2$ this is Theorem~B of Oliver~\cite{Oliver:2006a}, while for
odd primes it is Theorem~B of Oliver~\cite{Oliver:2004a}. 
\end{proof}

The following stronger theorem was proved later, and Oliver's
Theorem~\ref{th:Oliver} is a consequence.

\begin{theorem}\label{th:Chermak}
Given a fusion system $\cF$ on a finite $p$-group $D$, there exists a
unique centric linking system $\cL$ over $\cF$.
\end{theorem}
\begin{proof}
This was first proved by Chermak~\cite{Chermak:2013a} using his theory
of localities. The proof used the classification of finite simple
groups. Later, Oliver recast the proof in terms of obstruction theory.
His proof depends on the Meierfrankenfeld--Stellmacher classification of
quadratic best offenders~\cite{Meierfrankenfeld/Stellmacher:2012a}, 
which again relies on the classification of finite simple
groups. Finally, a classification free proof was given by Glauberman
and Lynd~\cite{Glauberman/Lynd:2016a}.
\end{proof}

\begin{remark}
Using Theorem~\ref{th:Chermak}, 
given a fusion system $\cF$ on a $p$-group $D$, it determines first 
a linking system $\cL$ over $\cF$, then the classifying space $|\cL|$,
then its $p$-completion $|\cL|\phat$, and finally the cochains
$C^*|\cL|\phat\simeq C^*|\cL|$. 
\end{remark}

\begin{corollary}\label{co:Chermak}
If a fusion system $\cF$ on a finite $p$-group $D$ occurs for a finite
group $G$ with Sylow $p$-subgroup $D$, then for any
linking system $\cL$ over $\cF$, we have $|\cL|\phat\simeq BG\phat$.\qed
\end{corollary}

\begin{remark}
The corresponding theorem for discrete $p$-toral groups is proved in
Levi and Libman~\cite{Levi/Libman:2015a}. This is relevant when trying
to understand $p$-completed classifying spaces of compact Lie
groups.\index{compact Lie group}\index{Lie group}
\end{remark}

\begin{remark}\label{rk:Linck}
Let $B$ be a block\index{block} of the group algebra $\kG$ of a finite group
$G$ over $\kk$. It follows from the work of Alperin and
Brou\'e~\cite{Alperin/Broue:1979a} that one can associate to $B$ a
fusion system $\cF_B$ describing the fusion of
subpairs\index{subpairs} associated to
the block.\index{fusion!system!of a block} 
This is spelled out in Section~3 of
Linckelmann~\cite{Linckelmann:2004a}. Thus using
Theorem~\ref{th:Chermak}, we can associate to $B$ a 
linking system\index{linking system!of a block}
$\cL_B$, and a classifying space\index{classifying space!of a block} 
$|\cL_B|$, whose cohomology
$H^*|\cL_B|$ is the cohomology of $B$ in the sense of 
Linckelmann~\cite{Linckelmann:1999b}. In the case of a principal
block, the defect groups are the Sylow $p$-subgroups of $G$, and the  
fusion system of the block is the same as the fusion system $\cF_G(D)$
of the group.
\end{remark}

Craven and Glesser~\cite{Craven/Glesser:2012a} studied fusion systems
on metacyclic groups.\index{metacyclic groups} 
Theorem~1.1 of that paper shows that for
dihedral, semidihedral and generalised quaternion $2$-groups, all
possible fusions systems are realised as fusion systems $\cF_D(G)$ of
some finite group $G$. It follows that when we discuss $2$-completed classifying
spaces of finite groups with these as Sylow $2$-subgroups, we are
really discussing the $2$-completed classifying spaces associated to
any fusion system on such a $2$-group, including classifying spaces of
blocks with these as defect groups.

\section{Hochschild cohomology of cochains}%
\index{Hochschild cohomology}

In this section, we compare Hochschild cohomology of $H^*BG$
with Hochschild cohomology of $H_*\Omega BG\phat$, both
being regarded as $A_\infty$ algebras. We show that these are 
isomorphic, see Theorem~\ref{th:HHOmega}, as a special case of a
general statement about pointed path connected 
spaces, see Theorem~\ref{th:coHH}. Recall that throughout, we are
working over a field $\kk$ of characteristic $p$, though
Theorem~\ref{th:coHH} is true over more general coefficients.

Let $A$ be an augmented DG algebra. Recall from Section~\ref{se:HH}
that we have the bar complex $\bbar_{*,*}(M,A,N)$ and the bar
construction $\bbar_{*,*}(A)=\bbar_{*,*}(\kk,A,\kk)$ 
whose total complex $\bbar_*(A)$ is used to
compute Hochschild homology and cohomology. Dually, if $C$ is a coaugmented
coalgebra, $M$ is a right comodule and $N$ is a left comodule, we have
the \emph{cobar complex}  
\[ \cobar_{*,*}(M,C,N) = \bigoplus_{n\ge 0}M\otimes
  (\Sigma^{-1}\overline C)^{\otimes n} \otimes N \]
with the appropriate differential dual to that of the bar complex, 
and the \emph{cobar construction}
$\cobar_{*,*}(C)=\cobar_{*,*}(\kk,C,\kk)$ that computes coHochschild homology
and cohomology. Here, $\overline C$ is the cokernel of the coaugmentation
$\kk\to C$.
The cobar construction was first introduced by Adams~\cite{Adams:1956a}
as a model for the loop space. CoHochschild theory is discussed in
Hess, Parent and Scott~\cite{Hess/Parent/Scott:2009a},
Hess~\cite{Hess:2016a}, Hess and Shipley~\cite{Hess/Shipley:2021a}, 
Rivera and Zeinalian~\cite{Rivera/Zeinalian:2018a}.
We shall need to be careful about two different versions of the total
complex of cobar, see Remark~\ref{rk:cobar}.

We first explain the simply connected case, which covers the case of 
$X=BG\phat$ with $G=O^p(G)$ (Notation~\ref{no:gt}).

\begin{theorem}
Let $X$ be a simply connected space. Then we have an isomorphism
\[ \HH^*C^*X\cong\HH^*C_*\Omega X \]
where $C^*X$ is the DG algebra\index{DG!algebra}
of singular cochains on $X$ and $C_*\Omega X$ is the DG algebra of singular
chains on the Moore loop space\index{Moore loop space} of $X$.
\end{theorem}
\begin{proof}
By Theorem~1 of F\'elix, Menichi and
Thomas~\cite{Felix/Menichi/Thomas:2005a}, we have
an isomorphism
\[ \HH^*C^*X\cong\HH^*\cobar_* C_*X. \]
Now $\cobar_* C_*X$  is quasi-isomorphic
to $C_*\Omega X$ (see the main theorem of~\cite{Adams:1956a}), so by
Proposition~3.3 
of~\cite{Felix/Menichi/Thomas:2005a}, this therefore induces an
isomorphism in Hochschild cohomology.
\end{proof}

For a general pointed path connected space, this is no longer true,
and the correct formulation is in terms of the 
coHochshild cohomology\index{coHochschild cohomology} of
the coalgebra of chains, see Theorem~\ref{th:coHH}. In the presence of
appropriate finiteness conditions, we can rewrite that as desired, see Theorem~\ref{th:HHloops}.

\begin{remark}\label{rk:cobar}
First we must discuss two versions of the cobar construction.
As a double complex, the cobar complex on a coaugmented DG
coalgebra,\index{DG!coalgebra}  $\cobar_{*,*}(C,C,C)$, 
and $\cobar_{*,*}(C)=\cobar_{*,*}(\kk,C,\kk)$ are
unambiguous, and are dual to the bar complex on an augmented DG algebra
$\bbar_{*,*}(A,A,A)$, and $\bbar_{*,*}(A)$. But for the total 
complex, we can either take direct sums or direct products, and they
give us different answers unless the DG coalgebra $C$ is $1$-reduced,
meaning that it is zero in negative degrees, $\kk$ in degree zero, and
zero in degree one. We shall write $\cobar_*(C,C,C)$ and $\cobar_*(C)$
for the direct sum versions, $\cobar^{\prod}_*(C,C,C)$ and 
$\cobar^{\prod}_*(C)$ for the direct product  
versions.\index{cobar complex!direct product version} 

We shall assume that all our DG coalgebras are coaugmented, and  
conilpotent, meaning that the coproduct on any given element of the  
coaugmentation coideal vanishes after enough iterations. 
Under these conditions, 
$\cobar_*$ is left adjoint to $\bbar_*$.
Conilpotence 
ensures that the counit map $C \to \bbar_*\cobar_*C$ is well defined.  
Beware that $\cobar_*$ 
does not send quasi-isomorphic DG coalgebras to quasi-isomorphic DG  
algebras, but see Lemma~\ref{le:cobar-qi} below.  
 
On the other hand, $\cobar^{\prod}_*$ is
right adjoint to $\bbar_*$, and is the
version involved in exotic convergence\index{exotic convergence} 
of the Eilenberg--Moore spectral sequence in
the non-simply connected case, see Dwyer~\cite{Dwyer:1975b}. 
\end{remark}

The following is a generalisation of Adams' theorem.

\begin{lemma}\label{le:RZ}
Let $X$ be a pointed path connected space. Then the reduced cobar 
complex on the singular chains on $X$
is quasi-isomorphic to the singular chains on $\Omega X$, 
\[ \cobar_* C_*X\simeq C_*\Omega X. \]
\end{lemma}
\begin{proof}
This is proved in Proposition~9.2 of~\cite{Rivera/Zeinalian:2018a}. 
\end{proof}

We saw in Section~\ref{se:HH} that
Hochschild homology of a DG algebra $A$ is defined as
$\HH_*A=H_*(\bbar_*A\otimes_\kk A,\partial)$, while Hochschild cohomology of $A$
is defined as $\HH^*A=H_*(\Hom_\kk(\bbar_*A,A),\delta)$. Dually, we use
the direct sum version of the cobar construction to define
coHochschild homology and coHochschild cohomology of a DG coalgebra
$C$. These are defined as 
\begin{align*}
\text{co}\HH_*C&=H_*(C\otimes_{C,C}\cobar_*(C,C,C))=H_*(C\otimes_\kk\cobar_*C,\partial)\\
\text{co}\HH^*C&=H_*\Hom_{C,C}(C,\cobar_*(C,C,C))=H_*(\Hom_\kk(C,\cobar_*C),\delta),
\end{align*}
see~\cite{Hess:2016a,Hess/Parent/Scott:2009a,Hess/Shipley:2021a}. 

\begin{remark}\label{rk:HomCA}
When
$A=\cobar_*C$, there is a smaller resolution of $A$ as an
$A$-$A$-bimodule than the bar resolution, namely we can put a
differential on $A\otimes_\kk C \otimes_\kk A$ as follows. 
If $c\in C$ then $\Sigma^{-1}c$ is naturally an element of $A$ of bar
length one (or is zero if $c$ is in the image of the
coaugmentation).\footnote{The role of $\Sigma^{-1}\colon C \to A$ here
  is that it is the universal twisting cochain out of $C$.\index{twisting cochain}}
Then the differential is
\begin{multline*} 
\partial=d_A\otimes 1\otimes 1 + 1\otimes d_C\otimes 1+ 1 \otimes 1\otimes d_A \\
-\left(\,(1\otimes 1 \otimes \mu)(1\otimes 1\otimes \Sigma^{-1}\otimes 1)
+(\mu\otimes 1\otimes 1)(1\otimes\Sigma^{-1}\otimes 1\otimes 1)\,\right)
(1\otimes \Delta\otimes 1) 
\end{multline*}
where $\mu$ is multiplication in $A$ and $\Delta$ is comultiplication
in $C$.

Then $\HH^*A=H_*(\Hom_{A,A}(A\otimes_\kk C\otimes_\kk A, A)) =
H_*(\Hom_\kk(C,A))$ with the induced differential. The latter is then
by definition the coHochschild cohomology $\text{co}\HH^*(C)$.
This is analogous
to Section~2.2 of~\cite{Hess:2016a}, where Hochschild homology is
computed as $H_*(C\otimes_\kk A)$ with the induced differential.
\end{remark}

\begin{theorem}\label{th:coHH}
Let $X$ be a pointed path connected space, and let $C_*X$ be the
singular chain complex on $X$. Then we have
\[ \HH^*C_*\Omega X\cong \text{\rm co}\HH^* C_*X. \]
Here, we are regarding $C_*\Omega X$ as an augmented DG algebra and
$C_*X$ as an augmented DG coalgebra.
\end{theorem}
\begin{proof}
Using Lemma~\ref{le:RZ}, and Remark~\ref{rk:HomCA},
we have
\begin{align*}
\HH^*C_*\Omega X&\cong\HH^*\cobar_*C_*X\\
&\cong H_*\Hom_\kk(C_*X,\cobar_*C_*X) \\
&\cong\text{co}\HH^*C_*X.
\qedhere
\end{align*}
\end{proof}

\begin{lemma}\label{le:minimal-Kan}
Suppose that $X$ is a connected Kan complex\index{Kan complex} 
and $\pi_i(X)$ is finite
for all $i>0$. Then the minimal Kan complex\index{minimal!Kan complex} 
equivalent to $X$ is $0$-reduced, and has
finitely many non-degenerate simplices in each degree.
\end{lemma}
\begin{proof}
See for example Section~I.10 of Goerss and
Jardine~\cite{Goerss/Jardine:1999a}. The non-degenerate simplices in
degree $i$ are
in one to one correspondence with the elements of $\pi_i(X)$.
\end{proof}

If $X$ is a simplicial set, we write 
$\Simp_*X$\index{Simp@$\Simp_*X$, $\Simp^*X$} for the simplicial
chains on $X$. This is the chain complex which in degree $n$ is the
$\kk$-vector space with the non-degenerate $n$-simplices of $X$ as
basis, and with differential coming from the alternating sum of the
face maps. We write $\Simp^*X$ for dual complex.

\begin{lemma}\label{le:cobar-qi}
If $f\colon X\to Y$ is a homotopy equivalence of $0$-reduced Kan complexes then
$f$ induces a quasi-isomorphism of cobar complexes on the simplicial
chains\index{simplicial chains}
\[ \cobar_*\Simp_*X\simeq\cobar_*\Simp_*Y . \] 
\end{lemma}
\begin{proof}
This follows from Propositions~8.1 and~8.2 of~\cite{Rivera/Zeinalian:2018a}.
\end{proof}

\begin{lemma}\label{le:coHH}
Let $X$ be a pointed connected space with $\pi_i(X)$ a finite
$p$-group for all $i>0$. 
Then $\text{\rm co}\HH^*C_*X\cong \HH^*C^*X$.
\end{lemma}
\begin{proof}
Let $Y$ be a minimal Kan complex homotopy equivalent to $X$. Then using
Lemma~\ref{le:minimal-Kan}, $Y$ is $0$-reduced, and has finitely many
non-degenerate simplices in each degree. So $C_*Y$ has finite type,
and $C^*Y$ is its dual. By Lemma~\ref{le:cobar-qi}
and Lemma~\ref{le:Cotor}, we have equivalences of doubly graded complexes
\[ \cobar_{*,*}C_*X\simeq\cobar_{*,*}\Simp_*Y\simeq\Hom_\kk(\bbar_{*,*}\,\Simp^*Y,\kk). \]
However, when we pass to total complexes,
$\Hom_\kk(\bbar_*\,\Simp^*Y,\kk)$ gives the direct product version
$\cobar^{\prod}_*\Simp_*Y$
rather than the direct sum version $\cobar_*\Simp_*Y$ of the cobar
construction, see Remark~\ref{rk:cobar}. 
Our hypotheses imply that $Y$ is $p$-complete
(Theorem~\ref{th:Bousfield-Kan}), and it is shown in
Dwyer~\cite{Dwyer:1975b} 
that if $Y$ is $p$-complete then the 
Eilenberg--Moore spectral sequence with $E^2$ page
$\Cotor_{*,*}^{H_*Y}(\kk,\kk)$ converges to $H_*\Omega Y$. Thus we have
\[ \cobar^{\prod}_*\Simp_*Y\simeq \Simp_*\Omega Y, \] 
which by Lemma~\ref{le:RZ} is quasi-isomorphic to
$\cobar_*\Simp_*Y$. 
Thus we have
\begin{align*}
\text{co}\HH^*C_*X&\cong \text{co}\HH^*\Simp_*Y\\
&= H_*\Hom_\kk(\Simp_*Y,\cobar_*\Simp_*Y)\\
&\cong H_*\Hom_\kk(\Simp_*Y,\cobar^{\prod}_*\Simp_*Y)\\
&\cong H_*\Hom_\kk(\bbar_*\Simp^*Y,\Simp^*Y)\\
&=\HH^*\Simp^*Y\\
&\cong \HH^*C^*X.
\end{align*}
Here, we have used both $p$-completeness of $X$ and
the finiteness hypothesis.
\end{proof}

\begin{theorem}\label{th:HHloops}
  If $X$ is a pointed connected space whose homotopy groups are finite
$p$-groups then $\HH^*C^*X\cong\HH^*C_*\Omega X$.
\end{theorem}
\begin{proof}
This follows by combining Theorem~\ref{th:coHH} with
Lemma~\ref{le:coHH}:
\begin{equation*}
\HH^*C^*X\simeq\text{co}\HH^*C_*X\simeq\HH^*C_*\Omega X.
\qedhere
\end{equation*}
\end{proof}

Recall from Section~\ref{se:HH} that for a DG algebra $A$, $\HH^*A$ is
isomorphic to $\HHinf^*H_*A$.

\begin{corollary}\label{co:HHloops}
If $X$ is a pointed connected space whose homotopy groups are finite $p$-groups,
then regarding $H^*X$ and $H_*\Omega X$ as $A_\infty$ algebras we have
$\HHinf^* H^*X\cong\HHinf^* H_*\Omega X$.
\end{corollary}
\begin{proof}
This follows by combining Theorem~\ref{th:HHloops} with
Corollary~\ref{co:AinfinityqiHH}. 
\end{proof}

\begin{theorem}\label{th:HHOmega}
Let $G$ be a finite group. Then
there is an isomorphism 
\[ \HHinf^*H^*BG\cong\HHinf^*H_*\Omega BG\phat. \]
\end{theorem}
\begin{proof}
By Theorem~\ref{th:Bousfield-Kan}\,(iv), the homotopy groups
$\pi_i(BG\phat)$ are finite $p$-groups. The theorem now follows from 
Corollary~\ref{co:HHloops}.
\end{proof}

\begin{remark}
Among other applications, in Section~\ref{se:HH*kD} we use
Theorem~\ref{th:HHOmega} to recompute 
the ring structure on Hochschild cohomology of the group algebras of
dihedral groups, first computed by Siegel and
Witherspoon~\cite{Siegel/Witherspoon:1999a} using more elementary
methods. 
\end{remark}

\section{Abelian Sylow subgroups}\index{abelian Sylow subgroups}%
\index{Sylow subgroup!abelian}

Let $G$ be a finite group with abelian Sylow $p$-subgroup $D$, and let
$\kk$ be a field of characteristic $p$. Then by a classical theorem of
Burnside, the normaliser $N_G(D)$ controls $G$-fusion in $D$. See for example
Theorem~7.1.1 of Gorenstein~\cite{Gorenstein:1968a}. This implies that
the inclusion $N_G(D) \hookrightarrow G$ induces an isomorphism of
fusion systems $\cF_{N_G(D)}(D) \to \cF_G(D)$. It follows that the
ring of stable elements\index{stable!elements} in $H^*BD$ (see Theorem~\ref{th:BLO2}) is just
the invariants of the normaliser. So we have the classical theorem
of Swan.\index{Swan's theorem}

\begin{theorem}\label{th:Swan}
Suppose that $G$ has an abelian Sylow $p$-subgroup $D$, and let $\kk$
be a field of characteristic $p$. Then the inclusion $N_G(D)\to G$ and
the quotient map $N_G(D) \to N_G(D)/O_{p'}N_G(D)$ (Notation~\ref{no:gt})
induce isomorphisms of rings
\[ H^*BG \cong H^*BN_G(D)\cong H^*B(N_G(D)/O_{p'}N_G(D)) \cong H^*BD^{N_G(D)/C_G(D)}. \] 
Here, the right hand side denotes the invariants of the action of
$N_G(D)/C_G(D)$ on $H^*BD$.
\end{theorem}
\begin{proof}
See Swan~\cite{Swan:1960c}.
\end{proof}

The consequence for the $p$-completed classifying space is the
following.

\begin{theorem}\label{th:abelian}
Suppose that $G$ has an abelian Sylow $p$-subgroup $D$, and let $\kk$
be a field of characteristic $p$. Then we have homotopy equivalences
\[ B(N_G(D)/O_{p'}N_G(D))\phat \leftarrow BN_G(D)\phat \rightarrow BG\phat. \]
\end{theorem}
\begin{proof}
By Theorem~\ref{th:Swan}, the inclusion of a Sylow
$p$-normaliser $N_G(D)\to G$ and the quotient map $N_G(D)\to
N_G(D)/O_{p'}N_G(D)$ induce mod $p$ cohomology equivalences
\[ B(N_G(D)/O_{p'}N_G(D)) \leftarrow BN_G(D) \rightarrow BG. \]
Hence by Theorem~\ref{th:Bousfield-Kan}\,(i) and Corollary~\ref{co:BGphat}, after $p$-completion
these give homotopy equivalences.
\end{proof}

\begin{remark}\label{rk:DH}
By the Schur--Zassenhaus Theorem\index{Schur--Zassenhaus theorem} 
(see for example Theorem~6.2.1 of Gorenstein~\cite{Gorenstein:1968a}), 
the group $N_G(D)/O_{p'}N_G(D)$ is isomorphic to a semidirect product
of $D$ by a $p'$-subgroup $H$ of $\Aut(D)$. This reduces the study of $BG\phat$
to the study of $B(D\rtimes H)\phat$.
\end{remark}

We shall discuss the case where $D$ is cyclic in
Section~\ref{se:cyclic} and the case where $D$ is an elementary
abelian $2$-group in Theorem~\ref{th:elemab2}.

\section{Singularity and cosingularity categories}\label{se:sing}

We begin this section with motivational background. 

\begin{definition}\label{def:thick}
Recall that a \emph{thick subcategory}\index{thick subcategory}  
of a triangulated category is a
full triangulated subcategory that is closed under taking direct
summands.
\end{definition}

\begin{definition}
If $R$ is a ring, we define the \emph{singularity
  category}\index{singularity category} to be
the Verdier quotient\index{Verdier quotient} 
 $\Db(R)/\Thick(R)$ of the bounded derived
category by the thick subcategory generated by the regular
representation of $R$. Objects in $\Thick(R)$ are \emph{perfect
  complexes}.\index{perfect complex}

If $R\to \kk$ is an augmented algebra, then regarding $\kk$ as an
$R$-module via the augmentation, we define the 
\emph{cosingularity category}\index{cosingularity category} to be the
Verdier quotient $\Db(R)/\Thick(\kk)$.
\end{definition}

\begin{example}
If $R$ is a commutative Noetherian ring, then the singularity category
$\Dsg(R)$ is generated by the modules $R/\fp$, with $\fp$ a prime
ideal in the singular locus\index{singular locus} of $R$. This follows
from the main theorem of Schoutens~\cite{Schoutens:2003a}. The ring
 $R$ is regular\index{regular ring} if and only if its
singularity category is zero. This can be seen as motivation for the
term ``singularity category''.
\end{example}

\begin{example}
If $R$ is a commutative 
Gorenstein ring\index{Gorenstein ring} 
then $\Dsg(R)$ is equivalent to the category $\MCM(R)$ of maximal
Cohen--Macaulay\index{maximal Cohen--Macaulay modules} $R$-modules
(Buchweitz~\cite{Buchweitz:1986a}). 
\end{example}

\begin{example}
For commutative non-Gorenstein rings, the singularity category is
often quite bad. As an example, consider the ring
$R=\kk[x,y]/(x^2,xy,y^2)$. The short exact sequence of $R$-modules
\[ 0 \to \kk\oplus \kk \to R \to \kk \to 0 \]
shows that in the category $\Dsg(R)$, the module $\kk$ decomposes as a
direct sum of two copies of a shift of $\kk$. It follows that $\kk$
decomposes into arbitrarily many summands, none of which is
indecomposable. This example is investigated in detail in
Chen~\cite{Chen:2011a}. 
\end{example}

\begin{example}
If $R$ is a finite dimensional graded
Koszul algebra\index{Koszul duality} 
with left Noetherian Koszul dual $R^!$ then there is a derived
equivalence $\Db(R)\simeq\Db(R^!)$ that swaps the roles of $R$ and
$\kk$. See for example Theorem~2.12.6
of Beilinson, Ginzburg and Soergel~\cite{Beilinson/Ginzburg/Soergel:1996a}. It follows that
$\Dsg(R)\simeq\Dcsg(R^!)$ and 
$\Dcsg(R)\simeq\Dsg(R^!)$.  
This equivalence can be seen as motivation for the  
term ``cosingularity category''.  
 For more general Koszul algebras, one must
change the definition of bounded, to accommodate the grading, as in the
discussion of $A_\infty$ algebras below.
\end{example}

Now let $\fa$ be an $A_\infty$ algebra. Recall from Section~\ref{se:D}
that the derived category $\sfD(\fa)$ has as its objects the
$A_\infty$ modules over $\fa$ and as morphisms the homotopy classes of
 morphisms of $A_\infty$ modules.
If $H_*\fa$ is commutative Noetherian, we define the 
\emph{bounded derived category}\index{bounded derived category} 
$\Db(\fa)$ to be the thick subcategory of 
$\sfD(\fa)$ whose objects are the modules whose homology is finitely
generated as a module over the ring $H_*\fa$.

We also need a suitable notion when $H_*\fa$ is
not commutative or Noetherian, in order to deal with the case of
$\fa=H_*\Omega BG\phat$. The appropriate condition there involves a
suitable notion of Noether normalisation
(Definition~3.7 of Greenlees and Stevenson~\cite{Greenlees/Stevenson:2020a}):

\begin{definition}\label{def:sing-cosing}
We say that $\fb \to \fa$
is a \emph{normalisation}\index{normalisation} of $\fa\to \kk$
if both $\fa$ and $\kk$ are in the thick subcategory
$\Thick(\fb)\subseteq\sfD(\fb)$ generated by $\fb$.
For example, if $H_*\fa$ is finitely presented then
the set of generators in a finite presentation leads to a normalisation 
(Theorem~3.13 of~\cite{Greenlees/Stevenson:2020a}). In the cases where
a normalisation exists, we make the following definitions.

If $\fb \to \fa$ is a normalisation, we define the 
\emph{bounded derived category} $\Db(\fa)$ to be 
full subcategory of $\sfD(\fa)$ consisting of those objects whose
restriction to $\fb$ are in $\Thick(\fb)\subseteq\sfD(\fb)$.
Under suitable hypotheses this is independent of the normalisation 
(Propositions~4.3 and~7.2 of~\cite{Greenlees/Stevenson:2020a}).

We define the 
\emph{singularity category}\index{singularity category}
$\Dsg(\fa)$ to be the Verdier quotient
of $\Db(\fa)$ by the thick
subcategory $\Thick(\fa)$ generated by $\fa$. We define the
\emph{cosingularity category}\index{cosingularity category}  $\Dcsg(\fa)$
to be the Verdier quotient of $\Db(\fa)$
by the thick subcatgory $\Thick(\kk)$ generated by the field $\kk$.
\end{definition}

We are interested in the cases of $H^*BG\phat$ and $H_*\Omega
BG\phat$. Recall from
Proposition~\ref{pr:fibration-phat}\ref{it:fibration-phat/completed} that
we have a fibration sequence
\[ (U(n)/G)\phat \to BG\phat \to BU(n)\phat. \]
This gives rise to maps of $A_\infty$ algebras
\[ H^*BU(n)\phat \to H^*BG\phat \to H^*(U(n)/G)\phat \]
and 
\[ H_*\Omega(U(n)/G)\phat \to H_*\Omega BG\phat \to H_*U(n)\phat. \]
These have the property that $H^*BU(n)\phat \to H^*BG\phat$ and
$H_*\Omega(U(n)/G)\phat \to H_*\Omega BG\phat$ are normalisations,
see Example~(10.6) of~\cite{Greenlees/Stevenson:2020a}, and so we may
define singularity and cosingularity categories for $H^*BG$ and
$H_*\Omega BG\phat$ using Definition~\ref{def:sing-cosing}. We
emphasise that these are singularity and cosingularity categories as
$A_\infty$-algebras, not just as algebras.

The following theorem expresses a version of Koszul
duality\index{Koszul duality} between $H^*BG$ and $H_*\Omega BG\phat$
(cf.\ Remark~\ref{rk:EMSS}).

\begin{theorem}\label{th:Dsg-Dcsg}
For a finite group $G$, the functor $\iHom_{H^*BG}(\kk,-)$ 
induces a triangulated equivalence of bounded
derived categories $\Db(H^*BG) \xrightarrow{\sim} 
\Db(H_*\Omega BG\phat)$ that sends $H^*BG$ to $\kk$ 
and sends $\kk$ to $H_*\Omega BG\phat$. 
It induces triangulated equivalences
\[ \Dsg(H^*BG) \xrightarrow{\sim} \Dcsg(H_*\Omega BG\phat), \qquad
\Dcsg(H^*BG) \xrightarrow{\sim} \Dsg(H_*\Omega BG\phat). \]
\end{theorem}
\begin{proof}
This follows from Theorem~9.1 and
Example~(10.6) of~\cite{Greenlees/Stevenson:2020a}.
\end{proof}

\section{Tame blocks}\index{tame!blocks}\label{se:tame}

The trichotomy theorem\index{trichotomy theorem} 
of Drozd~\cite{Drozd:1977a,Drozd:1979a,Drozd:1980a} (see also
Theorem~B on page 478 and
Corollary~C on page 480 of Crawley-Boevey~\cite{Crawley-Boevey:1988a})
for finite dimensional algebras over an algebraically closed field states 
that every finite dimensional algebra is of finite, tame or wild
representation type,\index{representation!type} 
and these types are mutually exclusive. Roughly
speaking, finite representation type\index{finite representation type} 
means that there are only
finitely many isomorphism classes of finitely generated indecomposable
modules. Tame representation type\index{tame!representation type} 
means that the finitely generated 
indecomposables of any particular dimension (over an infinite field) 
come in one parameter families with finitely many exceptions, and 
wild representation type\index{wild representation type} 
means that the module theory for a free algebra\index{free algebra} on
two generators can be encoded in the category of finite dimensional
modules for the given algebra. For details, see for example
Section~4.4 of~\cite{Benson:1998b}.

In the case of blocks of finite groups, the representation type only
depends on the defect group.\index{defect group} 

\begin{theorem}\label{th:tame}
Let $\kk$ be a field of characteristic $p$ and $G$ be a finite group.
Let $B$ be a block of $\kG$ with defect group $D$. Then
\begin{enumerate}[label={\rm(\roman*)}]
\item $B$ has finite representation type if and only if $D$
  is cyclic.
\item $B$ has tame representation type if and only if $p=2$,
  and $D$ is dihedral,\index{dihedral group} 
semidihedral,\index{semidihedral group} or 
generalised quaternion.\index{generalised quaternion group}
\item In all other cases, $B$ has wild representation type.
\end{enumerate}
\end{theorem}
\begin{proof}
It follows from Theorem~11 of Green~\cite{Green:1959a}, generalising a
theorem of Higman~\cite{Higman:1954b}, that the
representation type of $B$ only depends on the defect group $D$, which
is a finite $p$-subgroup of $G$ well defined up to conjugation in
$G$. For finite $p$-groups, 
the representation type was determined in Theorem~1 of
Bondarenko and Drozd~\cite{Bondarenko/Drozd:1982a}.  Note that
Ringel~\cite{Ringel:1974a} described the representation type of all
local algebras over an algebraically closed field, with the exception
of the generalised quaternion and semidihedral group algebras. The
representation type of these is discussed in
Crawley-Boevey~\cite{Crawley-Boevey:1989a}. There, the indecomposable
representations of the semidihedral groups are described in detail in
Theorem~1.9. 

One can deduce from this classification that the generalised quaternion
group algebras have tame representation type. This is because a
generalised quaternion group may be embedded in a semidihedral group
of twice the order. By Green's Indecomposability
Theorem~\cite[Theorem~8]{Green:1959a},%
\index{Green's indecomposability theorem} a finite dimensional
indecomposable module for 
the generalised quaternion group induces to an indecomposable for the
semidihedral group, and the restriction of that back to the
generalised quaternion group contains the original module as a direct
summand. So at most two isomorphism classes induce to a given
isomorphism class for the semidihedral group. On the other hand,
nobody has yet succeeded in giving a classification of the indecomposable
modules for the generalised quaternion groups.
\end{proof}

Blocks with cyclic defect group\index{cyclic!defect group}%
\index{defect group!cyclic} were completely described in the
work of Brauer~\cite{Brauer:1941a} for the prime cyclic case. Using
methods originating in a paper of Thompson~\cite{Thompson:1967a}, 
Dade~\cite{Dade:1966a} completed the general case of cyclic defect
group. 

The case of tame representation type was the subject of extensive work
of Erdmann~\cite{Erdmann:1977a,Erdmann:1979a,
Erdmann:1987a,Erdmann:1988a,Erdmann:1988b,
Erdmann:1988c,Erdmann:1988d,Erdmann:1990a,
Erdmann:1990b,Erdmann:1990c,Erdmann:1992a,
Erdmann:1994a,Erdmann/Michler:1977a,
Erdmann/Skowronski:2019a}, giving an almost complete description of
the Morita types of these blocks. 
Our work leans heavily on these papers. To make life easy, 
it follows from a case by case analysis, described in
Theorem~\ref{th:principal-block} below, that for each
isomorphism type of defect group of tame representation type and each
fusion system on it, there 
is a principal block of some finite group $G$ 
with the same fusion system, and all such have
equivalent classifying spaces by Oliver's
Theorem~\ref{th:Oliver}. Judicious choice of $G$ minimises the work
involved in understanding $H^*BG$ and $H_*\Omega BG\twohat$
as $A_\infty$ algebras. So the $2$-completed classifying spaces of blocks (as in
Remark~\ref{rk:Linck}) of tame representation type are all equivalent to
$2$-completed classifying spaces of finite groups of tame
representation type.

\begin{theorem}\label{th:principal-block}
For each isomorphism type of $2$-group $D$ of tame representation type,
and each fusion system $\cF$ on it, there is a 
finite group $G$ with the same fusion system. For any linking system
$\cL$ over $\cF$, and any such $G$, we have $B\cL\phat\simeq BG\phat$.
In particular, for any $2$-block of a finite group with defect group
$D$, there is a finite group with homotopy equivalent $2$-completed 
classifying space.
\end{theorem}
\begin{proof}
This is a case by case analysis. For blocks with dihedral defect
groups, Theorem~1.1 of Craven and Glesser~\cite{Craven/Glesser:2012a}
shows that  there are three possible fusion systems on a dihedral
group, and they have one, two and three conjugacy classes of
involutions. These all occur in
Sylow $2$-subgroups of finite groups. By Theorem~\ref{th:Chermak}, there is a
unique linking system $\cL$ for each possible fusion system, and by
Corollary~\ref{co:Chermak}, we then have $|\cL|\twohat\simeq BG\twohat$.
The reader should consult Section~\ref{se:dihedral-Sylow} for more
details of the finite groups with these fusion systems.

For blocks with semidihedral defect groups, Theorem~1.1
of~\cite{Craven/Glesser:2012a} shows that there are four possible
fusion systems, again all realised by finite groups. These groups are
described in Section~\ref{se:semidihedral-Sylow}. For blocks with
generalised quaternion defect groups, Theorem~1.1
of~\cite{Craven/Glesser:2012a} shows that there are three possible fusion
systems, again all realised by finite groups. These groups are described in
Section~\ref{se:generalised-quaternion-Sylow}. 
\end{proof}

\begin{remark}
Finite dimensional local symmetric algebras of tame representation
type are listed in Theorem~III.1 of Erdmann~\cite{Erdmann:1990a}. The
group algebras of finite $2$-groups among these are as given 
as follows. 
The dihedral groups have type III.1\,(c) with $k$ a power of two (or
type III.1\,(b) for the Klein four group), 
the semidihedral groups have type III.1\,(d$'$), and 
the generalised quaternion groups have type III.1\,(e$'$).
This is a slightly more precise statement than given
in~III.13 of~\cite{Erdmann:1990a}; see also
Sections~\ref{se:semidihedral-groups} 
and~\ref{se:Qu} for further comments.
\end{remark}

\section{Cohomology of complete intersections}\label{se:ci}

We shall need to compute
$\Ext^*_R(\kk,\kk)$\index{Ext@$\Ext^*_R(\kk,\kk)$} 
and $\HH^*(R)$ in the case where
$R=H^*BG$ is a complete intersection.\index{complete intersection}
For this reason, we give a brief review of
cohomology of complete intersections, following
Avramov~\cite{Avramov:1998a}, Sj\"odin~\cite{Sjodin:1976a} and
Buchweitz and Roberts~\cite{Buchweitz/Roberts:2015a}.

Let $R$ be a complete intersection of the form $Q/I$, 
where $Q=\kk[x_1,\dots,x_n]$ is a
positively graded polynomial ring and $I$ is generated by a 
homogeneous regular sequence $f_1,\dots,f_c$ in $\fm^2$, where $\fm$
is the ideal $(x_1,\dots,x_n)$.\footnote{By negating degrees, results
here apply equally well to negatively graded rings.}
We can take partial derivatives in the usual way to give polynomials
\[ b_{i,k}=\displaystyle \frac{\partial f_k}{\partial x_i}\in Q. \] 
We can then
take their images in $R$, which we denote $\bar b_{i,k}$.  But then there's a
problem when it comes to second partial derivatives,\index{partial derivatives} 
because in
characteristic two the second partial derivative of $x^2$ with respect
to $x$ vanishes. To remedy this, the 
\emph{divided partial derivatives},\index{divided partial derivative}
see Section~2.1 of~\cite{Buchweitz/Roberts:2015a}, are
defined to be the corresponding terms in the Taylor expansion of the
polynomial, so that for example $\displaystyle \frac{\partial^{(2)}(x^2)}{\partial x^2}
 = 1$. Thus we have
$\displaystyle \frac{\partial^2f_k}{\partial
  x_i^2}=2\frac{\partial^{(2)}f_k}{\partial x_i^2}$. 
More generally, if $\bfx=x_1,\dots,x_n$ then
\[ f(\bfx+d\bfx)=\sum_{\bfa=(a_1,\dots,a_n)}
\frac{\partial^{(|\bfa|)}f(\bfx)}{\partial\bfx^\bfa}d\bfx^\bfa \]
where $|a|=\sum_ia_i$ and the sum is over $n$-tuples of non-negative
integers $(a_1,\dots,a_n)$. The usual partial derivative is
$a_1!\cdots a_n!$ times the divided partial derivative.
So now set 
{\renewcommand{\arraystretch}{1.3}\[ 
a_{i,j,k}=\begin{cases}
\displaystyle
\frac{\partial^2 f_k}{\partial x_i\partial x_j} & i\ne j\vspace{2mm} \\ 
\displaystyle
\frac{\partial^{(2)}f_k}{\partial x_i^2} & i=j,
\end{cases} \]}%
as an element of $Q$, and write $\bar a_{i,j,k}$ for the image of
$a_{i,j,k}$ in $R$. These are the coefficients of the 
\emph{Hessian quadratic form} $\qq$ associated to the relations
defining $R$, see Section~2 of~\cite{Buchweitz/Roberts:2015a}.

\begin{definition}\label{def:Cliffq}
Following~\cite{Buchweitz/Roberts:2015a}, we 
define the \emph{Clifford algebra}\index{Clifford algebra} 
$\Cliffq$\index{Cliff@$\Cliffq$} to be the differential bigraded 
algebra\index{differential!bigraded algebra} over $R$ with generators 
$\xx_i$ dual to the $x_i$, in degree $(-1,-|x_i|)$ ($1\le
i \le n$) and $s_k$ dual to the $f_k$, in degree $(-2,-|f_k|)$ ($1\le j\le c$). The
multiplicative structure is given by making $s_j$ central, and
\[ \xx_i\xx_j+\xx_j\xx_i=\sum_{k=1}^c \bar a_{i,j,k}s_k\qquad (i\ne j),\qquad \qquad
\xx_i^2  = \sum_{k=1}^c\bar a_{i,i,k}s_k. \]
The differential $d\colon \Cliffq \to \Cliffq$ vanishes on $A$ and
on the $s_k$, and on the $\xx_i$ it is given by
\[ d(\xx_i)=\sum_{k=1}^c \bar b_{i,k} s_k. \]
\end{definition}

\begin{theorem}\label{th:ExtRkk}
We have $\Ext^{*}_R(\kk,\kk)=\kk\otimes_R \Cliffq$.
\end{theorem}
\begin{proof}
This is proved in Theorem~5 of Sj\"odin~\cite{Sjodin:1976a} when the characteristic
is not two. The general statement can be obtained from Avramov~\cite{Avramov:1998a}
by combining Theorem~10.2.1\,(5) and Example~10.2.2 there. 
\end{proof}

\begin{remark}
According to Theorem~\ref{th:ExtRkk}, $\Ext^*_R(\kk,\kk)$ is generated over $\kk$ by
elements $\xx_i$ in degree $(-1,-|x_i|)$ ($1\le i\le n$) and $s_k$ in degree
$(-2,-|f_k|)$ ($1\le k\le c$). The elements $s_k$ generate a central polynomial subring
over which $\Ext^*_R(\kk,\kk)$ is a free module of rank $2^n$. The
relations express $\xx_i\xx_j+\xx_j\xx_i$ and $\xx_i^2$ as linear combinations
of the $s_k$ with coefficients in $\kk$ given by the constant terms $a_{i,j,k}(0)$ of
the polynomials
$a_{i,j,k}$.  These only depend on the quadratic parts
of the polynomials $f_k$, so we have
$f_k=\sum_{i,j=1}^na_{i,j,k}(0)x_ix_j+{}$terms of degree at least
three.
\end{remark}

\begin{remark}
The algebra $\Ext^*_R(\kk,\kk)$ carries a Hopf algebra\index{Hopf algebra} structure for
which the elements $\xx_i$ and $s_k$ are primitive. This gives the
multiplication on the graded dual Hopf algebra
$\Tor_*^R(\kk,\kk)$. See Theorem~5 of~\cite{Sjodin:1976a}, Section~10 of~\cite{Avramov:1998a}.
\end{remark}

\begin{theorem}\label{th:HHR}
We have $\HH^{*}(R)= H^{*}(\Cliffq,d)$, the cohomology of $\Cliffq$
with respect to the differential 
$d$.\index{Hochschild cohomology!of complete intersection}
\end{theorem}
\begin{proof}
This is Theorem~2.11 of Buchweitz and Roberts~\cite{Buchweitz/Roberts:2015a}.
\end{proof}

\begin{remark}
The images 
of the $s_k$ in $\Ext^*_R(\kk,\kk)$ and in $\HH^*(R)$ are called the 
\emph{Eisenbud operators}\index{Eisenbud operator}~\cite{Eisenbud:1980a} 
for the relations $f_k$.  
\end{remark}

\begin{remark}
Since the relations $f_k$ are required to be in $\fm^2$, we have $b_{i,k}\in
\fm$, so $b_{i,k}(0)=0$, and the differential $d$ disappears on $\kk \otimes_R \Cliffq$. So
there is a natural map from $H^*(\Cliffq,d)$ to $\kk\otimes_R
\Cliffq$. This is the usual map $\HH^{*}(R) \to \Ext^{*}_R(\kk,\kk)$
obtained by applying $-\otimes_R \kk$ to a bimodule resolution of $R$
to obtain a module resolution of $\kk$.
\end{remark}

\begin{remark}
The interpretation of $\HH^1(R)$ given by Theorem~\ref{th:HHR} in
terms of derivations is as follows. Since $R$ is commutative, inner
derivations are zero. 
The element $\hat x_i$ is to be interpreted as the derivation $Q\to Q$
sending $x_i$ to $x_i$, and $x_j$ to zero if $j\ne i$. A linear
combination of $\hat x_i$ preserves the relations defining $R$, and
hence defines a derivation on $R$, if and
only if it is in the kernel of the differential $d$ of $\Cliffq$.
\end{remark}

\section{Koszul duality for graded algebras}\label{se:KoszulDual}%
\index{Koszul duality}

For many graded commutative rings whose relations are quadratic, 
Koszul duality provides a computation of 
both $\Ext$ and Hochschild cohomology. It is usual to have the
generators of a Koszul algebra in degree one. Because we want to apply
the theory to algebras which do not have this property, we deal with
bigraded algebras, where the generators are in degree $(n,1)$ for some
$n\in\bZ$. On the other hand, we shall only need to consider the case
where the degree zero part is $\kk$ and the degree one part is of
finite type with respect to the first degree.
So we use the following definition, adapted from the one in
Beilinson, Ginzburg and
Soergel~\cite{Beilinson/Ginzburg/Soergel:1996a}. They allow the degree
zero part to be semisimple, so our definition is more restrictive,
but suffices for our purposes.

\begin{definition}\label{def:Koszul}
A \emph{graded Koszul algebra}\index{graded Koszul algebra} is an
$\bZ\times\bZ$-graded $\kk$-algebra $R$ generated over $R_{0,0}=\kk$
by $R_{1,*}$, each $R_{1,i}$ is finite dimensional,
and with the property that
the minimal projective resolution of $\kk$ as a graded $R$-module is
\emph{linear}.\index{linear resolution}
In other words, the maps in the minimal resolution are given by
multiplication by linear combinations of the generators in degrees
$(1,j)$. Thus the $i$th projective in the resolution is generated by
its elements in degree $(i,*)$.
\end{definition}

The relations in a graded Koszul algebra are quadratic (Proposition~1.2.3
of~\cite{Beilinson/Ginzburg/Soergel:1996a}), but not every graded 
algebra with quadratic relations is Koszul. 

\begin{theorem}\label{th:Froberg}
If the generators in $R_{1,*}$ can be chosen in such a way that the
relations consist of monomials of first degree two, or consists of monomials
of first degree two together 
with for each $i<j$ a relation of the form $v_iv_j=c_{i,j}v_jv_i$ with
$c_{i,j}\in\kk$, then $R$ is Koszul.
\end{theorem}
\begin{proof}
This is proved in Fr\"oberg~\cite{Froberg:1975a}.
\end{proof}

But even for a
general graded commutative algebra with quadratic relations it is hard
to know whether the algebra is Koszul. Anick showed how to produce
commutative algebras with quadratic relations, for which the
generating function for the dimensions of $\Ext^i_R(\kk,\kk)$ is not a rational function, see
Example~7.1 of~\cite{Anick:1987a}. Such an algebra cannot be Koszul,
as we shall see below.

\begin{definition}
If $V$ is a graded vector space of finite type, we write $V^*$ for its
graded dual, with $(V^*)_i=(V_{-i})^*$. 
The pairing between $V$ and $V^*$ induces a pairing
between $V\otimes V$ and $V^*\otimes V^*\cong (V\otimes V)^*$.  If $S$
is a linear subspace of $V\otimes V$, we write $S^\perp$ for its
annihilator in $V^*\otimes V^*$. Thus $S^{\perp\perp}=S$.
If $v_i$, $i\in I$ form a basis for $V$ then we write $\hat v_i$ for
the dual basis of $V^*$.
\end{definition}

\begin{definition}
Let $V$ is a graded vector space of finite type and $S$ be a linear
subspace of $V\otimes V$ such that $R=\kk\langle V\rangle/(S)$ 
is a Koszul algebra. Then the \emph{Koszul dual} $R^!$ is defined to
be 
\[ R^!=\kk\langle V^*\rangle/(S^\perp).\]
Thus we have a natural isomorphism $R^{!!}\cong R$.
\end{definition}

\begin{theorem}
If $R$ is a graded Koszul algebra as defined in~\ref{def:Koszul}, then
$\Ext^{*,*}_R(\kk,\kk)\cong (R^!)^\op$. 
\end{theorem}
\begin{proof}
This is Theorem~2.10.1 of~\cite{Beilinson/Ginzburg/Soergel:1996a}.
\end{proof}

\begin{theorem}\label{th:Negron}
  Let $R=\kk\langle x_1,\dots,x_n\rangle /(S)$ be a graded Koszul algebra,
  with $S$ a set of quadratic relations, and let $R^!=\kk\langle
  \xx_1,\dots,\xx_n\rangle /(S^\perp)$ be the Koszul dual. Then as a
  $\kk$-algebra, the Hochschild cohomology $\HH^*R$ can be computed as
  $H^*(R\otimes R^!,d)$, where the differential
  $d$ given by $[e,-]$ where $e=x_1\otimes \xx_1 + \dots +
  x_n\otimes \xx_n$. Here, the variables $\xx_1,\dots,\xx_n$ are put
  in homological degree $-1$ in the complex, while $x_1,\dots,x_n$ are
  in homological degree zero.
\end{theorem}
\begin{proof}
In this form, this is proved in Theorem~1.2 of
Negron~\cite{Negron:2017a}, but see also Section~2 of Buchweitz,
Green, Snashall and
Solberg~\cite{Buchweitz/Green/Snashall/Solberg:2008a}, where this is
described in a more basis dependent way.
\end{proof}

\begin{remark}\label{rk:Negron}
In the case of a complete intersection $R$ with quadratic relations, of
course $S$ includes the commutativity relations, whereas the Koszul
dual $R^!$ is usually non-commutative.
In this case, the complex given in Theorem~\ref{th:Negron} is 
isomorphic to $(\Cliffq,d)$ appearing in Theorem~\ref{th:HHR}. The
advantage of this approach is that the same computation also computes
$\HH^*R^!$, but watching out for the change of degrees.
\end{remark}

\begin{remark}
If $R$ and $R^!$ are Koszul dual algebras, there is a relation
between the generating functions for the dimensions. 
Let
\[ p_R(s,t)=\sum_{i,j}s^it^j\dim_\kk R_{i,j}. \]
Then we have
\begin{equation}\label{eq:pR!st} 
p_{R^!}(s,t) = 1/p_R(-st^{-1},t^{-1}). 
\end{equation}
Without the internal grading, setting $t=1$ we recover the the more well  
known formula 
\begin{equation}\label{eq:pR!s} 
p_{R^!}(s)=1/p_R(-s) 
\end{equation}
proved in Theorem~2.11.1
of~\cite{Beilinson/Ginzburg/Soergel:1996a}. The proof in the doubly
graded case is the same. 
In some sense, we should also replace $s$ by
$s^{-1}$ for $R^!$ in order to make consistent use of homological
degrees; but then $R^!$ would be generated in degree $-1$ instead of
$1$. With this convention we have
\begin{equation}\label{eq:pR!st-1}
p_{R^!}(s,t)= 1/p_R(-s^{-1}t^{-1},t^{-1}).
\end{equation}
 \end{remark}

\begin{example}
If $R=\kk[x,y]$ with $x$ in degree $(-2,1)$ and $y$ in degree
$(-4,1)$ then $R^!=\Lambda(\hat x,\hat y)$ with $\hat x$ in degree $(1,-1)$ and
$\hat y$ in degree $(3,-1)$. Then we have
\begin{align*} 
p_R(s,t)&=1/(1-st^{-2})(1-st^{-4}), \\
p_{R^!}(s,t)&= (1-(-s^{-1}t^{-1})t^2)(1-(-s^{-1}t^{-1})t^4)=(1+s^{-1}t^{-1})(1+s^{-3}t^{-1}). 
\end{align*}
\end{example}

\section{The cyclic case}\label{se:cyclic}

In this section, we summarise the results on groups with cyclic Sylow
$p$-subgroups,\index{cyclic!Sylow subgroups}\index{Sylow subgroup!cyclic} 
from Benson and Greenlees~\cite{Benson/Greenlees:2021a,Benson/Greenlees:2023a}.

Let $G$ be a finite group with cyclic Sylow $p$-subgroups, and let
$\kk$ be a field of characteristic $p$. Then by
Theorem~\ref{th:abelian}, the inclusion of a Sylow
$p$-normaliser $N_G(D)\to G$ and the quotient map $N_G(D)\to
N_G(D)/O_{p'}N_G(D)$ (Notation~\ref{no:gt}) induce 
homotopy equivalences
\[ B(N_G(D)/O_{p'}N_G(D))\phat \leftarrow BN_G(D)\phat \rightarrow
  BG\phat. \]
The automorphism group of a cyclic group of order $p^n$ is abelian,
and its $p'$-part is cyclic of order $p-1$.
So using Remark~\ref{rk:DH}, it suffices to discuss the case
$\bZ/p^n\rtimes\bZ/q$, where $q\ge 2$ is 
a divisor of $p-1$ and $\bZ/q$ acts faithfully on $\bZ/p^n$. Indeed,
even in the case of a block with cyclic defect group of order $p^n$
and inertial index $q$, the $p$-completed classifying space 
(see Remark~\ref{rk:Linck}) has the homotopy type of $B(\bZ/p^n\rtimes\bZ/q)\phat$.

So set 
\[ G=\langle g,s\mid g^{p^n}=1, s^q=1, sgs^{-1}=g^\gamma\rangle
\cong\bZ/p^n\rtimes\bZ/q, \] 
where $\gamma$ is a primitive $q$th root of
unity modulo $p^n$. Setting
\[ U = \sum_{\substack{1\le j\le p^n-1,\\j^p\equiv j \pmod{p^n}}}g^j/j, \]
the group algebra is given by
\[ \kG=\kk\langle s,U\mid U^{p^n}=0,\ s^q=1,\ sU=\gamma Us \rangle \]
where $\gamma$ is a primitive $q$th root of unity. This
has a unique grading
up to scalar multiplication. It is convenient to use a
$\bZ[\frac{1}{q}]$-grading  and set $|s|=0$, $|U|=1/q$. With this
grading, the cohomology is the doubly graded ring given by
$H^*BG = \kk[x]\otimes\Lambda(t)$ with $|x|=(-2q,-p^n)$ and
$|t|=(-2q+1,-h)$ with $h=p^n-(p^n-1)/q$. The
$A_\infty$ structure is completely determined by
\[ m_{i}(t,\dots,t)=\begin{cases}
(-1)^{p^n(p^n-1)/2}x^h&i=p^n\\0&\text{otherwise}.
\end{cases} \]

The homology of the loop space on the $p$-completion 
$H_*\Omega BG\phat$ looks very similar. We have
\[ H_*\Omega BG\phat =k[\tau] \otimes \Lambda(\xi) \]
where $|\tau|=(2q-2,h)$ and $|\xi|=(2q-1,p^n)$. The $A_\infty$
structure is completely determined by
\[ m_i(\xi,\dots,\xi)=\begin{cases}
(-1)^{h(h-1)/2}& i=h \\
0 & \text{otherwise}.
\end{cases} \]
Thus the roles of $h$ and $p^n$ have been reversed. There is one
exceptional case. If $h=2$ then $q=2$ and $p^n=3$. In this case,
the formula above gives $m_2(\xi,\xi)=-\tau^3$. Thus $\xi$ is no
longer an exterior variable, but rather we have the
 formal $A_\infty$ algebra\index{formal $A_\infty$ algebra}
$H_*\Omega BG\phat=k[\tau,\xi]/(\xi^2+\tau^3)$ in this case. For
details, see Theorem~1.3 of~\cite{Benson/Greenlees:2021a}.

In Theorem~1.3 of~\cite{Benson/Greenlees:2023a}, it is proved that the category
$\Dsg(H^*BG)\simeq \Dcsg(H_*\Omega BG\phat)$ is equivalent 
to $\Db(H_*\Omega BG\phat[\tau^{-1}])$. Here, we are regarding these as singularity and
cosingularity categories of $A_\infty$ algebras.
This is a  finite
Krull--Schmidt triangulated category with $(q-1)(h-1)$ isomorphism
classes of indecomposable objects. The Auslander--Reiten
quiver\index{Auslander--Reiten quiver} of
this category is isomorphic to $\bZ A_{h-1}/T^{q-1}$, a cylinder of
height $h-1$ and circumference $q-1$. 
Here, $T$ is the
translation functor $\Sigma^{-2q}=\Sigma^{-2}$. The triangulated shift $\Sigma$ reverses the
ends of the cylinder, so that there are $[h/2]$ orbits of $\Sigma$ on
indecomposables. 

It is also shown in Theorem~1.3 of~\cite{Benson/Greenlees:2023a} that
the category $\Dcsg(H^*BG)\simeq \Dsg(H_*\Omega BG\phat)$ is
equivalent 
to $\Db(H^*BG[x^{-1}])$. This is a finite Krull--Schmidt triangulated
category with $q(p^n-1)$ isomorphism classes of indecomposable
objects. The Auslander--Reiten quiver is isomorphic to $\bZ
A_{p^n-1}/T^q$, a cylinder of height  $p^n-1$ and circumference $q$.
The translation functor this time is $T=\Sigma^{2(q-1)}$, and the
triangulated shift $\Sigma$ again reverses the ends of the cylinder,
so that there are $(p^n-1)/2$ orbits of $\Sigma$ on indecomposables.

\section{Squeezed resolutions}\index{squeezed resolution}\label{se:squeezed}

In one case, see Section~\ref{se:SD3}, the methods of extra gradings
used in the rest of this work do not give enough information, and we
resort to the method of squeezed resolutions. This is a method
developed in Benson~\cite{Benson:2009b} for computing 
$H_*\Omega BG\phat$ using modular representation theory. 

Let $\kk$ be a field of characteristic $p$, and let $G$ be a finite
group. If $M$ is a $\kG$-module, we define $[O^p(G),M]$ 
(Notation~\ref{no:gt}) to be the subspace of $M$ consisting of the 
elements $g(m)-m$ with $g\in O^p(G)$ and $m\in M$. This is the 
largest submodule of $M$ with no non-trivial map to $\kk$, the trivial
module, and is also the smallest with quotient filtered by copies of the
trivial module.

We begin with $P_0=N_0$, the projective cover of the trivial
$\kG$-module $\kk$. For $i\ge 1$, $M_{i-1}=[O^p(G),N_{i-1}]$, $P_i$
the projective cover of $M_{i-1}$, and $N_i=\Omega(M_{i-1})$, the
kernel of $P_i\to M_{i-1}$.
\[ \xymatrix@=4mm{
\cdots\ar[r]&P_2\ar[dr]\ar[rrr]&&&P_1\ar[dr]\ar[rrr]&&&P_0\\
M_3\ar[ur]&&[O^p(G),M_2]\ar@{^(->}[r]&M_2\ar[ur]&&[O^p(G),M_1]\ar@{^(->}[r]
&M_1\ar[ur]} \]
This is the \emph{minimal squeezed resolution} for $G$.
Finally, we define the \emph{squeezed homology}  
denoted $H_*^\Omega(G,\kk)$, to be the homology of this
complex. 

More generally, a complex $\cdots\to P_2\to P_1\to P_0$ is a (left)
\emph{squeezed resolution} if the following conditions are satisfied:
\begin{enumerate}
\item[\rm (i)] each $P_i$ is projective,
\item[\rm (ii)] $H_n(P_* \otimes_{\kG} \kk)\cong
\begin{cases}\kk(G/O^p(G))&n=0\\0&n>0,\end{cases}$
\item[\rm (iii)] $[O^p(G),H_n(P_*)]=0$ for all $n\ge 0$.
\end{enumerate}

\begin{theorem}
All squeezed resolutions are homotopy equivalent.
The homology of a squeezed resolution $H_*^\Omega(G,\kk)$ is
isomorphic to $H_*\Omega BG\phat$. 
\end{theorem}
\begin{proof}
This is Theorem~1.2 of~\cite{Benson:2009b}.
\end{proof}

\begin{remark}
There is an error\index{errors}
in~\cite{Benson:2009b}, which only shows up if $\Omega BG\phat$ is not
connected, namely when $G/O^p(G)$ (Notation~\ref{no:gt}) is not trivial. Namely, the
augmentation in Definition~3.2 and Theorem~3.4 of that paper should be to $\kG/O^p(G)$
rather than to $\kk$. This affects the computation of products in
Section~4.6 of the paper, which we shall be using.
\end{remark}

To multiply elements of $H_*^\Omega(G,\kk)$, we lift an element to a
map of complexes
\[ \xymatrix{\cdots\ar[r]&P_1\ar[r]\ar[d]&P_0\ar[r]\ar[d]&0\ar[d]\\
\cdots\ar[r]&P_{n+1}\ar[r]&P_n\ar[r]&P_{n-1}\ar[r]&\cdots\ar[r]&P_0\ar[r]&0} \]
and compose, see~\cite[4.6]{Benson:2009b}. This gives us the loop
product on $H_*\Omega BG\phat$.

\chapter{The dihedral case}\index{dihedral group}%
\index{Sylow subgroup!dihedral}

\section{Introduction}

In this chapter, we discuss the case of finite groups with dihedral Sylow $2$-subgroups.
The $A_\infty$ structure on the cohomology of dihedral groups 
is given in the following theorem. 

\begin{theorem}\label{th:H*D}
Let $\D$ be a dihedral group of order $4q$, where $q\ge 2$ is a power
of two, and let $\kk$ be a field of characteristic two. Then
as a ring, we have $H^*B\D=\kk[x,y,t]/(xy)$ with $|x|=|y|=-1$
and $|t|=-2$. Up to quasi-isomorphism, the $A_\infty$ structure on
$H^*B\D$ is determined by 
\[  m_{2q}(x,y,\dots,x,y) = m_{2q}(y,x,\dots,y,x) = t. \] 
\end{theorem}

We give an explicit $A_\infty$ structure within this quasi-isomorphism
class in Theorem~\ref{th:AinftyHBG}. 
It has $m_i\ne 0$ if and only if $i$ is congruent to $2$ modulo
$2q-2$. For completeness, we also describe the case $q=1$, which
behaves differently.

The idea of the proof is to put a double grading on the group algebra
$\kD$. This gives a triple grading on $H^*B\D$, which then restricts
the possibilities for the higher $m_i$. It is then easy to check that
$m_i=0$ unless $i$ is congruent to $2$ modulo $2q-2$. 
Then $m_{2q}$ is
interpreted as a Hochschild cocycle on $H^*B\D$, and 
quasi-isomorphism amounts to changing it by a Hochschild 
coboundary. We write down explicit formulas for all the $m_i$,
using some Hochschild cohomology computations involving the
circle product of Gerstenhaber.

Similar computations give the $A_\infty$ structure on the cohomology
of a finite group $G$ with dihedral Sylow $2$-subgroups of order
$4q$. 
These groups were classified by Gorenstein and
Walter~\cite{Gorenstein/Walter:1962a,Gorenstein/Walter:1965abc}. There
are three 
possible $2$-local structures, which are distinguished by the number
of conjugacy classes of involutions\index{involutions} 
(one, two or three). We examine the three
possibilities in detail, and determine the $A_\infty$ structures on
$H^*BG$ and $H_*\Omega BG\twohat$ in each case. In the case with three
conjugacy classes, $G$ has a normal $2$-complement, so $BG\twohat\simeq
B\D\twohat$. So we concentrate on the remaining two cases. 

Perhaps the most interesting case is the one where all involutions are
conjugate, as this happens if and only if $G$ has no subgroup of index
two. In this case, if $q\ge 2$ we have
\[ H^*BG=\kk[t,\xi,\eta]/(\xi\eta) \]
with $|t|=-2$ and $|\xi|=|\eta|=-3$ (homological grading).
If $q=1$ then $H^*BG=\kk[t,\xi,\eta]/(\xi\eta+t^2)$.
This time, the $A_\infty$ structure is determined up to
quasi-isomorphism by
\[ m_{2q}(\xi,\eta,\dots,\xi,\eta)=m_{2q}(\eta,\xi,\dots,\eta,\xi)
  = t^{2q+1}, \]
where the $\xi$ and $\eta$ alternate. Again the $m_i$ are non-zero for $i$
congruent to $2$ modulo $2q-2$, and zero otherwise, and we give an
explicit description of the non-zero ones.

The $A_\infty$ structure on $H_*\Omega BG\twohat$ in this case is easier to
describe than that of $H^*BG$. This is because there are only two
non-zero $m_i$. 

\begin{theorem}
Let $G$ be a finite group with dihedral Sylow $2$-subgroups of order
$4q$ with $q\ge 1$ a power of two, and with no normal subgroup of index
two, and let $\kk$ be a field of characteristic two. Then
the ring structure on the homology of $\Omega BG\twohat$ is given by
\[ H_*\Omega BG\twohat = \Lambda(\tau) \otimes \kk\langle
  \alpha,\beta \mid \alpha^2=\beta^2=0\rangle \]
where $|\tau|=1$, $|\alpha|=|\beta|=2$. The $A_\infty$ structure is
determined by 
\[ m_{2q+1}(\tau,\dots,\tau) =  s^q, \]
where $s=\alpha\beta+\beta\alpha$. 
\end{theorem}

See Theorems~\ref{th:HOmegaBG-D1} and~\ref{th:AinftyHOmegaBG-D1}
for details. 

We describe a DG algebra $Q$ which is quasi-isomorphic to the
$A_\infty$ algebra
$H_*\Omega BG\twohat$, and use it to show that the degree $4$ element
$s=\alpha\beta+\beta\alpha$ is central. It may then be inverted, to
obtain equivalences of categories
\[ \Db(Q[s^{-1}]) \simeq
\Db(H_*\Omega BG\twohat[s^{-1}])
\simeq \Dcsg(H_*\Omega BG\twohat) \simeq \Dsg(H^*BG) \]
(cf.\ Theorem~\ref{th:Dsg-Dcsg}).

Finally, we observe that there is a Morita equivalence between
$Q[s^{-1}]$ and one of the algebras discussed
in Benson and Greenlees~\cite{Benson/Greenlees:2023a}. This allows us to carry over the
classification theorem there, to
classify the indecomposable objects in $\Db(H_*\Omega
BG\twohat[s^{-1}])$, and hence also of $\Dsg(H^*BG)$.

\begin{theorem}
Let $G$ be a finite group with dihedral Sylow $2$-subgroups
of order $4q$ with $q\ge 1$ a power of two, and with no normal subgroup of index two, and let $\kk$ be a
field of characteristic two. Then $\Dsg(H^*BG)\simeq \Dcsg(H_*\Omega BG\twohat)$ is a finite
Krull--Schmidt category with
$4q$ isomorphism
classes of indecomposable objects, which come in 
$q$ orbits of the suspension $\Sigma$, all of length four.
The Auslander--Reiten quiver\index{Auslander--Reiten quiver} 
is isomorphic to $\bZ A_{2q}/T^2$, where
$T$ is the translation functor\index{translation functor} $\Sigma^{-2}$.
\end{theorem}

This theorem is proved in Section~\ref{se:classification} 
(Theorem~\ref{th:D1}). The
corresponding theorem in the case where $G$ has two conjugacy classes
of involutions, so that $G$ has a normal subgroup of index
two but no normal subgroup of index four, is given in
Theorem~\ref{th:D2}.  

In contrast with Theorems~\ref{th:D1} and~\ref{th:D2}, the
category $\Dcsg(H^*BG)\simeq\Dsg(H_*\Omega BG\twohat)$ has infinite
representation type.

\section{\texorpdfstring{Dihedral $2$-groups}
{Dihedral 2-groups}}\label{se:dihedral}

Let $\D=\langle g,h \mid g^2=h^2=(gh)^{2q}=1\rangle$, a dihedral
group of order $4q$, with $q\ge 1$ a power of two, and let $\kk$ be a
field of characteristic two. As elements of $\kD$, let $X=g-1$ and
$Y=h-1$. Then the group algebra can be rewritten as 
\begin{equation}\label{eq:kD}
\kD = \kk\langle X,Y \mid X^2=Y^2=0, (XY)^q=(YX)^q  \rangle. 
\end{equation}
This algebra has tame representation type, and the finitely generated
$\kG$-modules were classified by Ringel~\cite{Ringel:1975a}. 

We shall regard $\kD$  as a $\bZ\times\bZ\,$-graded algebra, with $|X|=(1,0)$ and
$|Y|=(0,1)$. With this bigrading, the relations are homogeneous.
It is easy to compute the minimal resolution of $\kk$ as a $\kG$-module,  
and hence the cohomology ring.  
Recall that we are using homological degrees 
throughout, so that cohomological degrees come out negative.
We list first the homological degree,   
and then the two internal degrees.   

\begin{remark}\label{rk:Z2xZ2formal}
The case $q=1$ behaves differently from $q\ge 2$, so we discuss this
case first. If $q=1$ then $\kD$ is an exterior algebra on two
generators, so
$H^*B\D\cong\Ext^*_{\kG}(\kk,\kk)$ is a formal $A_\infty$
algebra\index{formal $A_\infty$ algebra} 
$\kk[x,y]$ with $|x|=-(1,1,0)$
and $|y|=-(1,0,1)$. One way to see that it has to be formal is that
the values $m_i$ for $i>2$ on non-zero elements land in zero groups
for degree reasons. Namely, $m_i$ adds internal degrees, but it adds
homological degree and then adds $i-2$ (see
Section~\ref{se:Ainfinity}). See Theorem~\ref{th:elemab2} for a more
general statement.
\end{remark}

We now assume, for the rest of this section, that $q\ge 2$. We have
\[ H^*B\D\cong\Ext^*_{\kD}(\kk,\kk) \cong \kk[x,y,t]/(xy) \] 
where $|x|=-(1,1,0)$, $|y|=-(1,0,1)$ 
and $|t|=-(2,q,q)$. This can be computed using, for example, the
resolution given in Hamada~\cite{Hamada:1963a} and the diagonal
approximation given in Handel~\cite{Handel:1993a}, or alternatively
using the diagrammatic method described in Benson and
Carlson~\cite{Benson/Carlson:1987a}. The minimal resolution of $\kk$ is the
total complex of the following double complex that displays the gradings
on cohomology
\[ \xymatrix{
\vdots\ar[d]^{Y'}&\vdots\ar[d]^Y&\vdots\ar[d]^{Y'}
&\vdots\ar[d]^Y\\
\kD\ar[d]^Y&\kD\ar[l]_{X'}\ar[d]^{Y'}&\kD\ar[l]_X\ar[d]^Y
&\kD\ar[l]_{X'}\ar[d]^{Y'}&\cdots\ar[l]_X\\
\kD\ar[d]^{Y'}&\kD\ar[l]_X\ar[d]^Y&\kD\ar[l]_{X'}\ar[d]^{Y'}
&\kD\ar[l]_X\ar[d]^Y&\cdots\ar[l]_{X'}\\
\kD\ar[d]^Y&\kD\ar[l]_{X'}\ar[d]^{Y'}&\kD\ar[l]_X\ar[d]^Y
&\kD\ar[l]_{X'}\ar[d]^{Y'}&\cdots\ar[l]_X\\
\kD&\kD\ar[l]_X&\kD\ar[l]_{X'}&\kD\ar[l]_X
&\cdots\ar[l]_{X'}} \]
where $X'=(XY)^{q-1}X$, $Y'=(YX)^{q-1}Y$. 
The elements $x$ and $y$ in
$H^1B\D$ are dual to $X$ and $Y$ in $J(\kD)/J^2(\kD)$. The element 
\[ t \in H^2B\D \cong \Ext^1_{\kD}(\Omega \kk, \kk) \]
is represented by the short exact sequence of bigraded modules
\[ 0 \to \kk \to M \oplus N \to \Omega \kk \to 0, \]
where $M$ and $N$ are uniserial modules of length $2q$.
Examination of these
uniserial $\kD$-modules, using Theorem~2.4 of Dwyer~\cite{Dwyer:1975a}, 
shows that we have Massey products\index{Massey product}
\[ \langle x,y,\dots,x,y\rangle = \langle y,x,\dots,y,x \rangle= t \]
in $H^*B\D$.
Here, in both expressions
the arguments $x$ and $y$ alternate, and there are $q$ of each, for a
total of $2q$ terms.
Note that these Massey products are only well defined up to adding
multiples of $x^2$ and $y^2$. However, if we take the internal grading into account,
the Massey product is well defined, with no ambiguity.

\section{\texorpdfstring{$\HH^*H^*B\D$}{HH*H*BD}}\label{se:HH*H*BD}

Wishing to understand further the $A_\infty$ structure of the
cohomology of dihedral groups, 
it follows from Proposition~\ref{pr:HH} and Lemma~\ref{le:m=0} that
we should next compute the Hochschild cohomology $\HH^*H^*B\D$.
Here, $H^*B\D$ is just being regarded as a $k$-algebra without
reference to higher structure.
Since $H^*B\D$ is a
hypersurface\index{hypersurface} 
(i.e., a codimension one complete intersection),\index{complete intersection} 
we can compute Hochschild cohomology using Theorem~\ref{th:HHR}.
So first we compute
the algebra $\Cliffq$\index{Cliff@$\Cliffq$} for $H^*B\D$,
where $\D$ is a dihedral group of order $4q$ with $q\ge 2$.
This will also be useful for computing $\Ext^{*,*}_{H^*B\D}(\kk,\kk)$
using Theorem~\ref{th:ExtRkk}.
Recall that $H^*B\D=\kk[x,y,t]/(xy)$ with $|x|=-(1,1,0)$,
$|y|=-(1,0,1)$ and $|t|=-(2,q,q)$.

\begin{proposition}\label{pr:Cliffq-dihedralgroup}
The DG algebra $\Cliffq$ is equal to $H^*B\D\langle \hat x,
\hat y, \tau; s\rangle$, where $s$ is central, and $\hat x^2=0$, 
$\hat y^2=0$,  $\hat x\hat y+\hat y\hat x = s$, $\tau^2=0$,
$\hat x\tau=\tau\hat x$, $\hat y\tau=\tau\hat y$. The degrees are
given by 
$|x|=(0,-1,-1,0)$, $|y|=(0,-1,0,-1)$, $|t|=(0,-2,-q,-q)$,
$|\hat x|=(-1,1,1,0)$, $|\hat y|=(-1,1,0,1)$, $|\tau|=(-1,2,q,q)$,
$|s|=(-2,2,1,1)$. The differential is given by $d(\hat x)=ys$, $d(\hat
y)=xs$, $d(\tau)=0$, $d(s)=0$.
\end{proposition}
\begin{proof}
Let $f(x,y,t)=xy$. Then we have
\begin{gather*} 
\frac{\partial f}{\partial x}=y,\qquad 
\frac{\partial f}{\partial y}=x,\qquad 
\frac{\partial f}{\partial z}=0, \\
\frac{\partial^{(2)} f}{\partial x^2}=0, \qquad
\frac{\partial^{(2)} f}{\partial y^2}=0,\qquad
\frac{\partial^{(2)} f}{\partial z^2}=0,\\
\frac{\partial^2 f}{\partial x\partial y}=1, \qquad
\frac{\partial^2 f}{\partial x\partial z}=0,\qquad
\frac{\partial^2 f}{\partial y\partial z}=0.
\end{gather*}
Plugging these into Definition~\ref{def:Cliffq}, we get the given
relations and differential for $\Cliffq$.
\end{proof}

\begin{theorem}\label{th:HHHBD}
Let $G$ be a dihedral group of order $4q$ with $q\ge 2$ a power of two,
and let $\kk$ be a field of characteristic two.
The Hochschild cohomology\index{Hochschild cohomology} $\HH^*H^*B\D$
has generators $s$, $t$, $\tau$, $x$, $y$, $u$, $v$
with
\begin{align*}
|s|&=(-2,2,1,1) \\
|t|&= (0,-2,-q,-q) &
|\tau|&= (-1,2,q,q) \\
|x|&= (0,-1,-1,0)&
|y|&= (0,-1,0,-1) \\
|u|&= (-1,0,0,0) &
|v|&= (-1,0,0,0).
\end{align*}
The relations are given by $u^2=v^2=uv=\tau^2=xy=xv=yu=xs=ys=0$, 
and $us=vs$. The non-zero monomials and their degrees are as
follows, with $i_1,i_2\ge 0$, $\ep_1,\ep_2\in\{0,1\}$. The
first two cases overlap for $i_1>0$, the first and third, and the
second and fourth overlap for $i_1=0$.
\begin{align*}
|s^{i_1}t^{i_2}\tau^{\ep_1}u^{\ep_2}|&
=(-2i_1-\ep_1-\ep_2,2i_1-2i_2+2\ep_1,i_1+q(-i_2+\ep_1),i_1+q(-i_2+\ep_1)), \\
|s^{i_1}t^{i_2}\tau^{\ep_1}v^{\ep_2}|&
=(-2i_1-\ep_1-\ep_2,2i_1-2i_2+2\ep_1,i_1+q(-i_2+\ep_1),i_1+q(-i_2+\ep_1)), \\
|x^{i_1}t^{i_2}\tau^{\ep_1}u^{\ep_2}|&
=(-\ep_1-\ep_2,-i_1-2i_2+2\ep_1,-i_1+q(-i_2+\ep_1),q(-i_2+\ep_1)) \\
|y^{i_1}t^{i_2}\tau^{\ep_1}v^{\ep_2}|&
=(-\ep_1-\ep_2,-i_1-2i_2+2\ep_1,q(-i_2+\ep_1),-i_1+q(-i_2+\ep_1)) 
\end{align*}
(the top two coincide with the lower two when $i_1=0$, and are
otherwise disjoint).
There is only one monomial
in degree $(-i,i-2,0,0)$ with $i>2$, namely $s^qt$, with
\[ |s^qt|=(-2q,2q-2,0,0). \]
\end{theorem}
\begin{proof}
By Theorem~\ref{th:HHR}, $\HH^*H^*B\D$ is the cohomology of the
DG algebra $\Cliffq$. By Proposition~\ref{pr:Cliffq-dihedralgroup},
this is therefore 
as described in the theorem, with $u=x\hat x$ and $v=y\hat y$. Since
$d(\hat x\hat y) = (x\hat x+y\hat y)s$ in $\Cliffq$, we have $us=vs$ in $\HH^*H^*B\D$.

We also mention another approach to this computation, as this will
become relevant in the proof of
Proposition~\ref{pr:HHHOmegaBG-D1}. 
Namely, we can use Theorem~\ref{th:Negron}.
This gives rise to the same complex as above. Here, the $x_i$ are 
$x$, $y$ and $t$ and the $\xx_i$ are $\hat x$, $\hat y$ and $\tau$. 
The advantage of this approach is that it makes
it easy to compute $\HH^*A^!$ using the same computation, but watching
out for the changes of degrees, see Remark~\ref{rk:Negron}. This approach also makes it clear that
$\hat x$ and $\hat y$ are really just avatars for the elements $X$ and $Y$ of $\kD$.

For the last statement, we note that for the last two coordinates to
be zero, the monomial must be of one of the first two types. Then we
have
\begin{align*}
i_1+q(-i_2+\ep_1)&= 0, \\
(-2i_1-\ep_1-\ep_2)+(2i_1-2i_2+2\ep_1)&=-2.
\end{align*}
Twice the first equation minus $q$ times the second gives
$(\ep_1+\ep_2)q+2i_1=2q$, and so $i_1$ is either zero or $q$.
If $i_1=0$ then $\ep_1=\ep_2=1$, which then implies $i_2=1$, and the
resulting monomials have $i=2$. On the other hand, if $i_1=q$ then
$\ep_1=\ep_2=0$, and again we have $i_2=1$, and the resulting monomial
is $s^qt$.
\end{proof}

\section{\texorpdfstring{$A_\infty$ structure of $H^*B\D$}
{A∞ structure of H*BD}}\label{se:AinftyHBD}

In this section, we completely describe the 
$A_\infty$ structure\index{Ainfinity@$A_\infty$ algebra} of
$H^*B\D$. This makes extensive use of Section~\ref{se:circle}, describing
Gerstenhaber's circle product on Hochschild cochains and its relation
to the structure maps of an $A_\infty$ algebra. The multiple gradings
that we have developed on $\HH^*H^*B\D$ in Section~\ref{se:HH*H*BD}
play an essential role in the computation. 

Let $\D$ be a dihedral group of order $4q$ with $q\ge 2$ a power of  
two, and let $\kk$ be a field of characteristic two.

\begin{lemma}\label{le:m=0}
For any $A_\infty$ structure
on $H^*B\D$ that preserves internal degrees, 
we have $m_n=0$ unless $n-2$ is divisible by $2q-2$. In
particular for $2<n<2q$ we have $m_n=0$.
\end{lemma}
\begin{proof}
Looking at the degrees of the generators $x$, $y$ and $z$, 
for any monomial $\zeta$ in $H^*B\D$ of degree $(a,b,c)$ we have
$a\equiv b+c \pmod{2q-2}$. So for an $n$-tuple
$(\zeta_1,\dots,\zeta_n)$, the degree of
$m_n(\zeta_1,\dots,\zeta_i)$ satisfies $a\equiv b+c+ n -
2 \pmod{2q-2}$. It follows that for $m_n(\zeta_1,\dots,\zeta_n)$ to be non-zero we must have
$n - 2 \equiv 0 \pmod{2q-2}$.
\end{proof}

\begin{theorem}\label{th:AinftyHBG}
 The $A_\infty$ structure on $H^*B\D$ is given as follows.
The $m_n$ are $\kk[t]$-multilinear maps with $m_n=0$ 
for $n$ not congruent to $2$ modulo $2q-2$. For $i,j\ge 1$,
\[ m_{2q}(x^i,y,x,y,\dots,x,y^j)=m_{2q}(y^j,x,y,x,\dots,y,x^i)=x^{i-1}y^{j-1}t \]
where the arguments alternate between $x$ and $y$, and the right hand
side is zero unless either $i=1$ or $j=1$; $m_{2q}$ is zero on all
other tuples of monomials not involving $t$.
The maps $m_{\ell(2q-2)+2}$ with $\ell>1$ similarly vanish on all tuples of
monomials not involving $t$, except the ones which look as above,
but for some choice of indices:
\begin{multline*} 
1\le e_1\le e_2\le \cdots \le e_{\ell-1} 
< e_{\ell-1}+(2q-2)+1   
\le e_{\ell-2}+2(2q-2)+1 \\
\le \cdots \le e_1+(\ell-1)(2q-2)+1\le \ell(2q-2)+2. 
\end{multline*}
the exponents on the terms are increased by one
(or correspondingly more if an index is repeated).
The value on these tuples is $x^{i-1}y^{j-1}t^{\ell}$.
Thus 
\[
m_{\ell(2q-2)+2}(x^{i+\alpha_1},y^{\alpha_2},x^{\alpha_3},\dots,
x^{\alpha_{\ell(2q-2)+1}},y^{j+\alpha_{\ell(2q-2)+2}}) =
  x^{i-1}y^{j-1}t^\ell \]
where each $\alpha_\sigma$ is one plus the number of indices in the list above that are equal
to $\sigma$.
\end{theorem}

\begin{remark}\label{rk:HH}
To illustrate this rather complicated looking condition, suppose that
$q=4$. Then $m_8$ is given by
\[ m_8(x^i,y,x,y,x,y,x,y^j)=m_8(y^j,x,y,x,y,x,y,x^i)=x^{i-1}y^{j-1}t, \]
and then $m_{14}$ is the next non-zero $m_n$. The value of
each of the following seven expressions is $x^{i-1}y^{j-1}t^2$:
\begin{small}
\begin{gather*}
m_{14}(x^{i+1},y,x,y,x,y,x,y^2,x,y,x,y,x,y^j) \\
m_{14}(x^i,y^2,x,y,x,y,x,y,x^2,y,x,y,x,y^j) \\
m_{14}(x^i,y,x^2,y,x,y,x,y,x,y^2,x,y,x,y^j) \\
m_{14}(x^i,y,x,y^2,x,y,x,y,x,y,x^2,y,x,y^j) \\
m_{14}(x^i,y,x,y,x^2,y,x,y,x,y,x,y^2,x,y^j) \\
m_{14}(x^i,y,x,y,x,y^2,x,y,x,y,x,y,x^2,y^j) \\
m_{14}(x^i,y,x,y,x,y,x^2,y,x,y,x,y,x,y^{j+1})
\end{gather*}
\end{small}%
There are seven more such expressions with non-zero values of $m_{14}$, where $x$
and $y$ have been interchanged.
A typical non-zero value of $m_{20}$, which is the next non-zero $m_n$,
corresponding to $\ell=3$, is given by
\[  m_{20}(x^i,y,x,y^2,x,y^2,x,y,x,y,x,y,x^2,y,x,y,x^2,y,x,y^j)=x^{i-1}y^{j-1}t^3, \]
where the indices $4\le 6<13\le 17$ come from $e_1=4$, $e_2=6$. An example
with a repeated index is
\[  m_{20}(x^i,y,x,y,x,y^3,x,y,x,y,x,y,x^2,y,x,y,x,y,x^2,y^j)=x^{i-1}y^{j-1}t^3, \]
with indices $6\le 6 < 13 \le 19$ coming from $e_1=e_2=6$.
\end{remark}

\begin{proof}[Proof of Theorem~\ref{th:AinftyHBG}]
By Lemma~\ref{le:m=0}, for every $A_\infty$ structure preserving
degrees, $m_n=0$ for $2<n<2q$. So in order to determine
$m_{2q}$, we invoke Proposition~\ref{pr:HH}. We have Massey
products\index{Massey product}
\[ \langle x,y,\dots,x,y\rangle=\langle y,x,\dots,y,x\rangle=t, \]
well defined modulo the ideal generated by $x$ and $y$.
It follows from Theorem~\ref{th:Massey}
that $m_{2q}(x,y,\dots,x,y)$ and $m_{2q}(y,x,\dots,y,x)$ are non-zero.
So $m_{2q}$ represents a non-zero Hochschild cohomology class
in degree $(-2q,2q-2,0,0)$ in $\HH^*H^*B\D$. 
By Theorem~\ref{th:HHHBD},
up to scalar multiplication, there is only one non-zero possibility for
$m_{2q}$ up to quasi-isomorphism. It is easy to check that the given
formula for $m_{2q}$ is indeed a Hochschild cocycle. Replacing $t$
by a non-zero multiple of $t$ if necessary (or by working over
$\bF_2$) makes this the correct Hochschild cohomology class.

Again using Lemma~\ref{le:m=0}, we see that the next possible
$m_n$ after $m_{2q}$ is $m_{4q-2}$. Using~\eqref{eq:Ainfty},
this has to satisfy
\begin{gather*} 
m_2(\id \otimes m_{4q-2}) + \sum_{r=0}^{4q-3} m_{4q-2}(\id^{\otimes r} \otimes m_2 \otimes
  \id^{\otimes(4q-r-3)}) 
+  m_2(m_{4q-2}\otimes \id) \\
{}\qquad\qquad + \sum_{r=0}^{2q-1} m_{2q}(\id^{\otimes r} \otimes m_{2q}\otimes
\id^{\otimes 2q-r-1})=0.
\end{gather*}
Now the first three terms are the Hochschild coboundary of $m_{4q-2}$, 
while the last sum is the Gerstenhaber circle
product\index{Gerstenhaber!circle product}
$m_{2q}\circ m_{2q}$, see Section~\ref{se:circle}.
So as in Proposition~\ref{pr:circle}, we rewrite the above equation as 
\begin{equation}\label{eq:m4q-2} 
\HHd m_{4q-2}=m_{2q}\circ m_{2q}, 
\end{equation}
where $\HHd$ is the Hochschild coboundary.
Subject to this condition, $m_{4q-2}$ is well defined modulo
Hochschild coboundaries. But by Theorem~\ref{th:HHHBD}, the Hochschild cohomology
$\HH^*H^*B\D$ is zero in degree $(-4q+2,4q-4,0,0)$, so any choice of
$m_{4q-2}$ satisfying~\eqref{eq:m4q-2} is valid. The one we have
constructed satisfies this. 

We continue by induction on $\ell$. If we have constructed
$m_{2q},m_{4q-2}, \dots,m_{(\ell-1)(2q-2)+2}$, then the equation
satisfied by $m_{\ell(2q-2)+2}$ is
\[ \HHd m_{\ell(2q-2)+2}=\sum_{i+j=\ell}m_{i(2q-2)+2}\circ
  m_{j(2q-2)+2}. \]
Again $\HH^*H^*B\D$ is zero in degree $(-\ell(2q-2)-2,\ell(2q-2),0,0)$, and so
any choice of $m_{\ell(2q-2)+2}$ satisfying this equation is
valid. The one we have constructed satisfies this.
\end{proof}

\begin{remark}
Let us illustrate the way equation~\eqref{eq:m4q-2} works, with the
example of Remark~\ref{rk:HH}. We have
\begin{align*} 
(m_8 \circ m_8)(x^i,y,x,y^2,&x,y,x,y,x,y,x,x,y,x,y^j) \\ 
&= m_8(x^i,y,x,m_8(y^2,x,y,x,y,x,y,x),x,y,x,y^j) \\
&=m_8(x^i,y,x,yt,x,y,x,y^j) \\
&= x^{i-1}y^{j-1}t^2, 
\end{align*}
and correspondingly,
\begin{align*}
\HHd m_{14}(x^i,y,x,y^2,&x,y,x,y,x,y,x,x,y,x,y^j) \\
&=m_{14}(x^i,y,x,y^2,x,y,x,y,x,y,x^2,y,x,y^j) \\
&=x^{i-1}y^{j-1}t^2.
\end{align*}
\end{remark}

\section{\texorpdfstring{Loops on $B\D\twohat$}
{Loops on BDٛ₂}}\label{se:OmegaBDtwohat}

Throughout this section, $\D$ is a dihedral group of order $4q$, where  
$q$ is a power of two.  
Since $\D$ is a finite $2$-group, completing $B\D$ makes no difference
to its homotopy type, by Theorem~\ref{th:Bousfield-Kan}. So $\Omega B\D\twohat$ has contractible
connected components, and is homotopy equivalent to $D$ with the group
multiplication. It follows that $C_*\Omega B\D\twohat\simeq\kD$, and
we should expect to see the Eilenberg--Moore 
spectral sequence~\eqref{eq:Cotor} 
converging to $\kD$.\index{Eilenberg--Moore spectral sequence}

\begin{proposition}\label{pr:Ext*H*BDkk}
We have
\[ \Ext^{*,*}_{H^*B\D}(\kk,\kk) = \Lambda(\tau) \otimes \kk\langle
  \hat x,\hat y\mid
  \hat x^2=0,\ \hat y^2=0\rangle \]
with 
\[ |\tau|=(-1,2,q,q),\qquad |\hat x|=(-1,1,1,0),\qquad |\hat y|=(-1,1,0,1). \]
\end{proposition}
\begin{proof}
By Theorem~\ref{th:ExtRkk}, $\Ext^{*,*}_{H^*B\D}(\kk,\kk)\cong
\kk\otimes_{H^*B\D}\Cliffq$, where $\Cliffq$ is as described in 
Proposition~\ref{pr:Cliffq-dihedralgroup}. The positive degree
elements of $H^*B\D$ are killed by the tensor product with $\kk$. The
element $s$ there is 
redundant, as it is equal to $\hat x\hat y+\hat y\hat x$, so we are
left with $\tau$, $\hat x$ and $\hat y$. The element
$\tau$ squares to zero and commutes with $\hat x$ and $\hat y$. The
elements $\hat x$ and $\hat y$ do not commute, but they both square to
zero. There are no further relations. A basis for
$\Ext^{*,*}_{H^*B\D}(\kk,\kk)$ is given by the words in $\hat x$ and
$\hat y$ that alternate between them, together with $\tau$ times such words.
\end{proof}

The $E^2$ page of the spectral sequence
\[ \Ext^{*,*}_{H^*B\D}(\kk,\kk) \Rightarrow \kD \]
is given by Proposition~\ref{pr:Ext*H*BDkk}. Recall that we grade
everything homologically, so the first, homological degree in $\Ext^{*,*}$ is
negative, as indicated in the Proposition. The second, internal degree
is really a $\bZ^3$ grading consisting of the second, third and fourth
in the Proposition. The differential only changes the first and second
degree, with the third and fourth carried along for the ride.

By examining degrees, the
only page of the  
spectral sequence that can support a non-zero differential is
$E^{2q-1}$. Since the spectral sequence has to converge to something
finite, $d^{2q-1}$ has to be non-zero. The only possibility that has a
chance of giving rise to the structure of $\kD$ given in~\eqref{eq:kD}
is given by the equation
\[ d^{2q-1}(\tau)=(\hat x\hat y+\hat y\hat x)^q. \]
Then 
\[ E^{2q}=E^\infty\cong \kD=\kk\langle \hat x,\hat y\mid 
\hat x^2=0,\ \hat y^2=0,
(\hat x\hat y+\hat y\hat x)^q=0\rangle. \] 
Ungrading, $X$ represents
$\hat x$ and $Y$ represents $\hat y$ to give an
isomorphism with $\kD$.
Note that $\kD$ is isomorphic to its associated graded with respect to
the radical filtration, which is reflected in the fact that there is
no ungrading to be done in this case. In the generalised quaternion
and semidihedral situations, this will be more of an issue.

\section{\texorpdfstring{$\HH^*\kD$}{HH*kD}}\label{se:HH*kD}

Throughout this section, $\D$ is a dihedral group of order $4q$, where 
$q\ge 2$ is a power of two. 
In this section we use the spectral sequence~\eqref{eq:HHHA}
\begin{equation}\label{eq:ssHHkD} 
\HH^*H^*B\D \Rightarrow \HHinf^*H^*B\D \cong \HHinf^*H_*\Omega B\D\cong \HH^*\kD 
\end{equation}
to recompute $\HH^*\kD$ (see Remark~\ref{rk:HHkD}). Here, the first
isomorphism comes from 
Theorem~\ref{th:HHOmega} and the second comes from the
equivalence $\Omega B\D\simeq\D$.  This computation is not needed in
the rest of the paper, but is an illustration of the power of the internal
$\bZ\times\bZ$-grading on $\kD$. 
We shall see that the only differential comes from the
analysis of the map $m_{2q}$ in the $A_\infty$ structure on $H^*B\D$.
There is just one ungrading problem, which turns out to be the
only difficult part of the computation, and we do this by interpreting
$\HH^0\kD$ as the centre $Z(\kD)$, and $\HH^1\kD$ as derivations
$\kD\to\kD$ modulo inner derivations.

\begin{theorem}\label{th:d2q-1tau}
In the spectral sequence $\HH^*H^*B\D\Rightarrow\HHinf^*H^*B\D$, whose $E^2$ page
was computed in Theorem~\ref{th:HHHBD},  we have
$d^{2q-1}(\tau)=s^q$.
\end{theorem}
\begin{proof}
We use the description of the Hochschild complex given in Section~\ref{se:HH}.
The element $\tau$ on the $E^2$ page corresponds to the Hochschild
cochain $\tilde\tau\colon [t^i] \mapsto
it^{i-1}$, all other monomials going to zero. Applying the formula~\eqref{eq:AinfinityHH*} for
the differential, we have
\[ (\HHd\tilde\tau)[\underbrace{x,y,\dots,x,y}_{2q}] 
=\tilde\tau(m_{2q}(x,y,\dots,x,y)) = \tilde\tau(t)=1, \] 
and similarly 
\[ (\HHd\tilde\tau)[\underbrace{y,x,\dots,y,x}_{2q}] 
=\tilde\tau(m_{2q}(y,x,\dots,y,x)) = \tilde\tau(t)=1, \]
Since $s[x,y]=s[y,x]=1$, $\HHd\tilde\tau$ takes the same values as $s^q$,
and hence $\HHd\tilde\tau=s^q$. Examining the locations of these terms in
the filtration of the bar complex giving rise to the spectral
sequence, we deduce that this corresponds to the differential
$d^{2q-1}$ taking $\tau$ to $s^q$.
\end{proof}

\begin{theorem}\label{th:HHkD}
The algebra $\HHinf^*H^*B\D\cong \HHinf^*H_*\Omega B\D\twohat \cong \HH^*\kD$
has generators $s$, $t$, $x$, $y$, $u$, $v$, $w_1$, $w_2$, $w_3$ with
$|s|=(0,1,1)$, $|t|=(-2,-q,-q)$, $|x|= (-1,-1,0)$, $|y|=(-1,0,-1)$,
$|u|=|v|=(-1,0,0)$, $|w_1|=(0,q-1,q)$, $|w_2|=(0,q,q-1)$, $|w_3|=(0,q,q)$. These satisfy
the degree zero relations 
\[ w_1^2=w_2^2=w_3^2=w_1w_2=w_1w_3=w_2w_3=sw_1=sw_2=sw_3=s^q=0, \] 
the degree $-1$ relations
\begin{gather*}
vw_1=uw_2=uw_3=vw_3=xs=ys=0,\\ 
us=vs,\quad
xw_2=us^{q-1},\quad yw_1=vs^{q-1},\quad
xw_3=uw_1,\quad 
yw_3=vw_2, 
\end{gather*} 
and the degree $-2$ relations
\[ u^2=v^2=uv=xy=xv=yu=0. \] 
\end{theorem}
\begin{proof}
By the centraliser decomposition\index{centraliser decomposition} 
of Hochschild cohomology of a finite group algebra 
\[ \HH^*\kG\cong\displaystyle\bigoplus_{g\in G}H^*BC_G(g) \] 
(direct sum indexed by conjugacy classes of group elements, see for example
Theorem~2.11.2 of~\cite{Benson:1991b}), it is easy to compute that
\[ \dim_\kk\HH^n\kD=4n+q+3 \]
(see also Proposition~1.3 of Generalov~\cite{Generalov:2010b}).
The $E_2$ page of the spectral sequence 
\[ \HH^*H^*B\D \Rightarrow\HH^*\kD \] 
was computed in Theorem~\ref{th:HHHBD}. In this spectral
sequence, by Theorem~\ref{th:d2q-1tau} 
we have $d^{2q-1}(\tau)=s^q$. 
Let $w_1$, $w_2$ and $w_3$ be representatives in $E^{2q}$ of $x\tau$,  
$y\tau$ and $(u+v)\tau$ respectively. Then in $E^{2q}$ the relations
satisfied by the new generators $w_1$, $w_2$ and $w_3$ are
\[ vw_1=uw_2=uw_3=vw_3=xw_2=yw_1=0,\quad xw_3=uw_1,\quad
  yw_3=vw_2. \]
If $d^{2q-1}$ is the only differential then
the dimensions at
the $E^\infty$ page already match those for $\HH^n\kD$. 
This is because $\HH^0$ is spanned by
$s^i$ ($1\le i\le q$), $w_1$, $w_2$ and $w_3$, $\HH^1$ is spanned by
$u$, $v$, $us^i=vs^i$ ($1\le i\le q$), $x$, $xw_1$, $xw_3$, $y$,
$yw_2$, $yw_3$, and for $n\ge 2$, 
$\HH^n$ is spanned by
$t.\HH^{n-2}$ together with the eight elements $x^n$, $x^nw_1$,
$x^nw_3$, $y^n$, $y^nw_2$, $y^nw_3$,
 $x^{n-1}u$ and $y^{n-1}v$.
So $d^{2q-1}$ is the only differential, it's zero on all generators
except $\tau$, and we have $E^{2q}=E^\infty$.
The $E^\infty$ page is therefore as given in the statement of the
Theorem, but with $xw_2=yw_1=0$. 
It remains to ungrade the relations.

We begin with degree zero, and we work over $\bF_2$ so that the only
possible scalars are zero and one, before extending to $\kk$. The dimension of the algebra 
$\HH^0\kD=Z(\kD)$ is $q+3$, and it is spanned by 
$s^i=(XY+YX)^i$ with $0\le i\le q-1$, together with the elements
$w_1=(YX)^{q-1}Y$, 
$w_2=(XY)^{q-1}X$, 
and $w_3=(XY)^q=(YX)^q$.
These have the required internal degrees, and satisfy the degree zero
relations listed above. In particular, note that $s^q=(XY)^q+(YX)^q$ is equal to
zero and not to $w_3$, even though this has the right degree.

For the degree $-1$ and $-2$ relations, most have nothing
lower in the filtration, in the right internal degree so they ungrade
to the same relations. The exception is the relations $xw_2=yw_1=0$,
which ungrade to give that each of $xw_2$ and $yw_1$ is equal to some $\bF_2$-linear
combination of $us^{q-1}$ and 
$vs^{q-1}$. To determine what linear combinations they are, we look at
the interpretation of $\HH^1\kD$
as derivations of $\kD$. 

The only derivation of degree $(-1,0)$ is $X\mapsto 1$, $Y\mapsto 0$,
so this is the element $x\in \HH^1\kD$. Working in the same way, we have the
following table for the effects of the derivations $x$, $y$, $u$, $v$
on $X$ and $Y$. To distinguish the effects of $u$ and $v$, which have
the same degrees, we use the already
determined relations $xw_3=uw_1$, $yw_3=vw_1$.
\[ \begin{array}{|c|cccc|} \hline
&x&y&u&v\\ \hline
X&1&0&X&0\\
Y&0&1&0&Y\\ \hline
\end{array} \]
Thus the derivation $xw_2$ sends $X\mapsto (XY)^{q-1}X$, $Y\mapsto 0$,
as does $us^{q-1}$, while the derivation $yw_1$ sends $X\mapsto 0$,
$Y\mapsto (YX)^{q-1}Y$, as does $vs^{q-1}$. Thus we have
$xw_2=us^{q-1}$, $yw_1=vs^{q-1}$.
\end{proof}

\begin{remark}\label{rk:HHkD}
The algebra $\HH^*\kD$ was computed in Section~9 of Siegel and
Witherspoon~\cite{Siegel/Witherspoon:1999a}. They chose a different
basis, whose elements are not homogeneous with respect to our
grading, and which complicates their relations. See also
Generalov~\cite{Generalov:2010b}, where the Hochschild cohomology is
computed for algebras in Erdmann's class
III.1\,(c)~\cite{Erdmann:1990a} for any parameter $q$. The degree $-2$
relations depend on the parity of $q$, but are determined already on the $E^2$
page of the spectral sequence; for us, $q$ is always even. 
The same algebras in odd characteristic
are discussed in Generalov~\cite{Generalov:2010a}, where generators of
degree $-3$ and $-4$ also occur in the Hochschild cohomology.
\end{remark}

\section{\texorpdfstring{Groups with dihedral Sylow $2$-subgroups}
{Groups with dihedral Sylow 2-subgroups}}\label{se:dihedral-Sylow}

The computation for groups with a dihedral Sylow $2$-subgroup is
analogous to the dihedral group case described above.
These groups were classified by Gorenstein and
Walter~\cite{Gorenstein/Walter:1962a,Gorenstein/Walter:1965abc}, 
see also Bender and 
Glauberman~\cite{Bender:1981a,Bender/Glauberman:1981a}.
The representation theory
was investigated by Bleher~\cite{Bleher:2009a},
Brauer~\cite{Brauer:1966a,Brauer:1974a},
Cabanes and Picaronny~\cite{Cabanes/Picaronny:1992a},  
Donovan and Freislich~\cite{Donovan/Freislich:1978a},    
Donovan~\cite{Donovan:1979a},   
Erdmann~\cite{Erdmann:1977a,Erdmann:1987a,Erdmann:1990a,Erdmann:1990b},
Erdmann and Michler~\cite{Erdmann/Michler:1977a},
Holm~\cite{Holm:1997a},  
Holm and Zimmermann~\cite{Holm/Zimmermann:2008a},  
Kauer~\cite{Kauer:1998a},  
Koshitani~\cite{Koshitani:1982a},
Koshitani and Lassueur~\cite{Koshitani/Lassueur:2020a},
Landrock~\cite{Landrock:1976a},    
Linckelmann~\cite{Linckelmann:1994a}.
The cohomology rings were investigated by 
Martino and Priddy~\cite{Martino/Priddy:1991a},
Asai~\cite{Asai:1991a},
Asai and Sasaki~\cite{Asai/Sasaki:1993a}, 
Generalov et al.~\cite{
Balashov/Generalov:2002a,
Generalov:2002a, 
Generalov:2004a,
Generalov:2021a,
Generalov/Kosmatov:2005a,
Generalov/Kosmatov:2007a,  
Generalov/Kosmatov:2007b,  
Generalov/Osiyuk:2004a},
and the Hochschild cohomology in
Generalov et al.~\cite{
Filippov/Generalov:2019a,
Generalov:2005a,
Generalov:2010a,
Generalov:2010b,
Generalov:2024a,
Generalov:2024b,
Generalov/Kachalova/Mostovskij:2025a,
Generalov/Kosovskaya:2015a,
Generalov/Kosovskaya:2019a,
Generalov/Romanova:2016a,
Generalov/Romanova/Zilberbord:2016a,
Generalov/Zilberbord/Rogov:2023a},
Holm~\cite{Holm:2004a},
Kachalova~\cite{Kachalova:2025a},
Taillefer~\cite{Taillefer:2019a}.
The homology of $\Omega BG\twohat$ was computed 
by Levi~\cite{Levi:1995a}.

Let $G$ be a finite group with a 
dihedral Sylow $2$-subgroup $\D$ of order $4q$ 
with $q\ge 1$, and let $\kk$ be a field of 
characteristic two. Then by Theorem~1 and Lemmas~2.1 and~3.1 
of Gorenstein and Walter~\cite{Gorenstein/Walter:1965abc}, 
there are three mutually exclusive cases, according to the fusion on
the dihedral groups, described in the three cases below. 
By Theorem~1.1 of Craven and 
Glesser~\cite{Craven/Glesser:2012a}, these also represent the only possible
fusion systems\index{fusion!system} on dihedral $2$-groups.
For the number of isomorphism classes of simple modules in the
principal block of a finite group with dihedral Sylow $2$-subgroups, see
Section~VII of Brauer~\cite{Brauer:1966a}. For  the corresponding
result for blocks with dihedral defect, see Brauer~\cite{Brauer:1974a}.

\begin{case}\label{case:D1}
If $G$ has one class of involutions\index{involutions} then $G/O(G)$
(Notation~\ref{no:gt}) is isomorphic to 
either the alternating group $A_7$\index{alternating group $A_7$} or
a subgroup of $P\Gamma L(2,p^m)$\index{PGammaL@$P\Gamma L(2,p^m)$} 
with $p^m$ a power of an odd prime, 
containing $PSL(2,p^m)$\index{PSL2pm@$PSL(2,p^m)$} with odd index
(Notation~\ref{no:GLetc}). The principal 
block of $\kG$ has three isomorphism classes of simple modules.
\end{case}

\begin{case}\label{case:D2}
If $G$ has two classes of involutions then $G$ has a normal
subgroup of index two, but no normal subgroup of index four. 
In this case, $G/O(G)$ is 
a subgroup of $P\Gamma L(2,p^m)$ with $p^m$ a power of an odd prime,
containing $PGL(2,p^m)$\index{PGL@$PGL(2,p^m)$} 
with odd index. The principal block of $\kG$ has
two isomorphism classes of simple modules. In this case we have $q\ge 2$.
\end{case}

\begin{case}\label{case:D3}
If $G$ has three classes of involutions then $O(G)$ is a normal
complement in $G$ to a Sylow $2$-subgroup $\D$, so that
$G/O(G)\cong \D$ and $H^*BG\cong H^*B\D$. The principal block of $\kG$ is
isomorphic to $\kD$, and
has one isomorphism class of simple modules, namely the trivial module.
\end{case}

\begin{remark}
For  $p$ odd, we have
\begin{align*}
|P\Gamma L(2,p^m)|&=m(p^m-1)p^m(p^m+1),\\
|PGL(2,p^m)|&=(p^m-1)p^m(p^m+1),\\
|PSL(2,p^m)|&= (p^m-1)p^m(p^m+1)/2,
\end{align*} 
and $PSL(2,p^m)$ is simple for $p\ge 5$.
\end{remark}

\begin{proposition}\label{pr:BGtwohatD}
Suppose that $G$ has a dihedral Sylow $2$-subgroup $\D$. Then the homotopy
type of $BG\twohat$ is determined by $|\D|$ and the number of conjugacy
classes of involutions.
\end{proposition}
\begin{proof}
This follows from Theorem~\ref{th:Oliver} and the case analysis
described above.
\end{proof}

We shall deal with Cases~\ref{case:D1} and
\ref{case:D2} in turn. Case~\ref{case:D1} is the most interesting, because
this is the case where $G$ has no subgroup of index two, so 
$\Omega BG\twohat$ is connected. Case~\ref{case:D2}
is computationally quite similar, but $\Omega BG\twohat$ has two
connected components, and so we give the details anyway
for completeness. Case~\ref{case:D3} has already been dealt with in
Sections~\ref{se:dihedral}--\ref{se:HH*kD}, because $\Omega BG \twohat$ 
is homotopy equivalent to $\D$. Nonetheless, the Eilenberg--Moore
spectral sequence has an interesting differential in this case, as we
saw in Section~\ref{se:OmegaBDtwohat}.

We end this section with a table of the various cases of algebras of
dihedral type in characteristic two, in Erdmann's classification.\medskip
\begin{center}
\begin{tabular}{|c|c|c|c|c|} \hline 
Erdmann~\cite{Erdmann:1990a}&Case&Group&$H^*$&$\HH^*$ \\ \hline
III.I(a)&---&---&\cite{Generalov:2021a}&\\ 
III.I(b)&\ref{case:D3}&fours group&&
\cite{Holm:1996a,Cibils/Solotar:1997a}\\
III.I(b$'$)&---&---&&\\
III.I(c)&\ref{case:D3}&dihedral&\cite{Munkholm:1969a}&
\cite{Generalov:2010b,Siegel/Witherspoon:1999a}\\
III.I(c$'$)&---&---&&\\
D(2$\mathcal A$)&\ref{case:D2}&$PGL(2,q),\ q\equiv 1\pmod{4}$&
\cite{Martino/Priddy:1991a,Asai/Sasaki:1993a,Generalov/Osiyuk:2004a}&\\
D(2$\mathcal B$)&\ref{case:D2}&$PGL(2,q),\ q\equiv 3\pmod{4}$&
\cite{Martino/Priddy:1991a,Asai/Sasaki:1993a,Generalov:2004a}&
\cite{Generalov/Kosovskaya:2015a,Generalov/Kosovskaya:2019a,Generalov/Romanova:2016a}\\
D(3$\mathcal A$)$_1$&\ref{case:D1}&$PSL(2,q),\ q \equiv 1\pmod{4}$&
\cite{Martino/Priddy:1991a,Asai/Sasaki:1993a,Balashov/Generalov:1999a}&
\cite{Generalov:2005a}\\
D(3$\mathcal A$)$_2$&---&---&&\\
D(3$\mathcal B$)$_1$&\ref{case:D1}&Alternating group $A_7$&
\cite{Martino/Priddy:1991a,Asai/Sasaki:1993a,Generalov:2002a}&
\cite{Generalov:2005a}\\
D(3$\mathcal B$)$_2$&---&---&&\\
D(3$\mathcal D$)$_1$&---&---&&
\cite{Generalov:2005a}\\
D(3$\mathcal D$)$_2$&---&---&&\\
D(3$\mathcal K$)&\ref{case:D1}&$PSL(2,q),\ q\equiv 3\pmod{4}$&
\cite{Martino/Priddy:1991a,Asai/Sasaki:1993a,Balashov/Generalov:2002a}&
\cite{Generalov:2005a,Generalov/Kachalova/Mostovskij:2025a}\\
D(3$\mathcal L$)&---&---&\cite{Generalov/Kosmatov:2005a}&\\
D(3$\mathcal Q$)&---&---&\cite{Generalov/Kosmatov:2007b}&\\
D(3$\mathcal R$)&---&---&\cite{Generalov/Kosmatov:2007a}&
\cite{Filippov/Generalov:2019a}
\\ \hline
\end{tabular}
\end{center}

\section{\texorpdfstring
{Loops on $BG\twohat$: one class of involutions}
{Loops on BGٛ₂: one class of involutions}}\label{se:D1}

In this section we begin the examination of 
Case~\ref{case:D1}. This is the case
where $G$ has a dihedral Sylow $2$-subgroup $\D$, and one conjugacy
class of involution. In this case, $G$ has no subgroup of index two,
and it has three isomorphism classes of simple modules in the
principal block. 

\begin{remark}
As we have already mentioned, Theorem~\ref{th:Oliver} shows that up to
quasi-isomorphism, the $A_\infty$ algebra
$H^*BG$ only depends on the fusion, and according
to Proposition~\ref{pr:BGtwohatD}, for dihedral Sylow $2$-subgroups
this only depends on $|\D|$ and the number of simple modules.
In fact more is true. 
Linckelmann, in Theorem~1 of~\cite{Linckelmann:1994a}, has proved that all
blocks of finite groups with 
dihedral defect groups\index{defect group!dihedral} of a given order,
and three isomorphism classes of simple modules, are derived
equivalent. Explicit derived equivalences\index{derived!equivalence} 
are described in that
paper, and in the case of principal blocks, it can be checked that the
derived equivalence may be chosen to take the trivial module to the
trivial module. 
The endomorphism DGA of the trivial module is a derived invariant up 
to quasi-isomorphism, and is also quasi-isomorphic to the $A_\infty$
algebra $H^*BG$. 
\end{remark}

Let $G$ be a group with dihedral Sylow $2$-subgroup $\D$ of order $4q$,
$q\ge 1$, and
one conjugacy class of
involutions. By Proposition~\ref{pr:BGtwohatD}, 
for the purpose of studying $BG\twohat$, we may assume
that $G=PSL(2,p)$ for a suitable prime $p\equiv 1 \pmod{4}$.
In this case, by Theorem~2 of Erdmann~\cite{Erdmann:1977a} and
Lemma~6.1 of Erdmann~\cite{Erdmann:1987a},
 the principal block $B_0$ of $\kG$ has
three simple modules, $\kk$, $\MM$ and $\NN$, whose
$\Ext^1$ quiver\index{Ext@$\Ext^1$ quiver}\index{quiver} is as follows:
\[ \xymatrix{\MM\ar@/^/[r]^{e_2} & \kk \ar@/^/[r]^{e_3}\ar@/^/[l]^{e_1} & \NN\ar@/^/[l]^{e_4}}. \]
The relations are 
\[ e_1e_2=0,\qquad e_3e_4=0,\qquad
(e_4e_3e_2e_1)^q=(e_2e_1e_4e_3)^q. \]
We put an internal grading on the basic algebra 
in this case by assigning degree $(\half,0)$ to $e_1$ and $e_2$
and degree $(0,\half)$ to $e_3$ and $e_4$. Thus we assign degree
$\half(n_1,n_2)$ to a path
involving $n_1$ arrows of type $e_1$ or $e_2$, and $n_2$ arrows of
type $e_3$ or $e_4$. 
This choice is appropriate, because the internal grading
it induces in cohomology is compatible with restriction to the Sylow
$2$-subgroup, as we shall see below.

\begin{remark}
It is not clear \emph{a priori} that there exists a grading on the
principal block compatible with the restriction map in cohomology. 
This explains the need for the computation. For a further
discussion of gradings in this context, see Bogdani\'c~\cite{Bogdanic:2011a}.
\end{remark}

Let $\bar e_1$ be the element of $\Hom_B(P_\MM,P_\kk)$ dual to
$e_1$, and so on. Then
the minimal resolution of $\kk$ as a $\kG$-module takes the form
\begin{multline*}  
\cdots \xrightarrow{\quad}
P_\MM \oplus P_\kk \oplus P_\kk \oplus P_\NN 
\xrightarrow{\left(\begin{smallmatrix}
\bar e_1 &
v& 0 & 0 \\ 0 &
\bar e_1\bar e_2 & \bar e_3\bar e_4 & 0 \\
0 & 0 & u& \bar e_3
\end{smallmatrix}\right)}
P_\kk \oplus P_\kk \oplus P_\kk 
\xrightarrow{\left(\begin{smallmatrix}
\bar e_1\bar e_2 & 
v&0 \\ 0 & u&\bar e_3\bar e_4 
\end{smallmatrix}\right)}
P_\kk \oplus P_\kk \\
\xrightarrow{\left(\begin{smallmatrix}
\bar e_2v&0 \\
\bar e_1\bar e_2 & \bar e_3\bar e_4\\
0&\bar e_4u
\end{smallmatrix}\right)}
 P_\MM \oplus P_\kk \oplus P_\NN 
\xrightarrow{\left(\begin{smallmatrix}
\bar e_1 & v& 0 \\
0 & u & \bar e_3
\end{smallmatrix}\right)} 
P_\kk \oplus P_\kk \xrightarrow{\left(\begin{smallmatrix}
\bar e_1\bar e_2 & \bar e_3\bar e_4
\end{smallmatrix}\right)}
P_\kk\xrightarrow{\left(\begin{smallmatrix}
\bar e_2v\\ \bar e_4u
\end{smallmatrix}\right)} 
P_\MM \oplus P_\NN \xrightarrow{(\bar e_1,\bar e_3)} P_\kk 
\end{multline*}
where $u=\bar e_1\bar e_2(\bar e_3\bar e_4\bar e_1\bar e_2)^{q-1}$ and  
$v=\bar e_3\bar e_4(\bar e_1\bar e_2\bar e_3\bar e_4)^{q-1}$.  
This is the total complex of the following double complex that
displays the gradings on cohomology.\vspace{-4mm}
\[ \xymatrix@=6mm{
&&&P_\kk\ar[d]^{\bar e_3\bar e_4} & \vdots\\
&&P_\NN\ar[d]^{\bar e_3} & 
P_\kk\ar[l]_{\bar e_4u}\ar[d]^{\bar e_3\bar e_4} 
& P_\kk \ar[l]_u\ar[d]^(.45){\bar e_3\bar e_4} \\
&&P_\kk\ar[d]^{\bar e_3\bar e_4} & 
P_\kk\ar[l]_u\ar[d]^(.45){\bar e_3\bar e_4} & 
P_\kk\ar[l]_u\ar[d]^v & 
P_\kk\ar[l]_{\bar e_1\bar e_2}\ar[d]^v & \cdots\\
&P_\NN\ar[d]^{\bar e_3} & 
P_\kk\ar[l]_{\bar e_4u}\ar[d]^(.45){\bar e_3\bar e_4}
& P_\kk\ar[l]_u\ar[d]^v & 
P_\kk\ar[l]_{\bar e_1\bar e_2}\ar[d]^v 
& P_\kk\ar[l]_{\bar e_1\bar e_2}\ar[d]^{\bar e_2v} 
& P_\kk\ar[l]_{\bar e_1\bar e_2}\\
&P_\kk\ar[d]^(.45){\bar e_3\bar e_4} & 
P_\kk\ar[l]_u\ar[d]^v & 
P_\kk\ar[l]_{\bar e_1\bar e_2}\ar[d]^{\bar e_2v}
& P_\kk\ar[l]_{\bar e_1\bar e_2} & P_\MM\ar[l]_{\bar e_1} \\
P_\NN \ar[d]^{\bar e_3} & 
P_\kk \ar[l]_{\bar e_4u} \ar[d]^{\bar e_2v} & 
P_\kk \ar[l]_{\bar e_1\bar e_2} &
P_\MM \ar[l]_{\bar e_1} \\
P_\kk & \ar[l]_{\bar e_1} P_\MM} \]
In this picture, the blank spots represent the zero module. Each
column and each row is a complex, and the vertical and horizontal maps
anticommute, so that the sum of the vertical and horizontal maps
squares to zero. For further information on double complexes for
minimal resolutions, see Section~7 of Benson and
Carlson~\cite{Benson/Carlson:1987b}. 

So with this grading, if $q\ge 2$, the cohomology ring is given by
\begin{equation}\label{eq:HPSL2p}
H^*(BG,\kk)=\kk[\xi,\eta,t]/(\xi\eta),\quad
 |\xi|=-(3,q+1,q),\ 
|\eta|=-(3,q,q+1),\  |t|=-(2,q,q)
\end{equation}
while if $q=1$, we have
\begin{equation}\label{eq:HPSL2p2}
H^*(BG,\kk)=\kk[\xi,\eta,t]/(\xi\eta+t^3),\quad
 |\xi|=-(3,2,1),\ 
|\eta|=-(3,1,2),\  |t|=-(2,1,1). 
\end{equation}

\begin{remark}
If $q=1$ then $\D$ is abelian, so by Theorem~\ref{th:abelian},
$BG\twohat\simeq BA_4\,\twohat$.
In this case, we shall see in Theorem~\ref{th:elemab2} that the $A_\infty$
structure on $H^*BG$ is intrinsically formal.\index{intrinsically formal}
\end{remark}

For $q\ge 2$ we have Massey products\index{Massey product}
\begin{equation}\label{eq:Massey-PSL2p}
\langle  \xi,\eta,\dots,\xi,\eta\rangle = 
\langle \eta,\xi,\dots,\eta,\xi\rangle =  
t^{2q+1}. 
\end{equation}
In both expressions the arguments $\xi$ and $\eta$ alternate, and
there are $2q$ of them. These Massey products are only well defined
up to adding elements of the ideal generated by $\xi$ and $\eta$, but
taking the grading into account, they are well defined with no ambiguity.

\begin{theorem}\label{th:HOmegaBG-D1}
In Case~\ref{case:D1} we have
\[ H_*\Omega BG\twohat = \Lambda(\tau) \otimes \kk\langle \alpha,\beta
\mid \alpha^2=0,\ \beta^2=0\rangle. \]
with 
\[ |\tau|=(1,q,q),\qquad |\alpha|=(2,q+1,q),\qquad |\beta|=(2,q,q+1). \]
In homological degree $4n$ we have monomials $(\alpha\beta)^n$ and
$(\beta\alpha)^n$, in degree $4n+2$ we have monomials
$(\alpha\beta)^n\alpha$ and $(\beta\alpha)^n\beta$, 
and in odd degrees we have $\tau$ times all of these. These monomials
form a basis in these degrees.
\end{theorem}
\begin{proof}
For $q\ge 1$, 
by~\eqref{eq:HPSL2p}, \eqref{eq:HPSL2p2} and Theorem~\ref{th:ExtRkk} we have
\[ \Ext^{*,*}_{H^*BG}(\kk,\kk)=\Lambda(\tau) \otimes 
\kk\langle\alpha,\beta\mid \alpha^2=0,\ \beta^2=0\rangle \]
where the generators have degrees 
\[ |\tau|=(-1,2,q,q), \qquad
|\alpha|=(-1,3,q+1,q), \qquad
|\beta|=(-1,3,q,q+1).  \]
The four degrees are first homological, then internal to
$H^*BG$, and finally the two gradings internal to $\kG$.
The elements $\tau$, $\alpha$ and $\beta$ come from the
generators $t$, $\eta$ and $\xi$, while the element
$s=\alpha\beta+\beta\alpha$ in degree $(-2,6,2q+1,2q+1)$
is the Eisenbud operator\index{Eisenbud operator} 
for the relation $\xi\eta=0$ (or $\xi\eta=t^3$) in $H^*BG$.

The Eilenberg--Moore spectral sequence~\eqref{eq:Cotor}%
\index{Eilenberg--Moore spectral sequence} converging to  
$H_*\Omega BG\twohat$ has this as its $E_2$ page.
There is no room for non-zero differentials, and there are no ungrading
problems, so we have
$E_2=E_\infty=H_*\Omega BG\twohat$. 
\end{proof}

\begin{remark}\index{errors}
Proposition~II.4.1.5 of Levi~\cite{Levi:1995a} gets the correct
additive structure for $H_*\Omega BG\twohat$ but it is incorrectly claimed there
that the ring structure is a polynomial tensor exterior algebra.
\end{remark}

The algebras $H^*BG$ and $H_*\Omega BG\twohat$, ignoring higher structure, are Koszul
dual to each other (see Section~\ref{se:KoszulDual}). This will play a role in the computation of
Hochschild cohomology.

\begin{lemma}\label{le:m=0-D1}
In Case~\ref{case:D1},
for any $A_\infty$ structure on $H^*BG$ that preserves internal
degrees, we have $m_i=0$ unless $i-2$ is divisible by 
$2q-2$. In particular, for $2<i<2q$ we have $m_i=0$.
\end{lemma}
\begin{proof}
The proof is the same as the proof of Lemma~\ref{le:m=0}.
\end{proof}

\begin{proposition}\label{pr:HHHBG-D1}
In Case~\ref{case:D1} with $q\ge 2$,
the Hochschild cohomology\index{Hochschild cohomology} $\HH^*H^*BG$
has generators $s$, $t$, $\tau$, $\xi$, $\eta$, $u$, $v$
with
\begin{align*}
|s|&=(-2,6,2q+1,2q+1) \\
|t|&= -(0,2,q,q) &
|\tau|&= (-1,2,q,q) \\
|\xi|&= -(0,3,q+1,q)&
|\eta|&= -(0,3,q,q+1) \\
|u|&= -(1,0,0,0) &
|v|&= -(1,0,0,0).
\end{align*}
The relations are given by $u^2=v^2=uv=\tau^2=0$, 
$\eta u=\xi v=0$, $\xi s=\eta s=0$, and $us=vs$. The non-zero monomials and their degrees are as
follows, with $i_1,i_2\ge 0$, $\ep_1,\ep_2\in\{0,1\}$.
\begin{align*}
|s^{i_1}t^{i_2}\tau^{\ep_1}u^{\ep_2}|&
=(-2i_1-\ep_1-\ep_2,6i_1-2i_2+2\ep_1,i_1+q(2i_1-i_2+\ep_1),i_1+q(2i_1-i_2+\ep_1)), \\
|s^{i_1}t^{i_2}\tau^{\ep_1}v^{\ep_2}|&
=(-2i_1-\ep_1-\ep_2,6i_1-2i_2+2\ep_1,i_1+q(2i_1-i_2+\ep_1),i_1+q(2i_1-i_2+\ep_1)), \\
|\xi^{i_1}t^{i_2}\tau^{\ep_1}u^{\ep_2}|&
=(-\ep_1-\ep_2,-3i_1-2i_2+2\ep_1,-i_1+q(-i_1-i_2+\ep_1),q(-i_1-i_2+\ep_1)) \\
|\eta^{i_1}t^{i_2}\tau^{\ep_1}v^{\ep_2}|&
=(-\ep_1-\ep_2,-3i_1-2i_2+2\ep_1,q(-i_1-i_2+\ep_1),-i_1+q(-i_1-i_2+\ep_1)) 
\end{align*}
There is only one monomial
in degree $(-i,i-2,0,0)$ with $i>2$, namely $s^qt^{2q+1}$, with
\[ |s^qt^{2q+1}|=(-2q,2q-2,0,0). \]
\end{proposition}
\begin{proof}
As in Theorem~\ref{th:HHHBD}, we use the approach of
Theorems~\ref{th:HHR} and~\ref{th:Negron}. 
Thus $\HH^*H^*BG$ is the homology of
the complex 
\[ (H^*BG \otimes H_*\Omega BG\twohat,d). \]
which we can read off from~\eqref{eq:HPSL2p} and Theorem~\ref{th:HOmegaBG-D1}.
The generators $t$, $\xi$ and $\eta$ are in homological degree
zero, the generators $\tau$, $\alpha$ and $\beta$ are in homological
degree $-1$, and the differential is given by $d=[e,-]$ where
$e = t \otimes \tau + \xi \otimes \alpha + \eta\otimes \beta$. 
Thus setting $s=\alpha\beta+\beta\alpha$, we have
$d(t)=d(\xi)=d(\eta)=d(\tau)=0$, $d(\alpha)=\eta s$,  
$d(\beta)=\xi s$. The generators and relations for the homology
of this complex are therefore as given, with $u=\xi\alpha$ and
$v=\eta\beta$.

For the last statement, the computation is similar to the
corresponding part of the proof of Theorem~\ref{th:HHHBD}.
\end{proof}

\begin{theorem}\label{th:AinftyHBG-D1}
In Case~\ref{case:D1} with $q\ge 2$,
the $A_\infty$ structure on $H^*BG=\kk[\xi,\eta,t]/(\xi\eta)$,
see~\eqref{eq:HPSL2p}, is given as follows. 
The $m_n$ are $\kk[t]$-multilinear maps with $m_n=0$ 
for $n$ not congruent to $2$ modulo $2q-2$, and for $i,j\ge 1$
\[ m_{2q}(\xi^i,\eta,\xi,\eta,\dots,\xi,\eta^j)=
m_{2q}(\eta^j,\xi,\eta,\xi,\dots,\eta,\xi^i)=\xi^{i-1}\eta^{j-1}t^{2q+1} \]
where the arguments alternate between $\xi$ and $\eta$, and the right hand
side is zero unless either $i=1$ or $j=1$; $m_{2q}$ is zero on all
other tuples of monomials not involving $t$.
The maps $m_{\ell(2q-2)+2}$ with $\ell>1$ similarly vanish on all tuples of
monomials not involving $t$, except the ones which look as above,
but for some choice of indices in the tuple:
\begin{multline*} 
1\le e_1\le e_2\le \cdots \le e_{\ell-1} 
< e_{\ell-1}+(2q-2)+1   
\le e_{\ell-2}+2(2q-2)+1 \\
\le \cdots \le e_1+(\ell-1)(2q-2)+1\le \ell(2q-2)+2. 
\end{multline*}
the exponents on the terms are increased by one
(or correspondingly more if an index is repeated).
The value on these tuples is $\xi^{i-1}\eta^{j-1}t^{\ell(2q+1)}$.
Thus 
\[
m_{\ell(2q-2)+2}(x^{i+\alpha_1},y^{\alpha_2},x^{\alpha_3},\dots,
x^{\alpha_{\ell(2q-2)+1}},y^{j+\alpha_{\ell(2q-2)+2}}) =
  x^{i-1}y^{j-1}t^{\ell(2q+1)} \]
where each $\alpha_\sigma$ is one plus the number of indices in the list above that are equal
to $\sigma$.
\end{theorem}
\begin{proof}
The proof is the same as the proof of Theorem~\ref{th:AinftyHBG}, but
using Lemma~\ref{le:m=0-D1} and Proposition~\ref{pr:HHHBG-D1} in
place of Lemma~\ref{le:m=0} and Theorem~\ref{th:HHHBD}.
\end{proof}

Recall from Theorem~\ref{th:HOmegaBG-D1} that in Case~\ref{case:D1} we have
\begin{align*} 
H_*\Omega BG\twohat &= \Lambda(\tau) \otimes \kk\langle \alpha,\beta
\mid \alpha^2=0,\ \beta^2=0\rangle \\
 |\tau|=(1,q,q),\qquad& |\alpha|=(2,q+1,q),\qquad |\beta|=(2,q,q+1).
\end{align*}
We now turn to the computation of the $A_\infty$ structure on this.
It turns out to be easier to describe than the $A_\infty$ structure on $H^*BG$.

\begin{lemma}\label{le:m=0v2-D1}
In Case~\ref{case:D1}, for any $A_\infty$ structure on
$H_*\Omega BG\twohat$ that preserves 
internal degrees, we have $m_i=0$ unless $i - 2$ is
divisible by $2q-1$. In particular,  for $2<i<2q+1$ we have $m_i=0$. 
\end{lemma}
\begin{proof}
The proof is similar to the proof of Lemma~\ref{le:m=0}.
Looking at the degrees of the generators $\tau$, $\alpha$ and $\beta$,
for any monomial $\zeta$ in $H_*\Omega BG\twohat$ we have $a \equiv
b+c\pmod{2q-1}$. So for any $i$-tuple $(\zeta_1,\dots,\zeta_i)$,
the degree of $m_i(\zeta_1,\dots,\zeta_i)$ satisfies
$a\equiv b+c+i-2 \equiv 0 \pmod{2q-1}$. So for $m_i(\zeta_1,\dots,\zeta_i)$ to be
non-zero we must have $i -2 \equiv 0 \pmod{2q-1}$.
\end{proof}

\begin{proposition}\label{pr:HHHOmegaBG-D1}
In Case~\ref{case:D1}, the Hochschild
cohomology\index{Hochschild cohomology}  
$\HH^*H_*\Omega BG\twohat$ has 
generators $s$, $t$, $\tau$, $\xi$, $\eta$, $u$, and $v$ in degrees
\begin{align*}
|s|&=(0,4,2q+1,2q+1),\\
|t|&=-(1,1,q,q),&
|\tau|&=(0,1,q,q), \\
|\xi|&=-(1,2,q+1,q),&
|\eta|&=-(1,2,q,q+1), \\
|u|&=-(1,0,0,0),&
|v|&=-(1,0,0,0).
\end{align*}
The relations are given by $\xi\eta=0$, $u^2=v^2=uv=\tau^2=0$, 
$\eta u=\xi v = 0$, $\xi s = \eta s =0$, and $us=vs$.
The non-zero monomials and their degrees are given as follows, with
$i_1,i_2\ge 0$, $\ep_1,\ep_2\in\{0,1\}$.
\begin{align*}
|s^{i_1}t^{i_2}\tau^{\ep_1}u^{\ep_2}|
&=(-i_2-\ep_2,4i_1-i_2+\ep_1,(2i_1-i_2+\ep_1)q+i_1,(2i_1-i_2+\ep_1)q+i_1),\\
|s^{i_1}t^{i_2}\tau^{\ep_1}v^{\ep_2}|
&=(-i_2-\ep_2,4i_1-i_2+\ep_1,(2i_1-i_2+\ep_1)q+i_1,(2i_1-i_2+\ep_1)q+i_1),\\
|\xi^{i_1}t^{i_2}\tau^{\ep_1}u^{\ep_2}|&
=(-i_1-i_2-\ep_2,-2i_1-i_2+\ep_1,-i_1+(-i_1-i_2+\ep_1)q,(-i_1-i_2+\ep_1)q),\\
|\eta^{i_1}t^{i_2}\tau^{\ep_1}v^{\ep_2}|&
=(-i_1-i_2-\ep_2,-2i_1-i_2+\ep_1,(-i_1-i_2+\ep_1)q,-i_1+(-i_1-i_2+\ep_1)q).
\end{align*}
Thus there is only one monomial with degree $(-i,i-2,0,0)$ with $i>2$, namely
\[ |s^qt^{2q+1}|=(-2q-1,2q-1,0,0). \]
\end{proposition}
\begin{proof}
Again, as in Theorem~\ref{th:HHHBD}, we use the approach of 
Theorems~\ref{th:HHR} and~\ref{th:Negron}. 
As in the proof of Proposition~\ref{pr:HHHBG-D1}, 
$\HH^*H_*\Omega BG\twohat$ is
the homology of the complex 
\[ (H_*\Omega BG\twohat \otimes H^*BG,d), \]
but this time the generators $\tau$, $\alpha$ and $\beta$ of 
$H_*\Omega BG\twohat$ are in homological degree zero, 
the generators $t$, $\xi$, $\eta$ of $H^*BG$ are in homological degree $-1$, and
the differential is given by $d=[e,-]$ where $e=\tau\otimes t +
\alpha \otimes \xi + \beta\otimes \eta$. 
So the answer is the same as in Proposition~\ref{pr:HHHBG-D1} but with the
degrees changed.

For the last statement, the computation is again similar to the
corresponding part of the proof of Theorem~\ref{th:HHHBD}.
\end{proof}

\begin{theorem}\label{th:AinftyHOmegaBG-D1}
In Case~\ref{case:D1}, 
the $A_\infty$ structure on $H_*\Omega BG\twohat$ is determined by
\[ m_{2q+1}(\tau,\tau,\dots,\tau) = s^q, \]
where $s=\alpha\beta+\beta\alpha$.
This implies that
\begin{equation}\label{eq:m2q+1}
m_{2q+1}(f_1(\alpha,\beta)\tau,f_2(\alpha,\beta)\tau,\dots,f_{2q+1}(\alpha,\beta)\tau)
= f_1(\alpha,\beta)\dots
f_{2q+1}(\alpha,\beta)s^q, 
\end{equation}
and all $m_n$ for $n>2$ on all other $n$-tuples of monomials give
zero.
\end{theorem}
\begin{proof}
By Lemma~\ref{le:m=0v2-D1}, we have $m_n=0$ for $2<n<2q+1$. So
in order to determine $m_{2q+1}$, we invoke
Proposition~\ref{pr:HH}. This shows that $m_{2q+1}$ has to be a
Hochschild cocycle, well defined up to adding Hochschild coboundaries.
By Proposition~\ref{pr:HHHOmegaBG-D1}, the dimension of $\HH^*H_*\Omega
BG\twohat$ is one dimensional in degree $(-2q-1, 2q-1,0,0)$. A
representative for a non-zero cohomology class is given by \eqref{eq:m2q+1}.
It is easy to check that this is a cocycle but not a coboundary. So by
rescaling $\tau$ if necessary (or by working over $\bF_2$) we may
assume that either $m_{2q+1}$ is either zero or as given in the
theorem. In both cases we can check that the Gerstenhaber circle
product $m_{2q+1}\circ m_{2q+1}$ is
equal to the zero cocycle in degree $-4q$. 

As in the proof of Theorem~\ref{th:AinftyHBG}, we can rewrite
Equation~\ref{eq:Ainfty} in degree $-4q$ as
\[ \HHd m_{4q}=m_{2q+1}\circ m_{2q+1}, \]
which as we just saw, is zero.
Now by Proposition~\ref{pr:HHHOmegaBG-D1} again, $\HH^*H_*\Omega BG\twohat$ is
zero in degree $(-4q,4q-2,0,0)$. So $m_{4q}$ is a Hochschild
coboundary, and we can therefore take $m_{4q}=0$, as it is only well
defined modulo Hochschild coboundaries. 
At this point, for $\ell > 2$, the equation we obtain for
$m_{\ell(2q-1)+2}$ is $\HHd m_{\ell(2q-1)+2}=0$. Again,
$\HH^*H_*\Omega BG\twohat$ is zero in degree
$(-\ell(2q-1)-2,\ell(2q-1),0,0)$, and so we may take
$m_{\ell(2q-1)+2}=0$.

This argument shows that there are two possibilities for the
$A_\infty$ structure up to isomorphism, namely the one given and the
formal\index{formal $A_\infty$ algebra} one with $m_n=0$ for all
$n>2$. 
Now as algebras without higher structure, $H^*BG$ and $H_*\Omega
BG\phat$ are Koszul algebras that are Koszul dual (see
Section~\ref{se:KoszulDual}). So if $H_*\Omega BG\phat$ were formal, then
$H^*BG$ would also be formal, which it is not
by~\eqref{eq:Massey-PSL2p}. So the $A_\infty$ structure is as given. 
\end{proof}

\begin{remark}\label{rk:D1-HHCBG}
In the spectral sequence~\eqref{eq:HHHA} 
\[ \HH^*H_*\Omega BG\twohat \Rightarrow
\HHinf^*H_*\Omega BG\twohat, \] 
we have
$d^{2q}(\tau)=s^qt^{2q}$. This implies that after inverting $s$ (we
discuss this later), we have
\[ \HHinf^*H_*\Omega BG\twohat[s^{-1}] =
  \kk[s,s^{-1}][u,v,t]/(u^2,v^2,uv,t^{2q}). \]
Since by Theorem~\ref{th:HHOmega} we have $\HHinf^*H^*BG\cong
\HHinf^*H_*\Omega BG\twohat$,  
this is also the structure of $\HHinf^*H^*BG[s^{-1}]$.
\end{remark}

\section{A differential graded model}\label{se:Q}%
\index{differential!graded algebra}\index{DG!algebra}

Throughout this section, we work in Case~\ref{case:D1},
where $G$ has dihedral Sylow $2$-subgroups and one conjugacy class of involutions.
As in Section~3 of Benson and Greenlees~\cite{Benson/Greenlees:2023a},
we produce a differential graded 
model $Q$ for the $A_\infty$ algebra $H_*\Omega BG\twohat$. The proofs are similar
to the ones in that paper, but we spell out the details because there
are some differences. One is that we are in characteristic two,
so we don't need to be careful about signs; another is that a
polynomial ring in one variable has been replaced by the
noncommutative ring $\kk\langle\alpha,\beta\rangle/(\alpha^2,\beta^2)$. 

Recall from Theorems~\ref{th:HOmegaBG-D1} 
and~\ref{th:AinftyHOmegaBG-D1} that 
\[ H_*\Omega BG\twohat \cong \Lambda(\tau) \otimes \kk\langle
  \alpha,\beta\rangle/(\alpha^2,\beta^2) \]
with $m_{2q+1}$ determined by 
$m_{2q+1}(\tau,\dots,\tau)=s^q$, where $s=\alpha\beta+\beta\alpha$, and
with all other $m_i$ zero for $i>2$.

The generators of $Q$ are
elements $\tau_1,\dots,\tau_{2q},\alpha,\beta$, where $\tau_1$ will
eventually be seen to correspond to the element
$\tau\in H_*\Omega BG\twohat$. The relations and
differential are as follows:
\begin{align*}
\alpha\tau_i&=\tau_i\alpha\\
\beta\tau_i&=\tau_i\beta\\
\alpha^2&=\beta^2=0\\
d\alpha&=d\beta=0\\
\sum_{j+k=i}\tau_j\tau_k &=\begin{cases}
d\tau_i & 1\le i\le 2q \\
s^q& i=2q+1 \\
0 & 2q+2\le i\le 4q.
\end{cases}
\end{align*}
where $s=\alpha\beta+\beta\alpha$.
The antipode is the algebra anti-automorphism given by
$S(\tau_i)=\tau_i$, $S(\alpha)=\alpha$, $S(\beta)=\beta$
(we are in characteristic two, so there are no signs), and the
comultiplication is given by making the generators primitive:
\begin{equation}\label{eq:Delta}
 \Delta(\tau_i)=\tau_i\otimes 1  + 1 \otimes \tau_i,\qquad
\Delta(\alpha)=\alpha\otimes 1 + 1 \otimes \alpha,\qquad
\Delta(\beta)=\beta\otimes 1 + 1 \otimes\beta. 
\end{equation}
Then $s=[\alpha,\beta]$ is also primitive, hence so is $s^q$ because
$q$ is a power of the characteristic.
The degrees are given by $|\tau_i|=(2i-1,iq,iq)$,
$|\alpha|=(2,q+1,q)$, $|\beta|=(2,q,q+1)$, and
$|s|=(4,2q+1,2q+1)$.
We shall see that this algebra $Q$ is quasi-isomorphic to 
$C_*\Omega BG\twohat$.

\begin{example}
If $q=1$, the algebra $Q$ is generated by $\tau_1,\tau_2,\alpha,\beta$
with
\begin{align*}
d(\alpha)&=0 & \alpha^2&=0& \tau_1\tau_2+\tau_2\tau_1&=s=\alpha\beta+\beta\alpha\\
d(\beta)&=0 & \beta^2&=0 &\tau_2^2&=0\\
d(\tau_1)&=0 &\alpha \tau_i&=\tau_i\alpha \\
d(\tau_2)&=\tau_1^2 & \beta\tau_i&=\tau_i\beta 
\end{align*}
with $|\tau_1|=(1,1,1)$, $|\tau_2|=(3,2,2)$, $|\alpha|=(2,2,1)$,
$|\beta|=(2,1,2)$ and $|s|=(4,3,3)$.
\end{example}

\begin{lemma} 
In the algebra $Q$, every element has a unique expression of the form  
\[ f(\tau_1,\dots,\tau_{2q-1}) +
  \tau_{2q}g(\tau_1,\dots,\tau_{2q-1}) \] 
where $f$ and $g$ are polynomials in the non-commuting variables
$\tau_1,\dots,\tau_{2q-1}$ 
with coefficients in $\kk\langle\alpha,\beta\rangle/(\alpha^2,\beta^2)$.  
\end{lemma}
\begin{proof}
The algebra relations (ignoring the differential) say first that the
elements $\tau_1,\dots,\tau_{2q}$ commute with $\alpha$ and $\beta$;
and the remaining relations can be rewritten in the form
\[ \tau_i\tau_{2q}=\tau_{2q}\phi_i(\tau_1,\dots,\tau_{2q-1})\]
with $1\le i \le 2q$ (note that $\phi_{2q}=0$). Thus all occurrences
of $\tau_{2q}$ may be moved to the beginning, and
$\tau_{2q}^2=0$. There are no relations among $\tau_1,\dots,\tau_{2q-1}$.
\end{proof}

\begin{definition}\label{def:standard}
We shall refer to a monomial in $\tau_1,\dots,\tau_{2q-1}$, or
$\tau_{2q}$ times such a monomial, as a \emph{standard monomial} in
the variables $\tau_1,\dots,\tau_{2q}$. By the lemma, these monomials
form a basis for $Q$ over $\kk\langle
\alpha,\beta\rangle/(\alpha^2,\beta^2)$.
\end{definition}

\begin{lemma}
In the algebra $Q$, we have $d^2=0$.
\end{lemma}
\begin{proof}
The differential is given by
\[ d(f+\tau_{2q}g)=(df+(\tau_1\tau_{2q-1}+\dots+\tau_{2q-1}\tau_1)g) +
  \tau_{2q}dg. \]
For $1\le i\le 2q-1$, se see that $dd(\tau_i)$ has two terms for each
way of writing $i$ as a sum of three positive integers, and they
cancel. So we have $d^2=0$ on the subalgebra they generate. Thus we
have 
\begin{align*}
d^2(f+\tau_{2q}g)&=d(df+(\tau_1\tau_{2q-1}+\dots
+\tau_{2q-1}\tau_1)g)+\tau_{2q}dg) \\
&=d^2f + (\tau_1\tau_{2q-1}+\dots+\tau_{2q-1}\tau_1)dg +
(\tau_1\tau_{2q-1}+\dots+\tau_{2q-1}\tau_1)dg 
=0.
\qedhere
\end{align*}
\end{proof}

\begin{proposition}
The definitions above make $Q$ into a cocommutative DG Hopf
algebra.\index{Hopf algebra}\index{DG!Hopf algebra}
\end{proposition}
\begin{proof}
The above lemmas show that $Q$ is a DG bialgebra. It is easy to check
that the antipode satisfies the identity
$S(x_{(1)})x_{(2)}=x_{(1)}S(x_{(2)})=0$ in Sweedler notation, for
elements of non-zero degree this only needs checking on the
generators, where it is clear. Cocommutativity also only needs
checking on generators, where it is true by~\eqref{eq:Delta}.
\end{proof}

\begin{theorem}\label{th:Qqi}
There is a quasi-isomorphism from the $A_\infty$ algebra 
$H_*\Omega BG\twohat$ to the DG algebra $Q$, sending $\alpha$ to $\alpha$, $\beta$
to $\beta$, and $\tau$ to $\tau_1$.
\end{theorem}
\begin{proof}
First, we show that $H_*Q$ is isomorphic to $H_*\Omega BG\twohat$ as
an algebra over the noncommutative ring 
$\kk\langle \alpha,\beta\rangle/(\alpha^2,\beta^2)$. We
define a $\kk\langle \alpha,\beta\rangle/(\alpha^2,\beta^2)$-module
homomorphism $\delta\colon Q \to Q$ sending a monomial of the form
$\tau_1\tau_i f$ to $\tau_{i+1}f$ for $1\le i \le 2q-1$, and all other
standard monomials to zero. Thus
$\delta(f+\tau_{2q}g)=\delta(f)$. Then we have
\begin{align*}
\delta d(\tau_1\tau_i f)
&=\delta(\tau_1(\tau_1\tau_{i-1}+\dots+\tau_{i-1}\tau_1)f 
+ \tau_1\tau_idf) \\
&=(\tau_2\tau_{i-1}+\dots+\tau_i\tau_1)f + \tau_{i+1}df \\
d\delta(\tau_1\tau_i f)
&=d(\tau_{i+1}f)=(\tau_1\tau_i+\dots+\tau_i\tau_1)f
+\tau_{i+1}df \\
(\delta d + d\delta)(\tau_1\tau_if)&=\tau_1\tau_if
\end{align*}
while for $j>1$ we have
\begin{align*}
\delta d(\tau_j f)
&=\delta((\tau_1\tau_{j-1}+\dots+\tau_{j-1}\tau_1)f+\tau_jdf)=\tau_jf \\
d\delta(\tau_jf)&=d(0)=0 \\
(\delta d + d\delta)(\tau_jf)&=\tau_jf.
\end{align*}
Thus $\delta d+d\delta$ is the identity on all monomials except those
in the $\kk\langle \alpha,\beta\rangle/(\alpha^2,\beta^2)$-submodule
spanned by $1$ and $\tau_1$, where it is zero. So $\delta$ defines a
homotopy from the identity map of $Q$ to the projection onto this
submodule. It follows that $H_*Q$ is isomorphic to $H_*\Omega
BG\twohat$ as an algebra over 
$\kk\langle\alpha,\beta\rangle/(\alpha^2,\beta^2)$,
with $\tau_1$ corresponding to $\tau$.

We have an $A_\infty$ morphism $f\colon A \to Q$ given by
$f_1(\alpha)=\alpha$, $f_1(\beta)=\beta$, and
\[ f_i(\tau,\dots,\tau)=
\tau_i,\qquad 1\le i \le 2q. \] 
The computation above shows that $f_1$ is
a quasi-isomorphism, and hence by definition so is $f$.
This computation is a practical illustration of Kadeishvili's
Theorem~\ref{th:Kadeishvili}.\index{Kadeishvili's Theorem}
\end{proof}

\begin{corollary}
The bounded derived categories $\Db(Q)$, $\Db(H_*\Omega
BG\twohat)$ and $\Db(H^*BG)$ are equivalent as triangulated categories.
\end{corollary}
\begin{proof}
This follows from Theorem~\ref{th:Dsg-Dcsg}, together with
Theorem~\ref{th:Qqi} above.
\end{proof}

The element $s=\alpha\beta+\beta\alpha$ is central in $Q$, so it makes 
sense to invert it in the $A_\infty$ algebra $H_*\Omega BG\twohat$. 

\begin{corollary}
We have equivalences of triangulated categories 
\[ \Db(Q[s^{-1}]) \simeq\Db(H_*\Omega BG\twohat[s^{-1}])\simeq 
\Dcsg(H_*\Omega BG\twohat) \simeq 
\Dsg(H^*BG). \]
\end{corollary}
\begin{proof}
Since $H_*\Omega BG\twohat$ is periodic, with periodicity
generator\index{periodicity generator}
$s$, the effect on $\Db(Q)$ of inverting $s$ is to quotient out the thick
subcategory generated by $\kk$.
So this corollary again follows from Theorem~\ref{th:Dsg-Dcsg}.
\end{proof}

\section{\texorpdfstring{Duality for $Q[s^{-1}]$-modules}
{Duality for Q[s⁻¹]-modules}}\label{se:D-dual}\index{duality}

In this section, we continue to work in Case~\ref{case:D1},
where $G$ has dihedral Sylow $2$-subgroups and one conjugacy class of involutions.

\begin{definition}
We write $K$ for
$\kk\langle\alpha,\beta\rangle/(\alpha^2,\beta^2)[s^{-1}]$,
where $s=\alpha\beta+\beta\alpha$.
\end{definition}

\begin{lemma}
The graded algebra $K$ is simple. The
trace form $K \otimes_{\kk[s,s^{-1}]} K \to \kk[s,s^{-1}]$ 
induces an isomorphism of $K$-modules
\begin{equation}\label{eq:Kdual} 
K\cong \Hom_{\kk[s,s^{-1}]}(K,\kk[s,s^{-1}]). 
\end{equation}
\end{lemma}
\begin{proof}
This is the algebra of endomorphisms of a graded vector space of
dimension two over the graded field $\kk[s,s^{-1}]$, with a basis element $u$ in degree
zero and a basis element $v$ in degree one. The element $\alpha$ sends
$u$ to $v$ and $v$ to zero, while $\beta$ sends $v$ to $su$ and $u$ to
zero. Thinking in terms of matrices over $\kk[s,s^{-1}]$
this can be visualised as
\begin{equation*}
\alpha\mapsto \begin{pmatrix}0&0\\1&0\end{pmatrix}, 
\qquad
\beta\mapsto \begin{pmatrix}0&s\\0&0\end{pmatrix}, 
\qquad
s\mapsto\begin{pmatrix} s&0\\0&s\end{pmatrix},
\end{equation*}
giving an isomorphism
\[ K \cong \Mat_2(\kk[s,s^{-1}]). \]
The trace form is given by multiplying matrices and taking the
trace. It takes $\alpha\otimes\beta$ and $\beta\otimes\alpha$ both to
$s$. It therefore induces an
isomorphism  of $K$-modules $K\cong \Hom_{\kk[s,s^{-1}]}(K,\kk[s,s^{-1}])$
sending $\alpha$ to the homomorphism sending $\alpha$ 
to zero and $\beta$ to $s$ and
sending $\beta$ to the homomorphism sending $\alpha$ to $s$ and
$\beta$ to zero.
\end{proof}

If $X$ is any $K$-module, we write $X^*=\Hom_{\kk[s,s^{-1}]}(X,\kk[s,s^{-1}])$.
Then using the lemma, the isomorphism given by~\eqref{eq:Kdual} gives
us isomorphisms
\begin{align*} 
X^*&=\Hom_{\kk[s,s^{-1}]}(X,\kk[s,s^{-1}]) \\
&\cong\Hom_{\kk[s,s^{-1}]}(K\otimes_K X,\kk[s,s^{-1}]) \\
& \cong \Hom_K(X,\Hom_{\kk[s,s^{-1}]}(K,\kk[s,s^{-1}])) \\
&\cong \Hom_K(X,K),
\end{align*}
and so we can just as well regard $X^*$ as $\Hom_K(X,K)$.

\begin{proposition}\label{pr:Qs-1-qi-dual}
There is a quasi-isomorphism of $Q[s^{-1}]$-bimodules 
\[ Q[s^{-1}]\to\Sigma Q[s^{-1}]^*. \]
\end{proposition}
\begin{proof}
The standard monomials form a free basis for $Q[s^{-1}]$ as a
$K$-module. We construct a $K$-module homomorphism 
$Q[s^{-1}]\to \Sigma^{|\tau|}Q[s^{-1}]^*$ as follows. It takes all
standard monomials (Definition~\ref{def:standard}) to zero except $1$
and $\tau_1$. It takes $1$ to 
the element of $Q[s^{-1}]^*$ taking value $1$ on $\tau_1$ and zero on
all other monomials, and it takes $\tau_1$ to the element of
$Q[s^{-1}]^*$ taking value $1$ on $1$ and value zero on all other
standard monomials. It is easy to check that this is a map of
$Q[s^{-1}]$-bimodules, and a quasi-isomorphism.
\end{proof}

\begin{proposition}\label{pr:hom-tensor}
If $X$ is a left $Q[s^{-1}]$-module and $Y$ is a right
$Q[s^{-1}]$-module,  then there is a natural
isomorphism of $K$-modules
\[ \Hom_{Q[s^{-1}]}(X,\Hom_K(Y,K)) \cong
  \Hom_K(Y\otimes_{Q[s^{-1}]} X,K). \]
If $Y$ is a $Q[s^{-1}]$-bimodule, this is an isomorphism of left
$Q[s^{-1}]$-modules. 
\end{proposition}
\begin{proof}
This is the usual Hom tensor adjunction.
\end{proof}

\begin{corollary}\label{co:duals}
If $X$ is a homotopically projective\index{homotopically projective} 
$Q[s^{-1}]$-module then we have a quasi-isomorphism 
\[ \Hom_{Q[s^{-1}]}(X,Q[s^{-1}]) \simeq \Sigma\,\Hom_K(X,K). \]
\end{corollary}
\begin{proof}
We have
\begin{align*}
\Hom_{Q[\tau^{-1}]}(X,Q[s^{-1}]) 
&\simeq\Hom_{Q[s^{-1}]}(X,\Sigma Q[s^{-1}]^*) \\
&\cong \Sigma\,\Hom_{Q[s^{-1}]}(X,\Hom_K(Q[s^{-1}],K)) \\
&\cong \Sigma\,\Hom_K(Q[s^{-1}]\otimes_{Q[s^{-1}]} X,K) \\
&\cong\Sigma\,\Hom_K(X,K).
\qedhere
\end{align*}
\end{proof}

\begin{theorem}\label{th:Qduality}
Let $X$ and $Y$ be $Q[s^{-1}]$-modules, such that $X$ is homotopically
projective,
and its image in $\Db(Q[s^{-1}])$ is compact.
Then we have a duality
\[    \Hom_{Q[s^{-1}]}(X,Y)^* \cong
\Hom_{Q[s^{-1}]}(Y, \Sigma^{-1}X). \]
\end{theorem}
\begin{proof}
Since $X$ is homotopically projective with compact image in $\Db(Q[s^{-1}])$,
we have quasi-isomorphisms
\[ \Hom_{Q[s^{-1}]}(X,Y) \simeq 
\Hom_{Q[s^{-1}]}(X,Q[s^{-1}]) \otimes_{Q[s^{-1}]}Y \]
and
\[ \Hom_{Q[s^{-1}]}(\Hom_{Q[s^{-1}]}(X,Q[s^{-1}]),Q[s^{-1}]) \simeq X. \]
Combining the second of these with Corollary~\ref{co:duals}, we have
\[ \Hom_{Q[s^{-1}]}(X,Q[s^{-1}])^* \simeq \Sigma^{-1}X. \]
Hence using Proposition~\ref{pr:hom-tensor}, we have
\begin{align*}
\Hom_{Q[s^{-1}]}(X,Y)^*
&=\Hom_K(\Hom_{Q[s^{-1}]}(X,Y),K) \\
&\simeq\Hom_K(\Hom_{Q[s^{-1}]}(X,Q[s^{-1}])\otimes_{Q[s^{-1}]} Y,K) \\
&\cong\Hom_{Q[s^{-1}]}(Y,\Hom_K(\Hom_{Q[s^{-1}]}(X,Q[s^{-1}]),K)) \\
&\simeq\Hom_{Q[s^{-1}]}(Y,\Sigma^{-1}X).
\qedhere
\end{align*}
\end{proof}

\section{Some indecomposables}

Let $G$ be
a finite group with dihedral Sylow $2$-subgroups and a single conjugacy class of
involutions. In this section we describe some indecomposable objects
in the categories $\Dsg(H^*BG)\cong \Dcsg(H_*\Omega BG\phat)$.

Consider first $A_\infty$ modules over the $A_\infty$ algebra
$B=H^*BG$. We refer to \eqref{eq:HPSL2p} and Theorem~\ref{th:AinftyHBG-D1} for its structure.
The quotient $B/(t^{2q+1})=\kk[\xi,\eta,t]/(\xi\eta,t^{2q+1})$ is 
formal,\index{formal $A_\infty$ algebra}
so ordinary modules over this ring pull back to $A_\infty$ modules
over $B$. For $1\le i \le 2q$, let $X_i$ be the module 
$B/(\eta,t^i)$ and $X'_i$ be the module $B/(\xi,t^i)$. Thus $X_i$ has periodic resolution
\[ \cdots
\xrightarrow{\left(\begin{smallmatrix} \xi & t^i \\ 0 &
\eta\end{smallmatrix} \right)} 
B\oplus B
\xrightarrow{\left(\begin{smallmatrix} \eta & t^i \\ 0 &
\xi\end{smallmatrix} \right)} 
B \oplus B \xrightarrow{\left(\begin{smallmatrix} \xi & t^i \\ 0 &
\eta\end{smallmatrix} \right)} 
B \oplus B \xrightarrow{(\eta,t^i)} B \to X_i \to 0 \]
and swapping $\eta$ and $\xi$ gives a resolution of $X'_i$.
In $\Db(B)$, the residue field $\kk$ 
sits in a triangle
\[  B/(t) \to B/(\eta,t) \oplus B/(\xi,t) \to  \kk. \]
Furthermore, $B/(t)$ sits in a triangle 
\[ \Sigma^{-2}B \xrightarrow{t} B \to B/(t)\to\Sigma^{-1}B. \]
So in $\Dsg(B)$, $B/(t)$ is isomorphic to zero, 
and $\kk$ decomposes as $B/(\eta,t)
\oplus B/(\xi,t)= X_1\oplus X'_1$.
Resolutions of $X_i$ and $X'_i$, keeping track of the shifts, are as follows.
\begin{gather*} 
\cdots\!\to\!\Sigma^{-9}B \oplus \Sigma^{-6-2i}B
\xrightarrow{\!\left(\begin{smallmatrix} \eta & t^i \\
0 & \xi\end{smallmatrix}\right)\!}
\Sigma^{-6}B \oplus \Sigma^{-3-2i}B
\xrightarrow{\!\left(\begin{smallmatrix} \xi & t^i \\ 
0 & \eta\end{smallmatrix}\right)\!}
\Sigma^{-3} B\oplus \Sigma^{-2i} B
\xrightarrow{\!(\eta,t^i)\!} B \to X_i \to 0, \\
\cdots\!\to\!\Sigma^{-9}B \oplus \Sigma^{-6-2i}B
\xrightarrow{\!\left(\begin{smallmatrix} \xi & t^i \\
0 & \eta\end{smallmatrix}\right)\!}
\Sigma^{-6}B \oplus \Sigma^{-3-2i}B
\xrightarrow{\!\left(\begin{smallmatrix} \eta & t^i \\ 
0 & \xi\end{smallmatrix}\right)\!}
\Sigma^{-3} B\oplus \Sigma^{-2i} B 
\xrightarrow{\!(\xi,t^i)\!} B \to X'_i \to 0. 
\end{gather*}
It follows that
$\Sigma^2X_i\cong X'_i$ and $\Sigma^2X'_i\cong X_i$
in $\Dsg(B)$.
The category $\Dsg(B)$ is periodic of period four, with periodicity
generator\index{periodicity generator} $s=\alpha\beta+\beta\alpha$.\medskip 

Let $A=B^!$ be the $A_\infty$ algebra $H_*\Omega BG\twohat$. 
For $1\le i\le 2q$, let $Y_i=\Ext^*_B(\kk,X_i)$, the indecomposable $A$-module with generators
$u$ and $v$ satisfying $\alpha u=0$, $\alpha v=0$, and
\begin{align*}
m_{i+1}(\tau,\dots,\tau,u)&=v, \\
m_{2q+2-i}(\tau,\dots,\tau,v)&=(\alpha\beta)^qu.
\end{align*}
Then in $\Dcsg(A)$, we have $\Sigma^{2i-1}Y_i \cong Y_{2q+1-i}$, so this gives 
$q$ isomorphism classes up to shift, all periodic with period four,
for a total of $4q$ isomorphism classes.
Note that the ring $A$ itself, as an
object in $\Dcsg(A)$, decomposes as $Y_1 \oplus \Sigma^2 Y_1$.

Here they are for $q=2$. In these diagrams, vertices represent basis elements.
Degree is given by
horizontal displacement, increasing to the right. Multiplication by
$\alpha$ and $\beta$ is given by the arrows, or by zero if there is no
arrow out of a vertex. The operations $m_j(\tau,\dots,\tau,-)$ are represented by
dotted lines:
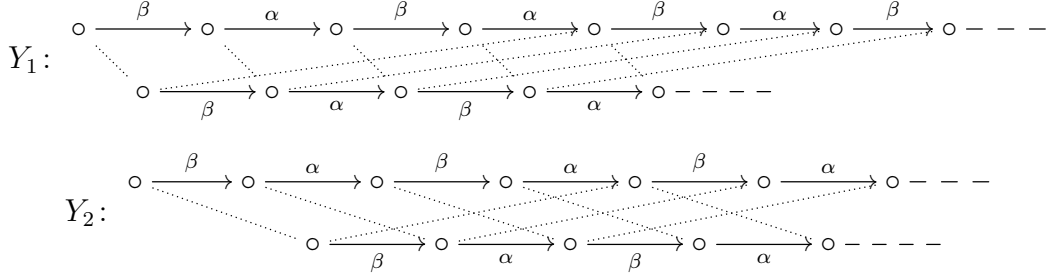
\[ Y_1\colon\vcenter{
\xymatrix@=4mm{\circ\ar[rr]^\beta\ar@{.}[dr]&&
\circ \ar[rr]^\alpha\ar@{.}[dr] && \circ\ar[rr]^\beta\ar@{.}[dr]&&
\circ\ar[rr]^\alpha\ar@{.}[dr]&&\circ\ar[rr]^\beta\ar@{.}[dr]&&
\circ\ar[rr]^\alpha&&\circ\ar[rr]^\beta&&\circ\ar@{--}[rr]&&\\ 
&\circ\ar[rr]_\beta\ar@{.}[urrrrrrr]&& \circ\ar[rr]_\alpha 
\ar@{.}[urrrrrrr]&&\circ\ar[rr]_\beta\ar@{.}[urrrrrrr]&&
\circ\ar[rr]_\alpha\ar@{.}[urrrrrrr]&&\circ\ar@{--}[rr]&&}} \]
\[ Y_2\colon\vcenter{
\xymatrix@=4mm{\circ\ar[rr]^\beta\ar@{.}[drrr]&&
\circ \ar[rr]^\alpha\ar@{.}[drrr] && \circ\ar[rr]^\beta\ar@{.}[drrr]&&
\circ\ar[rr]^\alpha\ar@{.}[drrr]&&\circ\ar[rr]^\beta\ar@{.}[drrr]&&
\circ\ar[rr]^\alpha&&\circ\ar@{--}[rr]&&\\ 
&&&\circ\ar[rr]_\beta\ar@{.}[urrrrr]&& \circ\ar[rr]_\alpha 
\ar@{.}[urrrrr]&&\circ\ar[rr]_\beta\ar@{.}[urrrrr]&&
\circ\ar[rr]_\alpha&&\circ\ar@{--}[rr]&&}} \]
Removing a finite number of nodes from the beginning of one of these
diagrams does not alter the isomorphism class in $\Dcsg(A)$.

\section{Classification  of indecomposables}\label{se:classification}

We continue to work in Case~\ref{case:D1} with $q\ge 1$, and write $B$ for the 
$A_\infty$ algebra $H^*BG$ and $A=B^!$ for the Koszul dual
$A_\infty$ algebra $H_*\Omega BG\twohat$. The way we classify the
indecomposables in $\Dcsg(A)\cong\Dsg(B)$ is via Morita equivalence,
reducing to the classification theorem of~\cite{Benson/Greenlees:2023a}.

Let $Y_i$, $1\le i\le 2q$, be the modules described in the previous
section. Then the regular representation of $A$ decomposes as
$Y_1\oplus \Sigma^2 Y_1$. 

Let $E$ be the $A_\infty$ algebra  $\Hom^*_A(Y_1,Y_1)$.
This is the algebra with $m_i=0$ for $i\ne 2,2q+1$,
defined as follows. The multiplication $m_2$ defines
the $\kk$-algebra structure as $\kk[s]\otimes\Lambda(\tau)$, with
generators $s$ and $\tau$ satisfying $|s|=(4,2q+1,2q+1)$,
$|\tau|=(1,q,q)$. We have 
\[ m_{2q+1}(s^{i_1}\tau,\dots,s^{i_{2q+1}}\tau)=s^{i_1+\dots+i_{2q+1}+q}, \] 
and $m_{2q+1}$ vanishes on all other tuples of monomials.

There is a right action of $E$ on $Y_1$ as an $A$-module given by 
\[ m_2(u,\tau)=v,\quad 
m_{2q+1}(v,\tau,\dots,\tau)=m_2(u,s^q). \] 
This makes $Y_1$ into an
$A$-$E$-bimodule, and $\Hom^*_A(Y_1,-)$  induces an equivalence of derived categories
$\Db(A)\simeq \Db(E)$ that sends $A$ to $E\oplus \Sigma^2E$ and $Y_1$
to $E$. 
It therefore also induces equivalences
$\Dcsg(A)\simeq \Dcsg(E)\simeq\Db(E[s^{-1}])$. Theorem~1.1
of~\cite{Benson/Greenlees:2023a} (with $a=1$, $b=2$, $h=2q+1$, $\ell=q$) 
therefore gives the following.

\begin{theorem}\label{th:D1}
The triangulated categories 
\[ \Dsg(B) \simeq \Dcsg(A)\simeq\Db(A[s^{-1}])\simeq \Dcsg(E) \simeq \Db(E[s^{-1}]) \]
satisfy the Krull--Schmidt theorem, and have $4q$ isomorphism classes
of indecomposable objects, in $q$ orbits of the shift functor
$\Sigma$. The Auslander--Reiten quiver\index{Auslander--Reiten quiver} 
of any of these triangulated categories
is isomorphic to $\bZ A_{2q}/T^2$, where $T$ is the 
translation functor\index{translation functor} $\Sigma^{-2}$. This
is a cylinder of height $2q$ and circumference $2$. The functor
$\Sigma$ switches the two ends of the cylinder.\qed
\end{theorem}

Here is a picture of the Auslander--Reiten
quiver\index{Auslander--Reiten quiver}
in the case $q=4$; the left and right side should be
identified to form a cylinder:
\[ \xymatrix@=3mm{Y_1 \ar[dr]\ar@{.}[dd] && \Sigma^2Y_1\ar[dr]&&Y_1\ar@{.}[dd] \\ 
& Y_2 \ar[ur]\ar[dr] &&\Sigma^2Y_2\ar[ur]\ar[dr]& \\
\Sigma^2Y_3\ar[ur]\ar[dr]\ar@{.}[dd]&& Y_3\ar[ur]\ar[dr]&& \Sigma^2Y_3\ar@{.}[dd] \\
& \Sigma^2Y_4\ar[ur]\ar[dr]&& Y_4\ar[ur]\ar[dr] & \\
\Sigma^3Y_4\ar[ur]\ar[dr]\ar@{.}[dd]&&\Sigma Y_4\ar[ur]\ar[dr]&&\Sigma^3Y_4\ar@{.}[dd] \\
& \Sigma Y_3\ar[ur]\ar[dr] && \Sigma^3Y_3\ar[ur]\ar[dr] &\\
\Sigma Y_2\ar[ur]\ar[dr]\ar@{.}[d]&&\Sigma^3Y_2\ar[ur]\ar[dr]&&\Sigma Y_2\ar@{.}[d] \\
&\Sigma^3 Y_1\ar[ur]&&\Sigma Y_1\ar[ur]& \\} \]

\begin{remark}\label{rk:D1}
In contrast with Theorem~\ref{th:D1}, the category
$\Dsg(A)\simeq\Dcsg(B)$ has infinite representation type. This can be
seen by examining the quotient $H_*\Omega BG\twohat/(\tau,s^q)$. By
Theorem~\ref{th:AinftyHOmegaBG-D1}, this is the formal $A_\infty$
algebra\index{formal $A_\infty$ algebra}
\[  \kk\langle\alpha,\beta|\alpha^2=0,\beta^2=0,(\alpha\beta)^q+(\beta\alpha)^q=0\rangle, \] 
which has tame representation type (Ringel~\cite{Ringel:1975a}). It
would be interesting to know whether $\Dsg(A)$ also has tame
representation type.
\end{remark}

\section{\texorpdfstring
{Loops on $BG\twohat$: two classes of involutions}
{Loops on BGٛ₂: two classes of involutions}}\label{se:D2}

We now turn to Case~\ref{case:D2}. 
This is the case where $G$ has a
dihedral Sylow $2$-subgroup $\D$ of order $4q$ with $q\ge 2$, 
and two conjugacy classes of involution.
In this case, $G$ has exactly one subgroup of index two, and it has
two isomorphism classes of simple modules in the principal block.

\begin{remark}\label{rk:c=0}
It follows from the work of Holm~\cite{Holm:1997a} that the
derived equivalence classes of algebras of dihedral type with two
isomorphism classes of simple modules are
determined by two parameters, namely a positive integer 
$k\ge 1$ and a field element $c\in\{0,1\}$.
For a block of a finite group with dihedral defect group of order
$4q$, the parameter $k$ is equal to $q$. 
Theorem~6.8 of Eisele~\cite{Eisele:2012a} shows that the case $c=1$
cannot occur for a block of a finite group, so we have $c=0$.
Note that by Corollary~2.3 of Generalov and
Romanova~\cite{Generalov/Romanova:2016a}, the cases $c=0$ and $c=1$
have different Hochschild cohomology rings, even in degree one.
\end{remark}

By Holm~\cite{Holm:1997a} and 
Proposition~\ref{pr:BGtwohatD}, for the purposes of studying
$BG\twohat$ we may assume that $G=PGL(2,p)$ for a suitable prime
$p\equiv 1 \pmod{4}$. In this case, the principal block $B_0$ of $\kG$ 
belongs to Erdmann's class D(2${\mathcal A}$). It has two simple
modules $\kk$ and $\MM$, whose 
$\Ext^1$ quiver\index{Ext@$\Ext^1$ quiver}\index{quiver} is as follows:
\[ \xymatrix{\kk\ar@/^/[r]^{e_2}\ar@(dl,ul)^{e_3} & \MM\ar@/^/[l]^{e_1}}. \]
Using Remark~\ref{rk:c=0}, the relations are
\[ e_2e_1=0,\qquad e_3^2=0, \qquad (e_1e_2e_3)^{q}=(e_3e_1e_2)^{q}. \]
We put an internal grading on the basic algebra in this case by
assigning degree $(\half,0)$ to $e_1$ and $e_2$ and $(0,1)$ to $e_3$.

We have
$H^*BG=\kk[\xi,y,t]/(\xi y)$ where  
\[ |\xi|=-(3,q+1,q),\qquad |y|=-(1,0,1),\qquad|t|=-(2,q,q). \] 
The restrictions to $\D$ are given by 
$\res_{\D}(\xi)=xt$, $\res_{\D}(y)=y$, and $\res_{\D}(t)=t$. Massey
products\index{Massey product}
are determined by
\[ \langle \xi,y,\dots,\xi,y\rangle = \langle
  y,\xi,\dots,y,\xi\rangle = t^{q+1}. \]
The computation of Hochschild cohomology is again very similar, and we
omit the details. The $A_\infty$ structure on $H^*BG$ again follows
the same lines as in Theorem~\ref{th:AinftyHBG}. This time, we only
replace $x$ by $\xi$, and again adjust the powers of $t$. So we have
\[ m_{2q}(\xi^i,y,\xi,y,\dots,\xi,y^j) 
=m_{2q}(y^j,\xi,y,\xi,\dots,y,\xi^i) 
=\xi^{i-1}y^{j-1}t^{q+1}. \]
The value of $m_{\ell(2q-2)+2}$ on the tuples at the end of the theorem is replaced by 
$\xi^{i-1}y^{j-1}t^{\ell(q+1)}$. 

\begin{theorem}
Let $G$ be a finite group with dihedral Sylow $2$-subgroups of order
$4q$ with $q\ge 2$ a power of two, and two classes of involutions, and
let $\kk$ be a field of characteristic two. Then
we have
\[ H_*\Omega BG\twohat = \Lambda(\tau) \otimes \kk\langle
  \alpha,Y\mid\alpha^2=0,\ Y^2=0\rangle \]
with
\[ |\tau|=(1,q,q),\qquad |\alpha|=(2,q+1,q),\qquad |Y|=(0,0,1). \]
In homological degree $2n$ we have monomials $(\alpha Y)^n$,
$(Y\alpha)^n$, $(\alpha Y)^{n-1}\alpha$ and $(Y\alpha)^nY$, and in odd
degrees we have $\tau$ times all of these.
\end{theorem}
\begin{proof}
Using Theorem~\ref{th:ExtRkk}, we have
\[ \Ext^{*,*}_{H^*BG}(\kk,\kk) = \Lambda(\tau) \otimes \kk\langle \alpha,
  Y\mid \alpha^2=0,\ Y^2=0\rangle \]
where the generators have degrees
\[ |\tau| = (-1,2,q,q), \qquad |\alpha|=(-1,3,q+1,q), \qquad |Y| =
  (-1,1,0,1). \]
The Eilenberg--Moore spectral  
sequence~\eqref{eq:Cotor}\index{Eilenberg--Moore spectral sequence}  
converging to $H_*\Omega BG\twohat$ has this as its $E_2$ page.
Here, the Eisenbud operator\index{Eisenbud operator} 
for the relation $\eta y=0$ is $s=\alpha Y + Y
\alpha$ in degree $(-2,4,q+1,q+1)$. Again there is no room for
non-zero differentials, and no ungrading problems, so we have $E_2=E_\infty=
H_*\Omega BG\twohat$. 
\end{proof}

\begin{lemma}\label{le:m=0-D2}
For any $A_\infty$ structure on $H^*BG$ that preserves internal
degrees, we have $m_n=0$ unless $n-2$ is divisible by 
$q-2$. In particular, for $2<n<q$ we have $m_n=0$.
\end{lemma}
\begin{proof}
The proof is essentially the same as the proof of Lemma~\ref{le:m=0}.
\end{proof}

\begin{proposition}\label{pr:HHHBG-D2}
Let $G$ be a group with a dihedral Sylow $2$-subgroup $\D$ of order
$4q$ with $q\ge 2$ a power of two,
and two conjugacy classes of involutions,
and let $\kk$ be a field of characteristic two.
The Hochschild cohomology\index{Hochschild cohomology} $\HH^*H^*BG$
has generators $s$, $t$, $\tau$, $\xi$, $y$, $u$, $v$
with
\begin{align*}
|s|&=(-2,4,q+1,q+1) \\
|t|&= -(0,2,q,q) &
|\tau|&= (-1,2,q,q) \\
|\xi|&= -(0,3,q+1,q)&
|y|&= -(0,1,0,1) \\
|u|&= -(1,0,0,0) &
|v|&= -(1,0,0,0).
\end{align*}
The relations are given by $u^2=v^2=uv=\tau^2=0$, 
$y u=\xi v=0$, $\xi s=y s=0$, and $us=vs$. The non-zero monomials and their degrees are as
follows, with $i_1,i_2\ge 0$, $\ep_1,\ep_2\in\{0,1\}$.
\begin{align*}
|s^{i_1}t^{i_2}\tau^{\ep_1}u^{\ep_2}|&
=(-2i_1-\ep_1-\ep_2,4i_1-2i_2+2\ep_1,i_1+q(i_1-i_2+\ep_1),i_1+q(i_1-i_2+\ep_1)), \\
|s^{i_1}t^{i_2}\tau^{\ep_1}v^{\ep_2}|&
=(-2i_1-\ep_1-\ep_2,4i_1-2i_2+2\ep_1,i_1+q(i_1-i_2+\ep_1),i_1+q(i_1-i_2+\ep_1)), \\
|\xi^{i_1}t^{i_2}\tau^{\ep_1}u^{\ep_2}|&
=(-\ep_1-\ep_2,-3i_1-2i_2+2\ep_1,-i_1+q(-i_1-i_2+\ep_1),q(-i_1-i_2+\ep_1)) \\
|y^{i_1}t^{i_2}\tau^{\ep_1}v^{\ep_2}|&
=(-\ep_1-\ep_2,-i_1-2i_2+2\ep_1,q(-i_2+\ep_1),-i_1+q(-i_2+\ep_1)) 
\end{align*}
There is only one monomial
in degree $(-i,i-2,0,0)$ with $i>2$, namely
\[ |s^qt^{q+1}|=(-q,q-2,0,0). \]
\end{proposition}
\begin{proof}
As in Theorem~\ref{th:HHHBD}, we use the approach of
Theorems~\ref{th:HHR} and~\ref{th:Negron}.
Thus $\HH^*H^*BG$ is the homology of
the complex 
\[ (H^*BG \otimes H_*\Omega BG\twohat,\partial), \]
where the generators $t$, $\xi$ and $\eta$ are in homological degree
zero, the generators $\tau$, $\alpha$ and $Y$ are in homological
degree $-1$, and the differential is given by $\partial=[e,-]$ where
$e = t \otimes \tau + \xi \otimes \alpha + y\otimes Y$. 
Thus setting $s=\alpha Y+Y\alpha$, we have $\partial(t)=0$, $\partial(\xi)=0$, 
$\partial(y)=0$, $\partial(\alpha)=y s$, $\partial(\tau)=0$,
$\partial(Y)=\xi s$. The generators and relations for the homology
of this complex are therefore as given, with $u=\xi\alpha$ and
$v=yY$.

For the last statement, the computation is similar to the
corresponding part of the proof of Theorem~\ref{th:HHHBD}.
\end{proof}

\begin{theorem}\label{th:AinftyHBG-D2}
The $A_\infty$ structure on $H^*BG$ is given as follows.
The $m_n$ are $\kk[t]$-multilinear maps with $m_n=0$ 
for $n$ not congruent to $2$ modulo $2q-2$, and for $i,j\ge 1$
\[ m_{2q}(\xi^i,y,\xi,y,\dots,\xi,y^j)=m_{2q}(y^j,\xi,y,\xi,\dots,y,\xi^i)=\xi^{i-1}y^{j-1}t^{q+1} \]
where the arguments alternate between $\xi$ and $y$, and the right hand
side is zero unless either $i=1$ or $j=1$; $m_{2q}$ is zero on all
other tuples of monomials not involving $t$.
The maps $m_{\ell(2q-2)+2}$ with $\ell>1$ similarly vanish on all tuples of
monomials not involving $t$, except the ones which look as above,
but for some choice of indices in the tuple:
\begin{multline*} 
1\le e_1\le e_2\le \cdots \le e_{\ell-1} 
< e_{\ell-1}+(2q-2)+1   
\le e_{\ell-2}+2(2q-2)+1 \\
\le \cdots \le e_1+(\ell-1)(2q-2)+1\le \ell(2q-2)+2. 
\end{multline*}
the exponents on the terms are increased by one
(or correspondingly more if an index is repeated).
The value on these tuples is $\xi^{i-1}y^{j-1}t^{\ell(q+1)}$.
Thus 
\[
m_{\ell(2q-2)+2}(x^{i+\alpha_1},y^{\alpha_2},x^{\alpha_3},\dots,
x^{\alpha_{\ell(2q-2)+1}},y^{j+\alpha_{\ell(2q-2)+2}}) =
  x^{i-1}y^{j-1}t^{\ell(q+1)} \]
where each $\alpha_\sigma$ is one plus the number of indices in the list above that are equal
to $\sigma$.
\end{theorem}
\begin{proof}
This is similar to the proof of Theorem~\ref{th:AinftyHBG}, but
using Lemma~\ref{le:m=0-D2} and
Proposition~\ref{pr:HHHBG-D2} instead of Lemma~\ref{le:m=0} 
and Theorem~\ref{th:HHHBD}.
\end{proof}

\begin{lemma}\label{le:m=0v2-D2}
For any $A_\infty$ structure on $H_*\Omega BG\twohat$ that preserves
internal degrees, we have $m_n=0$ unless $n - 2$ is
divisible by $q-1$. In particular,  for $2<n<q+1$ we have $m_n=0$. 
\end{lemma}
\begin{proof}
Looking at the degrees of the generators $\tau$, $\alpha$ and $\beta$,
for any monomial $\zeta$ in $H_*\Omega BG\twohat$ we have $a\equiv b
\pmod{q-1}$. So for any $n$-tuple $(\zeta_1,\dots,\zeta_n)$,
the degree of $m_n(\zeta_1,\dots,\zeta_n)$ satisfies
$a\equiv b+n-2 \pmod{q-1}$. So for this expression to be
non-zero we must have $n -2 \equiv 0 \pmod{q-1}$.
\end{proof}

\begin{proposition}\label{pr:HHHOmegaBG-D2}
  The Hochschild cohomology\index{Hochschild cohomology}
  $\HH^*H_*\Omega BG\twohat$ has 
generators $s$, $t$, $\tau$, $\xi$, $\eta$, $u$, and $v$ in degrees
\begin{align*}
|s|&=(0,2,q+1,q+1),\\
|t|&=-(1,1,q,q),&
|\tau|&=(0,1,q,q), \\
|\xi|&=-(1,2,q+1,q),&
|y|&=-(1,1,0,1), \\
|u|&=-(1,0,0,0),&
|v|&=-(1,0,0,0).
\end{align*}
The relations are given by $\xi\eta=0$, $u^2=v^2=uv=\tau^2=0$, 
$\eta u=\xi v = 0$, $\xi s = y s =0$, and $us=vs$.
The non-zero monomials and their degrees are given as follows, with
$i_1,i_2\ge 0$, $\ep_1,\ep_2\in\{0,1\}$.
\begin{align*}
|s^{i_1}t^{i_2}\tau^{\ep_1}u^{\ep_2}|
&=(-i_2-\ep_2,2i_1-i_2+\ep_1,(i_1-i_2+\ep_1)q+i_1,(i_1-i_2+\ep_1)q+i_1),\\
|s^{i_1}t^{i_2}\tau^{\ep_1}v^{\ep_2}|
&=(-i_2-\ep_2,2i_1-i_2+\ep_1,(i_1-i_2+\ep_1)q+i_1,(i_1-i_2+\ep_1)q+i_1),\\
|\xi^{i_1}t^{i_2}\tau^{\ep_1}u^{\ep_2}|&
=(-i_1-i_2-\ep_2,-2i_1-i_2+\ep_1,-i_1+(-i_1-i_2+\ep_1)q,(-i_1-i_2+\ep_1)q),\\
|y^{i_1}t^{i_2}\tau^{\ep_1}v^{\ep_2}|&
=(-i_1-i_2-\ep_2,-i_1-i_2+\ep_1,(-i_2+\ep_1)q,-i_1+(-i_2+\ep_1)q).
\end{align*}
Thus there is only one monomial with degree $(-i,i-2,0,0)$ with $i>2$, namely
\[ |s^qt^{q+1}|=(-q-1,q-1,0,0). \]
\end{proposition}
\begin{proof}
Again we use the approach of Theorem~\ref{th:Negron}.
This time, $\HH^*H_*\Omega BG\twohat$ is
the homology of the complex 
\[ (H_*\Omega BG\twohat \otimes H^*BG,\partial), \]
where the generators $\tau$, $\alpha$ and $Y$ of 
$H_*\Omega BG\twohat$ are in homological degree zero, 
the generators $t$, $\xi$, $y$ of $H^*BG$ are in homological degree $-1$, and
the differential is given by $\partial=[e,-]$ where $e=\tau\otimes t +
\alpha \otimes \xi + Y\otimes y$. 
So the answer is the same as in Proposition~\ref{pr:HHHBG-D2} but with the
degrees changed.

For the last statement, the computation is again similar to the
corresponding part of the proof of Theorem~\ref{th:HHHBD}.
\end{proof}

\begin{theorem}\label{th:AinftyHOmegaBG-D2}
In Case~\ref{case:D2}, 
the $A_\infty$ structure on $H_*\Omega BG\twohat$ is determined by
\[ m_{q+1}(\tau,\tau,\dots,\tau) = s^q, \]
where $s=\alpha Y+ Y\alpha$.
This implies that
\begin{equation*}
m_{q+1}(f_1(\alpha,Y)\tau,f_2(\alpha,Y)\tau,\dots,f_{q+1}(\alpha,Y)\tau)
= f_1(\alpha,Y)\dots f_{q+1}(\alpha,Y)s^q, 
\end{equation*}
and all $m_n$ for $n>2$ on all other $n$-tuples of monomials give
zero.
\end{theorem}
\begin{proof}
This is similar to the proof of
Theorem~\ref{th:AinftyHOmegaBG-D1}, but using
Lemma~\ref{le:m=0v2-D2} and
Proposition~\ref{pr:HHHOmegaBG-D2} in place of
Lemma~\ref{le:m=0v2-D1}  and
Proposition~\ref{pr:HHHOmegaBG-D1}.
\end{proof}

Everything from this point on is very similar to
Case~\ref{case:D1}, so we simply state the relevant results.

\begin{remark}\label{rk:D2-HHCBG}
In the spectral sequence $\HH^*H_*\Omega BG\twohat \Rightarrow
\HHinf^*H_*\Omega BG\twohat$, we have $d_q(\tau)=s^qt^q$. This implies
that after inverting $s$ we have
\[ \HHinf^*H^*BG[s^{-1}]\cong\HHinf^*H_*\Omega BG\twohat[s^{-1}] 
\cong \kk[s,s^{-1}][u,v,t]/(u^2,v^2,uv,t^q). \]
\end{remark}

The differential graded model $Q$ for the $A_\infty$ algebra $H_*\Omega BG\twohat$ is
essentially the same as that described in Section~\ref{se:Q},
except that $\beta$ is replaced by $Y$ in degree $(0,0,1)$, and
the element $s=\alpha Y+Y\alpha$ is in degree $(2,q+1,q+1)$.
So the generators for $Q$ are $\tau_1,\dots,\tau_q,\alpha,Y$, and the
relations between the $\tau_i$ are given by
\[ \sum_{j+k=i}\tau_j\tau_k = \begin{cases}
d\tau_i & 1\le i \le q \\ s^q & i = q+1 \\ 0 & q+2 \le i \le 2q.
\end{cases} \]
The final theorem in Case~\ref{case:D2} is as follows.

\begin{theorem}\label{th:D2}
The triangulated categories
\[ \Dsg(H^*BG) \simeq \Dcsg(H_*\Omega BG\twohat) \simeq \Db(H_*\Omega BG\twohat[s^{-1}]) \]
 satisfy the Krull--Schmidt theorem, and have $2q$ isomorphism classes
 of indecomposable objects, in $q$ orbits of the shift functor
 $\Sigma$. The Auslander--Reiten quiver\index{Auslander--Reiten quiver} 
is isomorphic to $\bZ A_{2q}/T$, where $T$ is the 
translation functor\index{translation functor} 
$\Sigma^{-2}$. This
 is a cylinder of height $2q$ and circumference one. The functor
 $\Sigma$ switches the two ends of the cylinder.\qed
\end{theorem}

\begin{remark}
Again, and for the same reason as in Remark~\ref{rk:D1}, 
in contrast with Theorem~\ref{th:D2} the category
$\Dsg(H_*\Omega BG\twohat)\simeq\Dcsg(H^*BG)$ 
has infinite representation type. This time, the 
formal\index{formal $A_\infty$ algebra} quotient is
\[ H_*\Omega BG\twohat/(\tau,s^q)=\kk\langle\alpha,Y\mid
\alpha^2=Y^2=(\alpha Y)^q+(Y\alpha)^q=0\rangle. \]
\end{remark}

\section{A related symmetric tensor category}

The symmetric tensor category $\cC_3$ in characteristic
two of non-finite tame representation type is discussed in
Section~5.2.3 in Benson and 
Etingof~\cite{Benson/Etingof:2019a}. This
has a basic algebra that is of dihedral type D(2$\mathcal A$) with $c=0$ and $k=1$ in Erdmann's
classification~\cite{Erdmann:1990a}, given by a
quiver and relations 
\begin{equation*}
\xymatrix{\kk\ar@(ul,dl)_a\ar@/^/[r]^b 
& \MM\ar@/^/[l]^c}
\end{equation*}
with relations
\[ a^2=0,\qquad bc=0,\qquad cba=acb. \]
This is not equivalent to a block of the group algebra of any finite
group, but it is quite similar in behaviour.

This algebra admits a $\bZ\times\bZ$-grading with $|a|=(1,0)$,
$|b|=|c|=(0,\half)$.
The minimal resolution over this algebra is the total complex of the
following double complex:
\[ \xymatrix{&&&\vdots\ar[d]&\vdots\ar[d]\\
&&P_\MM\ar[d]_{\bar b}&
P_\kk\ar[l]_{\bar a}\ar[d]_{\bar b\bar c}& 
 P_\kk \ar[l]_{\bar a}\ar[d]_{\bar b\bar c}&\cdots\ar[l]\\
&& P_\kk \ar[d]_{\bar b\bar c}&
P_\kk\ar[l]_{\bar a}\ar[d]_{\bar b\bar c} &
 P_\kk \ar[l]_{\bar a}\ar[d]_{\bar b\bar c}&\cdots\ar[l]\\
& P_\MM\ar[d]_{\bar b}&P_\kk\ar[l]_{\bar a}\ar[d]_{\bar b\bar c}
&P_\kk\ar[l]_{\bar a}\ar[d]_{\bar b\bar c} & 
P_\kk \ar[l]_{\bar a}\ar[d]_{\bar b\bar c}&\cdots\ar[l]\\
&P_\kk\ar[d]_{\bar b\bar c} & P_\kk\ar[l]_{\bar a}\ar[d]_{\bar b\bar c} 
&P_\kk\ar[l]_{\bar a}\ar[d]_{\bar b\bar c}&
P_k\ar[l]_{\bar a}\ar[d]_{\bar b\bar c}&\cdots\ar[l]  \\
P_\MM\ar[d]_{\bar b} & P_\kk\ar[l]_{\bar c\bar a}\ar[d]_{\bar b\bar c}
& P_\kk \ar[l]_{\bar a}\ar[d]_{\bar b\bar c}&
P_\kk\ar[l]_{\bar a}\ar[d]_{\bar b\bar c} 
&P_\kk\ar[l]_{\bar a}\ar[d]_{\bar b\bar c}&\cdots\ar[l] \\
P_\kk&P_\kk\ar[l]_{\bar a} &P_\kk\ar[l]_{\bar a} 
&P_\kk\ar[l]_{\bar a} 
&P_\kk\ar[l]_{\bar a}&\cdots\ar[l]} \]
The cohomology is 
\[ H^*\cC_3=\Ext^*_{\cC_3}(\kk,\kk)\cong\kk[x,y,z]/(xz+y^2) \] 
with $|x|=(-1,-1,0)$, $|y|=(-2,-1,-1)$, $|z|=(-3,-1,-2)$. 
The elements $x$, $y$ and $z$ are given by
shifts of degrees $(-1,0)$, $(-1,-1)$ and $(-1,-2)$ in this diagram,
killing the copies of $P_\MM$ and given by the identity map on all
copies of $P_\kk$. These maps commute, not just up to homotopy, 
and the relation $xz+y^2=0$ holds at
the level of cocycles. So this $\Ext$ algebra is formal as an
$A_\infty$ algebra. It is a Koszul algebra with Koszul dual
\[ {H^*\cC_3}^!\cong\kk[\eta]\langle
  \xi,\zeta\rangle/(\xi^2,\zeta^2,\xi\zeta+\zeta\xi+\eta^2), \]
with $\eta$ central, degrees $|\xi|=(0,1,0)$, $|\eta|=(1,1,1)$, $|\zeta|=(2,1,2)$,
and again formal as an $A_\infty$ algebra. As a module over $k[\eta]$
it is free of rank four, with basis $1$, $\xi$, $\zeta$, $\xi\zeta$.

Using Theorem~\ref{th:HHR},  and setting $u=x\xi+z\zeta$, we have
\[ \HH^*H^*\cC_3=k[x,y,z,\eta,u]/(x\eta^2,z\eta^2,u^2), \]
with $|x|=(0,-1,-1,0)$, $|y|=(0,-2,-1,-1)$, $|z|=(0,-3,-1,-2)$,
$|\eta|=(-1,2,1,1)$, $|u|=(-1,0,0,0)$. Then $\HHinf^*H^*\cC_3$ is
the same ring, but with the first two degrees added, so
$|x|=(-1,-1,0)$, $|y|=(-2,-1,-1)$, $|z|=(-3,-1,-2)$, $|\eta|=(1,1,1)$,
$|u|=(-1,0,0)$. 

\addtocontents{toc}{\protect\pagebreak}
\chapter{The semidihedral case}\index{semidihedral group}%
\index{Sylow subgroup!semidihedral}

\section{Introduction}

In this chapter,
we study the $A_\infty$ algebras $H^*BG$ and $H_*\Omega BG\twohat$
with coefficients in a field $\kk$ of characteristic two, in
the case where $G$ is a finite group with semidihedral Sylow
$2$-subgroup.
These groups were classified by Alperin, Brauer and
Gorenstein~\cite{Alperin/Brauer/Gorenstein:1970a}.
The simple groups of this type are the projective special 
linear groups $PSL(3,p^m)$\index{PSL3pm@$PSL(3,p^m)$} 
with $p^m\equiv 3\pmod{4}$, 
the projective special unitary groups 
$PSU(3,p^m)$\index{PSU3pm@$PSU(3,p^m)$} with 
$p^m\equiv 1 \pmod{4}$ (Notation~\ref{no:GLetc}), and the sporadic 
Mathieu group $M_{11}$\index{Mathieu group $M_{11}$} of
order $7920$.

We begin with the semidihedral group $\SD$ of order 
$8q$ itself. The group algebra in this
case was analysed by Bondarenko and
Drozd~\cite{Bondarenko/Drozd:1982a}, who gave a presentation as a
quiver with relations, but with a socle ambiguity. We resolve that
ambiguity in Theorem~\ref{th:kSD8q}, where we prove that for suitable
radical generators $X$ and $Y$ we have
\[ \kSD =\langle X, Y\mid X^2=0,\
  Y^2=X(YX)^{2q-1}+(YX)^{2q}\rangle. \]
We then recall the structure of $H^*B\SD$ and compute
$\Ext^*_{H^*B\SD}(\kk,\kk)$, and show how the Eilenberg--Moore
sequence~\eqref{eq:Cotor}\index{Eilenberg--Moore spectral sequence} 
with this as $E^2$ page converges to $\kSD$.

There are four cases for the possible fusion in $\SD$, leading to four
types for cochains on the classifying space of a finite group with
this fusion. Probably the most interesting is
the case where $G$ has no normal subgroup of index two. In that case,
it turns out that the basic algebra of the principal block admits a
grading, that endows the cohomology with a second, internal
grading. Also, the cohomology rings of these groups have the structure
of a complete intersection,\index{complete intersection} 
which allows for easy computation of the
Hochschild cohomology $\HH^*H^*BG$.
These facts together are what allows us to analyse the $A_\infty$
structure. The following theorem is proved in Sections~\ref{se:SD1}
to~\ref{se:SD1-dual}.

\begin{theorem}\label{th:SD1-main}
Let $G$ be a finite group with semidihedral Sylow $2$-subgroups of
order $8q$ ($q\ge 2$ a power of two), and
with no normal subgroup of index two, 
and let $\kk$ be a field of characteristic two. Then
the principal block $B$ of $\kG$ has an essentially unique grading.
This makes the cohomology ring 
\[ H^*BG=\kk[x,y,z]/(x^2y+z^2) \]
doubly graded, with $|x|=(-3,-q-1)$, $|y|=(-4,-4q)$ and $|z|=(-5,-3q-1)$.
This cohomology ring is formal as an $A_\infty$ 
algebra.\index{formal $A_\infty$ algebra} We have
\[ H_*\Omega BG\twohat = \Lambda(\xx,\yy) \otimes \kk[\zz] \]
with $|\xx|=(2,q+1)$, $|\yy|=(3,4q)$ and $|\zz|=(4,3q+1)$.
This is not formal, but the $A_\infty$ structure is given up to
quasi-isomorphism by the $\kk[\hat z]$-multilinear maps
\[
m_3(\hat x,\hat y,\hat x)=\hat z^2,\qquad
m_3(\hat x,\hat x\hat y,\hat x)=\hat x\hat z^2,\qquad
m_3(\hat y,\hat x,\hat x\hat y)=
m_3(\hat x\hat y,\hat x,\hat y)=\hat y\hat z^2, 
\]
and all $m_i$ with $i\ge 3$ vanish on all other tuples of monomials not involving
$\hat z$.
\end{theorem}

As part of this computation, we also compute Hochschild cohomology.

\begin{theorem}
Let $G$ be a finite group with semidihedral Sylow $2$-subgroups of
order $8q$ ($q\ge 2$ a power of two), and with no normal subgroups of
index two, and let $\kk$ be a field of characteristic two. Then
\[ \HH^*H^*BG
= H^*BG[\hat x,\hat z]/(\hat x^2+y\hat z^2,x^2\hat z^2), \]
with 
\begin{gather*}
|x|=(0,-3,-q-1),\qquad |y|=(0,-4,-4q),\qquad |z|=(0,-5,-3q-1), \\
|\hat x|=(-1,3,q+1),\qquad |\hat z|=(-1,5,3q+1).  
\end{gather*} 
The algebra $\HHinf^*H^*BG=\HHinf^*H_*\Omega BG\twohat$ is the same,
but with 
\begin{gather*}
|x|=(-3,-q-1),\qquad |y|=(-4,-4q),\qquad |z|=(-5,-3q-1), \\
|\hat x|=(2,q+1),\qquad |\hat z|=(4,3q+1). 
\end{gather*}
\end{theorem}

The second case in which we are able to make essentially complete
computations is where the Sylow $2$-subgroups of $G$ are semidihedral, 
$G$ has a normal subgroup $K$ of index two with generalised quaternion
Sylow $2$-subgroups, and $K$ has no normal subgroups of index two.
This case is very similar to the case discussed above. In particular,
again it turns out that the basic algebra of the principal block
admits a grading, that endows the cohomology with a second, internal
grading. The computations are similar, except that the degrees of
various elements have changed. Again the cohomology ring $H^*BG$ is
formal as an $A_\infty$ algebra. The corresponding theorems can be
found in Sections~\ref{se:SD2} to~\ref{se:Ainfty-SD2}.\medskip

The remaining case is the one where the Sylow $2$-subgroups of $G$ are semidihedral,
$G$ has a normal subgroup $K$ of index two with dihedral Sylow
$2$-subgroups, and $K$ has no normal subgroups of index two.
In this case, $H^*BG$ is not formal, but 
we compute $H_*(\Omega BG\twohat)$ 
(Theorem~\ref{th:SD3-HOmegaBGtwohat}) using the method of squeezed
resolutions from Section~\ref{se:squeezed}, since the Eilenberg--Moore
spectral sequence is difficult to ungrade directly. The information in this case
remains rather incomplete.

\section{Semidihedral groups}\label{se:semidihedral-groups}

The semidihedral group of order $8q$, $q\ge 2$ a power of two, is given by
the presentation
\[ \SD=\langle g,h\mid g^{4q}=1, h^2=1, hgh^{-1}=g^{2q-1}\rangle. \]
Let $\kk$ be a field of characteristic two. A modified version of the formulas of
Bondarenko and Drozd~\cite{Bondarenko/Drozd:1982a} describes the group
algebra $\kSD$ as follows.
Set
\[ X= 1+h,\qquad
Y= (1+h)\left(\sum_{i=0}^{\frac{q}{2}-1}g^{4i+1}+
\sum_{i=\frac{q}{2}+1}^q g^{4i-1}\right) + g^{2q}+hg^{4q-1}. \]

\begin{theorem}\label{th:kSD8q}
With this choice for $X$ and $Y$, the group algebra $\kSD$ of the
semidihedral group of order $8q$ has the presentation
\[ \kSD =\langle X, Y\mid X^2=0,\
  Y^2=X(YX)^{2q-1}+(YX)^{2q}\rangle. \]
\end{theorem}
\begin{proof}
The elements $X$ and $Y$ are in $J(\kSD)$, are independent modulo 
$J^2(\kSD)$, and the radical modulo its square is two dimensional,
so they generate $\kSD$. We have $X^2=(1+h)^2=0$,
so we must check the other relation. Set
\[ u = \sum_{i=0}^{\frac{q}{2}-1}g^{4i+1}+
\sum_{i=\frac{q}{2}+1}^q g^{4i-1} \]
so that $Y=(1+h)u+g^{2q}+hg^{-1}$. Write $N_1$, $N_2$ and $N_4$ for the norm
elements for $\langle g\rangle$, $\langle g^2\rangle$ and $\langle
g^4\rangle$ respectively (the \emph{norm element}\index{norm element}
of a group is the sum of the group elements, as an element of the
group algebra). Then we
have $u^2=N_4g^2$, $uh+hu=N_2gh$, $(1+h)u(1+h)=N_2g(1+h)$, $((1+h)u)^2=0$,
$ug=gu$, $g^{2q}h=hg^{2q}$,
and $u(g^{-1}+g^{2q+1})=1+g^{2q}$.
So in the expression for $Y$, the first and second terms commute, as
do the second and third. So squaring $Y$, we have square terms and cross terms
between the first and third term:
\begin{align*} 
Y^2&=0 + 1 + g^{2q}+(1+h)uhg^{-1}+hg^{-1}(1+h)u\\
&=1+g^{2q}+uhg^{-1}+ug^{-1}+N_2+hg^{-1}u+g^{2q+1}u\\
&=1+g^{2q}+u(g^{-1}+g^{2q+1})+(uh+hu)g^{-1}+N_2\\
&=N_2(1+h).
\end{align*}
On the other hand, we have
\begin{align*} 
YX&=(1+h)u(1+h)+(g^{2q}+hg^{-1})(1+h)\\
&=(N_2g+g^{2q}+g^{2q+1})(1+h).
\end{align*}
Since 
\begin{align*} 
(1+h)(N_2g+g^{2q}+g^{2q+1})(1+h)&=(1+h)g^{2q+1}(1+h)\\
&=(g^{-1}+g^{2q+1})(1+h), 
\end{align*}
by induction on $m\ge 1$ we have
\[ (YX)^m=(N_2g+g^{2q}+g^{2q+1})(g^{-1}+g^{2q+1})^{m-1}(1+h). \]
We have $N_2g(g^{-1}+g^{2q+1})=0$, so this simplifies for $m\ge 2$ to
\[ (YX)^m=(g^{2q}+g^{2q+1})(g^{-1}+g^{2q+1})^{m-1}(1+h). \]
We also have
$ (g^{-1}+g^{2q+1})^{2q-2}=(g^{-2}+g^2)^{q-1}=g^2N_4$, 
and so 
\[ (YX)^{2q-1}=(g^{2q}+g^{2q+1})g^2N_4(1+h)=(g^2+g^3)N_4(1+h), \] 
and 
\begin{align*}
X(YX)^{2q-1}&=(1+h)(g^2+g^3)N_4(1+h)\\
&=(g^2+g^3+g^{2q-2}+g^{2q-3})N_4(1+h)\\
&=gN_2(1+h). 
\end{align*}
Similarly, we have
\[ (g^{2q}+g^{2q+1})(g^{-1}+g^{2q+1})^{2q-1}=(g^{2q}+g^{2q+1})(g^{-1}+g^{2q+1})g^2N_4=N_1, \]
and so $(YX)^{2q}=N_1(1+h)$. Thus
\begin{align*} 
Y^2 &= N_2(1+h)\\
&=(gN_2+N_1)(1+h)\\
&=X(YX)^{2q-1}+(YX)^{2q}. 
\end{align*}

We thus have a surjective map from the algebra with the given
presentation to $\kSD$.
The relations
$X^2=0$ and $Y^2=X(YX)^{2q-1}+(YX)^{2q}$
imply that $Y^2X=XY^2=0$, and that the element
$Y^3=(YX)^{2q}=(XY)^{2q}$ is killed by $X$ and $Y$, and is therefore
in the socle. Thus the $8q$ alternating words in $X$ and $Y$, beginning
with $1$, $X$, $Y$, and ending with $(XY)^{2q}=(YX)^{2q}$ span 
the algebra with the given presentation. The surjective map to $\kSD$ is
therefore an isomorphism, and these alternating words form a basis.
\end{proof}

\begin{remark}
The reference~\cite{Bondarenko/Drozd:1982a}  uses a more complicated
choice of generators, and
gets the same relations, but only modulo the socle element $(XY)^{2q}=(YX)^{2q}$. 
It is erroneously stated\index{errors}
without proof in Section~15 of Benson and
Carlson~\cite{Benson/Carlson:1987a},
and in the papers of Generalov (page~530 of~\cite{Generalov:2009a},
page~164 of~\cite{Generalov:2010c},
page~279 of~\cite{Generalov:2012a}, and page~507 of~\cite{Generalov/Nikulin:2020a})
that the group algebra of the semidihedral
group is as given here, but without the extra term $(YX)^{2q}$ in the
expression for $Y^2$. 
See also Theorem~VIII.3 of Erdmann~\cite{Erdmann:1990a},  
where these two possibilities are given, labelled III.1\,(d) and
III.1\,(d$'$), but without deciding which is true.
In Corollary~7.2 of Erdmann~\cite{Erdmann:1988a},  and the tables at  
the back of~\cite{Erdmann:1990a} the incorrect choice is given.   
Theorem~\ref{th:kSD8q} shows that the correct answer  
is III.1\,(d$'$), whereas these sources state it as III.1\,(d).  
It is shown in Proposition~5.1 
of Bia{\l}kowski, Erdmann, Hajduk, Skowro\'nski and 
Yamagata~\cite{Bialkowski/Erdmann/Hajduk/Skowronski/Yamagata:2022a} that
these two algebras are not isomorphic.
\end{remark}

\iftrue
\begin{remark}
For a particular value of $q$
the following {\sc Magma}
code~\cite{Bosma/Cannon/Playoust:1997a}\index{Magma@{\sc Magma}}
(change the first line if necessary) checks the second relation in
Theorem~\ref{th:kSD8q}:
{\small
\begin{verbatim}
q:=8;
F<g,h>:=FreeGroup(2);
G:=PermutationGroup(quo<F|g^(4*q),h^2,h*g*h*g^(2*q+1)>);
kG:=GroupAlgebra(GF(2),G);
g:=kG!G.1; h:=kG!G.2; u:=kG!0;
for i in [0..(q/2-1)] do u:=u+g^(4*i+1); end for; 
for i in [(q/2+1)..q] do u:=u+g^(4*i-1); end for;
X:=(1+h)*u+g^(2*q)+h*g^(4*q-1); Y:=(1+h);
X^2+Y*(X*Y)^(2*q-1)+(X*Y)^(2*q) eq 0;
\end{verbatim}
}
\end{remark}
\fi

Here is a diagram of the case $q=2$ (only accurate modulo the extra
socle term in the expression for $X^2$). In this diagram, the single
lines represent multiplication by $Y$, and the double lines represent
multiplication by $X$.
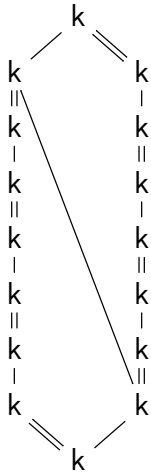
\[ \xymatrix@C=4mm@R=2mm{&\kk \ar@{-}[dl]\ar@{=}[dr] \\
\kk\ar@{=}[d]\ar@{-}[ddddddrr] &&\kk\ar@{-}[d] \\
\kk\ar@{-}[d] &&\kk\ar@{=}[d] \\
\kk\ar@{=}[d] &&\kk\ar@{-}[d] \\
\kk\ar@{-}[d] &&\kk\ar@{=}[d] \\
\kk\ar@{=}[d] &&\kk\ar@{-}[d] \\
\kk\ar@{-}[d] &&\kk\ar@{=}[d] \\
\kk\ar@{=}[dr] &&\kk\ar@{-}[dl] \\
&\kk} \]
This algebra has tame representation type, and its modules were classified by Bondarenko and 
Drozd~\cite{Bondarenko/Drozd:1982a}, 
Crawley-Boevey~\cite{Crawley-Boevey:1989b,Crawley-Boevey:1989a}. 
The cohomology ring was computed first by
Munkholm~\cite{Munkholm:1969a} and later also 
by Evens and Priddy~\cite{Evens/Priddy:1985a}, and is as follows.
\begin{equation}\label{eq:HBSD} 
H^*B\SD=\kk[x,y,z,w]/(xy,y^3,yz,z^2+x^2w), 
\end{equation}
with $|x|=|y|=-1$, $|z|=-3$ and $|w|=-4$. Here, $x$ and $y$ are dual
to $X$ and $Y$.

This is not formal as an $A_\infty$ algebra (see
Theorem~\ref{th:p-group}). In the next section we compute a few of the
higher multiplications.

\begin{remark}
The subalgebra $\cA_1$ of the Steenrod algebra generated by $\Sq^1$
and $\Sq^2$ is closely related to $\kSD$, with presentation
\[ \kk\langle \Sq^1,\Sq^2\mid (\Sq^1)^2=0,
  (\Sq^2)^2=\Sq^1\Sq^2\Sq^1\rangle.\]
This is like a (nonexistent) semidihedral group of order eight, but
without the socle element in the second relation. So it has 
type~III.(d) rather than~III.(d)$'$ in Erdmann's 
classification~\cite{Erdmann:1990a}. The cohomology
is the same ring as above~\eqref{eq:HBSD}, but with a different
$A_\infty$ structure.
\end{remark}

\section{\texorpdfstring{Resolutions for $\kSD$}{Resolutions for kSD}}

In this section, we write out the minimal resolution of $\kk$ over
$\kSD$. 
Since it is no extra work, we compute the minimal resolution of $\kk$
over an algebra of type III.I(d) or III.I(d$'$) in Erdmann's
classification~\cite{Erdmann:1990a} of algebras of semidihedral type
in characteristic two, given by the presentation
\[ \Lambda = \kk\langle X,Y\mid X^2=0,\ Y^2=X(YX)^{k-1}+\lambda
  (YX)^k\rangle \]
with $\lambda\in\kk$ and $k\ge 2$. The case of the group algebra $\kSD$ of a  
semidihedral group of order $8q$ with $q$ a power of two is then 
recovered by setting $\lambda=1$ and $k=2q$.

The minimal resolution of $\kk$ is the total complex of the following
double complex, where we have written $v$ for 
$\bar Y(\bar X\bar Y)^{k-1}$ and $w$ for $(\bar Y\bar
X)^{k-1}(1+\lambda \bar Y)$.
\[ \xymatrix@C=12mm@R=12mm{
&&&&\Lambda\ar[d]^{\bar Y} & 
\Lambda\ar[l]_{\bar Y}\ar[d]^w & \cdots\ar[l] \\
&&&&\Lambda\ar[d]^v & 
\Lambda\ar[l]_{\bar X}\ar[d]^v &\cdots\ar[l] \\
&&\Lambda\ar[d]^{\bar Y} & 
\Lambda\ar[l]_{\bar Y}\ar[d]^w &
\Lambda\ar[l]_{\bar Y\bar X}\ar[d]^v & 
\Lambda\ar[l]_{\bar X}\ar[d]^v & \cdots\ar[l] \\
&&\Lambda\ar[d]^v & 
\Lambda\ar[l]_{\bar X}\ar[d]^v&
\Lambda\ar[l]_{\bar X}\ar[d]^v & 
\Lambda\ar[l]_{\bar X}\ar[d]^v & \cdots\ar[l] \\
\Lambda\ar[d]^{\bar Y} & 
\Lambda\ar[l]_{\bar Y}\ar[d]^w &
\Lambda\ar[l]_{\bar Y\bar X}\ar[d]^v&
\Lambda\ar[l]_{\bar X}\ar[d]^v & 
\Lambda\ar[l]_{\bar X}\ar[d]^v&
\Lambda\ar[l]_{\bar X}\ar[d]^v&\cdots\ar[l] \\
\Lambda & \Lambda \ar[l]_{\bar X}&\Lambda\ar[l]_{\bar X}&
\Lambda\ar[l]_{\bar X} & \Lambda\ar[l]_{\bar X} & 
\Lambda\ar[l]_{\bar X}& \cdots\ar[l]} \]
Here $\bar X$ and $\bar Y$ are the elements of 
$\End_\Lambda(\Lambda)\cong \Lambda^\op$
corresponding to $X$ and $Y$ in $\Lambda$.

The cohomology element $x$ is represented by the map of double
complexes given by left shift 
followed by multiplication by the following elements of $\kSD$ in the
corresponding places in the double complex. It is easy to check that
this is a map of double complexes that induces $x$ in cohomology.
{\tiny
\[ \xymatrix@R=3mm@C=3mm{
&&&&\vdots\\
&&
(\bar X\bar Y)^{k-2}\bar X(1+\lambda\bar Y)&
\bar Y(1+\lambda \bar Y) & 
\cdots \\
&&1&1\\
(\bar X\bar Y)^{k-2}\bar X(1+\lambda\bar Y)&
\bar Y(1+\lambda \bar Y) & 
1&1\\
1&1&1&1&\cdots} \]
}
The element $y$ is represented by a map which is zero on most of the
copies of $\Lambda$, and only non-zero on the upper boundary. On this upper
boundary it is given by the following maps, followed by multiplication
by the labelled elements of $\kSD$ in the following diagram. Again, it
is easy to check that this is a map of double complexes.
{\tiny
\[ \xymatrix@R=10mm@C=10mm{
&&&&\bullet\ar[d]^1&\cdots\ar[l]^1\\
&&&&
\bullet\ar[d]^{(\bar Y\bar X)^{k-1}} \\
&&\bullet\ar[d]^1 & \bullet\ar[l]^1&
\bullet\ar[l]^{\bar X}\\
&&\bullet\ar[d]^{(\bar Y\bar X)^{k-1}}\\
\bullet\ar[d]^1&\bullet\ar[l]^1&
\bullet\ar[l]^{\bar X} \\
\bullet} \]
}
The element $z$ is represented by a shift two to the left and one
down, followed by multiplication by the following elements of $\kSD$.
{\tiny
\[ \xymatrix@R=3mm@C=3mm{&&(\bar X\bar Y)^{k-1}&
\bar Y+\lambda (\bar X \bar Y)^{k-1}\bar X & \cdots \\
&&1&1 \\
(\bar X\bar Y)^{k-1}&
\bar Y+\lambda (\bar X \bar Y)^{k-1}\bar X&1&1\\
1&1&1&1& \cdots} \]
}
Finally, the element $w$ is represented by a shift two to the left and
two down. This strictly commutes with $x$, $y$ and $z$.

In particular, we can
read off from the structure and minimal resolution of $\kSD$ that part of
the $A_\infty$ structure on $H^*B\SD$ is given over the
central subalgebra $\kk[w]$ by 
\begin{gather*}
m_4(y,x,y,z)=w,\qquad
m_{2k-1}(x,y,x,\dots,y,x)=y^2.
\end{gather*}

\section{\texorpdfstring{Loops on $B\SD\twohat$}
{Loops on BSDٛ₂}}\label{se:OmegaBSDtwohat}

Since $\SD$ is a finite $2$-group, we have $\Omega B\SD\twohat\simeq
\SD$. So we should expect to see the Eilenberg--Moore spectral
sequence\index{Eilenberg--Moore spectral sequence}~\eqref{eq:Cotor}
converging to $\kSD$. 

\begin{theorem}
We have
\[ \Ext^{*,*}_{H^*B\SD}(\kk,\kk)=\Lambda(\hat w) \otimes
\kk\langle \hat x,\hat y,\hat z,\eta\mid 
\hat x^2=\hat y^2=0, \hat x\hat z=\hat z\hat x,\eta\hat y=\hat y \eta
 \rangle \]
where $|\hat x|=(-1,1)$, 
$|\hat y|=(-1,1)$, $|\hat z|=(-1,3)$, $|\hat w|=(-1,4)$,
$|\eta|=(-2,3)$, and
$\eta$ is the Massey triple product\index{Massey product} $\langle\hat y,\hat y,\hat
y\rangle$. The Poincar\'e series is
\[ \sum_{i,j=0}^\infty t^{i}u^{j}\dim_\kk\Ext^{i,-j}_{H^*B\SD}(\kk,\kk)
=\frac{(1+tu^{4})(1+tu)}{1-tu-tu^{3}-t^2u^{3}}. \]
Note that $\Ext^{i,-j}$ is homologically indexed $(-i,j)$, so that
the coefficient of $t^iu^j$ is the dimension of the space of elements
of degree $(-i,j)$.
\end{theorem}
\begin{proof}
The element $w\in H^*B\SD$, see~\eqref{eq:HBSD}, is a regular element. 
Its appearance in the relations is in terms that are at least cubic,
so we have an algebra isomorphism 
\[ \Ext^{*,*}_{H^*B\SD}(\kk,\kk) \cong \Lambda(\hat w) \otimes
\Ext^{*,*}_{R}(\kk,\kk). \]
where 
\[ R=H^*B\SD/(w)=\kk[x,y,z]/(xy,y^3,yz,z^2). \] 
This algebra $R$ is the fibre product of 
$\kk[x,z]/(z^2)\to \kk$ and $\kk[y]/(y^3) \to \kk$. 
So by Theorem A of Moore~\cite{Moore:2009a},
$\Ext^{*,*}_R(\kk,\kk)$ is the coproduct over $\kk$ of the algebras
\begin{align*}
\Ext^{*,*}_{\kk[x,z]/(z^2)}(\kk,\kk)&=\kk[\hat x,\hat z]/(\hat x^2) \\
\Ext^{*,*}_{\kk[y]/(y^3)}(\kk,\kk)&=\kk[\hat y,\eta]/(\hat y^2),
\end{align*}
where $\eta$ is the Massey triple product $\langle\hat y,\hat y,\hat
y\rangle=m_3(\hat y,\hat y,\hat y)$ in degree $(-2,3)$.
So we have
\[\Ext^{*,*}_R(\kk,\kk)= \kk\langle \hat x,\hat y,\hat z,\eta\mid \hat x^2=\hat y^2=0,\hat
  x\hat z = \hat z\hat x, \eta\hat y= \hat y \eta \rangle \]
which has Poincar\'e series
\[ \sum_{i,j=0}^\infty t^iu^j\dim_\kk\Ext^{i,-j}_{R}(\kk,\kk)
=\frac{1+tu}{1-tu-tu^3-t^2u^3}. \]
Finally, tensoring with $\Lambda(\hat w)$ multiplies the Poincar\'e
series by $(1+tu^4)$.
\end{proof}

The differentials in the Eilenberg--Moore spectral 
sequence~\eqref{eq:Cotor}\index{Eilenberg--Moore spectral sequence}
\[ \Ext^{*,*}_{H^*B\SD}(\kk,\kk)\Rightarrow \kSD \] 
are determined by the fact that it has to converge to the associated
graded of $\kSD$. They
are given by $d^2(\hat z) =\eta\hat x+\hat x \eta$,
\[ E^3 = \Lambda(\hat w) \otimes \kk[\eta] \otimes 
\kk\langle \hat x,\hat y \mid \hat x^2 =\hat y^2 = 0\rangle, \]
then $d^3(\hat w)=\eta^2$,
\[ E^4 = E^{4q-2}= \Lambda(\eta) \otimes
\kk\langle \hat x,\hat y \mid \hat x^2 =\hat y^2 = 0\rangle, \]
and finally
$d^{4q-2}(\eta)=(\hat x \hat y)^{2q}+(\hat y\hat x)^{2q}$. So
\[ E^{4q-1}= E^\infty = 
\kk\langle \hat x,\hat y \mid \hat x^2 =\hat y^2 = 0, (\hat x\hat
y)^{2q}=(\hat y\hat x)^{2q} \rangle, \]
which is the associated graded of the group algebra $\kSD$.

\section{\texorpdfstring{Groups with semidihedral Sylow $2$-subgroups}
{Groups with semidihedral Sylow 2-subgroups}}\label{se:semidihedral-Sylow}

Groups with semidihedral Sylow $2$-subgroups were classified by 
Alperin, Brauer and
Gorenstein~\cite{Alperin/Brauer/Gorenstein:1970a},
see also Wong~\cite{Wong:1964a,Wong:1966a}.
By Section~VIII of Brauer~\cite{Brauer:1966a}, or Proposition~1.1
of~\cite{Alperin/Brauer/Gorenstein:1970a},  
there are four possibilities for the $2$-fusion in a finite group $G$
with semidihedral Sylow $2$-subgroups, which are distinguished by the
numbers of conjugacy classes of involutions and of elements of order
four. By Theorem~1.1 of Craven and
Glesser~\cite{Craven/Glesser:2012a}, these represent the only possible
fusion systems\index{fusion!system} on semidihedral $2$-groups.

To describe these, we first describe some particular finite groups
with semidihedral Sylow $2$-subgroups. 
First, we describe the groups
$SL^\pm(2,p^m)$\index{SL@$SL^\pm(2,p^m)$}  
and $SU^\pm(2,p^m)$ (cf.\ Notation~\ref{no:GLetc}).\index{SU@$SU^\pm(2,p^m)$}
These are the subgroups of $GL(2,p^m)$, respectively $GU(2,p^m)$, consisting
of elements of determinant $\pm 1$.
If $p^m \equiv 3 \pmod{4}$ then $SL^\pm(2,p^m)$ has semidihedral Sylow
$2$-subgroups, while if $p^m\equiv 1\pmod{4}$ then $SU^\pm(2,p^m)$ has
semidihedral Sylow $2$-subgroups. We remark that
$SL(2,p^m)$ and $SU(2,p^m)$ are isomorphic.

Next, we describe the group denoted
$PGL^*(2,p^{2m})$\index{PGL@$PGL^*(2,p^{2m})$} 
in Section~II.2 of~\cite{Alperin/Brauer/Gorenstein:1970a}. For $p$ odd,
the group $P\Gamma L(2,p^{2m})$ is a semidirect product of $PGL(2,p^{2m})$ by
a cyclic group of order $2m$ acting as Galois automorphisms. The group
$PGL(2,p^{2m})$ has $PSL(2,p^{2m})$ as a normal subgroup of index two.
Thus $P\Gamma L(2,p^{2m})$ contains three distinct subgroups, each
having $PSL(2,p^{2m})$ as a subgroup of index two. One of these is
$PGL(2,p^{2m})$, one is a semidirect product of $PSL(2,p^{2m})$ by the
Galois automorphism of order two, and the third one is the group we
denote by $PGL^*(2,p^{2m})$. For example, $PGL^*(2,9)$ is isomorphic
to the stabiliser $M_{10}$ of a point in the Mathieu group
$M_{11}$.\index{Mathieu group $M_{11}$} It is
proved in Lemma~2.3 of Gorenstein~\cite{Gorenstein:1969a} that the
Sylow $2$-subgroups of $PGL^*(2,p^{2m})$ are semidihedral.

The following four cases correspond to the four parts of Proposition~1
in Section~II, page~10 of~\cite{Alperin/Brauer/Gorenstein:1970a}. In
that paper, the first three cases are called QD-groups, Q-groups and
D-groups. The number of simple modules in the principal block in these cases, and the
decomposition matrices, were determined
by Erdmann~\cite{Erdmann:1988a,Erdmann:1988b}.

\begin{case}\label{case:SD1}
 $G$ has one class of involutions and one class of elements of
order four. In this case, $G$ has no normal subgroup of index
two. The group $G/O(G)$ (Notation~\ref{no:gt}) has a simple normal
subgroup with odd index, 
isomorphic to $PSL(3,p^m)$ with $p^m\equiv 3\pmod{4}$, 
$PSU(3,p^m)$ with $p^m\equiv 1\pmod{4}$, or the Mathieu group $M_{11}$.
The principal block of $\kG$ has three isomorphism classes of
simple modules.
\end{case}

\begin{case}\label{case:SD2}
 $G$ has two classes of involutions\index{involutions} and one 
class of elements of order four. In this case, $G$
has a normal subgroup $K$ of index two with generalised
quaternion Sylow $2$-subgroups, and $K$ has no normal subgroups of
index two. The group $G/O(G)$ is either
isomorphic to a subroup of 
$\Gamma L(2,p^m)$ containing $SL^\pm(2,p^m)$ with odd index, for some
prime power $p^m\equiv 3\pmod{4}$, or it is isomorphic to a subgroup
of $\Gamma U(2,p^m)$ containing $SU^\pm(2,p^m)$ with odd index, for some
prime power $p^m\equiv 1 \pmod{4}$.
The principal block of $\kG$ has two isomorphism classes of
simple modules.
\end{case}

\begin{case}\label{case:SD3}
 $G$ has one class of involutions and two classes of elements of
order four. In this case, $G$
a normal subgroup $K$ of index two with dihedral Sylow
$2$-subgroups, and $K$ has no normal subgroups of index two.
The group $G/O(G)$ is isomorphic to a subgroup of $P\Gamma
L(2,p^{2m})$ containing $PGL^*(2,p^{2m})$ with odd index, for some odd
prime $p$ and positive integer $m$.
The principal block of $\kG$ has two isomorphism classes of simple
modules. 
\end{case}

\begin{case}\label{case:SD4}
$G$ has two classes of involutions and two classes of elements of
order four. In this case, $O(G)$ is a normal complement to a Sylow
$2$-subgroup $\SD$, so that $G/O(G)\cong \SD$ and $H^*BG\cong H^*B\SD$. The
principal block of $\kG$ is isomorphic to $\kSD$, and has one
isomorphism class of simple modules, namely the trivial module.
\end{case}

Representation theory and cohomology of groups with 
semidihedral Sylow $2$-subgroups, and more generally,
blocks with semidihedral defect groups and finite dimensional
algebras of semidihedral type, are discussed in 
Erdmann~\cite{Erdmann:1979a,  
Erdmann:1988a,Erdmann:1990c,Erdmann:1990a,  
Erdmann:1990b,Erdmann:1994a},   as well as
Benson and Carlson~\cite{Benson/Carlson:1987a},
Bogdani\'c~\cite{Bogdanic:2015b,Bogdanic:2015a},  
Brauer~\cite{Brauer:1966a} (Section~VIII),
Carlson, Mazza and Th\'evenaz~\cite{Carlson/Mazza/Thevenaz:2013a},
Chin~\cite{Chin:1995a},
Evens and Priddy~\cite{Evens/Priddy:1985a},  
Generalov et~al.~\cite{Antipov/Generalov:2004a,
Generalov:2002b,Generalov:2005b,Generalov:2006a,
Generalov:2007a,Generalov:2007b,Generalov:2009a,
Generalov:2010c,Generalov:2012a,Generalov:2012b,
Generalov:2013a,Generalov:2013b,Generalov:2014a,
Generalov:2017a,Generalov:2018a,
Generalov:2023a,Generalov/Nikulin:2020a,
Generalov/Zaykovskiy:2018a,Generalov/Zaykovskiy:2019a,
Generalov/Zilberbord:2016a},
Hayami~\cite{Hayami:2011a,Hayami:2018a},  
Holm~\cite{Holm:1997a,Holm:1999a},  
Holm and Zimmermann~\cite{Holm/Zimmermann:2008a},
Kawai and Sasaki~\cite{Kawai/Sasaki:2006a},
Koshitani, Lassueur, and
Sambale~\cite{Koshitani/Lassueur/Sambale:2022a},
Macgregor~\cite{Macgregor:2022a},
Martino and Priddy~\cite{Martino/Priddy:1991a},
M\"uller~\cite{Muller:1974a},
Olsson~\cite{Olsson:1975a},
Sasaki~\cite{Sasaki:1994a},
Taillefer~\cite{Taillefer:2019a},
Zhou and Zimmermann~\cite{Zhou/Zimmermann:2011a}.
The homology of $\Omega BG\twohat$
was computed by Levi~\cite{Levi:1995a}.

\begin{proposition}\label{pr:BGtwohatSD}
Suppose that $G$ has a semidihedral Sylow $2$-subgroup $\SD$. Then the homotopy
type of $BG\twohat$ is determined by $|\SD|$ and the number of 
classes of involutions and of elements
of order four. In particular, if $G$ has no normal subgroup of index
two, then the homotopy type of $BG\twohat$ is determined by $|\SD|$.
\end{proposition}
\begin{proof}
This follows from Theorem~\ref{th:Oliver} and the main theorem 
of~\cite{Alperin/Brauer/Gorenstein:1970a} described above.
\end{proof}

We end this section with a table of the various cases of algebras of
semidihedral type in characteristic two. Note that the definition of
semidihedral type in~\cite{Erdmann:1990a} is slightly broader than
in~\cite{Erdmann:1988a,Erdmann:1990c}, leading to the extra case
$SD(3{\mathcal C})_2$. 
In each case except
$SD(3{\mathcal K})$, there
is a positive integer parameter $k$, which in our context is equal to
$2q$, and in some cases there are also further parameters. In the case
of $SD(3{\mathcal K})$ there are three integer parameters $a\ge b\ge
c$, $a\ge 2$.\bigskip

\begin{center}
\begin{tabular}{|c|l|c|c|c|c|} \hline
\!Erdmann~\cite{Erdmann:1990a}\! &
\cite{Erdmann:1988a,Erdmann:1990c} & 
Case & Group & $H^*$&$\HH^*$ \\ \hline
III.I(d) &&---&---&\cite{Generalov:2009a}&\cite{Generalov:2010c} \\
III.I(d$'$)  &&\ref{case:SD4}& semidihedral&
\!\!\!\cite{Munkholm:1969a,Evens/Priddy:1985a,Generalov:2009a}\!\!&
\cite{Generalov/Nikulin:2020a,Holm:2004a} \\
$SD(2{\mathcal A})_1$&\cite{Erdmann:1988a} II&
\ref{case:SD3}&$SU^\pm(2,p^m)$
&\cite{Chin:1995a,Generalov:2006a}& \\
&&&$p^m\equiv 1 \pmod{4}$&&\\
$SD(2{\mathcal A})_2$&\cite{Erdmann:1988a} III&
\ref{case:SD2}&$PGL^*(2,p^{2m}),$&\cite{Chin:1995a,Generalov:2006a}& \\
$SD(2{\mathcal B})_1$&\cite{Erdmann:1988a} IV&
\ref{case:SD2}&\quad ---
\blue{$[B_1(3M_{10})]$}
&\cite{Antipov/Generalov:2004a}
&\cite{Generalov/Zaykovskiy:2019a} \\
$SD(2{\mathcal B})_2$&\cite{Erdmann:1988a} I&
\ref{case:SD3}&$SL^\pm(2,p^m)$
&\cite{Chin:1995a,Generalov:2013a}&
\!\!\!\cite{Generalov:2013b,Generalov:2014a,Generalov:2018a}\!\! \\
&&&$p^m\equiv 3\pmod{4}$&&\\
$SD(2{\mathcal B})_3$&\cite{Erdmann:1988a} V&
\ref{case:SD2}&---&\cite{Antipov/Generalov:2004a}& \\
$SD(3{\mathcal A})_1$&\cite{Erdmann:1990c} II, \S5 &
\ref{case:SD1}&$PSU(3,p^m),$&\cite{Generalov:2002b}&\cite{Holm:2004a} \\
&&&$p^m\equiv 1\pmod{4}$&&\\
$SD(3{\mathcal A})_2$&\cite{Erdmann:1990c} VII, \S3 &---&---&\cite{Generalov:2006a}& \\
$SD(3{\mathcal B})_1$&\cite{Erdmann:1990c} IV, \S7&&---&\cite{Generalov:2007a}& \\
$SD(3{\mathcal B})_2$&\cite{Erdmann:1990c} I, \S7 &&---&\cite{Generalov:2007b}& \\
$SD(3{\mathcal C})_1$&\cite{Erdmann:1990c} VI, \S3 &---&---&& \\
$SD(3{\mathcal C})_2$& (excluded) &---&---&& \\
$SD(3{\mathcal D})$&\cite{Erdmann:1990c} III, \S6 &
\ref{case:SD1}&$PSL(3,p^m),\,p^m\equiv 3$&\cite{Generalov:2002b}&\cite{Holm:2004a} \\
&&&$\pmod{4},\,M_{11}$&&\\
$SD(3{\mathcal F})$&\cite{Erdmann:1990c} VIII, \S10 &---&---&& \\
$SD(3{\mathcal H})$&\cite{Erdmann:1990c} IX, \S10 &&
---
&& \\
$SD(3{\mathcal K})$&\cite{Erdmann:1990c} V, \S9 &---
&---&\cite{Generalov:2005b}&\cite{Generalov/Zilberbord:2016a} \\
\hline
\end{tabular}\bigskip
\end{center}

\begin{remarks}
The types with three simple modules are all 
derived equivalent\index{derived!equivalence} to an
algebra in the family $SD(3{\mathcal K})$ with uniquely determined
values of $a\ge b \ge c$, by 
Theorem~4.8 of Holm~\cite{Holm:1999a}. For blocks with semidihedral
defect group\index{defect group!semidihedral} 
of order $8q$ and three simple modules, 
these parameters are $2q\ge 2 \ge 1$, so
they are all derived equivalent.

Note that by 
Rickard~\cite{Rickard:1989a,Rickard:1991a},
for self-injective algebras,\index{self-injective algebra} 
a derived equivalence induces a stable
equivalence of Morita type.\index{stable!equivalence of Morita type} 
By a theorem of Happel (see for example 
Proposition~2.21.9 of Linckelmann~\cite{Linckelmann:2018a}),
for symmetric algebras,\index{symmetric algebra} 
derived equivalence also induces an
isomorphism in Hochschild cohomology. Therefore
reference~\cite{Generalov/Zilberbord:2016a} effectively computes
the dimensions of $\HH^n$ for all cases with three simple modules,
see also~\cite{Holm:2004a}.

Unfortunately,\index{errors} there 
are are some copying errors in~\cite{Erdmann:1990c,Erdmann:1990a}, and
an incorrect correction in~\cite{Craven:2019a}. 
It is erroneously
reported in statement (11.15)\,(c) of~\cite{Erdmann:1990c}
(incorrectly labelled (11.5)\,(c))
that the principal block of $M_{11}$ belongs to
family IV. In Table~1 of~\cite{Erdmann:1990c}, for family IV, $P_2$
should be ``as in I'' and not ``as in III''; the conditions for it to 
be a block should be $t=1$ and $k=2^{n-2}$, not the other way round. 
In the tables at the back of~\cite{Erdmann:1990a},
the principal blocks of $PSL(3,p^m)$ with $p^m\equiv 3 \pmod{4}$ 
are incorrectly assigned to
$SD(3\mathcal{B})_1$ rather than $SD(3\mathcal{D})$. In case
$SD(3{\mathcal K})$, the parameters should be $a\ge b \ge c\ge 1$, $a\ge 2$
rather than $a\ge b\ge c \ge 2$. 
On pages~143--144 of \cite{Craven:2019a}, the correction there incorrectly states that both
$M_{11}$ and $PSL(3,p^m)$ with $p^m \equiv 3 \pmod{4}$ belong to
family $SD(3\mathcal{B})_1$, and that there is only one simple module
with a non-trivial self-extension; in fact, the family is
$SD(3{\mathcal D})$, and there are two such simple modules.
\end{remarks}

\section{One class of involutions, one of order four}\label{se:SD1}

We begin with Case~\ref{case:SD1}, where $G$ has one class of
involutions, and one class of elements of order four.
In this case, $G$ has no normal subgroup of
index two, and Proposition~2.2
of~\cite{Alperin/Brauer/Gorenstein:1970a} implies that $G/O(G)$ contains a
simple normal subgroup with odd index.
By the main theorem of that paper, the simple groups with 
semidihedral Sylow $2$-subgroups are as follows.\medskip

(a) The projective special linear 
groups $PSL(3,p^m)$\index{PSL3pm@$PSL(3,p^m)$} 
with $p^m\equiv 3 \pmod{4}$.\medskip

(b) The projective special 
unitary groups $PSU(3,p^m)$\index{PSU3pm@$PSU(3,p^m)$} 
with $p^m\equiv 1 \pmod{4}$.\medskip

(c) The sporadic Mathieu group $M_{11}$\index{Mathieu group $M_{11}$} 
of order $7920$.\bigskip

Let $G$ be a finite group with semidihedral Sylow $2$-subgroups of
order $8q$ and no normal subgroups of index two, and let $\kk$ be 
a field of characteristic two. 
Let $B$ be the principal block of $\kG$.
The structure of the projective indecomposable
$B$-modules was determined by Erdmann~\cite{Erdmann:1979a}.

\begin{remark}
The one case not treated in~\cite{Erdmann:1979a} is $G=M_{11}$, which
was treated in the thesis of Schneider~\cite{Schneider:1979a}, and also in
unpublished work of Alperin.

The principal blocks of $M_{11}$ and $PSL(3,p^m)$ with $p^m\equiv 3 \pmod{4}$
are in family III of~\cite{Erdmann:1990c}, which is $SD(3\mathcal{D})$
of~\cite{Erdmann:1990a}. 
The principal blocks $PSU(3,p^m)$ with 
$p^m \equiv 1 \pmod{4}$ are in family II of~\cite{Erdmann:1990c},
which is $SD(3\mathcal{A})_1$ of~\cite{Erdmann:1990a}.
\end{remark}

\begin{remark}
Fortunately, Proposition~\ref{pr:BGtwohatSD} allows us to do the analysis
for just one group for each size $8q$ of semidihedral Sylow
$2$-subgroup. We choose to examine $PSL(3,p^m)$, where the $2$-part of
$p^m+1$ is $2q$.
\end{remark}

Let us look first at the cases of $PSL(3,3)$ and $M_{11}$, whose
principal blocks are Morita equivalent.
There are three isomorphism classes of simple $B$-modules, 
all self-dual, denoted
$\kk$, $\MM$ and $\NN$. These have dimensions $1$, $12$ and $26$ in
the case of $PSL(3,3)$, and dimensions
$1$, $44$ and $10$ in the case of $M_{11}$. 
Their projective covers are given by the following diagrams.
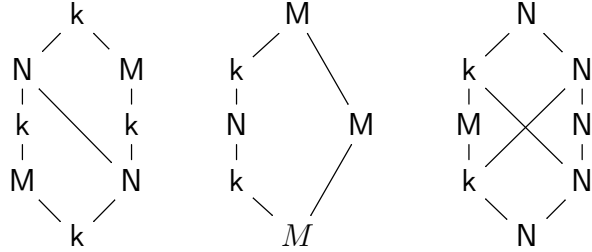
\[  \xymatrix@=2mm{&\kk\ar@{-}[dl]\ar@{-}[dr]\\
\NN\ar@{-}[d]\ar@{-}[ddrr]&&\MM\ar@{-}[d]\\
\kk\ar@{-}[d] &&\kk\ar@{-}[d]\\
\MM\ar@{-}[dr]&&\NN\ar@{-}[dl]\\
&\kk}\qquad
\xymatrix@=2mm{&\MM\ar@{-}[dl]\ar@{-}[ddr]\\
\kk\ar@{-}[d]\\
\NN\ar@{-}[d]&&\MM\ar@{-}[ddl]\\
\kk\ar@{-}[dr]\\
&M}\qquad
\xymatrix@=2mm{&\NN\ar@{-}[dl]\ar@{-}[dr]\\
\kk\ar@{-}[d]\ar@{-}[ddrr]&&\NN\ar@{-}[d]\ar@{-}[ddll]\\
\MM\ar@{-}[d]&&\NN\ar@{-}[d]\\
\kk\ar@{-}[dr]&&\NN\ar@{-}[dl]\\
&\NN} \]
Note that $\NN$ is periodic with period four, while $\kk$ and $\MM$ are not periodic.
The quiver\index{quiver} for $B$ is
\begin{equation}\label{eq:SD1-quiver} 
\xymatrix{\NN\ar@(ul,dl)_d\ar@/^/[r]^c & \kk \ar@/^/[r]^b
 \ar@/^/[l]^a& \MM\ar@/^/[l]^e\ar@(ur,dr)^f} 
\end{equation}
with relations 
\[ ef=0, \qquad be=0, \qquad fb=0,\qquad da=aeb,
\qquad cd=ebc, \qquad f^2=bcae,\qquad ac=d^3. \] 
This gives a presentation for the basic algebra of $B$.
This corresponds to the case discussed in Theorem~VIII.9.12
(with $k=1$, $s=4$, $t=2$)
and Proposition~IX.6.6\,(ii) (with $n=4$) of
Erdmann~\cite{Erdmann:1990a},

The unique self-dual grading (up to scalar multiples) 
on this quiver algebra is given by 
\begin{equation*} 
\textstyle |a|=|c|=\frac{3}{2},\quad |b|=|e|=\half, \quad
  |d|=1,\quad |f|=2.
\end{equation*}
We choose not to double these degrees, as the choice above 
makes the degrees in $H^*BG\cong \Ext^*_B(\kk,\kk)$ 
into integers with no common factor.

The principal blocks of the simple groups $PSL(3,p^m)$ with $p^m \equiv 3 \pmod{4}$ are
very similar, see Erdmann~\cite{Erdmann:1979a}. The only difference is
that if the $2$-part of $p^m+1$ is $4q$ (with $q$ a power of two) then
there are more repetitions of the simple module $\NN$ in its projective
cover. The case treated above is $q=2$, and in the general case there
are $2q-1$ copies of $\NN$ in the unserial module 
on the right hand side of the diagram instead of three.
So the relation $ac=d^3$ is replaced by $ac=d^{2q-1}$. The
Morita type of the principal block only depends on $q$, and not on
$p^m$. So for example, the principal blocks of $M_{11}$, $PSL(3,3)$,
$PSL(3,11)$ and $PSL(3,19)$ are Morita equivalent, with Sylow
$2$-subgroups of order $16$,  and the principal blocks of
$PSL(3,7)$ and $PSL(3,23)$ and $PSL(3,71)$ are all Morita equivalent,
with Sylow $2$-subgroups of order $32$. The grading also needs to be
adjusted, as follows.

\begin{theorem}
Let $G=PSL(3,p^m)$, where the $2$-part of $p^m+1$ is $2q$ ($q\ge 2$), 
or $G=M_{11}$ with $q=2$, and let
$\kk$ be a field of characteristic two. Then the basic algebra of the
principal block is given by the quiver~\eqref{eq:SD1-quiver},\index{quiver} with
relations
\[ ef=0, \qquad be=0, \qquad fb=0,\qquad da=aeb,
\qquad cd=ebc, \qquad f^2=bcae,\qquad ac=d^{2q-1}. \] 
The unique self dual grading, up to scalar multiples, on this algebra is given by
\begin{equation*} 
\textstyle |a|=|c|=q-\half,\qquad |b|=|e|=\half, 
\qquad |d|=1,\qquad |f|=q. 
\end{equation*}
\end{theorem}
\begin{proof}
The quiver with relations follows from the work of
Erdmann~\cite{Erdmann:1979a}. Given the relations, 
the uniqueness of the grading up to
scalar multiples is easy linear algebra. Self duality just means that $|a|=|c|$
and $|b|=|e|$.
\end{proof}

\begin{remark}
The given relations imply that
\begin{align*} 
f^3&=bcaef=0,\\
cac&=cd^{2q-1}=ebcd^{2q-2}=ebebcd^{2q-3}=0,\\
aca&=d^{2q-1}a=d^{2q-2}aeb=d^{2q-3}aebeb=0,\\
d^{2q+1}&=acd^2=aebcd=aebebc=0.
\end{align*}
\end{remark}

Let $\alpha$, $\beta$, $\gamma$, $\delta$, $\ep$, $\phi$ 
in $\Ext^1_B(\kk\oplus\MM\oplus\NN,\kk\oplus\MM\oplus\NN)$ be the
elements dual to $a$, $b$, $c$, $d$, $e$, $f$. These have degrees
\[ \textstyle |\alpha|=|\gamma|=(-1,-q+\half),\qquad
|\beta|=|\ep|=(-1,-\half),\qquad |\delta|=(-1,-1),\qquad
|\phi|=(-1,-q), \]

We can compute minimal resolutions of the simple modules as in~\cite{Benson/Carlson:1987a}, 
and the result is as follows when $q=2$. For larger values of $q$, the
only difference is that the chains of copies of $\NN$ in the
resolutions of $\kk$ and $\NN$ are longer.
\begin{gather*} 
\Omega(\kk)=\vcenter{
\xymatrix@=1mm{&\alpha&&\beta\\
&\NN\ar@{-}[dl]\ar@{-}[ddr]&&\MM\ar@{-}[d]\\
\kk\ar@{-}[d]&& &\kk\ar@{-}[dl]\\
\MM\ar@{-}[dr]&&\NN\ar@{-}[dl]\\ 
&\kk}}\quad
\Omega^2(\kk)=\vcenter{
\xymatrix@=1mm{&\delta\alpha&&\\
&\NN\ar@{-}[dl]\ar@{-}[ddr]\\ 
\NN\ar@{-}[d]&&&&\phi\beta\\
\NN\ar@{-}[dr]&&\kk\ar@{-}[dl]\ar@{-}[dr]&&\MM\ar@{-}[dl]\\
&\NN&&\MM}} \quad
\Omega^3(\kk)=\vcenter{
\xymatrix@=1mm{&&&\ep\phi\beta\\
&&&\kk\ar@{-}[d]\\
&\phi^2\beta&&\NN\ar@{-}[d]\\
&\MM\ar@{-}[dl]\ar@{-}[dr]&&\kk\ar@{-}[dl]\\
\kk\ar@{-}[d]&&\MM\\
\NN}}\\[-4mm]
\qquad
\Omega^4(\kk)=\vcenter{\xymatrix@=1mm{
&&&&\!\!\beta\ep\phi\beta\!\!\\
&&&&\MM\ar@{-}[d]\\
&&\phi^3\beta&&\kk\ar@{-}[d]\\
\kk\ar@{-}[dr]&&\MM\ar@{-}[dl]\ar@{-}[dr]&&\NN\ar@{-}[dl]\\
&\MM&&\kk}}\qquad
\Omega^5(\kk)=\vcenter{\xymatrix@=1mm{
&&&&&&&\!\!\!\ep\phi^3\beta\!\!\!&&\!\!\!\!\!\phi\beta\ep\phi\beta\!\!\!\!\!\\
&&&&&&&\kk\ar@{-}[dl]\ar@{-}[dr]&&\MM\ar@{-}[dl]\\
&&&&\!\!\phi^4\beta\!\!&&\NN\ar@{-}[d]&&\MM\\
&\NN\ar@{-}[dl]\ar@{-}[ddr]&&&\MM\ar@{-}[dl]\ar@{-}[dr]&&\kk\ar@{-}[dl]\\
\kk\ar@{-}[d]&&&\kk\ar@{-}[dl]&&\MM\\
\MM\ar@{-}[dr]&&\NN\ar@{-}[dl]\\
&\kk}}\\[-12mm]
\Omega^6(\kk)=\hspace{-2mm}\vcenter{\xymatrix@=1mm{
&&&&&&&&&\!\!\!\!(\ep\phi\beta)^2\!\!\!\!\\
&&&&&&&&&\kk\ar@{-}[d]\\
&&&&&&&\!\!\!\!\!\!\!\beta\ep\phi^3\beta\!\!\!\!\!\!\!&&\NN\ar@{-}[d]\\
&\NN\ar@{-}[dl]\ar@{-}[ddr]&&&&&&\MM\ar@{-}[dl]\ar@{-}[dr]&&\kk\ar@{-}[dl]\\
\NN\ar@{-}[d]&&&&\!\!\phi^5\beta\!\!&&\kk\ar@{-}[d]&&\MM\\
\NN\ar@{-}[dr]&&\kk\ar@{-}[dl]\ar@{-}[dr]&&\MM\ar@{-}[dl]\ar@{-}[dr]&&\NN\ar@{-}[dl]\\
&\NN&&\MM&&\kk
}}\qquad 
\Omega^7(\kk)=\hspace{-3mm}\vcenter{\xymatrix@=1mm{
&&&&&&&&\!\!\!\!\!\!\ep(\phi\beta\ep)^2\!\!\!\!\!\!\\
&&&&&&&&\MM\ar@{-}[d]\\
&&&&\!\!\!\!\ep\phi^5\beta\!\!\!\!&&\!\!\!\!\!\!\beta\ep\phi^3\beta\!\!\!\!\!\!&&\kk\ar@{-}[d]\\
&&&&\kk\ar@{-}[dl]\ar@{-}[dr]&&\MM\ar@{-}[dl]\ar@{-}[dr]&&\NN\ar@{-}[dl]\\
&\!\!\!\!\phi^6\beta\!\!\!\!&&\NN\ar@{-}[d]&&\MM&&\kk\\
&\MM\ar@{-}[dl]\ar@{-}[dr]&&\kk\ar@{-}[dl]\\
\kk\ar@{-}[d]&&\MM\\
\NN
}}\\ \hline
\end{gather*}\vspace{-10mm}
\begin{gather*}
\Omega(\MM)=\vcenter{\xymatrix@=1mm{\ep\\
\kk\ar@{-}[d]\\
\NN\ar@{-}[d]&&\phi\\
\kk\ar@{-}[dr]&&\MM\ar@{-}[dl]\\
&\MM}}\quad
\Omega^2(\MM)=\vcenter{\xymatrix@=1mm{
\beta\ep&&&&\ep\phi\\
\MM\ar@{-}[d]&&&&\kk\ar@{-}[d]\\
\kk\ar@{-}[d]&&\phi^2&&\NN\ar@{-}[d]\\
\NN\ar@{-}[dr]&&\MM\ar@{-}[dl]\ar@{-}[dr]&&\kk\ar@{-}[dl]\\
&\kk&&\MM}}\quad
\Omega^3(\MM)=\hspace{-5mm}\vcenter{\xymatrix@=1mm{
\!\phi\beta\ep\!&&\!\!\ep\phi^2\!\!&&&&&\!\beta\ep\phi\!\\
\MM\ar@{-}[dr]&&\kk\ar@{-}[dl]\ar@{-}[dr]&&&&&\MM\ar@{-}[d]\\
&\MM&&\NN\ar@{-}[d]&&\phi^3&&\kk\ar@{-}[d]\\
&&&\kk\ar@{-}[dr]&&\MM\ar@{-}[dl]\ar@{-}[dr]&&\NN\ar@{-}[dl]\\
&&&&\MM&&\kk}}\\[-8mm]
\Omega^4(\MM)=\vcenter{\xymatrix@=1mm{
\!\!\ep\phi\beta\ep\!\!\\
\kk\ar@{-}[d]\\
\NN\ar@{-}[d]&&\!\!\beta\ep\phi^2\!\!&&&&&&\ep\phi^3&&\!\!\!\phi\beta\ep\phi\!\!\!\\
\kk\ar@{-}[dr]&&\MM\ar@{-}[dl]\ar@{-}[dr]&&&&&&\kk\ar@{-}[dl]\ar@{-}[dr]&&\MM\ar@{-}[dl]\\
&\MM&&\kk\ar@{-}[d]&&\phi^4&&\NN\ar@{-}[d]&&\MM\\
&&&\NN\ar@{-}[dr]&&\MM\ar@{-}[dl]\ar@{-}[dr]&&\kk\ar@{-}[dl]\\
&&&&\kk&&\MM}}\qquad\quad
\\ \hline
\end{gather*}\vspace{-10mm}
\begin{gather*}
\Omega(\NN)=\vcenter{\xymatrix@=1mm{\gamma&&\delta\\
\kk\ar@{-}[d]\ar@{-}[ddrr]&&\NN\ar@{-}[d]\ar@{-}[ddll]\\
\MM\ar@{-}[d]&&\NN\ar@{-}[d]\\
\kk\ar@{-}[dr]&&\NN\ar@{-}[dl]\\
&\NN}}\quad\ 
\Omega^2(\NN)=\vcenter{\xymatrix@=1mm{
&\alpha\gamma&&\gamma\delta\\
&\NN\ar@{-}[dl]\ar@{-}[ddr]&&\kk\ar@{-}[ddl]\ar@{-}[dr]\\
\kk\ar@{-}[d]&&&&\MM\ar@{-}[d] \\
\MM\ar@{-}[dr]&&\NN\ar@{-}[dl]\ar@{-}[dr]&&\kk\ar@{-}[dl]\\
&\kk\ar@{-}&&\NN}}\quad\ 
\Omega^3(\NN)=\vcenter{\xymatrix@=1mm{&\alpha\gamma\delta\\
&\NN\ar@{-}[dl]\ar@{-}[dr]\\
\kk\ar@{-}[d]\ar@{-}[ddrr]&&\NN\ar@{-}[d]\ar@{-}[ddll]\\
\MM\ar@{-}[d]&&\NN\ar@{-}[d]\\
\kk&&\NN
}}\quad\ 
\Omega^4(\NN)=\NN
\end{gather*}

For all values of $q$, the minimal resolution of $\kk$ takes the form
{\small
\begin{multline*} 
\dots \to P_\MM \oplus P_\kk \oplus P_\MM \oplus P_\kk 
\xrightarrow{\left(\begin{smallmatrix}
\bar c\bar b & 0 & 0 & 0 \\
\bar f & \bar e & 0 & 0 \\
0 & \bar e\bar a\bar c & \bar f & 0 \\
0 & 0 & \bar a \bar c\bar b & \bar b
\end{smallmatrix}\right)}
P_\NN\oplus P_\MM \oplus P_\MM \oplus P_\kk 
\xrightarrow{\left(\begin{smallmatrix}
\bar d & 0 & 0 & 0 \\
\bar e\bar a & \bar f & 0 & 0 \\
0 & \bar a\bar c\bar b & \bar b & 0 \\
0 & 0 & \bar f & \bar e
\end{smallmatrix}\right)}
P_\NN \oplus P_\MM \oplus P_\kk \oplus P_\MM \qquad \\
\xrightarrow{\left(\begin{smallmatrix}
\bar a &\bar b & 0 & 0 \\
0 & \bar f & \bar e & 0 \\
0 & 0 & \bar e \bar a\bar c & \bar f
\end{smallmatrix}\right)}
P_\kk \oplus P_\MM \oplus P_\MM
\xrightarrow{\left(\begin{smallmatrix}
\bar e\bar a\bar c & \bar f & 0 \\ 0 & \bar a\bar c\bar b & \bar b
\end{smallmatrix}\right)}
P_\MM \oplus P_\kk
\xrightarrow{\left(\begin{smallmatrix}
\bar c\bar b & 0 \\ \bar f & \bar e 
\end{smallmatrix}\right)}
P_\NN \oplus P_\MM
\xrightarrow{\left(\begin{smallmatrix} 
\bar d & 0 \\ \bar e\bar a & \bar f
\end{smallmatrix}\right)} 
P_\NN \oplus P_\MM
\xrightarrow{(\bar a\ \bar b)} P_\kk 
\end{multline*}
}
This is the total complex of the following double complex:
{\small
\begin{equation}\label{eq:SD1-resolution} 
\vcenter{\xymatrix@=5mm{
&&&&P_\NN \ar[d]^{\bar a} &
P_\NN\ar[l]_{\bar d}\ar[d]^{\bar e\bar a} \\
&&&&P_\kk \ar[d]^{\bar e\bar a\bar c}&
P_\MM\ar[l]_{\bar b}\ar[d]^{\bar f} &
\cdots \\
&&P_\NN \ar[d]^{\bar a} 
&P_\NN \ar[l]_{\bar d}\ar[d]^{\bar e\bar a} &
P_\MM \ar[l]_{\bar c\bar b}\ar[d]^{\bar f}&
P_\MM \ar[l]_(.45){\bar f}\ar[d]^{\bar a\bar c\bar b}\\
&&P_\kk\ar[d]^{\bar e\bar a\bar c}& 
P_\MM\ar[l]_{\bar b}\ar[d]^{\bar f}&
P_\MM\ar[l]_(.45){\bar f}\ar[d]^{\bar a\bar c\bar b}&
P_\kk\ar[l]_{\bar e}\ar[d]^{\bar e\bar a\bar c} &
\cdots  \\
P_\NN\ar[d]^{\bar a} & 
P_\NN \ar[l]_{\bar d}\ar[d]^{\bar e\bar a} &
P_\MM\ar[l]_{\bar c\bar b}\ar[d]^{\bar f} & 
P_\MM\ar[l]_(.45){\bar f}\ar[d]^{\bar a\bar c\bar b} &
P_\kk \ar[l]_{\bar e} \ar[d]^{\bar e\bar a\bar c} &
P_\MM\ar[l]_(.45){\bar b}\ar[d]^{\bar f} \\
P_\kk &P_\MM\ar[l]_{\bar b} & 
P_\MM\ar[l]_(.45){\bar f} & P_\kk\ar[l]_{\bar e} &
P_\MM \ar[l]_(.45){\bar b}& P_\MM\ar[l]_(.45){\bar f} &
\cdots
} }
\end{equation}
}
Here, $\bar a$ is the element of $\Hom_B(P_\NN,P_\kk)$ opposite to $a$,
and so on, so that the barred variables satisfy the reverse of the relations
in the quiver.

The extensions $\alpha,\dots,\phi$ satisfy the following relations, which are easy to
verify using the grading and the minimal resolutions above:
\begin{gather}
\alpha\ep=0, \qquad
\beta\gamma=0,\qquad 
\gamma\alpha=0,\qquad 
\delta^2=0,\qquad 
\ep\beta=0,\notag\\ 
\beta\ep\phi^2=\phi^2\beta\ep,\qquad 
\ep\phi^2\beta=0,\qquad
\alpha\gamma\delta=\delta\alpha\gamma,\qquad
\gamma\delta\alpha=0, \notag\\ 
m_3(\alpha,\ep,\beta)=\delta\alpha,\qquad 
m_3(\beta,\gamma,\alpha)=0,\qquad
m_3(\gamma,\alpha,\ep)=0,\label{eq:H-rels}\\
m_{2q-1}(\delta,\dots,\delta)=\alpha\gamma,\qquad 
m_3(\ep,\beta,\gamma)=\gamma\delta,\qquad 
m_3(\gamma,\delta\alpha,\ep)=\ep\phi^2,\notag\\
m_3(\beta,\gamma\delta,\alpha)=\phi^2\beta,\qquad
m_4(\beta,\gamma,\alpha,\ep)=\phi^2,\notag\\
m_4(\gamma,\alpha,\ep,\phi^2\beta)=m_4(\ep\phi^2,\beta,\gamma,\alpha).\notag
\end{gather} 

\begin{remark}\label{rk:consistent}
The relation $\gamma\delta\alpha=0$ follows from the remaining
relations in two ways:
\begin{align*} 
\gamma\delta\alpha
&=\gamma m_3(\alpha,\ep,\beta)
=m_3(\gamma,\alpha,\ep)\beta=0,\\
\gamma\delta\alpha
&=m_3(\ep,\beta,\gamma)\alpha
=\ep m_3(\beta,\gamma,\alpha)=0.
\end{align*}
The last relation describes the unlabelled copy of $\kk$ at the top of the left
end of $\Omega^4(\kk)$. When postmultiplied by
$\ep$ or premultiplied by $\beta$, this relation follows from the remaining relations:
\begin{align*}
 m_4(\gamma,\alpha,\ep,\phi^2\beta)\ep
&=m_4(\gamma,\alpha,\ep,\phi^2\beta\ep)
=m_4(\gamma,\alpha,\ep,\beta\ep\phi^2)
=m_3(\gamma,m_3(\alpha,\ep,\beta),\ep\phi^2)\\
&=m_3(\gamma,\delta\alpha,\ep\phi^2)
=m_3(\gamma,\delta\alpha,\ep)\phi^2
=\ep\phi^4\\
&=\ep\phi^2m_4(\beta,\gamma,\alpha,\ep)
=m_4(\ep\phi^2,\beta,\gamma,\alpha)\ep,\\
\beta m_4(\gamma,\alpha,\ep,\phi^2\beta)
&=m_4(\beta,\gamma,\alpha,\ep)\phi^2\beta
=\phi^4\beta
=\phi^2m_3(\beta,\gamma\delta,\alpha)\\
&=m_3(\phi^2\beta,\gamma\delta,\alpha)
=m_3(\phi^2\beta,m_3(\ep,\beta,\gamma),\alpha)
=m_4(\phi^2\beta\ep,\beta,\gamma,\alpha)\\
&=m_4(\beta\ep\phi^2,\beta,\gamma,\alpha)
=\beta m_4(\ep\phi^2,\beta,\gamma,\alpha).
\end{align*}
\end{remark}

\begin{theorem}\label{th:H*M11}
Let $G=PSL(3,p^m)$, where the $2$-part of $p^m+1$ is $2q$ ($q\ge 2$), 
 or $G=M_{11}$ with $q=2$. The cohomology ring 
$H^*BG=\Ext^*_B(\kk,\kk)$ is generated by the commuting elements
\[ x=\ep\phi\beta,\qquad 
y=m_4(\gamma,\alpha,\ep,\phi^2\beta)
=m_4(\ep\phi^2,\beta,\gamma,\alpha),\qquad
z=\ep\phi^3\beta, \] 
subject to one relation:
\[ H^*BG =\Ext^*_B(\kk,\kk)= \kk[x,y,z]/(x^2y+z^2) \]
where $|x|=(-3,-q-1)$, $|y|=(-4,-4q)$ and $|z|=(-5,-3q-1)$. 
\end{theorem}
\begin{proof}
The structure of the cohomology ring of $M_{11}$ was
computed in~\cite{Benson/Carlson:1987a}, and is as above, if we ignore
the internal degrees. The principal blocks of $PSL(3,p^m)$ with
$p^m\equiv 3 \pmod{8}$ are Morita equivalent to that of $M_{11}$, and
therefore give the same answer. 
The analogous computation with possibly larger values of
$q$ gives exactly the same answer for $PSL(3,p^m)$ with $p^m\equiv 3
\pmod{4}$.  We show that the given
elements satisfy these relations, using the relations~\eqref{eq:H-rels}.
 We begin by observing (as in Remark~\ref{rk:consistent}) that
\begin{align*}
\beta y &= \beta m_4(\gamma,\alpha,\ep,\phi^2\beta)
=m_4(\beta,\gamma,\alpha,\ep)\phi^2\beta 
=\phi^4\beta,\\
y\ep &= m_4(\ep\phi^2,\beta,\gamma,\alpha)\ep
=\ep\phi^2 m_4(\beta,\gamma,\alpha,\ep)
=\ep\phi^4,
\end{align*}
and so 
\[ x^2y=(\ep\phi\beta\ep\phi)(\beta y)
=(\ep\phi\beta\ep\phi)(\phi^4\beta) =\ep\phi(\beta\ep\phi^2)\phi^3\beta
= \ep\phi(\phi^2\beta\ep)\phi^3\beta=(\ep\phi^3\beta)(\ep\phi^3\beta)
=z^2. \]
Commutativity is automatic for elements of $H^*BG$, but also follows
from the relations above:
\begin{align*}
yx&=(y\ep)(\phi\beta)=(\ep\phi^4)(\phi\beta)=(\ep\phi)(\phi^4\beta)
=(\ep\phi)(\beta y)=xy \\
zx&=(\ep\phi^3\beta)(\ep\phi\beta)=(\ep\phi)(\phi^2\beta\ep)(\phi\beta)
=(\ep\phi)(\beta\ep\phi^2)(\phi\beta)=(\ep\phi\beta)(\ep\phi^3\beta)=xz,\\
zy&=(\ep\phi^3)(\beta y)=(\ep\phi^3)(\phi^4\beta)=(\ep\phi^4)(\phi^3\beta)
=(y\ep)(\phi^3\beta)=yz.
\qedhere
\end{align*}
\end{proof}

\begin{remark}\label{rk:HPSU}
Since the homotopy type of $BG\twohat$ only depends on the Sylow $2$-subgroup and
the fusion, the cohomology ring is the same for $G=PSU(3,p^m)$ where
the $2$-part of $p^m-1$ is $2q$ ($q\ge 2$).

It would be possible, but not necessary for the currrent purposes, 
to do a similar analysis for $PSU(3,p^m)$ to that contained in this
section. The quiver\index{quiver}  in that case is as follows
\[ \xymatrix{\NN\ar@/^/[r]^c&\kk\ar@/^/[l]^a\ar@/^/[r]^b&\MM\ar@/^/[l]^e} \]
with relations
\[ be=0, \qquad aca=a(ebca)^{2q-1}eb,\qquad cac=(ebca)^{2q-1}ebc,
  \qquad acaca=0,\qquad cacac=0 \]
see Erdmann~\cite{Erdmann:1979a,Erdmann:1988a,Erdmann:1990c}. 
In particular in~\cite{Erdmann:1990c}, (11.15) refers us family~II of
that paper, which is described in
Theorems~(5.1) and where the quiver and relations are given in Theorem~(5.2).

This algebra admits a self-dual grading given by $|a|=|c|=q$,
$|b|=|e|=1-q$. The problem here, though, is that the method
of~\cite{Benson/Carlson:1987a} for computing with projective
resolutions doesn't really apply, and this makes the details of the
computations quite tedious.
\end{remark}
\section{\texorpdfstring{Ext and Hochschild cohomology over $H^*BG$}
{Ext and Hochschild cohomology over H*BG}}\label{se:Ext-HH-SD1}

Throughout this section, we are still working in Case~\ref{case:SD1}. 
So we let $G$ be a finite group with a semidihedral Sylow
$2$-subgroup of order $8q$ and no normal subgroup of index two,
and $\kk$ a field of characteristic two.
Our next task is to compute $\Ext^{**}_{H^*BG}(\kk,\kk)$ and $\HH^*H^*BG$ 
by applying Theorems~\ref{th:ExtRkk} and~\ref{th:HHR}. 
Recall that by Theorem~\ref{th:H*M11} and
Remark~\ref{rk:HPSU} we have
$H^*BG=\kk[x,y,z]/(x^2y+z^2)$ with $|x|=(-3,-q-1)$, $|y|=(-4,-4q)$ and
$|z|=(-5,-3q-1)$. Let $f=x^2y+z^2\in \kk[x,y,z]$. Then we have
\begin{gather*} 
\frac{\partial f}{\partial x} = 0, \qquad 
\frac{\partial f}{\partial y}=x^2,\qquad
\frac{\partial f}{\partial z}=0,\\
\frac{\partial^{(2)}f}{\partial x^2}=y,\qquad
\frac{\partial^{(2)}f}{\partial y^2}=0,\qquad
\frac{\partial^{(2)}f}{\partial z^2}=1,\\
\frac{\partial^2f}{\partial x \partial y}=0,\qquad
\frac{\partial^2f}{\partial x \partial z}=0,\qquad
\frac{\partial^2f}{\partial y \partial z}=0. 
\end{gather*}
Plugging these into Definition~\ref{def:Cliffq}, for the algebra
$\Cliffq$\index{Cliff@$\Cliffq$} we have variables
$\xx$, $\yy$, $\zz$ dual to $x$, $y$ and $z$ and $s$ dual to $f$. 
These have degrees $|\xx|=(-1,3,q+1)$, $|\yy|=(-1,4,4q)$,
$|\zz|=(-1,5,3q+1)$, $|s|=(-2,10,6q+2)$. 
Here, the first is the $\Ext$ degree, the second 
comes from the homological degree in $H^*BG$, and the third is the 
internal degree coming from the grading on the algebra $B$. So the
degrees of the generators of $H^*BG$ come out as $|x|=(0,-3,-q-1)$,
$|y|=(0,-4,-4q)$ and $|z|=(0,-5,-3q-1)$.
Then $s$ is central, and we have relations $\xx^2=ys$, $\yy^2=0$,
$\zz^2=s$, $\xx\yy+\yy\xx=0$,
 $\xx\zz+\zz\xx=0$, $\yy\zz+\zz\yy=0$.
The relation $\zz^2=s$ makes $s$ a redundant generator, and we end
up with
\begin{equation}\label{eq:CliffM11}
\Cliffq = H^*BG[\xx,\yy,\zz]/(\xx^2+y\zz^2,\yy^2).
\end{equation}
The differential is given by 
\begin{equation}\label{eq:Cliff-differential}
 d\xx=0,\qquad 
d\yy=x^2\zz^2,\qquad 
d\zz=0. 
\end{equation}

\begin{theorem}\label{th:ExtHBG}
We have
\[ \Ext^{**}_{H^*BG}(\kk,\kk) = \Lambda(\xx,\yy) \otimes \kk[\zz]. \] 
with degrees given by $|\xx|=(-1,3,q+1)$, $|\yy|=(-1,4,4q)$ and  
$|\zz|=(-1,5,3q+1)$. 
\end{theorem}
\begin{proof}
This follows from
Theorem~\ref{th:ExtRkk} and the computation~\eqref{eq:CliffM11} of $\Cliffq$.
\end{proof}

\begin{theorem}\label{th:HOmegaBGtwohat-SD1}
We have
\[ H_*\Omega BG\twohat = \Lambda(\xx,\yy)\otimes \kk[\zz] \]
with $|\xx|=(2,q+1)$, $|\yy|=(3,4q)$ and $|\zz|=(4,3q+1)$.
\end{theorem}
\begin{proof}
\!Theorem~\ref{th:ExtHBG} gives the $E_2$ page of
the Eilenberg--Moore spectral sequence~\eqref{eq:Cotor}%
\index{Eilenberg--Moore spectral sequence} 
\[ \Ext^{**}_{H^*BG}(\kk,\kk) \Rightarrow H_*\Omega BG\twohat. \]
There is no room for differentials, and there are no ungrading problems.
\end{proof}

\begin{remark}
This agrees with the answer given in Proposition~II.4.2.6 of 
Levi~\cite{Levi:1995a}.
\end{remark}

\begin{remark}
When we compute the spectral sequence
\[ \Ext^{**}_{H_*\Omega BG\twohat}(k,k) \Rightarrow H^*BG \] 
we get $E_2=\kk[x,y] \otimes \Lambda(z)$ with $|x|=(-1,-2,-q-1)$,
$|y|=(-1,-3,-4q)$ and $|z|=(-1,-4,-3q-1)$. There are no differentials, but
the relation $z^2=0$ then ungrades to give $z^2=x^2y$.
\end{remark}

\begin{theorem}\label{th:HHHBG-SD1}
The Hochschild cohomology\index{Hochschild cohomology} 
ring of $H^*BG$ is given by
\[ \HH^*H^*BG=H^*BG[\xx,\zz]/(\xx^2+y\zz^2,x^2\zz^2), \]
with
\begin{gather*}  
|x|=(0,-3,-q-1),\qquad |y|=(0,-4,-4q),\qquad |z|=(0,-5,-3q-1),\\
|\xx|=(-1,3,q+1),\qquad |\zz|=(-1,5,3q+1). 
\end{gather*}
\end{theorem}
\begin{proof}
This follows from~\eqref{eq:CliffM11}
and~\eqref{eq:Cliff-differential}, 
using Theorem~\ref{th:HHR}.
\end{proof}

\begin{proposition}\label{pr:-n,n-2,0}
There are no non-zero elements
of degree $(-n,n-2,0)$ in the Hochschild cohomology $\HH^*H^*BG$ with $n>2$.
\end{proposition}
\begin{proof}
By Theorem~\ref{th:HOmegaBGtwohat-SD1}, 
we have a $\kk$-basis for $\HH^*H^*BG$ consisting of the monomials
$x^{i_1}y^{i_2}z^{\ep_3}\xx^{\ep_1}\zz^{i_3}$ with either $i_1\le 1$
or $i_3\le 1$. 
Suppose that such a monomial has degree $(-n,n-2,0)$.
Comparing degrees, we have
\begin{align}
-n&=-\ep_1-i_3\label{eq:SD1-HH1} \\
n-2&=-3i_1-4i_2-5\ep_3+3\ep_1+5i_3,\label{eq:SD1-HH2}\\
0&=-(q+1)i_1-4qi_2-(3q+1)\ep_3+(q+1)\ep_1+(3q+1)i_3.\label{eq:SD1-HH3}
\end{align}
We shall show that there are no solutions in non-negative integers with $n>2$.

First we deal with the case $q=2$. In this case, equation~\eqref{eq:SD1-HH3} becomes
\begin{equation} 
0=-3i_1-8i_2-7\ep_3+3\ep_1+7i_3.\label{eq:SD1-HH4}
\end{equation}
Adding equations~\eqref{eq:SD1-HH2} and~\eqref{eq:SD1-HH4}, we get
\begin{equation} 
n-2=-6i_1-12i_2-12\ep_3+6\ep_1+12i_3,\label{eq:SD1-HH5}
\end{equation}
and so 
\[ n \equiv 2 \pmod{6}. \]
If instead, we add equations~\eqref{eq:SD1-HH1} and~\eqref{eq:SD1-HH2}
and subtract equation~\eqref{eq:SD1-HH4}, we get
$-2 = 4i_2+2\ep_3-\ep_1-3i_3$, 
or
\begin{equation} 
4i_2+2\ep_3=\ep_1+3i_3-2.\label{eq:SD1-HH6}
\end{equation}
So $\ep_1$ and $i_3$ determine $i_2$ and $\ep_3$, and then $i_1$.
Let $n=6a+2$, so that equation~\eqref{eq:SD1-HH5} becomes
\begin{equation} 
a=-i_1-2i_2+\ep_1+2i_3\ge 1.\label{eq:SD1-HH7}
\end{equation}
From equation~\eqref{eq:SD1-HH1}, we have 
$i_3=6a+2-\ep_1$. Then equation~\eqref{eq:SD1-HH6}
gives $4i_2+2\ep_3=\ep_1+18a+6-3\ep_1-2$, so 
\begin{equation*}
2i_2=9a+2-\ep_1-\ep_3.
\end{equation*}
Finally, plugging these values of $i_2$ and $i_3$ into equation~\eqref{eq:SD1-HH7} gives
\begin{align*}
i_1&=a-2i_2-2\ep_3+\ep_1+2i_3\\
&=a-9a-2+\ep_1+\ep_3-2\ep_3+\ep_1+12a+4-2\ep_1\\
&=4a+2-\ep_3.
\end{align*}
Since $a\ge 1$, we see that both $i_1$ and $i_3$ are greater than one,
which is a contradiction. This completes the case $q=2$.

Now suppose that $q> 2$. Reading 
equations~\eqref{eq:SD1-HH1}, \eqref{eq:SD1-HH2}, and~\eqref{eq:SD1-HH3} modulo four, 
we see that $\ep_3+i_1$ and $n=\ep_1+i_3$ are both even, and are
congruent modulo four. So if $n>2$ then $n\ge 4$, $i_3\ge 3$, and
hence $i_1\le 1$. So either $i_1=\ep_3=0$ or $i_1=\ep_3=1$.
Adding equations~\eqref{eq:SD1-HH1} and~\eqref{eq:SD1-HH2}, we get
\[ -2 = -3i_1-4i_2-5\ep_3+2\ep_1+4i_3. \]
Since $-3i_1-5\ep_3$ is divisible by four, we deduce that $\ep_1=1$.

In the case $i_1=\ep_3=0$, $\ep_1=1$ we get $i_3=i_2-1$, $n=i_2$.
Equation~\eqref{eq:SD1-HH3} becomes
\begin{align*}
0&=-4qi_2+(q+1)+(3q+1)(i_2-1) \\
&= (-q+1)i_2-2q
\end{align*}
so $i_2$ is not an integer, which is a contradiction.

In the case $i_1=\ep_3=1$, $\ep_1=1$ we get $i_3=i_2+1$,
$n=i_2+2$. Equation~\eqref{eq:SD1-HH3} becomes
\begin{align*}
0&=-(q+1)-4qi_2-(3q+1)+(q+1)+(3q+1)(i_2+1) \\
&=(-q+1)i_2
\end{align*}
and so $i_2=0$, $n=2$, again a contradiction. So for $q>2$ there
are no monomials of this form.
\end{proof}

\begin{theorem}\label{th:SD1-formal}
In Case~\ref{case:SD1}, 
with the grading inherited from the internal grading on the basic
algebra of $\kG$, the $A_\infty$ structure of $H^*BG$ is 
intrinsically formal.\index{intrinsically formal}
\end{theorem}
\begin{proof}
This follows from Propositions~\ref{pr:HH} and~\ref{pr:-n,n-2,0}.
\end{proof}

\begin{remark}\label{rk:SD1-formal}
Another proof of formality, but which does not give intrinsic formality, 
in Theorem~\ref{th:SD1-formal} 
is to notice that there are endomorphisms of the 
resolution~\eqref{eq:SD1-resolution}
representing $x$, $y$ and $z$, and strictly satisfying the
relation $x^2y=z^2$. The endomorphism representing $y$
just moves the whole diagram two places down and two places to the
left. For $x$, we move three places to the left, but then we have to
multiply by the following elements of the basic algebra in the
corresponding places in the double complex. It is easy to check that
this composite defines a map of double complexes.
\[ \xymatrix@=1mm{
&&&&\bar c\bar b & \bar c\\
&&&&1&1&\cdots\\
&&\bar c\bar b&\bar c&1&1\\
&&1&1&1&1&\cdots\\
\bar c\bar b&\bar c&1&1&1&1\\
1&1&1&1&1&1&\cdots} \]
Similarly, for $z$ we move one place down and four to the left, 
and compose with the same maps. This defines a quasi-isomorphism from 
the cohomology ring $\Ext^*_{\kG}(\kk,\kk)$ to the DG algebra 
$\End^*_{\kG}(P_\kk)$, which in turn is
quasi-isomorphic to the $A_\infty$ algebra $H^*BG$.
\end{remark}

\begin{corollary}\label{th:SD1-HHCBG}
In Case~\ref{case:SD1}, we have 
\[ \HH^*H^*BG\cong\HHinf^*H^*BG\cong\HHinf^*H_*\Omega BG\twohat. \]
\end{corollary}
\begin{proof}
The first isomorphism follows from Theorem~\ref{th:SD1-formal}, while the
second follows from Theorem~\ref{th:HHOmega}.
\end{proof}

\section{\texorpdfstring{$A_\infty$ structure of $H_*\Omega BG\twohat$}
{A∞ structure of H⁎ΩBGٛ₂}}\label{se:Ainfty-SD1}

We continue to work in Case~\ref{case:SD1}. So $G$ is a finite group
with a semidihedral Sylow $2$-subgroup of order $8q$ and no normal
subgroup of index two, and $\kk$ is a field of characteristic two.
Recall from Theorem~\ref{th:HOmegaBGtwohat-SD1} that 
\[ H_*\Omega BG\twohat = \Lambda(\hat x,\hat y)\otimes\kk[\hat z] \]
with $|\hat x|=(2,q+1)$, $|\hat y|=(3,4q)$, $|\hat z|=(4,3q+1)$. 

\begin{theorem}\label{th:SD1-HHHOmega}
We have
\[ \HH^*H_*\Omega BG\twohat = \kk[x,y,\zz]\otimes
  \Lambda(\xx,\yy,z) \]
with 
\begin{gather*}
|x|=(-1,-2,-q-1),\qquad |y|=(-1,-3,-4q),\qquad |z|=(-1,-4,-3q-1),\\
|\xx|=(0,2,q+1),\qquad\qquad |\yy|=(0,3,4q),\qquad\qquad |\zz|=(0,4,3q+1).
\end{gather*}
\end{theorem}
\begin{proof}
The algebra $H_*\Omega BG\twohat$ is a Koszul algebra with Koszul dual
\[ \kk[x,y]\otimes \Lambda(z) \]
with $|x|=(-2,-q-1)$, $|y|=(-3,-4q)$, $|z|=(-4,-3q-1)$. 
Now apply Theorem~\ref{th:Negron}. The differential there is zero, so
$\HH^*H_*\Omega BG\twohat$ is the tensor product of the algebra in
homological degree zero and the Koszul dual in degree $-1$.
\end{proof}

\begin{theorem}\label{th:SD1-m3}
In Case~\ref{case:SD1}, 
up to quasi-isomorphism, the maps $m_i$ in the $A_\infty$ structure 
on $H_*\Omega BG\twohat$ may be
taken to be the $\kk[\hat z]$-multilinear maps determined by 
\[
m_3(\hat x,\hat y,\hat x)=\hat z^2,\qquad
m_3(\hat x,\hat x\hat y,\hat x)=\hat x\hat z^2,\qquad
m_3(\hat y,\hat x,\hat x\hat y)=
m_3(\hat x\hat y,\hat x,\hat y)=\hat y\hat z^2, 
\]
and all $m_i$ with $i\ge 3$ vanish on all other triples of monomials not involving
$\hat z$. We have 
$m_3\circ m_3=0$
(Gerstenhaber's circle product).\index{Gerstenhaber!circle product}

This is the unique $A_\infty$ algebra structure on this algebra, such 
that the map $m_3$ represents the class $x^2y\hat z^2$ of
degree $(-3,1,0)$ in the Hochschild 
cohomology $\HH^*H_*\Omega BG\twohat$. 
\end{theorem}
\begin{proof}
Comparing 
Theorem~\ref{th:HOmegaBGtwohat-SD1} with 
Theorem~\ref{th:SD1-HHHOmega},
we see that in the spectral sequence
\[ \HH^*H_*\Omega BG\twohat \Rightarrow \HHinf^*H_*\Omega BG\twohat \]
we have $d^2(\hat y)=x^2\hat z^2$, and no further
differentials, and
\begin{equation}\label{eq:SD1-E3} 
E^3=E^\infty=\kk[x,y,\zz]/(x^2\zz^2) \otimes \Lambda(\xx,z). 
\end{equation}
The relation $\hat x^2=0$ ungrades to $\hat x^2=y\hat
z^2$, while the relation $z^2=0$ ungrades to $z^2=x^2y$. It follows
that $m_3$ on $\hat x$, $\hat x$
and $\hat y$ in $H_*\Omega BG\twohat$ in some order is a non-zero multiple of $\hat z^2$.

In degree $(-3,1,0)$, the Hochschild
cohomology $\HH^*H_*\Omega BG\twohat$ is one dimensional, and is
spanned by $x^2y\hat z^2$.
Since the Hochschild cocycle $m_3$ is only
well defined modulo coboundaries, we examine the values of the
coboundary of a $2$-cochain $f_2$ on these elements. For degree
reasons, we have 
$f_2(\hat x,\hat x)=f_2(\hat x,\hat y)=f_2(\hat y,\hat x)=0$. Let
$f_2(\hat x,\hat x\hat y)=\lambda \hat z^2$ and $f_2(\hat x\hat y,\hat
x)=\mu \hat z^2$. Then we have
\[ \HHd f_2(\hat x,\hat x,\hat y)=\lambda \hat z^2, \qquad
\HHd f_2(\hat x,\hat y,\hat x)=(\lambda+\mu)\hat z^2, \qquad
\HHd f_2(\hat y,\hat x,\hat x)=\mu \hat z^2. \]
Now everything is defined over $\bF_2$.
So working modulo these coboundaries, any assignment with
\[ m_3(\hat x,\hat x,\hat y)+m_3(\hat x,\hat y,\hat x)
+m_3(\hat y,\hat  x,\hat x)=\hat z^2 \]
is valid. For symmetry we take
$m_3(\hat x,\hat y,\hat x)=\hat z^2$ and
$m_3(\hat x,\hat x,\hat y)=m_3(\hat y,\hat x,\hat x)=0$.

Using the fact that $m_3$ is a Hochschild cocycle, and $\hat x\hat
y=\hat y\hat x$, we then have 
\begin{gather*}
m_3(\hat x,\hat x,\hat x\hat y)=0, \qquad
m_3(\hat x\hat y,\hat x,\hat x)=0,\qquad
m_3(\hat x\hat y,\hat y,\hat y)=0, \qquad
m_3(\hat y,\hat y,\hat x\hat y)=0, \\
m_3(\hat x,\hat x\hat y,\hat y)=0, \qquad
m_3(\hat y,\hat x\hat y,\hat x)=0,\qquad
m_3(\hat y,\hat x,\hat x\hat y) =
m_3(\hat y\hat x,\hat x,\hat y),\\
m_3(\hat x,\hat y,\hat x\hat y)=
m_3(\hat x\hat y,\hat y,\hat x),\qquad\qquad
m_3(\hat x,\hat x\hat y,\hat x) = 
\hat x\,m_3(\hat x,\hat y,\hat  x)=\hat x\hat z^2, \\
m_3(\hat x,\hat y,\hat x\hat y)+
m_3(\hat x\hat y,\hat x,\hat y)=\hat y\hat z^2,\qquad
m_3(\hat y\hat x,\hat y,\hat x)+
m_3(\hat y,\hat x,\hat y\hat x)=\hat y\hat z^2,
\end{gather*}
The $2$-cochain $f_2$ with $f_2(\hat x\hat y,\hat
x\hat y)=\hat y\hat z^2$, and $f_2=0$ on other monomials, has coboundary
$\HHd f_2(\hat x,\hat y,\hat x\hat y)=\hat y\hat z^2$. So adding a
multiple of $\HHd f_2$ to $m_3$, we can assume that $m_3(\hat x,\hat
y,\hat x\hat y)=0$. It then follows that 
\begin{gather*} 
m_3(\hat x,\hat y,\hat x\hat y)= 
m_3(\hat x\hat y,\hat y,\hat x)=0,\qquad
m_3(\hat y,\hat x,\hat x\hat y)=
m_3(\hat x\hat y,\hat x,\hat y)=\hat y\hat z^2,\\
m_3(\hat x\hat y,\hat x\hat y,\hat x)=
m_3(\hat x,\hat x\hat y,\hat x\hat y)=0,\qquad
m_3(\hat x\hat y,\hat x\hat y,\hat y)=
m_3(\hat y,\hat x\hat y,\hat x\hat y)=0,\\
m_3(\hat x\hat y,\hat x,\hat x\hat y)=
m_3(\hat x\hat y,\hat y,\hat x\hat y)=0,\qquad
m_3(\hat x\hat y,\hat x\hat y,\hat x\hat y)=0.
\end{gather*}

Now, it is straightforward to compute directly that the Gerstenhaber circle 
product\index{Gerstenhaber!circle product}
$m_3\circ m_3$ is the zero Hochschild cochain. 
All terms are zero except possibly those involving the
$5$-tuples $(\hat x,\hat y,\hat x,\hat y,\hat x)$ and
$(\hat y,\hat x,\hat y,\hat x,\hat y)$, 
and those obtained from these
by replacing copies of $\hat y$ with $\hat x\hat y$.
These all result in cancellations that give zero. There are twelve
such computations to make, all very similar, and two of them are:
\begin{align*} 
(m_3\circ m_3)&(\hat x,\hat y,\hat x,\hat x\hat y,\hat x)\\
&=m_3(m_3(\hat x,\hat y,\hat x),\hat x\hat y,\hat x)
+m_3(\hat x,m_3(\hat y,\hat x,\hat x\hat y),\hat x)
+m_3(\hat x,\hat y,m_3(\hat x,\hat x\hat y,\hat x)) \\
&=m_3(\hat z^2,\hat x\hat y,\hat x)
+m_3(\hat x,\hat y\hat z^2,\hat x)
+m_3(\hat x,\hat y,\hat x\hat z^2)
= 0+\hat z^4+\hat z^4 = 0,\\
(m_3\circ m_3)&(\hat y,\hat x,\hat x\hat y,\hat x,\hat x\hat y)\\
&=m_3(m_3(\hat y,\hat x,\hat x\hat y),\hat x,\hat x\hat y)
+m_3(\hat y,m_3(\hat x,\hat x\hat y,\hat x),\hat x\hat y)
+m_3(\hat y,\hat x,m_3(\hat x\hat y,\hat x,\hat x\hat y))\\
&=m_3(\hat y\hat z^2,\hat x,\hat x\hat y)
+m_3(\hat y,\hat x\hat z^2,\hat x\hat y)
+m_3(\hat y,\hat x,0)
=\hat y\hat z^2+\hat y\hat z^2+0=0.
\end{align*}
By
Proposition~\ref{pr:circle}, we have $\HHd m_4=m_3\circ m_3$,
so $m_4$ is a Hochschild cocycle. Since
there are no non-zero Hochschild classes in degree $(-4,2,0)$, this 
makes $m_4$ a coboundary, and so we can take $m_4=0$.
Then $m_4\circ m_3+m_3\circ m_4$ vanishes, and so $m_5$ is a
Hochschild cocycle. Since there are no non-zero classes in degree
$(-5,3,0)$, it follows that $m_5$ is a coboundary, and may hence be
taken to be zero.
We could continue this way, but eventually there are non-zero elements
of Hochschild cohomology in degree $(-n,n-2,0)$.
So instead, define an
$A_\infty$ algebra $\fa$ with with 
$H_*\fa\cong
H_*\Omega BG\twohat$, the same
structure maps as 
$H_*\Omega BG\twohat$ up to $m_3$, and $m_i=0$ for $i\ge 4$.
Then the Koszul dual $A_\infty$ algebra $\fb=\Hom_{\Db(\fa)}(\kk,\kk)$ has
homology isomorphic to $H^*BG$ as an associative algebra.
To see this, we  compute the spectral sequence
$\Ext^*_{H_*\fa}(\kk,\kk) \Rightarrow H_*\fb$. The map $m_3$
determines the $d^2$ differential in this spectral sequence, and so
the $E^3$ page is given by~\eqref{eq:SD1-E3}. There is no room for further
non-zero differentials or for ungrading problems, so this is also $H_*\fb$.

 By Theorem~\ref{th:SD1-formal},
$H^*BG$ is intrinsically formal, and so $\fb$ is quasi-isomorphic to
$H^*BG$ as an $A_\infty$ algebra. This implies that 
\begin{equation*}
\fa\simeq \Hom_{\Db(\fb)}(\kk,\kk)
\simeq \Hom_{\Db(H^*BG)}(\kk,\kk)
\simeq H_*\Omega BG\twohat.
\qedhere
\end{equation*} 
\end{proof}

\section{A differential graded model}\label{se:DG-SD1}

We continue to work with Case~\ref{case:SD1}.
Theorem~\ref{th:SD1-m3} suggests that there may be a nice DG algebra
quasi-isomorphic to the $A_\infty$ algebra $H_*\Omega BG\twohat$.
Since $H^*BG$ is formal, in order to produce
such an algebra, we look at endomorphisms of the minimal resolution
of $\kk$ over $H^*BG$. This resolution is eventually periodic of
period one, and takes the following form.
\begin{multline*} 
\cdots\to
(H^*BG)^4 \xrightarrow{
\left(\begin{smallmatrix}
z&xy&&\\x&z&&\\y&&z&xy\\&y&x&z
\end{smallmatrix}\right)}
(H^*BG)^4 \xrightarrow{
\left(\begin{smallmatrix}
z&xy&&\\x&z&&\\y&&z&xy\\&y&x&z
\end{smallmatrix}\right)}
(H^*BG)^4\\
\xrightarrow{
\left(\begin{smallmatrix}
x&z&&\\y&&z&xy\\&y&x&z
\end{smallmatrix}\right)}
(H^*BG)^3
\xrightarrow{
\left(\begin{smallmatrix}
&y&x&z
\end{smallmatrix}\right)}
H^*BG \to \kk. 
\end{multline*}
The map $\hat z$ shifts one to the right by the $4\times 4$ identity matrix, except 
at the right hand end:
\[ \cdots,\qquad
\left(\begin{smallmatrix}1&0&0&0\\0&1&0&0\\0&0&1&0\\0&0&0&1
\end{smallmatrix}\right),\qquad
\left(\begin{smallmatrix}0&1&0&0\\0&0&1&0\\0&0&0&1
\end{smallmatrix}\right),\qquad
\left(\begin{smallmatrix} 0&0&1 \end{smallmatrix}\right). \]
This endomorphism is in the centre of the endomorphism ring of the
resolution, and so we can regard everything as defined over $\kk[\hat
z]$.

Similarly, $\hat y$ is given by shifting to the right and using the matrices
\[ \cdots,\qquad
\left(\begin{smallmatrix}0&0&0&0\\0&0&0&0\\1&0&0&0\\0&1&0&0
\end{smallmatrix}\right),\qquad
\left(\begin{smallmatrix}0&0&0&0\\1&0&0&0\\0&1&0&0
\end{smallmatrix}\right),\qquad
\left(\begin{smallmatrix}1&0&0
\end{smallmatrix}\right),\]
and $\hat x$ is given by shifting to the right and using the matrices
\[ \cdots,\qquad
\left(\begin{smallmatrix}0&y&0&0\\1&0&0&0\\0&0&0&y\\0&0&1&0
\end{smallmatrix}\right),\qquad
\left(\begin{smallmatrix}1&0&0&0\\0&0&0&y\\0&0&1&0
\end{smallmatrix}\right),\qquad
\left(\begin{smallmatrix}0&1&0
\end{smallmatrix}\right).\]
These matrices commute, and satisfy $\hat y^2=0$, but $\hat x^2$ is
not zero, but rather $y\hat z^2$. So we find a homotopy $\xi$ from $\hat x^2$ to zero:
\[ \cdots,\qquad
\left(\begin{smallmatrix}0&0&1&0\\0&0&0&1\\0&0&0&0\\0&0&0&0
\end{smallmatrix}\right),\qquad
\left(\begin{smallmatrix}0&0&0&1\\0&0&0&0\\0&0&0&0
\end{smallmatrix}\right),\qquad
\left(\begin{smallmatrix}0&0&0
\end{smallmatrix}\right).
 \]
Then we have $\xi^2=0$, $d\xi=\hat x^2$, $\xi\hat x=\hat x\xi$, 
and $\xi\hat y+\hat y\xi=\hat z^2$.

\begin{theorem}\label{th:SD1-Q}
Let $Q$ be the DG algebra over $\kk[\hat z]$ generated by elements
$\hat x$, $\hat y$ and $\xi$ with
\begin{gather*} 
d\hat x=0,\quad d\hat y=0,\quad 
\hat y^2=0,\quad  
 \hat x\hat y=\hat y\hat x,\\
d\xi=\hat x^2,\quad 
\xi\hat x=\hat x\xi,\quad 
\xi^2=0,\quad \xi\hat y+\hat y\xi = \hat z^2, 
\end{gather*}
and with degrees 
\[ |\hat x|=(2,q+1),\quad |\hat y|=(3,4q),\quad
|\hat z|=(4,3q+1),\quad |\xi|=(5,2q+2). \]
Then $Q$ is quasi-isomorphic to the $A_\infty$ algebra $H_*\Omega BG\twohat$.
\end{theorem}
\begin{proof}
The algebra relations imply that this has a free $\kk[\hat z]$-basis
consisting of the elements $\hat x^i\hat y^{\ep_1}\xi^{\ep_2}$ with
$i\ge 0$, $\ep_1,\ep_2\in\{0,1\}$. The differential sends the basis elements
with $\ep_2=1$ bijectively to the basis elements with $i\ge 2$ and
$\ep_2=0$. So $H_*Q$ is the algebra $\kk[\hat z] \otimes \Lambda(\hat
x,\hat y)$, which is isomorphic to $H_*\Omega BG\twohat$. 
The $A_\infty$ structure on $H_*Q$ is not formal. Indeed, it is easy
to check that $m_3$ represents the Hochschild class $x^2y\hat z^2$.
By Theorem~\ref{th:SD1-m3}, there is a unique
$A_\infty$ structure on this algebra such that $m_3$ represents this
class. It follows that $Q$ is quasi-isomorphic to $H_*\Omega BG\twohat$.
\end{proof}

\begin{remark}
We can give an explicit quasi-isomorphism $H_*\Omega BG\twohat \to Q$
as follows. 
The map $f_1$ is the $k[\hat z]$-module homomorphism which 
sends each monomial $1$, $\hat x$, $\hat y$, $\hat x\hat y$ 
in $H_*\Omega BG\twohat$ 
to the monomial with the same name in $Q$. The map $f_2$ is given
by 
\[ f_2(\hat x,\hat x)=\xi,\qquad 
f_2(\hat x,\hat x\hat y)=\xi\hat y,\qquad 
f_2(\hat x\hat y,\hat x)=\hat y\xi, \]
and $f_2$ is zero on all other pairs of monomials. All higher $f_i$
are the zero map. From this information, we can inductively compute the
higher multiplications $m_i$ on $H_*\Omega BG\twohat$, and they agree
with those given in Theorem~\ref{th:SD1-m3}. For example, we have
\begin{align*}
m_2(f_1\otimes f_2-f_2\otimes f_1)(\hat x,\hat y,\hat x)
&=0,\\
f_2(1\otimes m_2-m_2\otimes 1)(\hat x,\hat y,\hat x)
&=\xi\hat y +\hat y \xi = \hat z^2,
\end{align*}
and so $m_3(\hat x,\hat y,\hat x)=\hat z^2$. 
\end{remark}

\section{\texorpdfstring{Duality for $Q[\hat z^{-1}]$-modules}
{Duality for Q[ẑ⁻¹]-modules}}\label{se:SD1-dual}\index{duality}

Continuing with Case~\ref{case:SD1}, 
by Theorem~\ref{th:SD1-Q}, $\hat z$ is central in $Q$. It therefore
makes sense to invert it and examine $Q[\hat z^{-1}]$ as an
algebra over the graded field $\kk[\hat z,\hat z^{-1}]$. This
parallels Section~\ref{se:D-dual}, so we give fewer details.

If $X$ is any $\kk[\hat z,\hat z^{-1}]$-module, we write
\[ X^*=\Hom_{\kk[\hat z,\hat z^{-1}]}(X,\kk[\hat z,\hat z^{-1}])
\cong \Hom_\kk(X,\kk). \]

\begin{proposition}
There is a quasi-isomorphism of $Q[\hat z^{-1}]$-bimodules
\[ Q[\hat z^{-1}]\simeq \Sigma Q[\hat z^{-1}]^*. \]
\end{proposition}
\begin{proof}
The proof is similar to the proof of Proposition~\ref{pr:Qs-1-qi-dual}. Consider
the basis of $Q[\hat z^{-1}]$ as a $\kk[\hat z,\hat z^{-1}]$-module
given by the monomials $\hat x^i\hat y^{\ep_1}\xi^{\ep_2}$. The
map of $\kk[\hat z,\hat z^{-1}]$-modules $Q[\hat z^{-1}]
\to \Sigma^{|\hat x\hat y|}Q[\hat z^{-1}]^*$ sending all monomials to
zero except that
\[ 1 \mapsto (\hat x\hat y\mapsto 1), \qquad
\hat x \mapsto (\hat y \mapsto 1), \qquad
\hat y \mapsto (\hat x\mapsto 1),\quad
\hat x\hat y \mapsto(1\mapsto 1) \]
and these maps take all other monomials to zero,
is easily checked to be a quasi-isomorphism of 
$Q[\hat z^{-1}]$-bimodules. Now $|\hat x\hat y|=5$, and $\hat z$ is a
periodicity generator\index{periodicity generator} of degree four, so 
$\Sigma^{|\hat x\hat y|}Q[\hat z^{-1}]^*\cong \Sigma Q[\hat z^{-1}]^*$.
\end{proof}

\begin{corollary}\label{co:SD1-duals}
If $X$ is a homotopically projective\index{homotopically projective} 
$Q[\hat z^{-1}]$-module then
we have a quasi-isomorphism
\[ \Hom_{Q[\hat z^{-1}]}(X,Q[\hat z^{-1}]) \simeq
\Sigma\Hom_{\kk[\hat z,\hat z^{-1}]}(X,\kk[\hat z,\hat z^{-1}]). \]
\end{corollary}
\begin{proof}
The proof is essentially the same as that of Corollary~\ref{co:duals}.
\end{proof}

\begin{theorem}\label{th:SD1-Qduality}
Let $X$ and $Y$ be $Q[\hat z^{-1}]$-modules, such that $X$ is
homotopically projective, and its image in $\Db(Q[\hat z^{-1}])$ is
compact. Then we have a duality
\[ \Hom_{Q[\hat z^{-1}]}(X,Y)^*\cong 
\Hom_{Q[\hat z^{-1}]}(Y,\Sigma^{-1}X). \]
\end{theorem}
\begin{proof}
The proof is essentially the same as that of Theorem~\ref{th:Qduality}.
\end{proof}

We have $\Dsg(H^*BG)\simeq \Dcsg(H_*\Omega BG\twohat)$. Now the
element $\hat z\in H_*\Omega BG\twohat$ is the Eisenbud operator for
the relation $x^2y=z^2$ in $H^*BG$. It comes from an element of
$\HHinf^*H_*\Omega BG\twohat$ with the same name. It follows that $\hat
z$ is central, and we may invert it to obtain an equivalence
\[ \Dsg(H^*BG)\simeq \Dcsg(H_*\Omega BG\twohat) \simeq 
\Db(H_*\Omega BG\twohat[\hat z^{-1}]). \]

\section{Two classes of involutions, one of elements of order four}\label{se:SD2}

We now turn to the  Case~\ref{case:SD2} of a finite group $G$ with semidihedral Sylow
$2$-subgroup $\SD$ of order $8q$, $q\ge 2$, with two classes of involutions
and one class of elements of order four. In this case, $G$ has a
normal subgroup $K$ of index two with generalised quaternion Sylow
$2$-subgroups, and $K$ has no normal subgroups of index two.

The principal blocks of this type are all in Erdmann's classes
$SD(2{\mathcal A})_1$ and $SD(2{\mathcal B})_2$. Class 
$SD(2{\mathcal B})_2$ turns out to be the easier to deal with, so we take 
$G=SL^\pm(2,p^m)$ with $p^m\equiv 3\pmod{4}$. The basic algebra 
of the principal block is given by the quiver
\begin{equation}\label{eq:SD2-quiver} 
\xymatrix{\kk\ar@(ul,dl)_a\ar@/^/[r]^b 
& \MM\ar@/^/[l]^c\ar@(ur,dr)^d} 
\end{equation}
with relations
\[ db=bacba,\qquad
cd=acbac,\qquad
bc=d^{2q-1},\qquad
a^2=0. \]
These imply that
\begin{align*}
d^2b&=dbacba=bacba^2cba=0,\\
cd^2&=acbacd=acba^2cbac=0,\\
bcb&=d^{2q-1}b=d^{2q-3}(d^2b)=0,\\
cbc&=cd^{2q-1}=(cd^2)d^{2q-3}=0.
\end{align*}
This admits a $\bZ$-grading with $|a|=1-q$, $|b|=|c|=q-\half$,
$|d|=1$. 

The structures of the projective indecomposables are as follows:
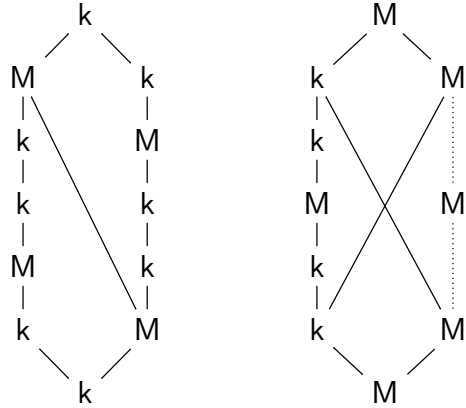
\[ \vcenter{\xymatrix@=3mm{
&\kk\ar@{-}[dl]\ar@{-}[dr] \\ 
\MM\ar@{-}[d]\ar@{-}[ddddrr]&&\kk\ar@{-}[d]\\
\kk\ar@{-}[d]&&\MM\ar@{-}[d]\\
\kk\ar@{-}[d]&&\kk\ar@{-}[d]\\
\MM\ar@{-}[d]&&\kk\ar@{-}[d]\\
\kk\ar@{-}[dr]&&\MM\ar@{-}[dl]\\
&\kk}}
\qquad\qquad
\vcenter{\xymatrix@=3mm{
&\MM\ar@{-}[dl]\ar@{-}[dr]\\
\kk\ar@{-}[d]\ar@{-}[ddddrr]&&\MM\ar@{.}[dd]\ar@{-}[ddddll]\\
\kk\ar@{-}[d]&&\\
\MM\ar@{-}[d]&&\MM\ar@{.}[dd]\\
\kk\ar@{-}[d]&&\\
\kk\ar@{-}[dr]&&\MM\ar@{-}[dl]\\
&\MM}} \]
Here, the case $q=2$ is as shown, and the dotted lines indicate that 
for $q>2$ there are more copies of $\MM$ in the right arm of $P_\MM$.

The minimal resolution of the trivial module may be computed
using the method of~\cite{Benson/Carlson:1987a}, and the result is
as follows.
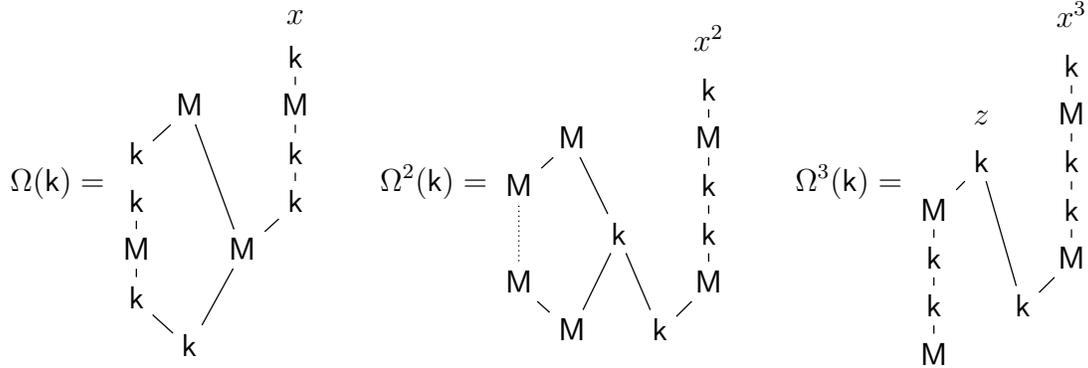
\begin{gather*} 
\Omega(\kk)=\vcenter{\xymatrix@=1mm{
&&&x\\
&&&\kk\ar@{-}[d]\\
&\MM\ar@{-}[dl]\ar@{-}[dddr]&&\MM\ar@{-}[d]\\
\kk&&&\kk\ar@{-}[d]\\
\kk\ar@{-}[d]&&&\kk\ar@{-}[dl]\\
\MM\ar@{-}[d]&&\MM\ar@{-}[ddl]\\
\kk\ar@{-}[dr]&&\\
&\kk&
}}\qquad
\Omega^2(\kk)=\vcenter{\xymatrix@=1mm{
&&&&x^2\\
&&&&\kk\ar@{-}[d]\\
&\MM\ar@{-}[dl]\ar@{-}[ddr]&&&\MM\ar@{-}[d]\\
\MM\ar@{.}[dd]&&&&\kk\ar@{-}[d]\\
&&\kk\ar@{-}[ddl]\ar@{-}[ddr]&&\kk\ar@{-}[d]\\
\MM\ar@{-}[dr]&&&&\MM\ar@{-}[dl]\\
&\MM&&\kk
}}\qquad
\Omega^3(\kk)=\vcenter{\xymatrix@=1mm{
&&&x^3\\
&&&\kk\ar@{-}[d]\\
&z&&\MM\ar@{-}[d]\\
&\kk\ar@{-}[dl]\ar@{-}[dddr]&&\kk\ar@{-}[d]\\
\MM\ar@{-}[d]&&&\kk\ar@{-}[d]\\
\kk\ar@{-}[d]&&&\MM\ar@{-}[dl]\\
\kk\ar@{-}[d]&&\kk\\
\MM
}}\vspace{-2cm} \\
\Omega^4(\kk)=\vcenter{\xymatrix@=1mm{
&&&&&x^4\\
&&&&&\kk\ar@{-}[d]\\
&&&xz&&\MM\ar@{-}[d]\\
&&&\kk\ar@{-}[dl]\ar@{-}[dddr]&&\kk\ar@{-}[d]\\
y&&\MM\ar@{-}[d]&&&\kk\ar@{-}[d]\\
\kk\ar@{-}[dddr]&&\kk\ar@{-}[d]&&&\MM\ar@{-}[dl]\\
&&\kk\ar@{-}[d]&&\kk\\
&&\MM\ar@{-}[dl]\\
&\kk
}}\qquad
\Omega^5(\kk)=\hspace{-3mm}
\vcenter{\xymatrix@=1mm{
&&&&&&&&&x^5\\
&&&&&&&&&\kk\ar@{-}[d]\\
&&&&&&&x^2z&&\MM\ar@{-}[d]\\
&&&&&&&\kk\ar@{-}[dddr]\ar@{-}[dl]&&\kk\ar@{-}[d]\\
&&&&xy&&\MM\ar@{-}[d]&&&\kk\ar@{-}[d]\\
&&&&\kk\ar@{-}[dl]\ar@{-}[dddr]&&\kk\ar@{-}[d]&&&\MM\ar@{-}[dl]\\
&\MM\ar@{-}[dl]\ar@{-}[dddr]&&\MM\ar@{-}[d]&&&\kk\ar@{-}[d]&&\kk\\
\kk\ar@{-}[d]&&&\kk\ar@{-}[d]&&&\MM\ar@{-}[dl]\\
\kk\ar@{-}[d]&&&\kk\ar@{-}[dl]&&\kk\\
\MM\ar@{-}[d]&&\MM\ar@{-}[ddl]\\
\kk\ar@{-}[dr]\\
&\MM
}}
\end{gather*}
The minimal resolution is the total complex of the following double
complex:
\begin{equation} 
\vcenter{\xymatrix@=6mm{
&&&&P_\MM\ar[d]^{\bar b}&
P_\MM\ar[l]_{\bar d}\ar[d]^{\bar b\bar c\bar a\bar b}\\
&&&&P_\kk\ar[d]^{\bar b\bar c\bar a\bar b\bar c}&
P_\kk\ar[l]_{\bar a}\ar[d]^{\bar b\bar c\bar a\bar b\bar c}&
\cdots\\
&&P_\MM\ar[d]^{\bar b}&
P_\MM\ar[l]_{\bar d}\ar[d]^{\bar b\bar c\bar a\bar b}&
P_\kk\ar[l]_{\bar c\bar a}\ar[d]^{\bar b\bar c\bar a\bar b\bar c}&
P_\kk\ar[l]_{\bar a}\ar[d]^{\bar b\bar c\bar a\bar b\bar c}\\
&&P_\kk\ar[d]^{\bar b\bar c\bar a\bar b\bar c}&
P_\kk\ar[l]_{\bar a}\ar[d]^{\bar b\bar c\bar a\bar b\bar c}&
P_\kk\ar[l]_{\bar a}\ar[d]^{\bar b\bar c\bar a\bar b\bar c}&
P_\kk\ar[l]_{\bar a}\ar[d]^{\bar b\bar c\bar a\bar b\bar c}&
\cdots\\
P_\MM\ar[d]^{\bar b}&
P_\MM\ar[l]_{\bar d}\ar[d]^{\bar b\bar c\bar a\bar b}&
P_\kk\ar[l]_{\bar c\bar a}\ar[d]^{\bar b\bar c\bar a\bar b\bar c}&
P_\kk\ar[l]_{\bar a}\ar[d]^{\bar b\bar c\bar a\bar b\bar c}&
P_\kk\ar[l]_{\bar a}\ar[d]^{\bar b\bar c\bar a\bar b\bar c}&
P_\kk\ar[l]_{\bar a}\ar[d]^{\bar b\bar c\bar a\bar b\bar c}\\
P_\kk&P_\kk\ar[l]_{\bar a}&P_\kk\ar[l]_{\bar a}&
P_\kk\ar[l]_{\bar a}&P_\kk\ar[l]_{\bar a}&
P_\kk\ar[l]_{\bar a}& \cdots}}\label{eq:SD2-resolution}
\end{equation}

The cohomology ring in this
case is therefore
\[ H^*BG = \kk[x,y,z]/(x^2y+z^2) \]
with 
\[ |x|=(-1,q-1),\qquad |y|=(-4,-4q), \qquad |z|=(-3,-q-1). \]
The situation is therefore very
similar to Case~\ref{case:SD1}. 

\section{Ext and Hochschild cohomology}\label{se:Ext-HH-SD2}

Continuing with Case~\ref{case:SD2}, 
the proofs of the following theorems are exactly as in 
the corresponding computations in Section~\ref{se:Ext-HH-SD1}
for Case~\ref{case:SD1}.

\begin{theorem}
We have 
\[ \Ext^{**}_{H^*BG}(\kk,\kk) = \Lambda(\hat x,\hat y) \otimes
\kk[\hat z] \]
with degrees given by $|\hat x|=(-1,1,1-q)$, 
$|\hat y|=(-1,4,4q)$ and $|\hat z|=(-1,3,q+1)$.\qed
\end{theorem}

\begin{theorem}\label{th:HOmegaBGtwohat-SD2}
We have 
\[ H_*\Omega BG\twohat = \Lambda(\hat x,\hat y) \otimes \kk[\hat z] \] 
with $|\hat x|=(0,1-q)$, $|\hat y|=(3,4q)$ and $|\hat z|=(2,q+1)$.\qed
\end{theorem}

\begin{theorem}\label{th:HHHBG-SD2}
We have
\[ \HH^*H^*BG=H^*BG[\hat x,\hat z]/(\hat x^2+y\hat z^2,x^2\hat z^2) \]
with $|x|=(0,-1,q-1)$, $|y|=(0,-4,-4q)$, $|z|=(0,-3,-q-1)$, $|\hat x|=(-1,1,1-q)$,
$|\hat z|=(-1,3,q+1)$.\qed
\end{theorem}

The proof of the following proposition follows along the same lines as the
proof of Proposition~\ref{pr:-n,n-2,0}, but the details are different, so we
spell them out.

\begin{proposition}\label{pr:-n,n-2,0again}
There are no non-zero elements of degree $(-n,n-2,0)$ in the
Hochschild cohomology $\HH^*H^*BG$ with $n>2$.
\end{proposition}
\begin{proof}
We have a $\kk$-basis for $\HH^*H^*BG$ consisting of the monomials
$x^{i_1}y^{i_2}z^{\ep_3}\hat x^{\ep_1}\hat z^{i_3}$ with either
$i_1\le 1$ or $i_3\le 1$. Suppose that such a monomial has degree
$(-n,n-2,0)$. Comparing degrees, we have
\begin{align}
-n&=-\ep_1-i_3\label{eq:SD2-HH1}\\
n-2&=-i_1-4i_2-3\ep_3+\ep_1+3i_3\label{eq:SD2-HH2}\\
0&=(q-1)i_1-4qi_2-(q+1)\ep_3-(q-1)\ep_1+(q+1)i_3.\label{eq:SD2-HH3}
\end{align}
We shall show that there are no solutions in non-negative integers
with $n>2$.

First we deal with the case $q=2$. In this case, equation~\eqref{eq:SD2-HH3}
becomes
\begin{equation} 
0=i_1-8i_2-3\ep_3-\ep_1+3i_3.\label{eq:SD2-HH4}
\end{equation}
Adding equations~\eqref{eq:SD2-HH2} and~\eqref{eq:SD2-HH4}, we get
\begin{equation} 
n-2=-12i_2-6\ep_3+6i_3,\label{eq:SD2-HH5}
\end{equation}
and so 
\[ n \equiv 2 \pmod{6}. \]
If instead, we add equations~\eqref{eq:SD2-HH1} 
and~\eqref{eq:SD2-HH2} and subtract equation~\eqref{eq:SD2-HH4},
we get $-2=4i_2+\ep_1-i_3$, or
\begin{equation} 
4i_2=-\ep_1+i_3-2.\label{eq:SD2-HH6}
\end{equation}
So $i_3$ determines $\ep_1$ and $i_2$.

Let $n=6a+2$, so that equation~\eqref{eq:SD2-HH5} gives
\begin{equation}
a=-2i_2-\ep_3+i_3 \ge 1.\label{eq:SD2-HH7}
\end{equation}
Equation~\eqref{eq:SD2-HH1} implies $i_3=6a+2-\ep_1$. 
Then equation~\eqref{eq:SD2-HH6} gives $i_2=(3a-\ep_1)/2$.
Plugging these values for $i_2$ and $i_3$ into
equation~\eqref{eq:SD2-HH7} gives
\[ a=-3a+\ep_1-\ep_3+6a+2-\ep_1=3a+2-\ep_3, \]
and so $\ep_3=2a+2$ is bigger than one. This contradiction completes
the case $q=2$.

Now suppose that $q> 2$. Reading 
equations~\eqref{eq:SD2-HH1}, \eqref{eq:SD2-HH2},
and~\eqref{eq:SD2-HH3} modulo four, we
see that $\ep_3+i_1$ and $n=\ep_1+i_3$ are both even, and are
congruent modulo four. Since $n>2$, we have $n\ge 4$, $i_3\ge 3$, and
hence $i_1\le 1$. So either $i_1=\ep_3=0$ or $i_1=\ep_3=1$. Adding the
equations~\eqref{eq:SD2-HH1} and~\eqref{eq:SD2-HH2}, we get
\begin{equation}
 -2 = -i_1-4i_2-3\ep_3+2i_3. \label{eq:SD2-HH8}
\end{equation}
Since $-i_1-3\ep_3$ is divisible by four, we deduce that $i_3$ is odd,
and hence $\ep_1=1$, and $n=1+i_3$.

Since $i_1=\ep_3$ and $\ep_1=1$, equation~\eqref{eq:SD2-HH3} becomes
\[ 0 = -4qi_2-2\ep_3-(q-1)+(q+1)i_3 \]
and equation~\eqref{eq:SD2-HH8} gives
\[ 2i_2=1-2\ep_3+i_3. \]
Substituting, we get
\begin{align*} 
0&=-2q(1-2\ep_3+i_3) -2\ep_3-(q-1)+(q+1)i_3\\
&=(1-q)i_3+(4q-2)\ep_3+1-3q
\end{align*}
and so 
\[ (4q-2)\ep_3=(q-1)i_3+(3q-1). \]
This is bigger than zero, so $\ep_3=1$, which then gives $i_3=1$. Then
by equation~\eqref{eq:SD2-HH1}, $n=\ep_1+i_3=2$, which is a contradiction.
\end{proof}

\begin{theorem}\label{th:SD2-formal}
In Case~\ref{case:SD2}, 
with the grading inherited from the internal grading on the basic
algebra of $\kG$, the $A_\infty$ structure of $H^*BG$ is formal.
\end{theorem}
\begin{proof}
This follows from Propositions~\ref{pr:HH} and~\ref{pr:-n,n-2,0again}.
\end{proof}

\begin{remark}\label{rk:SD2-formal}
Exactly as in Remark~\ref{rk:SD1-formal}, but with different degrees,
another proof of formality, but which does not give intrinsic formality, in
Theorem~\ref{th:SD2-formal} is to notice that there are endomorphisms
of the resolution~\eqref{eq:SD2-resolution} representing $x$, $y$ and
$z$, and strictly satisfying the relation $x^2y=z^2$. 
\end{remark}

\begin{corollary}\label{th:SD2-HHCBG}
In Case~\ref{case:SD2}, we have 
\[ \HH^*H^*BG\cong \HHinf^*H^*BG\cong \HHinf^*H_*\Omega BG\twohat. \]
\end{corollary}
\begin{proof}
The first isomorphism follows from Theorem~\ref{th:SD2-formal}, while
the second follows from Theorem~\ref{th:HHOmega}.
\end{proof}

\section{\texorpdfstring{$A_\infty$ structure, a DG model, and duality}
{A∞ structure, a DG model, and duality}}\label{se:Ainfty-SD2}\index{duality}

We continue to work in Case~\ref{case:SD2}, and because the details
are similar to those in Case~\ref{case:SD1}, we skip some details.
So $G$ is a finite group with a semidihedral Sylow $2$-subgroup of
order $8q$, and has a normal subgroup $K$ of index two with generalised
quaternion Sylow $2$-subgroups, and $K$ has no normal subgroup of index two.

\begin{theorem}\label{th:SD2-HHHOmega}
We have
\[ \HH^*H_*\Omega BG\twohat = \kk[x,y,\zz]\otimes
  \Lambda(\xx,\yy,z) \]
with 
\begin{gather*}
|x|=(-1,0,q-1),\qquad |y|=(-1,-3,-4q),\qquad |z|=(-1,-2,-q-1),\\
|\xx|=(0,0,1-q),\qquad\qquad |\yy|=(0,3,4q),\qquad\qquad |\zz|=(0,2,q+1).
\end{gather*}
\end{theorem}
\begin{proof}
This is a routine computation 
using Theorems~\ref{th:HHR} 
and~\ref{th:HOmegaBGtwohat-SD2}.
\end{proof}

\begin{theorem}\label{th:SD2-m3}
In Case~\ref{case:SD2}, up to quasi-isomorphism, the maps $m_i$ in the $A_\infty$ structure 
on $H_*\Omega BG\twohat$ may be
taken to be the $\kk[\hat z]$-multilinear maps determined by 
\[
m_3(\hat x,\hat y,\hat x)=\hat z^2,\qquad
m_3(\hat x,\hat x\hat y,\hat x)=\hat x\hat z^2,\qquad
m_3(\hat y,\hat x,\hat x\hat y)=
m_3(\hat x\hat y,\hat x,\hat y)=\hat y\hat z^2, 
\]
and all $m_i$ with $i\ge 3$ vanish on all other triples of monomials not involving
$\hat z$. We have 
$m_3\circ m_3=0$
(Gerstenhaber's circle product).\index{Gerstenhaber!circle product}

This is the unique $A_\infty$ algebra structure on this algebra, such 
that the map $m_3$ represents the class $x^2y\hat z^2$ of
degree $(-3,1,0)$ in the Hochschild 
cohomology $\HH^*H_*\Omega BG\twohat$. 
\end{theorem}
\begin{proof}
This is the same as the proof of Theorem~\ref{th:SD1-m3}.
\end{proof}

\begin{theorem}\label{th:SD2-Q}
Let $Q$ be the DG algebra over $\kk[\hat z]$ generated by elements
$\hat x$, $\hat y$ and $\xi$ with
\begin{gather*} 
d\hat x=0,\quad d\hat y=0,\quad 
\hat y^2=0,\quad  
 \hat x\hat y=\hat y\hat x,\\
d\xi=\hat x^2,\quad 
\xi\hat x=\hat x\xi,\quad 
\xi^2=0,\quad \xi\hat y+\hat y\xi = \hat z^2, 
\end{gather*}
and with degrees 
\[ |\hat x|=(0,q+1),\quad |\hat y|=(3,4q),\quad
|\hat z|=(2,3q+1),\quad |\xi|=(1,2q+2). \]
Then $Q$ is quasi-isomorphic to the $A_\infty$ algebra $H_*\Omega BG\twohat$.
\end{theorem}
\begin{proof}
This is proved in the same way as Theorem~\ref{th:SD1-Q}.
\end{proof}

Since $\hat z$ is central in $Q$, we may invert it. Let 
$Q[\hat z^{-1}]$ be the resulting DG algebra over 
$\kk[\hat z,\hat z^{-1}]$. 

\begin{proposition}
There is a quasi-isomorphism of $Q[\hat z^{-1}]$-bimodules
\[ Q[\hat z^{-1}]^*\cong \Sigma Q[\hat z^{-1}]. \]
\end{proposition}

\begin{corollary}
If $X$ is a homotopically projective $Q[\hat z^{-1}]$-module then we
have a quasi-isomorphism
\[ \Hom_{Q[\hat z^{-1}]}(X,Q[\hat z^{-1}])\simeq 
\Sigma \Hom_{\kk[\hat z,\hat z^{-1}]}(X,\kk[z,z^{-1}]). \]
\end{corollary}

\begin{theorem}
Let $X$ and $Y$ be $Q[\hat z^{-1}]$-modules, such that $X$ is
homotopically projective, and its image in $\Db(Q[\hat z^{-1}])$ is
compact. Then we have a duality
\[ \Hom_{Q[\hat z^{-1}]}(X,Y)^* \cong \Hom_{Q[\hat
    z^{-1}]}(Y,\Sigma^{-1}X). \]
\end{theorem}

\section{One class of involutions, two of elements of order four}\label{se:SD3}

Now we consider Case~\ref{case:SD3}, of a finite group $G$ with
semidihedral Sylow $2$-subgroup $\SD$ of order $8q$, $q\ge 2$, with one
class of involutions and two classes of elements of order four. In
this case, $G$ has a normal subgroup $K$ of index two with dihedral
Sylow $2$-subgroups, and $K$ has no normal subgroups of index two.
The group $K$ is therefore in Case~\ref{case:D1} of the
classification of groups with dihedral Sylow $2$-subgroups.

The principal blocks of this type are all in Erdmann's class
$SD(2{\mathcal A})_2$, see the tables at the back
of~\cite{Erdmann:1990a}. This causes a problem with socle relations. 
To see this, let us look at the principal
block $B_0$ of the group
$PGL^*(2,p^{2m})$.\index{PGL@$PGL^*(2,p^{2m})$} 
Let $\kk$ and $\MM$ be the two simple modules. 
\iftrue
In the
case $q=2$, their
projective covers are given by the following diagrams.
{\tiny
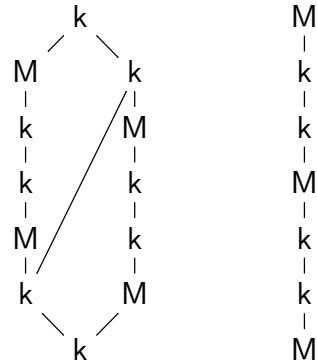
\[ \xymatrix@R=1mm@C=4mm{
&\kk\ar@{-}[dl]\ar@{-}[dr]\\
\MM\ar@{-}[d]&&\kk\ar@{-}[d]\ar@{-}[ddddddddddll]\\
\kk\ar@{-}[d]&&\MM\ar@{-}[d]\\
\kk\ar@{-}[d]&&\kk\ar@{-}[d]\\
\MM\ar@{-}[d]&&\kk\ar@{-}[d]\\
\kk\ar@{-}[d]&&\MM\ar@{-}[d]\\
\kk\ar@{-}[d]&&\kk\ar@{-}[d]\\
\MM\ar@{-}[d]&&\kk\ar@{-}[d]\\
\kk\ar@{-}[d]&&\MM\ar@{-}[d]\\
\kk\ar@{-}[d]&&\kk\ar@{-}[d]\\
\MM\ar@{-}[d]&&\kk\ar@{-}[d]\\
\kk\ar@{-}[dr]&&\MM\ar@{-}[dl]\\
&\kk}\qquad\qquad
\xymatrix@=1mm{
\MM\ar@{-}[d]\\
\kk\ar@{-}[d]\\
\kk\ar@{-}[d]\\
\MM\ar@{-}[d]\\
\kk\ar@{-}[d]\\
\kk\ar@{-}[d]\\
\MM\ar@{-}[d]\\
\kk\ar@{-}[d]\\
\kk\ar@{-}[d]\\
\MM\ar@{-}[d]\\
\kk\ar@{-}[d]\\
\kk\ar@{-}[d]\\
\MM} \]
}
\fi

The quiver\index{quiver} for $B_0$ is
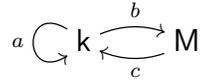
\begin{equation*}
\xymatrix{\kk\ar@(ul,dl)_a\ar@/^/[r]^b 
& \MM\ar@/^/[l]^c} 
\end{equation*}
with relations
\begin{equation}\label{eq:SD2A2}
bc=0,\qquad (cba)^{2q}=(acb)^{2q},\qquad a^2=cb(acb)^{2q-1} +
\lambda(cba)^{2q} 
\end{equation}
with $\lambda\in \kk$ unknown at this point. If $\lambda \ne 0$ then
any non-trivial grading on this algebra has $|a|=0$
and $|b|+|c|=0$, which then induces the trivial grading on cohomology.

We begin with the case $q=2$. In this case, we can choose 
$G=M_{10}=PGL^*(2,9)$. We can use 
{\sc Magma}~\cite{Bosma/Cannon/Playoust:1997a}\index{Magma@{\sc Magma}} to construct the basic
algebra as follows:
\begin{verbatim}
M11:=Group("M11");
M10:=Stabiliser(M11,1);
A:=BasicAlgebraOfPrincipalBlock(M10,GF(2));
\end{verbatim}\medskip
Beware that there is a random element in this construction, so
verifying the relations takes some care. In any case, 
we find that the socle constant $\lambda$ in the relations~\eqref{eq:SD2A2} for the algebra of
type $SD(2{\mathcal A})_2$ is equal to one in this case. It seems likely that this
is also true for higher values of $q$. 
It follows that there is no useful grading on the basic algebra, so we
are going to have to resort to other means.

Martino~\cite{Martino:1988a} computed the cohomology ring for groups in 
Case~\ref{case:SD3} to be
\[ H^*BG = \kk[y,z,w,v]/(y^3,vy,yz,v^2+z^2w) \]
with $|y|=-1$, $|z|=-3$, $|w|=-4$, $|v|=-5$. 
Part of the $A_\infty$ structure is given by
\begin{gather*}
m_3(z,y,y^2)=v, \qquad 
m_4(y^2,y,y^2,y)=w, \\
m_{4q-1}(v,y,v,y,\dots,y,v)= w^{2q}y^2.
\end{gather*}

The computation of the $\Ext$ ring is similar to the case of the
semidihedral group.

\begin{theorem}
In Case~\ref{case:SD3} we have
\[ \Ext^{*,*}_{H^*BG}(\kk,\kk)=\Lambda(\hat w) \otimes 
\kk\langle \hat y, \hat z,\hat v,\hat \eta\mid
\hat y^2=\hat z^2=0,\hat v\hat z = \hat z\hat v,
\eta\hat y=\hat y\eta\rangle \]
with $\eta=\langle\hat y,\hat y,\hat y\rangle$, $|\hat y|=(-1,1)$, 
$|\hat z|=(-1,3)$, $|\hat w|=(-1,4)$, $|\hat v|=(-1,5)$,
$|\eta|=(-2,3)$.
\end{theorem}

In the Eilenberg--Moore spectral sequence~\eqref{eq:Cotor}\index{Eilenberg--Moore spectral sequence}
\[ \Ext^{*,*}_{H^*BG}(\kk,\kk) \Rightarrow H_*\Omega BG\twohat \]
we have $d^2(\hat v) = \eta \hat z + \hat z \eta$,
\[ E^3=\Lambda(\hat w) \otimes \kk[\eta]\otimes \kk\langle\hat y,\hat z\mid
  \hat y^2=\hat z^2=0\rangle, \]
then $d^3(\hat w)=\eta^2$,
\[ E^4=E^\infty=\Lambda(\eta) \otimes \kk\langle\hat y,\hat z\mid \hat
  y^2=\hat z^2=0\rangle. \]
Ungrading, we have $|\hat y|=0$, $|\eta|=1$, $|\hat z|=2$,  and since
there are no lower terms in the filtration, we have
$\hat y^2=0$, and $\eta^2=0$. However, the relation $\hat z^2=0$
is harder to ungrade.

To compute the ring structure of $H_*\Omega BG\twohat$, 
and in particular the square of $\hat z$,
we resort to the method of 
squeezed resolutions\index{squeezed resolution} described
in Section~\ref{se:squeezed}. It is easy to compute the minimal (left) squeezed 
resolution for $G$, which is as follows.
\[ \cdots
\xrightarrow{\left(\begin{smallmatrix}
(\bar c\bar a\bar b)^{2q}&0\\0&(\bar c\bar a\bar b)^{2q}
\end{smallmatrix}\right)}
P_\MM \oplus P_\MM
\xrightarrow{\left(\begin{smallmatrix}
\bar c\bar a\bar b&0\\0&\bar c\bar a\bar b
\end{smallmatrix}\right)}
P_\MM \oplus P_\MM
\xrightarrow{\left(\begin{smallmatrix}
(\bar c\bar a\bar b)^{2q}&0\\0&(\bar c\bar a\bar b)^{2q}
\end{smallmatrix}\right)}
P_\MM \oplus P_\MM
\xrightarrow{\left(\begin{smallmatrix}\bar b&
\bar a\bar b\end{smallmatrix}\right)}
P_\kk \to 0. \]
After the first step, this repeats with period two. Indeed, after the
first step, it decomposes as a direct sum of two copies of the two-periodic complex
\[ \cdots \xrightarrow{(\bar c\bar a\bar b)^{2q}} P_\MM
\xrightarrow{\bar c\bar a \bar b} P_\MM
\xrightarrow{(\bar c\bar a\bar b)^{2q}} P_\MM. \]
The point here is that the $A_\infty$ algebra $H_*\Omega BG\twohat$ is quasi-isomorphic to
the $\kG$-endomorphism DG algebra of this squeezed resolution.

The element $\hat y\in H_*\Omega BG\twohat$ is 
represented by the map of complexes
\[ \xymatrix{\cdots\ar[r] &
P_\MM \oplus P_\MM \ar[r] 
\ar[d]^{\left(\begin{smallmatrix}0&0\\1&0\end{smallmatrix}\right)}& 
P_\MM \oplus P_\MM \ar[r] 
\ar[d]^{\left(\begin{smallmatrix}0&0\\1&0\end{smallmatrix}\right)}& 
P_\kk \ar[d]^{(\bar a)} \\
\cdots\ar[r] &
P_\MM\oplus P_\MM\ar[r] &
P_\MM\oplus P_\MM\ar[r] &P_\kk} \]
The square of this map is not zero, but is null homotopic,
with homotopy $u$ given by 
\[ \xymatrix@R=2.2cm@C=2.7cm{\cdots\ar[r] &
P_\MM \oplus P_\MM \ar[r] 
\ar[d]^{\left(\begin{smallmatrix}0&0\\0&0\end{smallmatrix}\right)}
\ar[dl]^{(\bar c\bar a\bar b)^{2q-1}\left(\begin{smallmatrix}0&1\\
0&\lambda\end{smallmatrix}\right)}& 
P_\MM \oplus P_\MM \ar[r] 
\ar[d]^{\left(\begin{smallmatrix}0&0\\0&0\end{smallmatrix}\right)}
\ar[dl]^{\left(\begin{smallmatrix}0&1\\
0&\lambda\end{smallmatrix}\right)}& 
P_\kk \ar[d]^{(\bar a^2)}
\ar[dl]^{(\bar c\bar a\bar b)^{2q-1}\bar c\left(\begin{smallmatrix} 1\\
\lambda\end{smallmatrix}\right)} \\
\cdots\ar[r] &
P_\MM\oplus P_\MM\ar[r] &
P_\MM\oplus P_\MM\ar[r] &P_\kk} \]
The element $\langle \hat y,\hat y,\hat y\rangle$ 
is represented by the map of complexes $u\hat y+\hat y u$:
\[ \xymatrix@C=1.2cm{\cdots \ar[r] &
P_\MM\oplus P_\MM \ar[r]
\ar[d]^{\left(\begin{smallmatrix}1&0\\
\lambda&1\end{smallmatrix}\right)} &
P_\MM\oplus P_\MM \ar[r]
\ar[d]^{(\bar c\bar a\bar b)^{2q-1}\left(\begin{smallmatrix}1&0\\
\lambda&1\end{smallmatrix}\right)} & 
P_\MM \oplus P_\MM \ar[r]
\ar[d]^{\left(\begin{smallmatrix}1&0\\
\lambda&1\end{smallmatrix}\right)} & P_\kk  
\ar[d]^{(\bar c\bar a\bar b)^{2q-1}\bar c\left(\begin{smallmatrix}
\bar a\\1+\lambda\bar a\end{smallmatrix}\right)} \\
\cdots \ar[r] &
P_\MM \oplus P_\MM \ar[r]
& P_\MM \oplus P_\MM \ar[r]
& P_\MM \oplus P_\MM \ar[r]
&P_\MM \oplus P_\MM \ar[r]
& P_\kk \oplus P_\kk
} \]
and the element $\hat z$ is represented by the map 
\[ \xymatrix{
\cdots\ar[r] &
P_\MM\oplus P_\MM\ar[r]
\ar[d]^{\left(\begin{smallmatrix}0&1\\0&0\end{smallmatrix}\right)} & 
P_\MM\oplus P_\MM\ar[r]
\ar[d]^{\left(\begin{smallmatrix}0&1\\0&0\end{smallmatrix}\right)} & 
P_\kk 
\ar[d]^{\left(\begin{smallmatrix}\bar c\\
0\end{smallmatrix}\right)}\\
\cdots\ar[r] &
P_\MM \oplus P_\MM \ar[r] &
P_\MM \oplus P_\MM \ar[r] &
P_\MM \oplus P_\MM\ar[r]& 
P_\MM\oplus P_\MM\ar[r]&P_\kk} \]
Now $u\hat y+\hat y u$ does not commute with $\hat z$, but $(u\hat
y+\hat y u)(1+\lambda \hat y)$ does, so it is more convenient to 
set 
\[ \eta = (u\hat y + \hat y u)(1+\lambda \hat y). \]
This is represented by the map
\[ \xymatrix@C=1.2cm{\cdots \ar[r] &
P_\MM\oplus P_\MM \ar[r]
\ar[d]^{\left(\begin{smallmatrix}1&0\\
0&1\end{smallmatrix}\right)} &
P_\MM\oplus P_\MM \ar[r]
\ar[d]^{(\bar c\bar a\bar b)^{2q-1}\left(\begin{smallmatrix}1&0\\
0&1\end{smallmatrix}\right)} & 
P_\MM \oplus P_\MM \ar[r]
\ar[d]^{\left(\begin{smallmatrix}1&0\\
0&1\end{smallmatrix}\right)} & P_\kk  
\ar[d]^{(\bar c\bar a\bar b)^{2q-1}\bar c\left(\begin{smallmatrix}
\bar a\\1\end{smallmatrix}\right)} \\
\cdots \ar[r] &
P_\MM \oplus P_\MM \ar[r]
& P_\MM \oplus P_\MM \ar[r]
& P_\MM \oplus P_\MM \ar[r]
&P_\MM \oplus P_\MM \ar[r]
& P_\kk \oplus P_\kk
} \]
Thus $\hat z^2=0$, and 
$\hat y\hat z+\hat z\hat y$ is the periodicity
generator\index{periodicity generator} of degree two,
central in the endomorphism algebra:
\[ \xymatrix{
\cdots\ar[r] &
P_\MM\oplus P_\MM\ar[r]
\ar[d]^{\left(\begin{smallmatrix}1&0\\0&1\end{smallmatrix}\right)} & 
P_\MM\oplus P_\MM\ar[r]
\ar[d]^{\left(\begin{smallmatrix}1&0\\0&1\end{smallmatrix}\right)} & 
P_\kk 
\ar[d]^{\left(\begin{smallmatrix}\bar c\bar a\\
\bar c\end{smallmatrix}\right)}\\
\cdots\ar[r] &
P_\MM \oplus P_\MM \ar[r] &
P_\MM \oplus P_\MM \ar[r] &
P_\MM \oplus P_\MM\ar[r]& 
P_\MM\oplus P_\MM\ar[r]&P_\kk} \]
To summarise, we have proved the following.

\begin{theorem}\label{th:SD3-HOmegaBGtwohat}
In Case~\ref{case:SD3}, we have
\[ H_*\Omega BG\twohat = \Lambda(\eta)\otimes 
\kk\langle \hat y,\hat z\mid \hat y^2=\hat z^2=0\rangle \]
with $|\eta|=1$, $|\hat y|=0$ and $|\hat z|=2$.\qed
\end{theorem}

It is interesting to note that the degree zero element $\hat y$ is not
central in $H_*\Omega BG\twohat$, while $\hat y\hat z+\hat z\hat y$ is
the central periodicity generator.\index{periodicity generator} 
This is very similar to what
happens for groups with dihedral Sylow $2$-subgroups in
Case~\ref{case:D2}.

Part of the $A_\infty$ structure is given by
$m_3(\hat y,\hat y,\hat y)= \eta(1+\lambda \hat y)$, 
$m_{2q}(\eta,\dots,\eta)=(\hat y\hat z+\hat z\hat y)^{2q}$.

\chapter{The generalised quaternion case}%
\index{quaternion group}\index{generalised quaternion group}%
\index{Sylow subgroup!generalised quaternion}

\section{Introduction}

In this chapter, we discuss the case of finite groups with generalised quaternion
Sylow $2$-subgroups. 
The theorem of Brauer and Suzuki~\cite{Brauer/Suzuki:1959a}
shows that there are no
simple groups of this type, but there are interesting perfect groups. 

We begin with the generalised quaternion group $\Q$ of order $8q$
itself. The group algebra in this case was analysed by Dade~\cite{Dade:1972a}, and
we describe a modified version of his presentation  as a quiver with
relations. If $q=1$ and $\kk$ contains $\bF_4$ then for suitable
radical generators we have
\[ \kQ = \kk\langle X,Y \mid X^2=YXY,\quad Y^2=XYX,\quad
  X^4=Y^4=0\rangle. \]
If $q\ge 2$, and $\kk$ is any field of characteristic two,
for suitable radical generators $X$ and $Y$ we have
\[ \kQ =\kk\langle X,Y\mid X^2=(YX)^{2q-1}Y+(XY)^{2q},\ 
Y^2=(XY)^{2q-1}X+(YX)^{2q},\
X^4=Y^4=0\rangle. \]
See Theorem~\ref{th:Q} for details.

There are three cases for the possible fusion in $\Q$, leading to
three types of cochains on the classifying space of a finite group
with this fusion. 
Probably the most interesting is the case where $G$
has no normal subgroup of index two. 

\begin{theorem}
Let $G$ be a non-2-nilpotent finite group with (generalised) quaternion Sylow
$2$-subgroup, and let $\kk$ be a field of
characteristic two. Then the following are equivalent:
\begin{enumerate}[label={\rm(\arabic*)}]
\item the $A_\infty$ algebra $H^*BG$  is 
formal,\index{formal $A_\infty$ algebra}
\item the $A_\infty$ algebra $H_*\Omega BG\twohat$ is formal,
\item $G$ has no normal subgroup of index two.
\end{enumerate}
\end{theorem}

The fact that (3) implies (1) and (2) follows from 
Theorem~\ref{th:Q1-formal} and
Corollary~\ref{co:Q1-formal}.
The converse follows from Corollary~\ref{co:HOmegaBG-quat-2classes}
and Theorem~\ref{th:Qu-non-formal}. 

Note also that if $G$ is
$2$-nilpotent with (generalised) quaternion Sylow $2$-subgroup
then by Theorem~\ref{th:p-group}, $H^*BG$ is not formal.

\section{Generalised quaternion groups}\label{se:Qu}

The generalised quaternion group
of order $8q$, $q$ a power of two, is given by the presentation
\[ \Q=\langle g,h\mid g^{2}=h^2=(g^{-1}h)^{2q} \rangle. \]
These relations imply that $g^2=h^2$ is central, and $g^4=h^4=1$.
If $q=1$, this is the quaternion group of order eight.

\begin{theorem}\label{th:Q}
We have the following presentations for $\kQ$.
\begin{enumerate}[label={\rm(\roman*)}]
\item
In the case $q=1$, suppose that
$\kk$ contains $\bF_4=\{0,1,\omega,\bar\omega\}$, with $1+\omega+\omega^2=0$.
Set 
\[ X=gh+\omega g+\bar\omega h,\qquad Y=gh+\bar\omega g+\omega h. \]
Then
\[ \kQ = \kk\langle X,Y \mid X^2=YXY,\quad Y^2=XYX,\quad
  X^4=Y^4=0\rangle. \]
The automorphism $g\mapsto h \mapsto gh\mapsto g$ of $\Q$ of order three sends
$X \mapsto \bar\omega X$ and $Y \mapsto \omega Y$.\label{th:Q/Q8}
\item
For $q\ge 2$, and any field $\kk$ of characteristic
two, set 
\[ u=g+h, \quad v=u^{4q-3}+\sum_{2^i=2}^qu^{2q-2^i}, \quad
 x=(g+1)+v, \quad y = (h+1)+v, \]
and finally, $X=x+(xy)^{2q-1}$, $Y=y+(yx)^{2q-1}$.
Then the group algebra has the following presentation:
\[ \kQ =\kk\langle X,Y\mid X^2=(YX)^{2q-1}Y+(XY)^{2q},\ 
Y^2=(XY)^{2q-1}X+(YX)^{2q},\
X^4=Y^4=0\rangle. \]
These relations imply that $(XY)^{2q}=X^3=(YX)^{2q}=Y^3$ is annihilated by $X$
and $Y$, and hence lie in $\Soc(\kQ)=J^{4q}(\kQ)$.\label{th:Q/Qge16}
\end{enumerate}
\end{theorem}
\begin{proof}
This follows Dade~\cite{Dade:1972a}, with a change of variables in the
second case.
To begin with, in both cases $X$ and $Y$ are in $J(\kQ)$ and are
linearly independent modulo $J^2(\kQ)$, and therefore generate $\kQ$.

\ref{th:Q/Q8} A somewhat long computation shows that 
$X^2=(1+g^2)(gh+\omega h+\bar\omega g)=YXY$, and hence $X^4=0$.
Applying the automorphism of $\bF_4$, we get 
$Y^2=XYX$ and $Y^4=0$. These relations imply that $\kQ$ is spanned by $1$, $X$,
$Y$, $XY$, $YX$, $XYX$, $YXY$ and $XYXY=YXYX$, so comparing
dimensions, these relations define $\kQ$.

\ref{th:Q/Qge16} It follows from Lemme~1.7 and Proposition~1.9
of~\cite{Dade:1972a} that the elements $x$ and $y$ in $J(\kQ)$ satisfy  
\[\kQ=\langle x,y\mid x^2=y^2=(xy)^{2q-1}x+(yx)^{2q-1}y+(xy)^{2q},\quad x^4=y^4=0\rangle. \]
These relations imply that $(xy)^{2q}=x^3=(yx)^{2q}=y^3$ 
spans $\Soc(\kQ)=J^{4q}(\kQ)$, $x^2=y^2$ is
central, and $J^{2q+1}(\kQ)=0$. Since $X$ and $Y$ are congruent to $x$
and $y$ modulo $J^{4q-2}(\kQ)$, it follows that monomials in $X$ and $Y$
of length at least three are equal to the corresponding monomials in
$x$ and $y$. So $X$ and $Y$ satisfy
\begin{align*}
X^2&=x^2+(xy)^{2q-1}x=(yx)^{2q-1}y+(xy)^{2q}=(YX)^{2q-1}Y+(XY)^{2q}\\
Y^2&=y^2+(yx)^{2q-1}x=(xy)^{2q-1}x+(yx)^{2q}=(XY)^{2q-1}Y+(YX)^{2q},
\end{align*}
and $X^4=Y^4=0$. Note that unlike $x^2$ and $y^2$, the elements 
$X^2$ and $Y^2$ are not central.
\end{proof}

\begin{remark}
It is erroneously stated on page~303 of~\cite{Erdmann:1990a}, 
page~38 of~\cite{Generalov:2007c}, and page~518
of~\cite{Generalov/Semenov:2020a} that the group
algebra of the generalised quaternion group is as given here, but
without the extra term $(XY)^{2q}$, $(YX)^{2q}$ in the expressions for
$X^2$ and $Y^2$.\index{errors}
\end{remark}

\begin{remark}\label{rk:Q8grading}
In the case of the quaternion group of order eight, over a field
containing $\bF_4$, there is a
$\bZ/3$-grading on the group algebra given by $|X|=1$ and $|Y|=-1$.
This is the grading induced by the automorphism of order three. 

In the case of the generalised quaternion groups of order at least
$16$, there is no non-trivial grading on the group algebra for which
the generators $X$ and $Y$ are homogeneous, because of the socle terms
in the relations.
\end{remark}

\begin{remark}
It is known that $\kQ$ has tame representation type, by embedding
$\Q$ into a semidihedral group of twice the order (cf.\
Section~\ref{se:isoclinism}). However, although the indecomposables 
for the semidihedral group are classified, nobody has yet been able to use
this to classify indecomposable modules for generalised quaternion
groups. For a further discussion of this, see the comments in the
proof of Theorem~\ref{th:tame}.
\end{remark}

The cohomology ring for $q=1$ is
\[ H^*BG=\kk[u,v,z]/(u^2+uv+v^2,u^2v+uv^2) \]
with $|u|=|v|=1$ and $|z|=4$, and with $u$ and $v$ dual to $U=\bar\omega
X+\omega Y=gh+h$ and $V=\omega X+\bar\omega Y=gh+g$.

For $q\ge 2$, we have
\[ H^*BG=\kk[x,y,z]/(xy,x^3+y^3), \]
again with $|x|=|y|=1$ and $|z|=4$, and with $x$ and $y$ 
dual to $X$ and $Y$. See for example
Rusin~\cite{Rusin:1987a} or Martino and Priddy~\cite{Martino/Priddy:1991a}.

Note that if $\kk$ contains $\bF_4$ then the cohomology of the
quaternion group of order eight can be made to fit the same pattern by
using the elements $x$, $y$ in $H^1BG$ dual to $X$ and $Y$ in
$J(\kQ)$. These are homogeneous with respect to the grading described
in Remark~\ref{rk:Q8grading}, so that the $\bZ\times\bZ/3$-grading is
given by $|x|=(-1,-1)$, $|y|=(-1,1)$ and $|z|=(-4,0)$.

\section{\texorpdfstring{$\HH^*H^*B\Q$}{HH*H*BQ}}

The cohomology ring $H^*B\Q=\kk[x,y,z]/(xy,x^3+y^3)$ 
is a complete intersection\index{complete intersection} of codimension
two, so we can calculate $\HH^*H^*B\Q$ and
$\Ext^{*,*}_{H^*B\Q}(\kk,\kk)$ using Theorems~\ref{th:HHR}
and~\ref{th:ExtRkk}. We first compute $\Cliffq$.\index{Cliff@$\Cliffq$}

\begin{proposition}\label{pr:Cliffq-quaterniongroup}
Let $\Q$ be a generalised quaternion group of order $8q$ with $q$ a
power of two. Let $\kk$ be a field of characteristic two, and if
$q=1$, we suppose that $\kk$ contains $\bF_4$. Then
the algebra $\Cliffq$ is equal to $H^*B\Q\langle 
\hat x,\hat y,\hat z;s_1,s_2\rangle$, where $s_1$ and $s_2$ are central, and 
\[ \hat x^2=\hat y^2=\hat z^2=0,\qquad 
\hat x\hat y+\hat y\hat x=s_1,\qquad 
\hat x\hat z=\hat z\hat x,\qquad  
\hat y\hat z=\hat z\hat y. \]
The degrees are given by $|\hat x|=|\hat y|=(-1,1)$, $|\hat
z|=(-1,4)$, $|s_1|=(-2,2)$, $|s_2|=(-2,3)$.
The differential $d$ on $\Cliffq$ is given by 
\[ d(\hat x)=ys_1+x^2s_2,\qquad
d(\hat y)=xs_1+y^2s_2,\qquad 
d(\hat z)=d(s_1)=
d(s_2)=0. \] 
\end{proposition}
\begin{proof}
Let $f_1(x,y,z)=xy$ and $f_2(x,y,z)=x^3+y^3$, so that
$H^*B\Q=\kk[x,y,z]/(f_1,f_2)$. 
Then we have
\[ \renewcommand{\arraystretch}{2.2}
\begin{array}{cccccc}
\displaystyle\frac{\partial f_1}{\partial x}=y,&
\displaystyle\frac{\partial f_1}{\partial y}=x, &
\displaystyle\frac{\partial f_1}{\partial z}=0,&
\displaystyle\frac{\partial f_2}{\partial x}=x^2,&
\displaystyle\frac{\partial f_2}{\partial y}=y^2,&
\displaystyle\frac{\partial f_2}{\partial z}=0,\\
\displaystyle\frac{\partial^{(2)}f_1}{\partial x^2}=0,&
\displaystyle\frac{\partial^{(2)}f_1}{\partial y^2}=0,&
\displaystyle\frac{\partial^{(2)}f_1}{\partial z^2}=0,&
\displaystyle\frac{\partial^{(2)}f_2}{\partial x^2}=x,&
\displaystyle\frac{\partial^{(2)}f_2}{\partial y^2}=y,&
\displaystyle\frac{\partial^{(2)}f_2}{\partial z^2}=0,\\
\displaystyle\frac{\partial^2 f_1}{\partial x\partial y}=1,&
\displaystyle\frac{\partial^2 f_1}{\partial x\partial z}=0,&
\displaystyle\frac{\partial^2 f_1}{\partial y\partial z}=0,&
\displaystyle\frac{\partial^2 f_2}{\partial x\partial y}=0,&
\displaystyle\frac{\partial^2 f_2}{\partial x\partial z}=0,&
\displaystyle\frac{\partial^2 f_2}{\partial y\partial z}=0.
\end{array} \]
Plugging these into Definition~\ref{def:Cliffq}, with $s_1$ and $s_2$
the degree $-2$ generators corresponding to the relations $f_1$ and
$f_2$, we get the given relations and differential for $\Cliffq$.
\end{proof}

\begin{remark}\label{rk:Cliffq-Q8}
In the case of the quaternion group of order eight without
assuming that $\kk$ contains $\bF_4$, the algebra $\Cliffq$ is
equal to $H^*B\Q\langle \hat u,\hat v,\hat z;s_1,s_2\rangle$,
where $s_1$ and $s_2$ are central, and
\[ \hat u^2=\hat v^2=\hat u\hat v+\hat v\hat u = s_1,\qquad
\hat u\hat z=\hat z\hat u, \qquad
\hat v\hat z=\hat z\hat v. \]
The degrees are given by 
$|\hat u|=|\hat v|=(-1,1)$, $|\hat
z|=(-1,4)$, $|s_1|=(-2,2)$, $|s_2|=(-2,3)$.
The differential $d$ on $\Cliffq$ is given by 
\[ d(\hat u)=vs_1+v^2s_2,\qquad
d(\hat v)=us_1+u^2s_2,\qquad 
d(\hat z)=d(s_1)=
d(s_2)=0. \] 
\end{remark}

\begin{theorem}\label{th:Q1-HHHBG}
The Hochschild cohomology\index{Hochschild cohomology}
$\HH^*H^*B\Q$ is 
\begin{multline*} 
H^*B\Q[\hat z,s_1,s_2,w_1,w_2]/
(x^2w_1+yw_2, 
y^2w_1+xw_2,\\
xs_1+y^2s_2, 
ys_1+x^2s_2, 
w_1s_1+w_2s_2,w_1^2,w_2^2,w_1w_2,\hat z^2). 
\end{multline*}
where
\[ w_1=x\hat x+y\hat y,\qquad 
w_2=y^2\hat x+x^2\hat y. \] 
The generators have degrees
$|\hat z|=(-1,4)$, $|s_1|=(-2,2)$, $|s_2|=(-2,3)$, 
$|w_1|=(-1,0)$, $|w_2|=(-1,-1)$. 
\end{theorem}
\begin{proof}
This follows from Theorem~\ref{th:HHR} and Proposition~\ref{pr:Cliffq-quaterniongroup}.
\end{proof}

\section{\texorpdfstring{Loops on $B\Q\twohat$}
{Loops on BQٛ₂}}\label{se:OmegaBQtwohat}

Since $\Q$ is a finite $2$-group, we have $\Omega B\Q\twohat\simeq
Q$. So we should expect to see the Eilenberg--Moore spectral 
sequence\index{Eilenberg--Moore spectral sequence}~\eqref{eq:Cotor}
converging to $\kQ$.

\begin{theorem}
If $\Q$ is a generalised quaternion group of order $8q$ with either  
$q\ge 2$ or $\kk$ containing $\bF_4$,  
then  
\[ \Ext^{*,*}_{H^*B\Q}(\kk,\kk) \cong \kk\langle \hat x,\hat y  
  \mid \hat x^2=\hat y^2=0\rangle \otimes \kk[\hat z,s]/(\hat z^2). \] 
The degrees are given by $|\hat x|=|\hat y| = (-1,1)$, 
$|\hat z|=(-1,4)$ and $|s|=(-2,3)$.  

 If $\Q$ is a quaternion group of order eight then
\[ \Ext^{*,*}_{H^*B\Q}(\kk,\kk) \cong \kk\langle \hat u,\hat v\mid
\hat u^2=\hat v^2 = \hat u\hat v + \hat v\hat u \rangle 
\otimes \kk[\hat z,s]/(\hat z^2). \]
The degrees are given by $|\hat u|=|\hat u| = (-1,1)$, $|\hat z|=(-1,4)$ and $|s|=(-2,3)$. 
\end{theorem}
\begin{proof}
In both cases, $H^*B\Q$ is a complete intersection, so we 
compute the $\Ext$ ring using Theorem~\ref{th:ExtRkk}. The algebra
$\Cliffq$ is given by Proposition~\ref{pr:Cliffq-quaterniongroup} and Remark~\ref{rk:Cliffq-Q8}.
The generator $s_1$ is redundant, so we eliminate it, and we write
$s$ for $s_2$.
\end{proof}

For $q=1$, the differentials in the Eilenberg--Moore spectral 
sequence~\eqref{eq:Cotor}\index{Eilenberg--Moore spectral sequence}
\[ \Ext^{*,*}_{H^*B\Q}(\kk,\kk) \Rightarrow \kQ \]
are given by $d^2(s)=\hat u^4=\hat v^4=(\hat u\hat v+\hat v\hat u)^2$ 
and $d^3(\hat z)=s^2$.

If $\kk$ contains $\bF_4$ then we can set $\hat x=\bar\omega\hat
u+\omega \hat v$
and $\hat y=\omega \hat u +\bar\omega \hat v$, 
so that $\Ext^{*,*}_{H^*B\Q}(\kk,\kk)$
becomes 
\[ \kk\langle \hat x,\hat y\mid \hat x^2=\hat y^2=0\rangle \otimes \kk[\hat
z,s]/(\hat z^2). \] 
We have $\hat u\hat v+\hat v\hat u=\hat x\hat y+\hat y\hat x$, 
and $d^2(s)=(\hat x\hat y+\hat y\hat x)^2$, and $d^2(s^2)=0$. Then
\[ E^3 = \kk\langle \hat x,\hat y \mid \hat x^2=\hat y^2=0,
(\hat x\hat y)^2=(\hat y\hat x)^2\rangle \otimes
  k[\hat z,s^2]/(\hat z^2). \]
The differential $d^3(\hat z)=s^2$ then gives
\[ E^4=E^\infty = \kk\langle \hat x,\hat y\mid
  \hat x^2=\hat y^2=0,(\hat x\hat y)^2=(\hat y\hat x)^2\rangle. \]
This is the associated graded of $\kQ$ with respect to the radical
filtration, with $\hat x$ representing $X$ and $\hat y$ representing $Y$.

The Eilenberg--Moore spectral sequence in the case $q\ge 2$ is
similar,  but the differentials happen in the opposite order. The
first differential is $d^3(\hat z) = s^2$, giving
\[ E^4= \kk\langle \hat x,\hat y\mid \hat x^2=\hat y^2=0\rangle \otimes k[s]/(s^2). \]
Then the next non-zero differential is 
$d^{4q-2}(s)=(\hat x\hat y+\hat y\hat x)^{2q}$, so
that
\[ E^{4q-1} = E^\infty = \kk\langle \hat x,\hat y \mid 
\hat x^2=\hat y^2=0, (\hat x\hat y)^{2q}=(\hat y\hat x)^{2q} \rangle. \]
This is again the associated graded of $\kQ$ with respect to the
radical filtration, with $\hat x$ representing $X$ and $\hat y$
representing $Y$.

\section{Isoclinism}\index{isoclinism}\label{se:isoclinism}

In order to describe the finite groups with generalised quaternion
Sylow $2$-subgroups in the next section, we first discuss
isoclinism and the groups $SL^\circ(2,p^m)$, which are isoclinic to
the groups $SL^\pm(2,p^m)$ described in
Section~\ref{se:semidihedral-Sylow}. We restrict ourselves to 
the situation we need, rather than describing isoclinism in general.
For the general definition, see Definition~4.28 in Suzuki~\cite{Suzuki:1986a}.

Suppose that $G$ is a finite group with a central subgroup of order
two, $Z=\{1,z\}$, contained in a normal subgroup $H$ of index 
$|G:H|=2$. Then we can make a new group of the same order as follows.
Consider the maps
\[ \bZ/2 \to G \times \bZ/4 \to \bZ/2, \]
where $\bZ/4=\langle \gamma\mid \gamma^4=1\rangle$, 
the first map sends the generator of $\bZ/2$ to $(z,\gamma^2)$,  
and the second map sends $G$ to $\bZ/2$ surjectively with kernel $H$,
and $\bZ/4$ to $\bZ/2$ with kernel $\langle\gamma^2\rangle$.
Let $\tilde G$ be the kernel of the second map modulo the image of the
first. Then $\tilde G$ has the same order as $G$, and is said to be
\emph{isoclinic} to $G$. The group $\tilde G$ has a normal subgroup of
index two isomorphic to $H$, and a central subgroup of order two
generated by $(z,1)$ with $\tilde G/\langle(z,1)\rangle\cong G/\langle
z\rangle$, but $G$ and $\tilde G$ are not in general isomorphic.

If $\rho\colon G \to GL(n,\bC)$ is a complex representation of $G$
with $z$ represented as minus the identity, 
then we can obtain a complex representation of $\tilde G$ by sending
elements of the subgroup of index two isomorphic to 
$H$ to the same matrices as before, but the elements outside are
multiplied by the complex number $i$. The character table of $\tilde
G$ therefore looks just like that of $G$ except that the character
values of elements outside $H$ 
on the representations with $z$ acting as minus the identity have been
multiplied by $i$. In particular, the character degrees are the same.

As an example, let $G$ be a semidihedral group of order $8q$, $q\ge 2$, with presentation
\[ G=\langle g,h\mid g^{4q}=1, h^2=1, hgh^{-1}=g^{2q-1} \rangle \]
as in Section~\ref{se:semidihedral-groups}. The element $z$ is
$g^{2q}$, and the normal subgroup $H$ of index two is the
(generalised) quaternion subgroup generated by
$g^2$ and $gh$. Then $\tilde G$ is
generated by $\tilde g=g\gamma$ and $\tilde h=h\gamma$,
with $g^{2q}$ identified with $\gamma^2=\tilde h^2$.
We have 
\[ \tilde h \tilde g\tilde h^{-1} = \gamma hgh^{-1}=\gamma
g^{2q-1}=\gamma^2\tilde g^{2q-1} = \tilde g^{-1}.\]
Thus 
\[ \tilde G = \langle \tilde g,\tilde h\mid \tilde g^{4q}=1,\tilde
  h^2=\tilde g^{2q}, \tilde h\tilde g\tilde h^{-1}=\tilde g^{-1} \rangle \]
which is a presentation for the generalised quaternion group of order
$8q$.
Thus the semidihedral group and generalised quaternion group of the
same order are isoclinic.

Now let $p^m$ be an odd prime power, and let us apply the same 
process to the groups $SL^\pm(2,p^m)$ and $SU^\pm(2,p^m)$ described in
Section~\ref{se:semidihedral-Sylow}. These have centres of order two,
and normal subgroups $SL(2,p^m)$ and $SU(2,p^m)$ of index two with
(generalised) quaternion Sylow $2$-subgroups. These data allow us to define
isoclinic groups, which we shall denote
$SL^\circ(2,p^m)$\index{SL@$SL^\circ(2,p^m)$}  and
$SU^\circ(2,p^m)$,\index{SU@$SU^\circ(2,p^m)$} 
of the same orders as $SL^\pm(2,p^m)$\index{SL@$SL^\pm(2,p^m)$} and
$SU^\pm(2,p^m)$.\index{SU@$SU^\pm(2,p^m)$} 
The groups $SL^\circ(2,p^m)$ with $p^m\equiv
3\pmod{4}$ and $SU^\circ(2,p^m)$ with $p^m\equiv 1\pmod{4}$ have
generalised quaternion Sylow $2$-subgroups, and will appear in
Case~\ref{case:Q2} in the next section.

\section{\texorpdfstring{Groups with generalised quaternion Sylow $2$-subgroups}
{Groups with generalised quaternion 
Sylow 2-subgroups}}\label{se:generalised-quaternion-Sylow}

Groups with generalised quaternion Sylow $2$-subgroups were classified
by Brauer and Suzuki~\cite{Brauer/Suzuki:1959a}, see also Section~VII
of Brauer~\cite{Brauer:1964a}, as well as Suzuki~\cite{Suzuki:1962a}, 
Glauberman~\cite{Glauberman:1974a}. The main theorem is that if $G$
has a generalised quaternion Sylow $2$-subgroup then the involution in
the centre of a Sylow $2$-subgroup has central image in $G/O(G)$ (Notation~\ref{no:gt}).
So the quotient of $G/O(G)$ by this central involution has dihedral Sylow
$2$-subgroups of order $4q$, and no odd order normal subgroups.
Such groups were analysed by Gorenstein and
Walter~\cite{Gorenstein/Walter:1965abc}. By Theorem~1.1 of Craven and
Glesser~\cite{Craven/Glesser:2012a}, these also represent the only
possible fusions systems on a generalised quaternion $2$-group.
As a consequence,  there are
three mutually exclusive possibilities for the fusion in $G$, listed
below.

The following case analysis comes from combining the Brauer--Suzuki
theorem with the Gorenstein--Walter classification of finite groups
with dihedral Sylow subgroups (see
Section~\ref{se:dihedral-Sylow}). The analysis of double covers of
those 
groups depends only on the fusion systems. The number of simple
modules in the principal block was analysed by
Erdmann~\cite{Erdmann:1988b,Erdmann:1988c,Erdmann:1988d}.\medskip 

\begin{case}\label{case:Q1}
If $G$ has one class of elements of order four then $G/O(G)$ is
isomorphic to either the double cover $2A_7$ of the alternating group
$A_7$,\index{alternating group $A_7$!double cover} or a subgroup of 
$\Gamma L(2,p^m)$\index{GammaL2pm@$\Gamma L(2,p^m)$} 
with $p^m$ a power of an odd
prime, containing $SL(2,p^m)$\index{SL@$SL(2,p^m)$} 
with odd index. The principal block of
$\kG$ has three isomorphism classes of simple modules.
\end{case}

\begin{case}\label{case:Q2}
If $G$ has two classes of elements of order four then $G$ has a normal
subgroup of index two, but no normal subgroup of index four. In this
case, $G/O(G)$ contains a normal subgroup of odd index
isomorphic to either $SL^\circ(2,p^m)$\index{SL@$SL^\circ(2,p^m)$} 
with $p^m\equiv 3\pmod{4}$ or
$SU^\circ(2,p^m)$\index{SU@$SU^\circ(2,p^m)$} with $p^m\equiv
1\pmod{4}$
(see Section~\ref{se:isoclinism}).
The principal block of $\kG$ has two isomorphism classes of
simple modules.
\end{case}

\begin{case}\label{case:Q3}
If $G$ has three classes of elements of order four then $O(G)$ is a
normal complement in $G$ to a Sylow $2$-subgroup $\Q$, so that
$G/O(G)\cong \Q$ and $H^*BG\cong H^*B\Q$. The principal block of $\kG$ is
isomorphic to $\kQ$, and has one isomorphism class of simple modules,
namely the trivial module.
\end{case}

\begin{proposition}
Suppose that $G$ has generalised quaternion Sylow $2$-subgroup $\Q$. Then the
homotopy type of $BG\twohat$ is determined by $|\Q|$ and the number of
conjugacy classes of elements of order four.
\end{proposition}
\begin{proof}
This follows from Theorem~\ref{th:Oliver} and the classification
theorem described above. 
\end{proof}

Representation theory and cohomology of groups with generalised quaternion Sylow
$2$-subgroups,\index{defect group!generalised quaternion} 
and more generally, of blocks with generalised quaternion defect
groups and finite dimensional algebras of quaternion type, are
discussed in Bleher~\cite{Bleher:2010a},
Bogdani\'c~\cite{Bogdanic:2021a,Bogdanic:2021b},
Cabanes and Picaronny~\cite{Cabanes/Picaronny:1992a},
Carlson, Mazza and Th\'evenaz~\cite{Carlson/Mazza/Thevenaz:2013a},
Eisele~\cite{Eisele:2016a},
Erdmann~\cite{Erdmann:1988b,Erdmann:1988c,   
Erdmann:1988d,Erdmann:1990a,Erdmann:1990b,  
Erdmann:1992a},   
Erdmann and Skowro\'nski~\cite{Erdmann/Skowronski:2019a},
Generalov et al.~\cite{Generalov:2007c,
Generalov:2009b,
Generalov/Ivanov/Ivanov:2008a,
Generalov/Semenov:2020a,
Generalov/Semenov:2024a},
Hayami~\cite{Hayami:2006a},
Holm~\cite{Holm:1999a},
Holm, Kessar and Linckelmann~\cite{Holm/Kessar/Linckelmann:2007a},
Ivanov et al.~\cite{Ivanov:2010a,Ivanov:2012a,Ivanov:2016a,
Ivanov/Ivanov/Volkov/Zhou:2015a}, 
Kawai and Sasaki~\cite{Kawai/Sasaki:2006a},
Kessar and Linckelmann~\cite{Kessar/Linckelmann:2002a},
Koshitani and Lassueur~\cite{Koshitani/Lassueur:2020b}, 
Langer~\cite{Langer:2014a},
Martino and Priddy~\cite{Martino/Priddy:1991a},
M\"uller~\cite{Muller:1974b},
Olsson~\cite{Olsson:1975a},
Semenov~\cite{Semenov:2025a},
Taillefer~\cite{Taillefer:2019a},
Zhou and Zimmermann~\cite{Zhou/Zimmermann:2011a}.
The homology of $\Omega BG\twohat$
was computed by Levi~\cite{Levi:1995a}.

\begin{remark}
Let $B$ be the principal block of $\kG$.
In Case~\ref{case:Q1}, one can put a $(\bZ\times\bZ)$-grading on
the basic algebra of $B$, and in Case~\ref{case:Q2}, one can
put a $\bZ$-grading on the basic algebra.
In Case~\ref{case:Q3}, there is no nontrivial grading. 
However, in all cases, these gradings are unhelpful, because
they induce the trivial grading on cohomology.
\end{remark}

We end this section with a table of the various cases of algebras of
quaternion type in characteristic two, in Erdmann's classification. We
note some minor misprints. In the appendix to~\cite{Erdmann:1988c}, the
entry $k+2$ in the Cartan matrix for type II should be $k+s$. In the
entry for type Q(3$\mathcal K$) in the tables at the end
of~\cite{Erdmann:1990a}, the last column should say $q\equiv 3$
mod~4 rather than $q\equiv 1$ mod~4. It seems unclear what happened to
type Q(2$\mathcal B$)$_2$ of~\cite{Erdmann:1990a} in the analysis
of~\cite{Erdmann:1988b}.\index{errors} \medskip

\begin{center}
\begin{tabular}{|c|c|c|c|c|c|} 
\hline
Erdmann~\cite{Erdmann:1990a}&\cite{Erdmann:1988b,Erdmann:1988c}&
Case&Group&$H^*$&$\HH^*$ \\ \hline
III.1(e)&&---&---&&\cite{Generalov:2007c} \\
III.1(e$'$)&&\ref{case:Q3}&$Q_{2^n}$&\cite{Munkholm:1969a}
&\cite{Generalov/Semenov:2020a,Hayami:2006a} \\
Q(2$\mathcal A$)&\cite{Erdmann:1988b} I&\ref{case:Q2}&
$SU^\circ(2,p^m),$&\cite{Martino/Priddy:1991a}& \\ 
&&&$p^m\equiv 1\pmod{4}$&&  \\
Q(2$\mathcal B$)$_1$&\cite{Erdmann:1988b} II& \ref{case:Q2}&
$SL^\circ(2,p^m),$&\cite{Martino/Priddy:1991a}&
\cite{Generalov:2009b,Generalov/Ivanov/Ivanov:2008a} \\
&&&$p^m\equiv 3\pmod{4}$&& \\
Q(2$\mathcal B$)$_2$&?&---&---&& \\
Q(2$\mathcal B$)$_3$&\cite{Erdmann:1988b} II ($k=1$)&---&---&& \\
Q(3$\mathcal A$)$_1$&\cite{Erdmann:1988c} II&---&---&& \\
Q(3$\mathcal A$)$_2$&\cite{Erdmann:1988c} III&\ref{case:Q1}&
$SL(2,p^m),$&\cite{Martino/Priddy:1991a}& \\
&&&$p^m\equiv 1\pmod{4}$&& \\
Q(3$\mathcal B$)&\cite{Erdmann:1988c} IV&\ref{case:Q1}&
$2A_7$&\cite{Martino/Priddy:1991a}& \\
Q(3$\mathcal C$)&\cite{Erdmann:1988c} I&---&---&& \\
Q(3$\mathcal D$)&\cite{Erdmann:1988c} V&---&---&& \\
Q(3$\mathcal K$)&\cite{Erdmann:1988c} VI&\ref{case:Q1}&
$SL(2,p^m),$&\cite{Martino/Priddy:1991a}& \\
&&&$p^m\equiv 3\pmod{4}$&& \\
\hline
\end{tabular}\medskip
\end{center}

\section{One class of elements of order four}\label{se:Q1}

Let $G$ be a finite group with quaternion or generalised 
quaternion Sylow $2$-subgroups of order
$8q$, and let $\kk$ be a field of
characteristic two. In this section 
we shall be interested in Case~\ref{case:Q1}, and our approach
will be to work directly with projective resolutions.\index{projective!resolution} Let us look
first at the case of $SL(2,3)\cong \Q_8\rtimes\bZ/3$, with $q=1$. 
There are three isomorphism classes of
simple $B$-modules, all of dimension one, which we shall denote $\kk$,
$\om$ and $\omb$. Their projective covers\index{projective!cover} 
are given by the following diagrams.
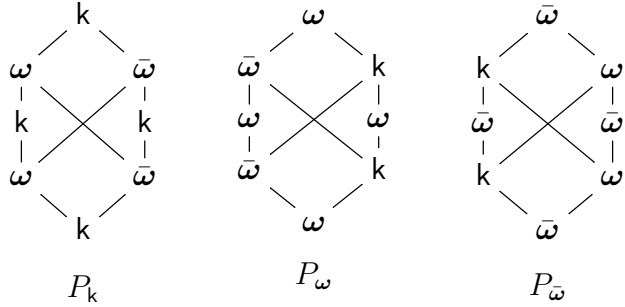
\[ \xymatrix@=2mm{
&\kk\ar@{-}[dl]\ar@{-}[dr] \\ 
\om\ar@{-}[d]\ar@{-}[ddrr]&&\omb\ar@{-}[d]\ar@{-}[ddll] \\
\kk\ar@{-}[d] && \kk\ar@{-}[d] \\
\om\ar@{-}[dr] &&\omb\ar@{-}[dl] \\
&\kk\\&P_\kk}
\qquad
\xymatrix@=2mm{
&\om\ar@{-}[dl]\ar@{-}[dr] \\ 
\omb\ar@{-}[d]\ar@{-}[ddrr]&&\kk\ar@{-}[d]\ar@{-}[ddll] \\
\om\ar@{-}[d] && \om\ar@{-}[d] \\
\omb\ar@{-}[dr] &&\kk\ar@{-}[dl] \\
&\om\\&P_\om}
\qquad
\xymatrix@=2mm{
&\omb\ar@{-}[dl]\ar@{-}[dr] \\ 
\kk\ar@{-}[d]\ar@{-}[ddrr]&&\om\ar@{-}[d]\ar@{-}[ddll] \\
\omb\ar@{-}[d] && \omb\ar@{-}[d] \\
\kk\ar@{-}[dr] &&\om\ar@{-}[dl] \\
&\omb\\&P_\omb} \]
Note that $\kk$, $\om$ and $\omb$ are all periodic with period four.
The quiver\index{quiver} for $B$ is
\begin{equation}\label{eq:2A4-quiver} 
\vcenter{\xymatrix@R=16mm@C=9mm{&\kk\ar@/^/[dr]^a\ar@/^/[dl]^d \\
\om\ar@/^/[ur]^c\ar@/^/[rr]^f&&\omb\ar@/^/[ll]^e\ar@/^/[ul]^b}}
\end{equation}
with relations
\begin{gather*} 
aba=fd,\qquad 
cdc=bf,\qquad
efe=db,\\
bab=ce,\qquad  
dcd=ea,\qquad 
fef=ac,\\
abf=0, \qquad
cdb=0,\qquad
efd=0.
\end{gather*}
These relations imply that 
\begin{equation}\label{eq:3-rels}
acd=bac=bfe=cef=dba=dce=eab=fdc=fea=0, 
\end{equation}
as well as 
\begin{equation}\label{eq:3-rels2}
cea=bfd,\qquad eac=dbf,\qquad ace=fdb,
\end{equation}
and that the composite of any
five arrows is zero.

For larger values of $q$, we can choose a prime power $p^m \equiv 3 \pmod{4}$
such that the $2$-part of $p^m+1$ is $4q$, and take $G=SL(2,p^m)$. In
this case, we label the two non-trivial simple modules $\MM$ and $\NN$
rather than $\om$ and $\omb$.
By~(1.3) of~\cite{Erdmann:1992a} 
(see also Theorem VII.8.8 of~\cite{Erdmann:1990a}), 
the structures of the projectives\index{projective!module}
are similar to the above, but longer. The quiver is the same, but the
relations are as follows:
\begin{gather*} 
(ab)^{2q-1}a=fd,\qquad 
(cd)^{2q-1}c=bf,\qquad
(ef)^{2q-1}e=db,\\
(ba)^{2q-1}b=ce,\qquad  
(dc)^{2q-1}d=ea,\qquad 
(fe)^{2q-1}f=ac,\\
abf=0, \qquad
cdb=0,\qquad
efd=0.
\end{gather*}
Again, these imply relations~\eqref{eq:3-rels} and~\eqref{eq:3-rels2}, and that 
the composite of any $4q+1$ arrows is zero.

We have 
\[ B \cong \End_B(B)^{\mathsf{op}}=\End_B(P_\kk\oplus P_\MM\oplus P_\NN)^{\mathsf{op}}, \] 
and we write $\bar a$ for the element of
$\Hom_B(P_\NN,P_\kk)$ opposite to $a$, and so on. The relations
satisfied by these are obtained by reversing those satisfied by the
original elements. 

\begin{theorem}\label{th:Q1-formal}
Let $G$ be a finite group with quaternion or generalised quaternion Sylow $2$-subgroups and no 
normal subgroup of index two. Then 
\[ H^*BG = \Lambda(y) \otimes k[z] \]
is a formal $A_\infty$ algebra.\index{formal $A_\infty$ algebra}
The degrees of the generators are $|y|=-3$ and $|z|=-4$.
\end{theorem}
\begin{proof}
The minimal resolution $P_*$ of the trivial module over $B$ is given by
\[ \cdots \xrightarrow{(\bar d,\bar a)}
P_\kk \xrightarrow{\bar a\bar e \bar c\,=\,\bar d\bar f\bar b}  
P_\kk \xrightarrow{\left(\begin{smallmatrix}
\bar c\\\bar b\end{smallmatrix}\right)}
P_\MM \oplus P_\NN 
\xrightarrow{\left(\begin{smallmatrix}
\bar c\bar d&\bar f\\\bar e&\bar b\bar a\end{smallmatrix}\right)}
P_\MM \oplus P_\NN\xrightarrow{(\bar d,\bar a)}
  P_\kk \]
The cohomology $H^*BG\cong\Ext^*_B(\kk,\kk)$ as an algebra is easily
read off from this, and is as given in the theorem.

One choice of a
map of minimal resolutions $\tilde y$ representing $y$ is as follows.
\[ \xymatrix@C=7.4mm@R=10mm{
P_\kk\ar[r]^{\left(\begin{smallmatrix} \bar c\\ \bar b\end{smallmatrix}\right)}
\ar[d]_{(1)}&P_\MM\oplus P_\NN
\ar[r]^{\left(\begin{smallmatrix}\bar c\bar d&\bar f\\\bar e&\bar b\bar a\end{smallmatrix}\right)} 
\ar[d]_{(\bar a\bar e,0)}
&P_\MM\oplus P_\NN\ar[r]^{\left(\begin{smallmatrix}\bar d& \bar a\end{smallmatrix}\right)}
\ar[d]_{\left(\begin{smallmatrix}\bar c\bar d&0\\0&0\end{smallmatrix}\right)}&
P_\kk\ar[r]^{\bar a\bar e \bar c\,=\,\bar d\bar f\bar b}
\ar[d]_(0.4){\left(\begin{smallmatrix}\bar f\bar b \\ 0\end{smallmatrix}\right)}
\ar[r]&P_\kk\ar[r]\ar[d]_{(1)}&
P_\MM\oplus P_\NN\ar[r] &P_\MM\oplus P_\NN\ar[r]&P_\kk\\
P_\kk\ar[r]_{\bar a\bar e\bar c\,=\,\bar d\bar f\bar b}&P_\kk
\ar[r]_(0.4){\left(\begin{smallmatrix} \bar c\\ \bar  b\end{smallmatrix}\right)}
&P_\MM\oplus P_\NN
\ar[r]_{\left(\begin{smallmatrix}\bar c\bar d&\bar f\\\bar e&\bar b\bar a\end{smallmatrix}\right)}
&P_\MM\oplus P_\NN\ar[r]_(.6){\left(\begin{smallmatrix}\bar d& \bar a\end{smallmatrix}\right)}&P_\kk} \]
The element $z$ lifts to 
the periodicity generator\index{periodicity generator} 
$\tilde z\colon P_* \to P_*$ of degree $-4$.
We have $\tilde y\tilde z=\tilde z\tilde y$, and we have
the following homotopy $u$ from $\tilde y \circ \tilde y$ to zero.
\[ \xymatrix@R=15mm@C=8mm{P_\kk
\ar[r]^{\left(\begin{smallmatrix} \bar c\\ \bar b\end{smallmatrix}\right)}
\ar[d]^(.4){\left(\begin{smallmatrix}\bar f\bar b \\ 0\end{smallmatrix}\right)}
& P_\MM \oplus P_\NN
\ar[r]^{\left(\begin{smallmatrix}\bar c\bar d&\bar f\\\bar e&\bar b\bar a\end{smallmatrix}\right)} 
\ar[d]^(.4){\left(\begin{smallmatrix}\bar a\bar e & 0\end{smallmatrix}\right)}
\ar[dl]^{\left(\begin{smallmatrix}\bar e & 0\\0&0\end{smallmatrix}\right)}
& P_\MM \oplus P_\NN
\ar[r]^{\left(\begin{smallmatrix}\bar d& \bar a\end{smallmatrix}\right)}
\ar[d]^(.4){\left(\begin{smallmatrix} 0 & 0\end{smallmatrix}\right)}
\ar[dl]^{\left(\begin{smallmatrix} 0 & 0\end{smallmatrix}\right)}
& P_\kk
\ar[r]^{\bar a\bar e\bar c\,=\,\bar d\bar f\bar b}
\ar[d]^(.4){\left(\begin{smallmatrix} 0 \\ 0\end{smallmatrix}\right)}
\ar[dl]^{(0)}
&P_\kk
\ar[r]^{\left(\begin{smallmatrix} \bar c\\ \bar b\end{smallmatrix}\right)}
\ar[d]^(.4){\left(\begin{smallmatrix} \bar f\bar b \\ 0\end{smallmatrix}\right)}
\ar[dl]^{\left(\begin{smallmatrix} 0 \\ \bar b\end{smallmatrix}\right)}
&P_\MM\oplus P_\NN\ar[r]
\ar[d]^(.4){\left(\begin{smallmatrix}\bar a\bar e & 0\end{smallmatrix}\right)}
\ar[dl]^{\left(\begin{smallmatrix}\bar e & 0\\0&0\end{smallmatrix}\right)}
&\cdots
\ar[dl]^{\left(\begin{smallmatrix} 0 & 0\end{smallmatrix}\right)}
\\
P_\MM \oplus P_\NN
\ar[r]_{\left(\begin{smallmatrix}\bar d& \bar a\end{smallmatrix}\right)}
&P_\kk\ar[r]_{\bar a\bar e\bar c\,=\,\bar d\bar f\bar b}
&P_\kk
\ar[r]_{\left(\begin{smallmatrix} \bar c\\ \bar b\end{smallmatrix}\right)}
& P_\MM\oplus P_\NN
\ar[r]_{\left(\begin{smallmatrix}\bar c\bar d&\bar f\\\bar e&\bar b\bar a\end{smallmatrix}\right)} 
& P_\MM\oplus P_\NN
\ar[r]_{\left(\begin{smallmatrix}\bar d& \bar a\end{smallmatrix}\right)}
& P_\kk\ar[r]
& 0
} \]
Thus $du = \tilde y \circ \tilde y$. Moreover, it is easy to check
that $u \tilde y = 0$
and $\tilde y u = 0$. Using the recipe of Kadeishvili given in
the proof of Theorem~\ref{th:Kadeishvili} for
computing the $A_\infty$ structure on $H^*BG$ from the differential
graded algebra structure on $\End_B(P_*)$, this implies that for all $n>2$ we have
$m_n(y,\dots,y)=0$.
Explicitly, let $A=\End_B(P_*)$ with $m_1$ the differential and $m_2$
the composition of endomorphisms. Then $A$ is quasi-isomorphic to the $A_\infty$
algebra $H^*BG$. The map $f_1$ takes $y^\ep z^i$ to $\tilde
y^\ep\tilde z^i$, and $f_2$ takes
$(yz^{i_1},yz^{i_2})$ to $uz^{i_1+i_2}$ and the remaining monomials to zero.
Then $u \tilde y=0$ implies that $m_2(f_2\otimes f_1)=0$ while
$\tilde y u=0$ implies that $m_2(f_1\otimes f_2)=0$. We also check that
$f_2(1\otimes m_2-m_2\otimes 1)=0$. Now applying
Remark~\ref{rk:formal}, we may take $f_i=0$ and $m_i=0$
for $i\ge 3$ to deduce that $H^*BG$ is formal.
\end{proof}

\begin{corollary}\label{co:Q1-formal}
Let $G$ be a finite group with quaternion or generalised quaternion Sylow
$2$-subgroup, and with no normal subgroup of index two.
\begin{enumerate}[label={\rm(\roman*)}]
\item
We have $\HH^*H^*BG\cong  \Lambda(y,\hat z) \otimes
\kk[z,\hat y]$ with $|y|=(0,-3)$, $|z|=(0,-4)$, $|\hat y|=(-1,3)$,
$|\hat z|=(-1,4)$. 
\item
The ring $\HHinf^*H^*BG\cong \HHinf^*H_*\Omega BG\twohat$, 
has the same generators and relations as in {\rm (i)}, but
the degrees are given by $|y|=-3$, $|z|=-4$, $|\hat y|=2$, $|\hat
z|=3$.
\item
We have $H_*\Omega BG\twohat\cong\Lambda(\hat z) \otimes
\kk[\hat y]$ with $|\hat z|=3$ and $|\hat y|=2$. This is formal as an
$A_\infty$ algebra.
\item
We have $\HH^*H_*\Omega BG\twohat\cong \Lambda(y,\hat z) \otimes
\kk[z,\hat y]$ with $|y|=(-1,-2)$, $|z|=(-1,-3)$, $|\hat y|=(0,2)$ and
$|\hat z| = (0,3)$.
\end{enumerate}
\end{corollary}
\begin{proof}
By Theorem~\ref{th:Q1-formal}, $H^*BG=\Lambda(y)\otimes \kk[z]$ is
formal. It is a complete intersection, and a Koszul
algebra, with Koszul dual $H_*\Omega BG\phat$, giving (iii).
Parts (i) and (iv) then follow using either
Theorem~\ref{th:HHR} or Theorem~\ref{th:Negron}. Part (ii) then follows
using Theorem~\ref{th:HHOmega} and Theorem~\ref{th:HHHA}. 
\end{proof}

\begin{theorem}
We have
\begin{align*} 
&\Dsg(H^*BG)\simeq
\Dcsg(H_*\Omega BG\twohat)\simeq\Db(\Lambda(\hat z)\otimes \kk[\hat
y,\hat y^{-1}]) \\
&\Dsg(H_*\Omega BG\twohat) 
\simeq
\Dcsg(H^*BG) 
\simeq \Db(\Lambda(y) \otimes \kk[z,z^{-1}]). 
\end{align*}
\end{theorem}
\begin{proof}
The equivalences of the singularity and cosingularity categories  follows
from Theorem~\ref{th:Dsg-Dcsg}. The equivalences with bounded derived
categories with $\hat y$ inverted in the first case and $z$
inverted in the second case follows from the fact that this inversion
precisely kills $\Thick(\kk)$ to form the cosingularity category.
\end{proof}

\section{Two classes of elements of order four}

In this section, we examine Case~\ref{case:Q2}, of a finite group $G$ with
generalised quaternion Sylow $2$-subgroups of order $8q$ 
and two classes of elements of order
four. This implies that $q\ge 2$.

\begin{theorem}\label{th:HBGquat-2classes}
Suppose that $G$ is a finite group with generalised quaternion Sylow
$2$-subgroups and two classes of elements of order four. Then
\[ H^*BG = \kk[y,z]/(y^4), \] 
with $|y|=-1$, $|z|=-4$.
\end{theorem}
\begin{proof}
Without loss of generality, assume that $O(G)=1$. Then $G$ has a
central involution $s$, and the class of the central extension of
$G/\langle s\rangle$ by $\langle s\rangle$, in the notation of
Section~\ref{se:D2}, is $t+y^2$. So in the spectral sequence of the
central extension, we have $d_2(w)=t+y^2$,
$d_3(w)=\Sq^1(t+y^2)=\xi+yt$, and $H^*BG=\kk[\xi,y,t]/(\xi
y,t+y^2,\xi+yt) \otimes \kk[w^4]$. In this ring, we have $t=y^2$,
$\xi=yt=y^3$, and $y^4=\xi y=0$. So letting $z$ be a representative of
$w^4$ in $H^*BG$, the structure is as given.
\end{proof}

\begin{theorem}\label{th:ExtHBQ}
The $\Ext$ ring of $H^*BG$ is given by
\[  \Ext^{*,*}_{H^*BG}(\kk,\kk) = \Lambda(\hat y,\hat z)\otimes \kk[\eta] \]
with $\eta = \langle \hat y,\hat y,\hat y,\hat y\rangle$, 
$|\hat y|=(-1,1)$, $|\hat z| =(-1,4)$, and $|\eta|=(-2,4)$.
\end{theorem}
\begin{proof}
The algebra $H^*BG$ is a 
complete intersection.\index{complete intersection} 
The second partial
derivatives of the relation all vanish, so this follows by an easy application
of Theorem~\ref{th:ExtRkk}. 
\end{proof}

\begin{corollary}\label{co:HOmegaBG-quat-2classes}
We have $H_*\Omega BG\twohat=\Lambda(\hat y,\hat z) \otimes \kk[\eta]$ 
with $\eta=\langle \hat y,\hat y,\hat y,\hat y\rangle$, $|\hat
y|=0$, $|\hat z|=3$ and $|\eta|=2$.
\end{corollary}
\begin{proof}
There is no room for non-zero differentials in the Eilenberg--Moore spectral
sequence~\eqref{eq:Cotor}
\[ \Ext^{*,*}_{H^*BG}(\kk,\kk) \Rightarrow H_*\Omega BG\twohat. \]
For the ungrading, the only issue is to choose the correct
representative for $\hat y$ so that it squares to zero. This is
possible, because the group of connected components is $\bZ/2$, and
the group algebra of this has a non-zero element that squares to zero.
\end{proof}

\begin{theorem}\label{th:HHHBQ2classes}
The Hochschild cohomology of $H^*BG$ is given by
\[ \HH^*H^*BG = \Lambda(\hat y,\hat z) \otimes \kk[y,z,\eta]/(y^4)  \]
with $|y|=(0,-1)$, $|z|=(0,-4)$, $|\hat y|=(-1,1)$, $|\hat z|=(-1,4)$
and $|\eta|=(-2,4)$.
\end{theorem}
\begin{proof}
By Theorem~\ref{th:HBGquat-2classes}, $H^*BG$ is a complete
intersection, so this follows from
Theorem~\ref{th:HHR}.
\end{proof}

\begin{theorem}
The Hochschild cohomology of $H_*\Omega BG\twohat$ is given by
\[ \HH^*H_*\Omega BG\twohat = \Lambda (\hat y,\hat z,\hat\eta) \otimes
  \kk[y,z,\eta] \]
with $|\hat y| = (0,0)$, $|\hat z|=(0,3)$, $|\hat\eta|=(-1,-2)$,
$|y|=(-1,0)$, $|z|=(-1,-3)$ and $|\eta|=(0,2)$.
\end{theorem}
\begin{proof}
By Corollary~\ref{co:HOmegaBG-quat-2classes}, $H_*\Omega BG\twohat$ is
a complete intersection, so this follows from Theorem~\ref{th:HHR}.
\end{proof}

The fact that in $\Ext^{*,*}_{H^*BG}(\kk,\kk)$ we have 
$\eta=\langle\hat y,\hat y,\hat y,\hat y\rangle$ 
implies that in the spectral sequence 
\[ \HH^*H_*\Omega BG\twohat \Rightarrow
\HH^*C_*\Omega BG\twohat \] 
we have a differential $d^3(\hat\eta)=y^4$.

\section{Non-formality}

Our goal in this section is to prove that in Case~\ref{case:Q2}, 
$H^*BG$ is not formal as an $A_\infty$ algebra. To do so, we shall
show that the Massey triple product\index{Massey product} 
$\langle y^2,y^2,y^2\rangle$
vanishes, but $\langle y^2,y^2,y^2,y^2\rangle$ is 
equal to $y^2z$. 

If $G$ is $SL^\circ(2,p^m)$\index{SL@$SL^\circ(2,p^m)$} 
with $p^m\equiv 3\pmod{4}$ then $G$ is an
example of Case~\ref{case:Q2}, and the principal block belongs to
Erdmann's~\cite{Erdmann:1990a} class $Q(2{\mathcal B})_1$, which is
labelled I in the Appendix to~\cite{Erdmann:1988b}.
If $G$ is $SU^\circ(2,p^m)$\index{SU@$SU^\circ(2,p^m)$} with
$p^m\equiv 1\pmod{4}$ then $G$ is also an example of Case~\ref{case:Q2},
and the principal block belongs to Erdmann's class
$Q(2{\mathcal A})$, which is labelled II in~\cite{Erdmann:1988b}.
The two types are derived equivalent, by the main theorem of
Holm~\cite{Holm:1997a}.\index{derived!equivalence} 

Type $Q(2{\mathcal B})_1$ is the easier to handle, so we assume that we
are in the case of $SL^\circ(2,p^m)$ with $p^m\equiv 3 \pmod{4}$. 
The quiver\index{quiver} for the principal block $B$ is
\begin{equation*}
\xymatrix{\kk\ar@(ul,dl)_a\ar@/^/[r]^b 
& \MM\ar@/^/[l]^c\ar@(ur,dr)^d} 
\end{equation*}
with relations
\[ bc=d^{2q-1},\quad
db=bacba, \quad
cd=acbac,\quad
a^2=cbacb+\lambda(acb)^{2},\quad
ba^2=0 \]
for some value of the parameter $\lambda \in \kk$ that has not been
determined. 
Note that these relations imply that
\begin{gather*} 
bcb = d^{2q-1}b=d^{2q-2}bacba=d^{2q-3}bacba^2cba=0, \\
cbc = cd^{2q-1}=acbacd^{2q-2}=acba^2cbacd^{2q-3}=0, \\
a^2c=cbacbc+\lambda acbacbc = 0.
\end{gather*}
The projective covers of the simple modules $\kk$ and $\MM$ have the
following diagrams (beware of the extra socle term in the
expression for $a^2$):
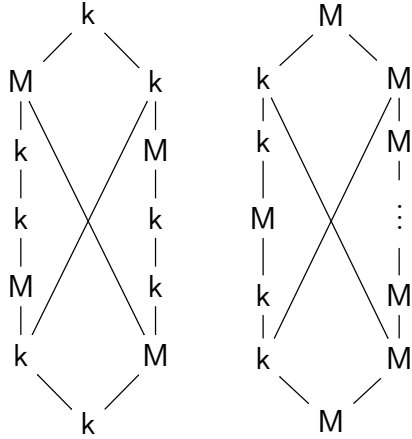
\[ \vcenter{\xymatrix@=3.7mm{&\kk\ar@{-}[dl]\ar@{-}[dr]\\
\MM\ar@{-}[d]\ar@{-}[ddddrr] && \kk\ar@{-}[d]\ar@{-}[ddddll] \\
\kk\ar@{-}[d]&&\MM\ar@{-}[d]\\
\kk\ar@{-}[d]&&\kk\ar@{-}[d]\\
\MM\ar@{-}[d]&&\kk\ar@{-}[d]\\
\kk\ar@{-}[dr]&&\MM\ar@{-}[dl]\\
&\kk
}}\qquad
\vcenter{\xymatrix@=3mm{&\MM\ar@{-}[dl]\ar@{-}[dr] \\ 
\kk\ar@{-}[d]\ar@{-}[ddddrr]&&\MM\ar@{-}[d]\ar@{-}[ddddll] \\
\kk\ar@{-}[d]&&\MM\ar@{-}[d]\\
\MM\ar@{-}[d]&&\vdots_{_{_{}}}\!\ar@{-}[d]\\
\kk\ar@{-}[d]&&\MM\ar@{-}[d]\\
\kk\ar@{-}[dr]&&\MM\ar@{-}[dl]\\
&\MM}} \]
where the number of copies of $\MM$ down the right hand side of the
projective cover of $\MM$ is $2q-1$.

The minimal resolution of $\kk$ is periodic with period four, and has
the following form:
\[ \cdots \to P_\kk \oplus P_\MM \xrightarrow{(\bar a,\bar b)}
P_\kk \xrightarrow{(\bar a^3)} P_\kk 
\xrightarrow{\left(\begin{smallmatrix}\bar a+\lambda\bar b
\bar c\bar a\bar b\bar c \\\bar c\end{smallmatrix}\right)} 
P_\kk \oplus P_\MM
\xrightarrow{\left(\begin{smallmatrix}\bar a & \bar b\bar c\bar a\bar b\\
\bar c \bar a \bar b \bar c + \lambda \bar c\bar a\bar b\bar c\bar a &
\bar d\end{smallmatrix}\right)} P_\kk \oplus P_\MM 
\xrightarrow{(\bar a,\bar b)} P_\kk  \]

We lift $y\in \Ext^1_\kQ(\kk,\kk)$ to a map $\tilde y$ of resolutions
as follows:
{\small
\[ \xymatrix@C=1.0cm@R=1.3cm{
P_\kk\ar[d]^{(1,0)} \oplus P_\MM\ar[r]^{(\bar a,\bar b)} &
P_\kk\ar[rr]^{(\bar a^3)}\ar[d]^{(\bar a^2)} &&
P_\kk\ar[rr]^-{\left(\begin{smallmatrix}\bar a+\lambda\bar b
\bar c\bar a\bar b\bar c \\\bar c\end{smallmatrix}\right)} 
\ar[d]^{\left(\begin{smallmatrix} 1\\\lambda \bar c\end{smallmatrix}\right)}&& 
P_\kk \oplus P_\MM
\ar[rr]^{\left(\begin{smallmatrix}\bar a & \bar b\bar c\bar a\bar b\\
\bar c \bar a \bar b \bar c + \lambda \bar c\bar a\bar b\bar c\bar a\!\!\! &
\bar d\end{smallmatrix}\right)}
\ar[d]^{\left(\begin{smallmatrix}1&0\\0&\bar c\bar a \bar b\end{smallmatrix}\right)} &&
P_\kk \oplus P_\MM\ar[r]^{(\bar a,\bar b)}\ar[d]^{(1,0)}&P_\kk\\
P_\kk \ar[r]_{(\bar a^3)} &
P_\kk \ar[rr]_-{\left(\begin{smallmatrix}\bar a+\lambda\bar b
\bar c\bar a\bar b\bar c \\\bar c\end{smallmatrix}\right)}  &&
P_\kk \oplus P_\MM\ar[rr]_{\left(\begin{smallmatrix}\bar a & \bar b\bar c\bar a\bar b\\
\bar c \bar a \bar b \bar c + \lambda \bar c\bar a\bar b\bar c\bar a\!\!\! &
\bar d\end{smallmatrix}\right)} && P_\kk \oplus P_\MM 
\ar[rr]_{(\bar a,\bar b)} && P_\kk} \]
}
The composite $\tilde y^8$ is zero, and we have the following homotopy $u_1$ from
 $\tilde y^4$ to zero:
{\small
\[ \xymatrix@C=1.0cm@R=1.8cm{
P_\kk\ar[r]^{(\bar a^3)}\ar[d]^{(\bar a^2)} &
P_\kk\ar[rr]^-{\left(\begin{smallmatrix}\bar a+\lambda\bar b
\bar c\bar a\bar b\bar c \\\bar c\end{smallmatrix}\right)} 
\ar[d]^{(\bar a^2)}\ar[dl]^{(0)}&& P_\kk \oplus P_\MM
\ar[rr]^-{\left(\begin{smallmatrix}\bar a & \bar b\bar c\bar a\bar b\\
  \bar c \bar a \bar b \bar c + \lambda \bar c\bar a\bar b\bar c\bar a\!\!\! &
\bar d\end{smallmatrix}\right)}
\ar[d]^-{\left(\begin{smallmatrix}\bar a^2&0\\0&0\end{smallmatrix}\right)}
\ar[dll]^{\substack{(0,\bar b\bar c\bar a\bar b\\ \quad+\lambda \bar a\bar b\bar
  c\bar a\bar b)}} &&
P_\kk \oplus P_\MM\ar[r]^-{(\bar a,\bar b)}
\ar[d]^(0.4){\left(\begin{smallmatrix}\bar a^2&0\\0&0\end{smallmatrix}\right)}
\ar[dll]^{\left(\begin{smallmatrix}0&\bar b \\ 0 & \lambda \bar c
      \bar b\end{smallmatrix}\right)}
&P_\kk\ar[r]^{(\bar a^3)}\ar[d]^{(\bar a^2)}
\ar[dl]^{\left(\begin{smallmatrix}\bar a\\0\end{smallmatrix}\right)}& 
P_\kk\ar[r]
& \cdots\\
P_\kk\ar[r]_{(\bar a^3)}&
P_\kk\ar[rr]_-{\left(\begin{smallmatrix}\bar a+\lambda\bar b
\bar c\bar a\bar b\bar c \\\bar c\end{smallmatrix}\right)} 
&&P_\kk\oplus P_\MM
\ar[rr]_-{\left(\begin{smallmatrix}\bar a & \bar b\bar c\bar a\bar b\\
\bar c \bar a \bar b \bar c + \lambda \bar c\bar a\bar b\bar c\bar a\!\!\! &
\bar d\end{smallmatrix}\right)}
&&P_\kk\oplus P_\MM\ar[r]_-{(\bar a,\bar b)}&P_\kk} \]
}
This satisfies $u_1^2=0$ and $u_1\tilde y^4=\tilde y^4 u_1$.

The composite $u_1\tilde y$ is given by the matrices
\[ (0),\qquad (\lambda \bar b\bar c\bar a\bar b\bar c + 
\lambda^2\bar a \bar b\bar c\bar a\bar b\bar c),\qquad 
\begin{pmatrix}0&\bar b\bar c\bar a\bar b \\ 0 & 0\end{pmatrix},\qquad
\begin{pmatrix}\bar a & 0 \\ 0 & 0\end{pmatrix} \]
while $\tilde y u_1$ is given by the matrices
\[ (\bar a),\qquad (0),\qquad 
\begin{pmatrix} 0 & \bar b\bar c\bar a\bar b + \lambda \bar a \bar
  b\bar c\bar a\bar b \\ 0 & \lambda^2\bar c\bar a\bar b\bar c\bar
  a\bar b\end{pmatrix},\qquad
\begin{pmatrix} 0 & \bar b \\ 0 & 0 \end{pmatrix} \]
Since $u_1\tilde y + \tilde y u_1$ is non-zero, we need to find a homotopy
from it to zero. The following is such a homotopy $u_2$:
\[ (0),\qquad (0,\lambda\bar b\bar c\bar a\bar b
+\lambda^2 \bar a\bar b\bar c\bar a\bar b),\qquad
\begin{pmatrix} 0 & 0\\ 0 & 0\end{pmatrix},\qquad
\begin{pmatrix} 1\\0\end{pmatrix} \]
Then $\tilde y u_2+u_2\tilde y$ is a homotopy from $\tilde y^2
u_1+u_1\tilde y^2$ to zero, $\tilde y^2 u_2 + \tilde y u_2 \tilde y +
u_2\tilde y^2$ is a homotopy from $\tilde y^3 u_1+u_1\tilde y^3$ to
zero, and $\tilde y^3u_2+\tilde y^2u_2\tilde y + \tilde yu_2\tilde y^2
+ u_2\tilde y^3$ is a homotopy from $\tilde y^4u_1+u_1\tilde y^4=0$ to
zero. 

At the next stage, the relevant composites are given by the matrices
\[ \begin{array}{ccccc}
\tilde y^3u_2\colon & \begin{pmatrix} 0 \\ 0 \end{pmatrix},&
(0,\lambda \bar b \bar c \bar a\bar b + 
\lambda^2 \bar a\bar b\bar c\bar a\bar b),&
(0,0),&
\begin{pmatrix}\bar a^2 \\ 0\end{pmatrix} \\
\tilde y^2u_2\tilde y\colon&
 \begin{pmatrix} 
\lambda^3\bar a\bar b\bar c\bar a \bar b\bar c \\ 0 
\end{pmatrix}, &
(0,0), &
(\bar a^2, 0), &
\begin{pmatrix} 0 \\ 0 \end{pmatrix} \\
\tilde y u_2\tilde y^2\colon&
\begin{pmatrix} 0 \\ 0 \end{pmatrix},&
(1,0),&
(0,0),&
\begin{pmatrix} 0 \\ 0 \end{pmatrix} \\
u_2\tilde y^3\colon & 
\begin{pmatrix} 1 \\ 0 \end{pmatrix},&
(0,0),&
(0,0),&
\begin{pmatrix} 0 \\ 0 \end{pmatrix}.
\end{array} \]
Thus $\tilde y^3u_2+\tilde y^2u_2\tilde y+\tilde yu_2\tilde
y^2+u_2\tilde y^3 + u_1^2$ is given by the matrices
\[ \begin{pmatrix} 
1+\lambda^3\bar a\bar b\bar c\bar a \bar b\bar c \\ 0 
\end{pmatrix},\qquad
(1,\lambda \bar b \bar c \bar a\bar b + 
\lambda^2 \bar a\bar b\bar c\bar a\bar b),\qquad
(\bar a^2,0),\qquad
\begin{pmatrix} \bar a^2 \\ 0 \end{pmatrix}. \]
This is homotopic to $\tilde y^2\tilde z$, which is given by the
matrices
\[ \begin{pmatrix} 
1 \\ \lambda \bar c\bar a\bar b\bar c 
\end{pmatrix},\qquad
(1,0),\qquad
(\bar a^2,0)\qquad
\begin{pmatrix} \bar a^2 \\ 0 \end{pmatrix}. \]
The following is such a homotopy $u_3$:
\[ \begin{pmatrix} 
\lambda^2\bar a & \lambda^2\bar b\bar c\bar a\bar b 
+ \lambda^3\bar a\bar b \bar c\bar a\bar b \\ 
\lambda^2\bar c\bar a\bar b\bar c\bar a &
\lambda \bar c\bar a\bar b\end{pmatrix},\qquad
(0,0),\qquad
(0),\qquad 
\begin{pmatrix} 0 \\ 0 \end{pmatrix} \]

So we set
\begin{align*}
f_1(y)&=\tilde y,\\
f_2(y,y^3)&=f_2(y^3,y)=f_2(y^2,y^2)=u_1,\\
f_3(y,y^3,y)&=u_2,\\
f_3(y^2,y^2,y^2)&=\tilde y u_2 + u_2 \tilde y,\\
f_3(y^3,y,y^3)&=\tilde y^2 u_2 + \tilde y u_2 \tilde y + u_2\tilde y^2
+ \tilde y^2\tilde z,\\
f_4(y^3,y,y^3,y)&=f_4(y^2,y^2,y^2,y^2)=f_4(y,y^3,y,y^3)=u_3
\end{align*}
Thus we have $m_3=0$,
$m_1f_3=m_2(f_2\otimes f_1+f_1\otimes f_2)$, and
\begin{equation}\label{eq:m4}
m_4(y^3,y,y^3,y)=m_4(y^2,y^2,y^2,y^2)=m_4(y,y^3,y,y^3)=y^2z.
\end{equation}
We can now check that the values of
$m_1f_4+f_1m_4$ and
$m_2(f_3\otimes f_1+f_2\otimes f_2 +f_1\otimes f_3)$ on the quadruples
$(y^3,y,y^3,y)$, $(y^2,y^2,y^2,y^2)$ and $(y,y^3,y,y^3)$ are all equal
to
\[ \tilde y^3 u_2 + \tilde y^2 u_2 \tilde y + \tilde y u_2\tilde
  y^2+u_2\tilde y^3+u_1^2. \]

Equation~\eqref{eq:m4} may now be interpreted in terms of Hochschild
cohomology, using Proposition~\ref{pr:HH}. Since $m_3=0$, $m_4$ is a Hochschild cocycle. It
represents the element $\eta^2y^2z$ in degree $(-4,2)$, which by
Theorem~\ref{th:HHHBQ2classes} is non-zero. 

At the next stage, the expression $m_2(f_4\otimes f_1+f_3\otimes f_2 +
f_2\otimes f_3+f_1\otimes f_4)$ sends the $5$-tuple
$(y^2,y^2,y^2,y^2,y^2)$ to $u_3\tilde y^2+(\tilde y u_2 +  u_2\tilde
y)u_1 + u_1(\tilde y u_2+u_2 \tilde y)+\tilde y^2 u_3$, which is
\[ (0),\qquad \begin{pmatrix} \lambda^2 \bar
    a^3&0 \end{pmatrix},\qquad
\begin{pmatrix} 0&0\\0&0\end{pmatrix},\qquad
\begin{pmatrix} \lambda^2\bar a^3 \\ 0\end{pmatrix}. \]
This is homotopic to zero, with homotopy $u_4$ given by
\[ (0),\qquad
\begin{pmatrix} \lambda^2&\lambda^3\bar b\bar c\bar a \bar
  b\end{pmatrix},\qquad
\begin{pmatrix} 0&0\\0&0\end{pmatrix},\qquad
\begin{pmatrix} \lambda^2 \\ 0\end{pmatrix}. \]
Thus we can take $m_5=0$ and 
$f_5(y^2,y^2,y^2,y^2,y^2)=u_4$, 
to obtain 
\[ m_1f_5+f_1m_5=m_2(f_4\otimes f_1+f_3\otimes f_2+f_2\otimes f_3
+ f_1\otimes f_4) .\]

\begin{theorem}\label{th:Qu-non-formal}
If $G$ is a finite group with generalised quaternion Sylow
$2$-subgroups and two classes of elements of order four, then $H^*BG$
and $H_*\Omega BG\twohat$ are not formal as $A_\infty$
algebras.\index{formal $A_\infty$ algebra}
\end{theorem}
\begin{proof}
We have just shown that $H^*BG$ is not formal, since we can choose
$m_3$ to be zero, but then $m_4$ cannot
be chosen to be zero. The fact that
$H_*\Omega BG\twohat$ is not formal follows from 
the relation $m_4(\tilde y,\tilde y,\tilde y,\tilde y)=\eta$. This in
turn follows from the fact that the Massey product\index{Massey product} 
$\langle \tilde y,\tilde y,\tilde y,\tilde y\rangle=\eta$ in
Corollary~\ref{co:HOmegaBG-quat-2classes} has no indeterminacy, see
Theorem~\ref{th:Massey}.
\end{proof}

This finally allows us to compute $\HHinf^*H^*BG\cong\HHinf^*H_*\Omega
BG\twohat$, see Theorem~\ref{th:HHOmega}.

\begin{theorem}\label{th:HHH-Qu-non-formal}
We have $\HHinf^*H^*BG \cong \HHinf^*H_*\Omega BG\twohat\cong 
\Lambda(\hat y) \otimes \kk[y,z,\eta]/(y^4,y^2\eta^2)$ with
$|y|=-1$, $|z|=-4$, $|\hat y|=0$, $\eta=2$.
\end{theorem}
\begin{proof}
The computation above of $m_3$ and $m_4$ in $H^*BG$, together with
Corollary~\ref{co:HOmegaBG-quat-2classes} show that in the spectral
sequence~\eqref{eq:HHHA}
\[ \HH^*H^*BG \Rightarrow \HHinf^*H^*BG \]
we have $d^2=0$ and  $d^3(\hat z)=y^2\eta^2$. This then implies that
$\eta$ is a universal cycle, and $E_4=E_\infty$. 

Ungrading $y^4=0$ in $E_\infty$, we see that $y^4$ has to be a linear
combination of the elements $y^2z\eta$ and $z^i\eta^{2i-2}$ with $i\ge 2$.
Since $\tilde y^8=0$, the relation ungrading $y^4=0$ has to satisfy
$y^8=0$, so it cannot involve the elements $z^i\eta^{2i-2}$. So $y^4$
is some multiple of $y^2z\eta$. But if it's a non-zero multiple then
$y^8$ is a non-zero multiple of $y^2z^3\eta^3$. This contradiction
shows that $y^4=0$ in $\HH^*H^*BG$.

Ungrading $y^2\eta^2=0$ in $E_\infty$, we see that $y^2\eta^2$ is a
linear combination of the elements $z^i\eta^{2i+1}$ with $i\ge 1$.
Again, since $\tilde y^8=0$, $y^2\eta^2$ is nilpotent, and so we
have $y^2\eta^2=0$ in $\HH^*H^*BG$.
\end{proof}

It follows from this, that in the spectral sequence
\[ \HH^*H_*\Omega BG\twohat \Rightarrow \HHinf^*H_*\Omega BG\twohat \cong
  \HHinf^*H^*BG, \]
see~\eqref{eq:HHHA} and Theorem~\ref{th:HHOmega},
we are forced to have $d^2(\hat z)=y^2\eta^2$, $d^3(\hat \eta)=y^4$,
to give the same answer for $\HHinf^*H^*BG$ as in Theorem~\ref{th:HHH-Qu-non-formal}.

\chapter{Beyond tame}

\section{Introduction}

In this chapter we discuss the $A_\infty$ algebras
$H^*BG$ and $H_*\Omega BG\phat$ beyond the
tame case (see Theorem~\ref{th:tame}). We've seen in our discussion of
tame representation type, 
that $H_*\Omega BG\phat$ always has polynomial growth. Furthermore, 
there is always a finitely generated central subalgebra over which
$H_*\Omega BG\phat$ is finitely generated as a module.

This is not always the case for finite groups.
There is a dichotomy, discovered by Ran Levi, between
polynomial and (semi-)exponential growth for $H_*\Omega BG\phat$, and we
discuss this in Section~\ref{se:poly-exp}. We give examples both of
exponential growth and of polynomial growth beyond tame
representation type. For example,  Chevalley and twisted Chevalley groups in
non-defining characteristic are always of polynomial growth, at least
for non-torsion primes, as we
shall explain in Section~\ref{se:non-defining}. 

The other aspect revealed by our discussion of tame representation
type is that in some unexpected cases it turns out that $H^*BG$ is
formal as an $A_\infty$ algebra. We begin with a discussion of this
phenomenon. 

\section{Formality}\index{formal $A_\infty$ algebra}

One of the surprising aspects of our work is the discovery 
that $H^*BG$ is formal in two of the cases with
semidihedral Sylow $2$-subgroups, see Theorems~\ref{th:SD1-formal}
and~\ref{th:SD2-formal}, and also 
when $G$ has generalised quaternion Sylow  
subgroups and no normal subgroup of index two, see  
Theorem~\ref{th:Q1-formal}. 
In this section, we discuss formality in general for the $A_\infty$ algebra
$H^*BG$. We begin with finite $p$-groups.

\begin{theorem}\label{th:p-group}
Let $G$ be a finite $p$-group.\index{finite $p$-group} 
Then the following are equivalent.
\begin{enumerate}[label={\rm(\roman*)}]
\item The comultiplication on a minimal resolution of $\kk$ as a $\kG$-module
can be made strictly coassociative.\label{it:p-group/coass}
\item The $A_\infty$ algebra $H^*BG$ is formal.\label{it:p-group/formal}
\item $H^*BG$ is a polynomial ring.\label{it:p-group/poly}
\item $p=2$ and $G$ is an 
elementary abelian $2$-group.\index{elementary abelian $2$-group}\label{it:p-group/elemab}
\end{enumerate}
\end{theorem}
\begin{proof}
\ref{it:p-group/coass} $\Rightarrow$ \ref{it:p-group/formal}:
Let $P_*$ be a minimal resolution of $\kk$ as
a $\kG$-module. Then the differential on
$\Hom_{\kG}(P_*,\kk)$ is zero, and so $H^*BG\cong \Hom_{\kG}(P_*,\kk)$
with the multiplication induced from the comultiplication on $P_*$.
Since the comultiplication is coassociative, $\Hom_{\kG}(P_*,\kk)$ is
a DG algebra with zero differential, and is therefore
formal.

\ref{it:p-group/formal} $\Rightarrow$ \ref{it:p-group/poly}:
If $H^*BG$ is formal then the 
Eilenberg--Moore spectral 
sequence\index{Eilenberg--Moore spectral sequence}~\eqref{eq:Cotor}
gives an isomorphism
$\Ext^{*,*}_{H^*BG}(\kk,\kk)\cong \kG$. In particular, 
$\Ext^{*,*}_{H^*BG}(\kk,\kk)$ has finite total dimension
over $\kk$, so $H^*BG$
has finite global dimension. It follows that it is regular, 
and hence a polynomial ring.

\ref{it:p-group/poly} $\Leftrightarrow$ \ref{it:p-group/elemab}
is proved in Corollary~6.6 
of Benson and Carlson~\cite{Benson/Carlson:1994a}. 

\ref{it:p-group/elemab} $\Rightarrow$ \ref{it:p-group/coass}:
If $G$ is cyclic of order two then the 
reduced bar resolution is minimal, so $G$ satisfies \ref{it:p-group/coass}. 
If $G$ is an elementary abelian $2$-group then we express 
$G$ as a direct product of cyclic groups of order two, and 
form the tensor product of their minimal resolutions. The 
comultiplication resulting from this is strictly coassociative. 
\end{proof}

As an illustration of the grading techniques, we prove the following,
which also gives another proof of
\ref{it:p-group/elemab} $\Rightarrow$ \ref{it:p-group/formal} in 
Theorem~\ref{th:p-group}.

\begin{theorem}\label{th:elemab2}
Suppose that $G$ is a finite group with elementary 
abelian Sylow $2$-subgroup $D$,
and $\kk$ is a field of characteristic two. Then there is a grading on
$\kG$ such that
the $A_\infty$ algebra $H^*BG$ is intrinsically
formal.\index{intrinsically formal} In particular,
without reference to grading, $H^*BG$ is formal.
\end{theorem}
\begin{proof}
By Theorem~\ref{th:abelian} and Remark~\ref{rk:DH}, 
we can suppose that $G$ is a semidirect product $D\rtimes H$, with $H$
a $p'$-subgroup of $\Aut(D)$.

Let $D=\langle g_1,\dots,g_r\rangle\cong (\bZ/2)^r$. The group $H$
acts on $kD$, and this gives a short exact sequene of $kH$-modules
\[ 0 \to J^2(kD) \to J(kD) \to J(kD)/J^2(kD) \to 0. \]
Since $p$ does not divide $|H|$, this sequence splits. Let $U$ be
an invariant complement to $J^2(kD)$ in $J(kD)$, and let
$X_1,\dots,X_r$ be a basis for $U$.
Then 
\[ \kk D = k[X_1,\dots,X_r]/(X_1^2,\dots,X_r^2). \]
We can put a grading on $\kk D$ by setting $|X_i|=1$. 
Putting elements of $H$ in degree zero then defines a grading on $kG$.  
This gives us an $H$-invariant 
grading on $H^*BD=k[x_1,\dots,x_r]$ with $|x_i|=(-1,-1)$.
Then the ring $H^*BG=(H^*BD)^H$ is doubly graded. The cohomological degrees
of elements are equal to their internal degrees. The $A_\infty$ maps
$m_i\colon H^*BG \to H^*BG$ have degrees $(i-2,0)$, see
Theorem~\ref{th:Kadeishvili}. So for $i>2$ we have $m_i=0$, since
either the source or the target is zero. 
\end{proof}

\begin{remark}\index{compact Lie group}\index{Lie group}
A discussion of formality for $H^*BG$ in the case of a compact Lie
group $G$ can be found in Benson and
Greenlees~\cite{Benson/Greenlees:2023b}. The last section of this
paper has a discussion of the literature.
\end{remark}

\section{Polynomial versus exponential growth}%
\index{polynomial growth}%
\index{exponential growth}%
\index{semi-exponential growth}%
\label{se:poly-exp}

\begin{definition}
Let $f$ be a real valued function on the non-negative integers.
We say that $f$ grows \emph{at most polynomially} if there exists
a polynomial function $p$ such that for all $n\ge 0$ we have
$|f(n)|\le p(n)$.
\end{definition}

In commutative algebra, we have the following theorem,
characterising complete intersections.\index{complete intersection}

\begin{theorem}[Gulliksen~\cite{Gulliksen:1980a}, Theorem~2.3]%
\label{th:Gulliksen}
Let $R$ be a commutative local ring with 
residue field\index{residue field} $\kk$. Then $R$ is a complete
intersection if and only if
$\Ext^*_R(\kk,\kk)$ has polynomial growth.
\end{theorem}

The corresponding theorem for loop space homology of finite
complexes is as follows.

\begin{definition}
Let $f$ be a real valued function on the non-negative integers.
We say that $f$ grows \emph{at least semi-exponentially} if there
exists a constant $C>1$ such that for $n$ large enough
$\sum_{i=0}^n|f(i)|\ge C^{\sqrt{n}}$.
\end{definition}

\begin{example}
The partition function\index{partition function} $p(n)$ 
satisfies $\log p(n) \sim \pi \sqrt{\frac{2n}{3}}$ as
  $n\to\infty$, so $p(n)$ has semi-exponential growth.
\end{example}

\begin{theorem}\label{th:FHT}
[Theorem~A of F\'elix, Halperin and Thomas~\cite{Felix/Halperin/Thomas:1993a}]
Let $X$ be a simply connected finite CW complex, and $p$ a prime.
Then $H_n(\Omega X;\bF_p)$ grows either at most polynomially, or at least
semi-exponentially.
\end{theorem}

\begin{remark}
Anick~\cite{Anick:1982a} found examples where the growth is 
semi-exponential but not exponential, in the contexts of
both Theorem~\ref{th:Gulliksen} and Theorem~\ref{th:FHT}.
\end{remark}

\begin{definition}
A finite CW complex is said to be 
\emph{elliptic}\index{elliptic space} at $p$ 
if $H_*(\Omega X;\bF_p)$ has polynomial growth.
\end{definition}

Using Theorem~\ref{th:FHT}, Levi~\cite{Levi:1997a} proved the
following, which shows that in this respect, $BG\phat$ resembles a
finite CW complex.

\begin{theorem}
For a finite group $G$, 
the loop space homology $H_*\Omega BG\phat$ grows either at most polynomially
or at least semi-exponentially.
\end{theorem}

Examples are given in Levi~\cite{Levi:1996a,Levi:1997a} of groups for which 
the homology contains a free algebra on two
variables, so that the growth is exponential. We shall display a
slightly easier example in Section~\ref{se:exp-finite}. No examples are currently
known where the growth is at least semi-exponential but not exponential.

\begin{remark}
A discussion of various derived notions of 
complete intersections,\index{derived!complete intersection} in the context of
polynomial versus semi-exponential growth, can be found in
Benson, Greenlees and Shamir~\cite{Benson/Greenlees/Shamir:2013a},
Greenlees, Hess and Shamir~\cite{Greenlees/Hess/Shamir:2013a}.
The hope is that for spaces of the form $BG\phat$ with $G$ a finite
group, these notions coincide, and describe when $H_*\Omega BG\phat$
has at most polynomial growth.
\end{remark}

\section{An exponential compact Lie example}\label{se:exp-Lie}

For non-connected compact Lie 
groups,\index{compact Lie group}\index{Lie group} 
it is not hard to cook up examples
of exponential growth. In this section, we give an example which is
not only of exponential growth, but also formal. In
Section~\ref{se:exp-finite} we give a finite group example based on
this one.

Let $\kk=\bQ$, let $T$ be an
$r$-dimensional torus,\index{torus} 
and let $G=T\rtimes \bZ/2$, where the
involution inverts every element of $T$. Then
$H^*BT=\bQ[x_1,\dots,x_r]$ with $|x_i|=2$ for $1\le i\le r$, and
$H^*BG=(H^*BT)^{\bZ/2}$, where the action of the generator of $\bZ/2$ sends each $x_i$
to $-x_i$. This subalgebra is generated by the elements
$x_{i,j}=x_ix_j$ with $1\le i\le j\le r$. The relations are
\[ x_{i,i}x_{j,j}=x_{i,j}^2,\qquad
x_{i,i}x_{j,k}=x_{i,j}x_{i,k},\qquad
x_{i,j}x_{k,\ell}=x_{i,k}x_{j,\ell}=x_{i,\ell}x_{j,k}. \]
Here, distinct letters in the subscripts represent distinct indices.
By Theorem~\ref{th:Froberg} this is a Koszul 
algebra,\index{Koszul duality} so $\Ext^{*,*}_{H^*BG}(\kk,\kk)$ is the Koszul 
dual, which is a non-commutative algebra generated by degree $(-1,2)$
elements $\hat x_{i,j}$ with relations
\begin{gather*}
\hat x_{i,i}^2=0,\qquad 
[\hat x_{i,i},\hat x_{i,j}]=0, \qquad
[\hat x_{i,i},\hat x_{j,j}]+\hat x_{i,j}^2=0,\\
[\hat x_{i,i},\hat x_{j,k}]+[\hat x_{i,j},\hat x_{i,k}]=0,\qquad
[\hat x_{i,j},\hat x_{k,\ell}]+[\hat x_{i,k},\hat x_{j,\ell}]
+[\hat x_{i,\ell},\hat x_{j,k}]=0. 
\end{gather*}
Note that here, for elements $x$, $x'$ of odd homological degree, $[x,x']$ means
$xx'+x'x$. 

The Eilenberg--Moore spectral 
sequence~\eqref{eq:Cotor}\index{Eilenberg--Moore spectral sequence}
\[ \Ext^{*,*}_{H^*BG}(\bQ,\bQ) \Rightarrow H_*\Omega BG\Qhat \]
has no room for differentials or
ungrading problems, because all the generators have degree $(-1,2)$.
So $H_*\Omega BG\Qhat$ is the same ring as
$\Ext^*_{H^*BG}(\bQ,\bQ)$, but where the
generators $\hat x_{i,j}$ are in degree one. For $r\ge 3$, this has exponential
growth. For example, when $r=3$ it is a free module on eight
generators 
over the free subalgebra $\bQ\langle
\hat x_{1,2},\hat x_{1,3},\hat x_{2,3}\rangle$, and the quotient by the ideal
generated by $\hat x_{1,2}$, $\hat x_{1,3}$, $\hat x_{2,3}$ is an
exterior algebra on $\hat x_{1,1}$, $\hat x_{2,2}$, $\hat x_{3,3}$. 
The Poincar\'e series for $r=3$ is
given by
\[ \sum_{n=0}^\infty t^n\dim_\bQ H_{n}\Omega BG\Qhat = \frac{(1+t)^3}{1-3t} 
= 1+6t + 21t^2+64t^3+192t^4+576t^5+\cdots \]

\begin{theorem}
For general $r$, we have
\[
\sum_{n=0}^\infty t^n\dim_\bQ H^{n}BG=
\frac{\sum_{i=0}^{\lfloor\frac{r}{2}\rfloor}
\binom{r}{2i}t^{2i}}{(1-t^2)^r}, \qquad
\sum_{n=0}^\infty t^n \dim_\bQ H_{n}\Omega BG\Qhat=
\frac{(1+t)^r}{\sum_{i=0}^{\lfloor\frac{r}{2}\rfloor}
\binom{r}{2i}(-t)^i}. 
\]
\end{theorem}
\begin{proof}
The computation of the Poincar\'e series for $H^*BG$ is an
application of Molien's theorem,\index{Molien's theorem} see for
example Theorem~2.5.2 of Benson~\cite{Benson:1993a}.
\begin{align*}
\sum_{n=0}^\infty t^n\dim_\bQ H^nBG&=\sum_{n=0}^\infty
t^n\dim_\bQ(H^nBT)^{\bZ/2}\\
&=\half\left(\frac{1}{(1-t)^r}+\frac{1}{(1+t)^r}\right)\\
&=\frac{\half((1+t)^r+(1-t)^r)}{(1-t)^r(1+t)^r} \\
&=\frac{\half\left(\sum_{i=0}^r\binom{r}{i}t^i+
\sum_{i=0}^r(-1)^i\binom{r}{i}t^i\right)}{(1-t^2)^r}\\
&=\frac{\sum_{i=0}^{\lfloor\frac{r}{2}\rfloor}
\binom{r}{2i}t^{2i}}{(1-t^2)^r}.
\end{align*}
The computation of the Poincar\'e
series for $H_*\Omega BG\Qhat$ then follows using the
general relation~\eqref{eq:pR!s} between the Poincar\'e
series of a Koszul algebra and its dual. Note that the generators of
$H^*BG$ are in degree two while the generators of the Koszul dual
$H_*\Omega BG\Qhat$ are in degree one, so $t^2$ must be replaced by
$t$ after applying~\eqref{eq:pR!s}.
\end{proof}

\begin{remark}
It follows from the main theorem of Benson and
Greenlees~\cite{Benson/Greenlees:2023b} that for this
family of examples, the
$A_\infty$ structure on $H^*BG$ is formal.\index{formal $A_\infty$ algebra}
Then since it is a Koszul algebra, it follows that the Koszul dual
$H_*\Omega BG\Qhat$ is also formal.
\end{remark}

\section{An exponential finite group example}\label{se:exp-finite}

 The loop space homology in the cases
discussed in Chapters~2--4 is of polynomial growth, 
and almost commutative, in the sense that
there is a central subring over which the whole ring is finitely
generated as a module. 
In this section, for contrast, we examine a finite group 
example where $H_*\Omega BG\phat$ has 
exponential growth. We take our cue from what happened in the compact
Lie example of Section~\ref{se:exp-Lie}.
This is related to Levi's example but is somewhat simpler to analyse
using our technique of introducing an internal grading on the group algebra.

Let $p$ be an odd prime, $\kk$ be a field of characteristic $p$, 
and let $G$ be the group 
\[ (\bZ/p\times \bZ/p)\rtimes \bZ/2 \]
given by the presentation
\[ \langle g,h,\sigma\mid g^p=h^p=\sigma^2=1, gh=hg, \sigma g=g^{-1}\sigma,
  \sigma h=h^{-1}\sigma\rangle. \]
Let $H$ be the subgroup of index two generated by $g$ and $h$, and
let $X=\sum_{i=1}^{p-1}g^i/i$ and $Y=\sum_{i=1}^{p-1} h^i/i$ 
as elements of $\kk H\le \kG$. Then since
$-X=\sum_{i=1}^{p-1}g^{-i}/i$ and $-Y=\sum_{i=1}^{p-1}h^{-i}/i$, we
have the following presentation 
for the group algebra:
\[ \kG=\kk\langle X,Y,\sigma \mid X^p=Y^p=0, XY=YX, \sigma X=-X\sigma,
  \sigma Y=-Y\sigma, \sigma^2=1\rangle. \]
We can put a double grading on this by setting $|X|=(1,0)$,
$|Y|=(0,1)$ and $|\sigma|=(0,0)$. Then 
\[ H^*BH=\kk[u,v] \otimes \Lambda(x,y) \]
with $|u|=(-2,-p,0)$, $|v|=(-2,0,-p)$, $|x|=(-1,-1,0)$ and
$|y|=(-1,0,-1)$. The cohomology ring $H^*BG$ is equal to the invariants of the action
of $\sigma$, which is the subring generated by 
$a=u^2$, $b=uv$, $c=v^2$, $\alpha=xu$, $\beta=xv$, $\gamma=yu$, 
$\delta=yv$, and $\ep=xy$. Regarding $a$, $b$ and $c$ as polynomial
generators, $\alpha$, $\beta$, $\gamma$ and $\delta$ as exterior
generators, and $\ep$ as a square zero commuting generator, the further relations
are: 
\begin{gather*} 
ac=b^2,\quad 
a\beta=b\alpha,\quad 
b\beta=c\alpha,\quad 
a\delta=b\gamma,\quad
b\delta=c\gamma,\quad 
a\ep=\alpha\gamma,\quad
b\ep=\alpha\delta=\beta\gamma,\\
c\ep=\beta\delta,\quad
\alpha\beta=0,\quad 
\gamma\delta=0,\quad 
\alpha\ep=0,\quad 
\beta\ep=0,\quad
\gamma\ep=0,\quad 
\delta\ep=0.
\end{gather*}
Ignoring the higher multiplications, by Theorem~\ref{th:Froberg} this is a Koszul algebra, and so
its $\Ext$ algebra is the Koszul dual, with eight generators and $16$ relations:
\begin{multline*} 
\Ext^{*,*}_{H^*BG}(\kk,\kk) = \kk\langle \hat a,\hat b,\hat c,
\hat\alpha,\hat\beta,\hat\gamma,\hat\delta,\hat\ep\mid
\hat a^2=\hat c^2=[\hat a,\hat b]=[\hat b,\hat c]=[\hat a,\hat c]+\hat b^2=0,\\
[\hat a,\hat\alpha]= [\hat a,\hat\gamma]=[\hat c,\hat\beta]=
[\hat c,\hat \delta]=0,\\
\qquad\qquad[\hat a,\hat\beta]+[\hat b,\hat\alpha]=
[\hat a,\hat\delta]+[\hat b,\hat\gamma]=
[\hat b,\hat \beta]+[\hat c,\hat \alpha]= 
[\hat b,\hat \delta]+[\hat c,\hat\gamma]=0,\\
\qquad[\hat a,\hat\ep]+[\hat\alpha,\hat\gamma]= 
[\hat b,\hat\ep]+[\hat\alpha,\hat\delta]
+[\hat\beta,\hat\gamma]= 
[\hat c,\hat\ep]+[\hat\beta,\hat\delta]=0\,
\rangle.\qquad
\end{multline*}
The degrees are 
\begin{gather*} 
|\hat a|=(-1,4,2p,0),\qquad |\hat b|=(-1,4,p,p),\qquad 
|\hat c|=(-1,4,0,2p),\qquad |\hat \alpha|=(-1,3,p+1,0), \\
|\hat \beta|=(-1,3,1,p),\qquad |\hat\gamma|=(-1,3,p,1),\qquad 
|\hat \delta|=(-1,3,0,p+1),\qquad |\hat \ep|=(-1,2,1,1).
\end{gather*}
Note that the sign in a commutator in this presentation is determined by the sum of the
first two degrees; it is the difference of two products if even, and
the sum if odd.

Since the differentials in the Eilenberg--Moore 
spectral sequence\index{Eilenberg--Moore spectral sequence}~\eqref{eq:Cotor}
preserve the two internal degrees, which are the third and fourth
degrees, it is easy to see that there is no
room for non-zero differentials. For example, $d^n(\hat a)$ has last
degree zero, so it can only involve $\hat a$ and $\hat\alpha$. No
monomial in these has the appropriate third degree. 

Similarly, when we ungrade the $E^\infty$ page of the spectral
sequence, there are no other monomials of the same internal degrees as
the quadratic terms in the list, so the ungraded relations are the
same as in $E^\infty$. It follows that $H_*\Omega BG\phat\cong
\Ext^{*,*}_{H^*BG}(\kk,\kk)$ is as described above, but where the
first two degrees have been added:
\begin{gather*} 
|\hat a|=(3,2p,0),\qquad |\hat b|=(3,p,p),\qquad 
|\hat c|=(3,0,2p),\qquad |\hat \alpha|=(2,p+1,0), \\
|\hat \beta|=(2,1,p),\qquad |\hat\gamma|=(2,p,1),\qquad 
|\hat \delta|=(2,0,p+1),\qquad |\hat \ep|=(1,1,1).
\end{gather*}
This algebra contains for example 
the free algebra\index{free algebra} $\kk\langle\hat\beta,
\hat\gamma\rangle$, and therefore has exponential
growth.

To obtain the Poincar\'e series for $R=H^*BG$, we use the version of Molien's
theorem\index{Molien's theorem} given in Theorem~5.2.2 of~\cite{Benson:1993a}. We need to be
careful about the Koszul grading, because the generators of $H^*BG$
are products of pairs of generators of $H^*BH$. To compensate for
this, we give the generators of $H^*BH$ Koszul degree $\half$. So
Molien's theorem gives us
\begin{align*}
p_R(s,t)&=\sum_{i,j}s^it^{-j}\dim_\kk H^{i,j}BG\\
&=\half\left(\left(\frac{1+s^\halfs t^{-1}}{1-s^\halfs t^{-2}}\right)^{\!2}+
\left(\frac{1-s^\halfs t^{-1}}{1+s^\halfs t^{-2}}\right)^{\!2}\right)\\
&=\frac{\half\left((1+s^\halfs t^{-1})^2(1+s^\halfs t^{-2})^2
+(1-s^\halfs t^{-1})^2(1-s^\halfs t^{-2})^2\right)}{(1-s^\halfs t^{-2})^2(1+s^\halfs t^{-2})^2}\\
&=\frac{1+st^{-2}+4st^{-3}+st^{-4}+s^2t^{-6}}{(1-st^{-4})^2}.
\end{align*} 
For $R^!=H_*\Omega BG\phat$ we use
formula~\eqref{eq:pR!st-1} to obtain
\begin{align*}
p_{R^!}(s,t)&=\sum_{i,j}s^it^j\dim_\kk H_{i,j}\Omega BG\phat
=\frac{1}{p_R(-s^{-1}t^{-1},t^{-1})}\\
&=\frac{(1+s^{-1}t^3)^2}{1-s^{-1}t-4s^{-1}t^2-s^{-1}t^3+s^{-2}t^4}. 
\end{align*}
Setting $s=1$ to ignore the Koszul grading, 
and cancelling $(1+t)^2$ from top and bottom gives the required Poincar\'e series
\begin{align*} 
\sum_{n=0}^\infty t^n\dim_\kk H_n\Omega BG\phat 
&=\frac{(1-t+t^2)^2}{1-3t+t^2}\\
&= 1 + t+5t^2+12t^3+32t^4+84t^5+220t^6+576t^7+\cdots 
\end{align*}
which agrees with the answer given in Section~2 of Benson~\cite{Benson:2009b}
in the case $p=3$. The zeros of the denominator occur at
$t=\frac{3\pm\sqrt 5}{2}$, and because $0<\frac{3-\sqrt 5}{2}<1$, the
radius of convergence is less than one, and the 
growth is exponential.

\begin{remark}
A similar analysis holds in larger rank. Let
\[ G=(\bZ/p)^r\rtimes\bZ/2=\langle g_1,\dots,g_r,\sigma\rangle \] 
with $g_i$
commuting elements of order $p$ and the involution $\sigma$ inverting each
$g_i$, $\sigma g_j=g_j^{-1}\sigma$.
Let $H$ be the normal subgroup of index two. 
The group algebra $\kG$ is $\bZ^r$-graded, with $s$ in degree zero and
$\sum_{i=1}^{p-1}g_j^i/i$ in degree one in the $j$th coordinate and
zero in the other coordinates. The cohomology ring $H^*BG$ is
equal to the invariants of the involution on $H^*BH$, and Molien's
theorem gives us the Poincar\'e
series
\[ p(s,t)=\half\left(\left(\frac{1+s^\halfs t^{-1}}
{1-s^\halfs t^{-2}}\right)^{\!r}+
\left(\frac{1-s^\halfs t^{-1}}
{1+s^\halfs t^{-2}}\right)^{\!r}\,\right). \]
 This is again a
Koszul algebra, with Koszul dual $\Ext^{*,*}_{H^*BG}(\kk,\kk)$. There
is no room for non-zero differentials in the Eilenberg--Moore spectral
sequence, and no ungrading problems, so we have $H_*\Omega
BG\phat=\Ext^{*,*}_{H^*BG}(\kk,\kk)$. This again has exponential
growth, with Poincar\'e series
\[ \frac{1}{p(-s^{-1}t^{-1},t^{-1})}=\frac{2{(1+s^{-1}t^3)}^r}
{{(1-(-st)^{-\halfs}t)}^r{(1-(-st)^{-\halfs}t^2)}^r+
{(1+(-st)^{-\halfs}t)}^r{(1+(-st)^{-\halfs}t^2)}^r} \]
via formula~\eqref{eq:pR!st-1}. The denominator is an even function of
$(-st)^{-\halfs}$, so this is a rational function of $s^{-1}$ and
$t$. Again we can set $s=1$ if we wish to ignore the Koszul grading,
but we have to do this after applying~\eqref{eq:pR!st-1} to obtain the
formula above.
\end{remark}

\section{Reflection groups}\index{reflection group}

We do not want to give the impression that polynomial growth for
$H_*\Omega BG\phat$ only
happens for finite or tame representation type. We therefore mention the
following. In the next two sections we give further examples of
polynomial growth.

\begin{theorem}
Suppose that $G$ is a semidirect product
$E\rtimes H$ with $E$ an elementary abelian $p$-group ($p$ odd), and 
$H$ a $p$-adic reflection group of order prime to $p$, acting on $E$ via the 
reduction modulo $p$ of the 
reflection representation.\index{p-adic@$p$-adic reflection group} 
Then $H^*BG$ is a polynomial tensor exterior algebra.
In this case, $H_*\Omega BG\phat$ is also usually polynomial tensor exterior,
and always has polynomial growth.
\end{theorem}
\begin{proof}
We have $H^*BG = (H^*BE)^H$, the invariants of $H$ on $H^*BE$.
It follows from a theorem of Solomon~\cite{Solomon:1963a} (see also
Section~7.3 of~\cite{Benson:1993a}) that 
$H^*BG$ is a polynomial algebra tensored with an exterior algebra. So
$\Ext^*_{H^*BG}(\kk,\kk)$ is also polynomial tensor exterior, and has
polynomial growth. Then the $E^2$ page of the Eilenberg--Moore spectral
sequence~\eqref{eq:Cotor}\index{Eilenberg--Moore spectral sequence} 
is also polynomial tensor exterior, and has polynomial growth. It follows
that $H_*\Omega BG\phat$ also has polynomial growth.
\end{proof}

\begin{remark}
In the theorem, if we use a single grading on $\kG$ by powers of the radical, the polynomial
generators for $H^*BG$ lie in degrees $(-2n_i,-pn_i)$ and the exterior
ones in degrees $(-2n_i+1,-p(n_i-1)-1)$, where $n_i$ runs over the degrees
of the fundamental invariants of the reflection group $H$. So the polynomial generators of 
$\Ext^*_{H^*BG}(\kk,\kk)$ are in degree $(-1,2n_i-1,p(n_i-1)+1)$ and the
exterior generators are in degrees $(-1,2n_i,pn_i)$. There is no room
for non-zero differentials, but it occasionally happens that the exterior
relations ungrade to have non-zero squares and commutators in the
polynomial part. An example of this is the symmetric group of degree
three at the prime three, with $E=\bZ/3$ and $H=\bZ/2$, see the last
part of Theorem~1.3 of Benson and Greenlees~\cite{Benson/Greenlees:2021a}.
\end{remark}

\section{Chevalley groups in non-defining characteristic}\label{se:non-defining}

In this section, we describe why, if $H$ is a finite Chevalley or
twisted Chevalley group in non-torsion 
non-defining characteristic $\ell$, $H_*\Omega BH\ellhat$ has polynomial
growth.\footnote{For the purposes of this section, we are working in
  characteristic $\ell$ rather than $p$, because we want to reserve $p$ for the
  characteristic of definition of the group. Also, $G$ denotes a
  compact Lie group, so the finite group is called $H$.} 
This is a consequence of a construction of
Quillen~\cite{Quillen:1972a}, elaborated in 
Friedlander \cite{Friedlander:1976a,Friedlander:1982a}, 
Fiedorowicz and Priddy~\cite{Fiedorowicz/Priddy:1978a},
Wilkerson~\cite{Wilkerson:1974a}, and
Kleinerman~\cite{Kleinerman:1982a}.

Let $G$ be a connected compact Lie 
group,\index{compact Lie group}\index{Lie group} and
let $H$ be a corresponding
finite group of Lie type over the finite field $\bF_{p^m}$. This is
a Chevalley group\index{Chevalley group} or twisted Chevalley
group\index{twisted Chevalley group}. 
Let
$\ell$ be a prime different from $p$, and let $\kk$ be a field of
characteristic $\ell$. Provided $\ell$ does not divide the order of
the Weyl group of $G$, Wilkerson~\cite{Wilkerson:1974a} showed that there
is an unstable Adams operation\index{Adams operation}%
\index{unstable Adams operation}
$\psi^{p^m}\colon BG\ellhat\to BG\ellhat$\index{psi@$\psi^{p^n}$} 
whose effect in cohomology in degree $2i$ is multiplication by $p^{mi}$.
Then there is a homotopy fibre square\index{homotopy fibre square}
\[ \xymatrix{BG(p^m)\ellhat \ar[r] \ar[d] & 
BG\ellhat \ar[d]^{\id  \times \id} \\
BG\ellhat \ar[r]^(0.4){\id \times\psi^{p^m}} & BG\ellhat \times BG\ellhat.} \]
In the case of a twisted group of Lie type, $\psi^{p^m}$ is replaced
by its composite with a diagram automorphism, and the same homotopy
fibre square results.

We wish to consider the Eilenberg--Moore spectral sequence of this fibre square 
with coefficients in $\kk$:
\begin{equation}\label{eq:EMBGp^m}
\Tor_{*,*}^{H^*BG\otimes H^*BG}
(H^*BG,H^*BG) 
\Rightarrow H^*BG(p^m) 
\end{equation}

Using \'etale homotopy theory,
Friedlander~\cite{Friedlander:1976a,Friedlander:1982a} showed that even in cases where $\ell$
does divide the order of the Weyl group, although the unstable Adams operations
do not literally exist, we can still behave as though such
a fibre square exists for the purposes of under\-standing
cohomology. The details are also summarised in Section~2 of 
Kleinerman~\cite{Kleinerman:1982a}. Recall that if $G$ is a compact
connected semisimple Lie group then there is a complex Lie group
$G_\bC$ for which $G$ is a maximal compact subgroup (see for example
Serre~\cite{Serre:1966a}), and $BG$ and $BG_\bC$ are homotopy
equivalent. By Kostant~\cite{Kostant:1966a}, there is a
reductive Chevalley group scheme $G_\bZ$ whose complex points are $G_\bC$.

\begin{theorem}
Let $G_\bZ$ be a reductive Chevalley group scheme and $\kk=\bar\bF_p$
and let $\phi\colon G_\kk\to G_\kk$ be a surjective endomorphism with
finite invariant group $G_\kk^\phi=H$. Let $\ell$ be a prime different
to $p$. Then 
there is a homotopy commutative square
\[ \xymatrix{BH\ellhat\ar[r]\ar[d]&
\displaystyle\lim_{\leftarrow} BG(\bC)\ellhat\ar[d]^\Delta\\
\displaystyle\lim_{\leftarrow}BG(\bC)\ellhat\ar[r]&
\displaystyle\lim_{\leftarrow}(BG(\bC)\times
BG(\bC))\ellhat} \]
This gives rise to a spectral sequence~\eqref{eq:EMBGp^m}.
\end{theorem}
\begin{proof}
The first part of this is Theorem~2.9 
of~\cite{Friedlander:1976a}. The limits here are associated to the result of applying a 
rigid \'etale type functor landing in pro-simplicial sets. Sections~3 
and~4 of~\cite{Kleinerman:1982a} explain how this gives rise to the
spectral sequence.
\end{proof}

Recall that if $G$ is a connected, simply connected compact Lie
group\index{compact Lie group}\index{simply connected}
and $p$ is a prime,
Borel~\cite{Borel:1953a,Borel:1954a,Borel:1960a,Borel:1960b,Borel:1967a} 
and Steinberg~\cite{Steinberg:1975a} have shown that the following are
equivalent (see also Section~5 of Benson and Greenlees~\cite{Benson/Greenlees:2023b}):
\begin{enumerate}
\item $H^*(G;\bZ)$ has no $p$-torsion.
\item $H^*(G;\bF_p)$ is an exterior algebra on odd degree classes.
\item $H^*(BG;\bF_p)$ is a polynomial algebra on even degree classes.
\item $H^*(BG;\bZ)$ has no $p$-torsion.
\item $H^*(G/T;\bF_p)$ is generated by elements of degree two, where
  $T$ is a maximal torus.
\item Every elementary abelian $p$-subgroup of $G$ is contained in a
  torus.
\item Every elementary abelian $p$-subgroup of $G$ of rank at most
  three is contained in a torus.
\item The multiplicity of some fundamental root in the dominant coroot
  is divisible by $p$.
\end{enumerate}

The primes $p$ for which this occurs are called the 
\emph{torsion primes}\index{torsion prime} for $G$.
They are a subset of the 
\emph{bad primes},\index{bad prime} which are those for which the
multiplicity of some fundamental root in the dominant root is
divisible by $p$.
Here is a table of the torsion primes and bad primes.
\begin{center}
\begin{tabular}{|ccc||ccc|} \hline 
type & bad & torsion & type & bad & torsion \\ \hline
$A_n$&$\varnothing$&$\varnothing$&$E_6$&$\{2,3\}$&$\{2,3\}$ \\
$B_n$&$\{2\}$&$\{2\}$&$E_7$&$\{2,3\}$&$\{2,3\}$ \\
$C_n$&$\{2\}$&$\varnothing$&$E_8$&$\{2,3,5\}$&$\{2,3,5\}$ \\
$D_n$&$\{2\}$&$\{2\}$&$F_4$&$\{2,3\}$&$\{2,3\}$\\
&&&$G_2$&$\{2,3\}$&$\{2\}$ \\ \hline 
\end{tabular}\medskip
\end{center}

\begin{theorem}
Let $G$ be a connected compact Lie group,  
$G_\bZ$ be the corresponding reductive Chevalley group scheme and $\kk=\bar\bF_p$, 
and let $\phi\colon G_\kk\to G_\kk$ be a surjective endomorphism with 
finite invariant group $G_\kk^\phi=H$. Let $\ell$ be a prime different from
$p$, and suppose that $H^*BG\ellhat$ is a polynomial ring. This happens, for
example, when $G$ is simply connected and $\ell$ is a non-torsion prime.

If $\ell$ is odd then $H^*BH$ is a polynomial tensor exterior algebra.
If $\ell=2$ then $H^*BH$ is a complete
intersection. In both cases, $H_*\Omega BG\ellhat$ has polynomial growth.
\end{theorem}
\begin{proof}
Since $\ell$ is not a torsion prime for $G$, $H^*BG$ is a polynomial
ring. In this case, the Eilenberg--Moore spectral
sequence~\eqref{eq:EMBGp^m} stops at 
the $E^2$ page, and gives a finite filtration on
$H^*BG(p^m)$ whose associated graded is a polynomial
algebra tensored with an exterior algebra. 

If $\ell$ is an odd prime
then there is no ungrading problem, and this gives the 
structure of the cohomology as a polynomial tensor exterior algebra. 
On the other hand, 
if $\ell=2$, it can happen that the
exterior generators ungrade to give elements whose square is not
necessarily zero, but is expressible in terms of the other generators.
The exact relations can be difficult to determine, but there is one relation
for the square of each exterior generator in the $E^\infty$
page. There is no problem with ungrading the commutation relations
between the exterior generators because the answer is supposed to be
graded commutative.
So independently of the
exact relations, the answer is a complete
intersection in characteristic two.

We apply another Eilenberg--Moore spectral sequence
(see Remark~\ref{rk:EMSS})
\[ \Ext^{*,*}_{H^*BG(p^m)}(\kk,\kk) 
\Rightarrow H_*\Omega BG(p^m)\ellhat. \]
By Theorem~\ref{th:ExtRkk}, the $E_2$ page of this spectral sequence
has polynomial growth. It is finite over its centre, and the centre is
finitely generated as a $\kk$-algebra. It follows that $E_\infty$ and
$H_*\Omega BG(p^m)$ have the same property.
\end{proof}

\begin{example}[Quillen~\cite{Quillen:1972a}, Friedlander~\cite{Friedlander:1976a}]
Let $G=U(n)$, of Lie type $A_{n-1}$. This is not simply connected, but
nonetheless $H^*BU(n)$ is a polynomial ring. We have
\[ H^*BU(n)=\kk[c_1,\dots,c_n], \] 
where the $c_i$ are the
Chern classes\index{Chern classes} of degree $2n$. Then
$G(p^m)$ is the general linear group
$GL(n,p^m)$\index{GL@$GL(n,p^m)$}. 

For $\ell$ odd, this gives 
\[ H^*BGL(n,p^m) = \kk[c_r,c_{2r},\dots,c_{tr}] \otimes
\Lambda(e_r,e_{2r},\dots,e_{tr}) \]
where $r$ is the order of $p^m$ modulo $\ell$, and $t$ is the integer
part of $n/r$. The degrees are $|c_{ir}|=-2ir$,
$|e_{ir}|=-2ir+1$. Then the associated graded of
$H_*\Omega BGL(n,p^m)\ellhat$ is
\[ \kk[\hat e_r,\hat e_{2r},\dots,\hat e_{tr}] \otimes
\Lambda(\hat c_r,\hat c_{2r},\dots,\hat c_{tr}) \]
with $|\hat e_{ir}|=2ir-2$, $|\hat c_{ir}|=2ir-1$.
Beware that there is no reason why the answer should be graded
commutative, so it is not obvious how to ungrade the square zero
relations for the $\hat c_{ir}$. For example, in characteristic three
we have 
\[ H_*\Omega BGL(2,2)\threehat=\kk[\hat e_2,\hat c_2]/
(\hat c_2^2+\hat e_2^3), \] 
see Section~\ref{se:cyclic}. This is commutative, but not graded
commutative. But in any case, the answer has 
polynomial growth.

For $\ell=2$, we have $r=1$ and $t=n$.  In this case, if 
$p^m\equiv 1\pmod{4}$ we get the same answer as above, but if
$p^m\equiv 3 \pmod{4}$
then we have 
\[ e_j^2=\sum_{a=0}^{j-1}c_{a}c_{2j-1-a}. \]
Here, $c_{2j-1-a}$ is interpreted as zero if $2j-1-a>n$. This gives a
complete intersection with Krull dimension $n$, with
$n+\lfloor\frac{n}{2}\rfloor$ generators $e_1,\dots,e_n,c_2,c_4,\dots$ and 
$\lfloor\frac{n}{2}\rfloor$ relations. So by Theorem~\ref{th:ExtRkk}, 
the $E^2$ term of the Eilenberg--Moore spectral sequence is finite as
a module over a central polynomial subring with
$\lfloor\frac{n}{2}\rfloor$ generators. So $H_*\Omega
BGL(n,p^m)\ellhat$ has polynomial growth.
\end{example}

\begin{example}[Kleinerman~\cite{Kleinerman:1982a}]
Let $G=G_2$ and $\ell=2$. Two is a torsion prime for $G$, and we have 
\[ H^*BG_2 = \kk[d_4,d_6,d_7]. \]
The Eilenberg--Moore spectral sequence gives
the associated graded of $H^*BG_2(p^m)$ ($p$ an odd prime) 
to be $\kk[d_4,d_6,d_7]
\otimes \Lambda(y_3,y_5,y_6)$. Ungrading the relations gives
$y_3^2=y_6$, $y_5^2=y_3d_7+y_6d_4$ and
$y_6^2=y_5d_7+y_6d_6$ (Grbi\'c~\cite{Grbic:2006a}). 
So $H^*BG_2(p^m)$ is the complete intersection
\[ \kk[d_4,d_6,d_7,y_3,y_5]/(y_5^2+y_3d_7+y_3^2d_4,\,
y_3^4+y_5d_7+y_3^2d_6). \]
Using Theorem~\ref{th:ExtRkk}, we see that 
the $E^2$ page of the Eilenberg--Moore
spectral sequence for $H_*\Omega BG_2(p^m)\twohat$ is generated 
over the central subalgebra $\kk[s_{10},s_{12}]$ by elements
$\hat d_4$, $\hat d_6$, $\hat d_7$, $\hat y_3$, $\hat y_5$.
The relations say that all squares and commutators of the latter elements are zero
except for
\[ \hat y_5^2 = [\hat d_7,\hat y_3] = s_{10},\qquad
[\hat d_7,\hat y_5]=s_{12}. \]
There is no room for differentials, but ungrading the $E^\infty$ page
requires some work. This is done in the paper of Levi and
Seeliger~\cite{Levi/Seeliger:2012a}, where they also compute the 
coproduct and action of the dual Steenrod algebra.
It turns out that $E^\infty$
as given above is isomorphic to $H_*\Omega BG_2(p^m)\twohat$.
The degrees are added, so that the elements $\hat d_i$ and $\hat y_i$
now have degree $i-1$, and the elements $s_j$ have degree $j-2$.

This is of polynomial growth, since it is finitely generated (actually free of
rank $2^5$) over the 
central polynomial subalgebra $k[s_{10},s_{12}]$
with Poincar\'e series
\[ \sum_{n=0}^\infty \dim_\kk H_n\Omega BG_2(p^m)\twohat =
\frac{(1+t^2)(1+t^3)(1+t^4)(1+t^5)(1+t^6)}{(1-t^8)(1-t^{10})}
=\frac{(1+t^3)(1+t^6)}{(1-t^2)(1-t^5)}. \]
\end{example}

\section{\texorpdfstring{An exotic example: $B\Sol(q)$}
{An exotic example: BSol(q)}}\index{BSol@$B\Sol(q)$}

Let $q$ be an odd prime power, and let 
$\Sol(q)$\index{Sol@$\Sol(q)$} be the exotic 
Benson--Solomon $2$-local finite group. This was
originally discussed as a configuration that was 
proved not to come from a finite group in
Solomon~\cite{Solomon:1974a}. 
Its classifying space $B\Sol(q)$ was then discussed in Benson~\cite{Benson:1998a}, 
and finally it was constructed as a fusion system and linking system
by Levi and Oliver~\cite{Levi/Oliver:2002a}. 
Using the fibre square like the one in the previous section,
the associated graded of
$H^*B\Sol(q)$ was computed in~\cite{Benson:1998a} to be a polynomial
ring on generators in degrees $8$, $12$, $14$ and $15$ tensored with
an exterior algebra on generators in degrees $7$, $11$, $13$ and $14$.
The ungrading was carried out 
by Grbi\'c~\cite{Grbic:2006a}, who computed it
to be the codimension three complete intersection
\[ H^*B\Sol(q)=\kk[u_8,u_{12},u_{14},u_{15},y_7,y_{11},y_{13}]/(f_{22},f_{26},f_{28}) \]
where the (homological) degrees are minus the subscripts, and where 
\begin{align*}
f_{22}&=y_{11}^2+u_8y_7^2+u_{15}y_7,\\
f_{26}&=y_{13}^2+u_{12}y_7^2+u_{15}y_{11},\\
f_{28}&=y_7^4+u_{14}y_7^2+u_{15}y_{13}.
\end{align*}

Applying Theorem~\ref{th:ExtRkk}, we find that
$\Ext_{H^*B\Sol(q)}^{**}(\kk,\kk)$ is generated over a central subalgebra
$\kk[s_{22},s_{26},s_{28}]$ by elements $\hat u_8,\hat u_{12},\hat u_{14},\hat
u_{15},\hat y_7,\hat y_{11},\hat y_{13}$. The degrees of the elements
$\hat u_i$ and $\hat y_i$ are $(-1,i)$, while the degrees of the $s_j$
are $(-2,j)$. The relations say that all
squares and commutators of the latter elements are zero except for
\begin{equation}\label{eq:Sol-rel}
  \hat y_{11}^2=[\hat u_{15},\hat y_7]=s_{22},\qquad
\hat y_{13}^2=[\hat u_{15},\hat y_{11}]=s_{26},\qquad
[\hat u_{15},\hat y_{13}]=s_{28}.
\end{equation}
In the Eilenberg--Moore spectral sequence
\[ \Ext_{H^*B\Sol(q)}^{**}(\kk,\kk) \Rightarrow H_*\Omega B\Sol(q)\twohat \]
there is no room for non-zero differentials, but ungrading the
$E^\infty$ page takes more work. This is done in the paper of 
Levi and Seeliger~\cite{Levi/Seeliger:2012a}, where they also compute
the coproduct and action of the dual Steenrod algebra. It turns out that
$E^\infty$ as given above is isomorphic to
$H_*\Omega B\Sol(q)\twohat$. The degrees are added, so that the
$\hat u_i$ and $\hat y_i$ now have degree $i-1$ and the $s_j$ have
degree $j-2$.

\begin{theorem}
$H_*\Omega B\Sol(q)\twohat$ has a central subalgebra $\kk[s_{22},s_{26},s_{28}]$,
where the degrees are two less than the subscripts. Over this it is generated
by elements
$\hat u_8$,
$\hat u_{12}$, $\hat u_{14}$, $\hat u_{15}$, $\hat y_7$,
$\hat y_{11}$, $\hat y_{12}$, where the degrees are one less than the
subscripts. The relations say that the squares and commutators of
the latter elements are zero except for the
relations~\eqref{eq:Sol-rel} above.
\end{theorem}

This algebra $H_*\Omega B\Sol(q)\twohat$ is of polynomial growth,
since it is finitely generated (actually free of rank $2^7$) 
as a module over the central polynomial 
subalgebra $\kk[s_{22},s_{26},s_{28}]$, with Poincar\'e series
\begin{align*} 
\sum_{n=0}^\infty t^n\dim_\kk H_n\Omega B\Sol(q)\twohat &=
\frac{(1+t^7)(1+t^{11})(1+t^{13})(1+t^{14})(1+t^{6})
(1+t^{10})(1+t^{12})}{(1-t^{20})(1-t^{24})(1-t^{26})} \\
&=\frac{(1+t^7)(1+t^{11})(1+t^{14})
}{(1-t^{6})(1-t^{10})(1-t^{13})}.
\end{align*}
The polynomial growth conforms to the idea that we should expect 
$\Sol(q)$ to behave like a finite group of Lie type in non-defining characteristic.

\bibliographystyle{amsplain}
\bibliography{../repcoh}
\printindex

\end{document}